\documentclass[a4paper,12pt,oldfontcommands]{memoir}

\usepackage{coursmonoeng}
\usepackage[backend=biber, style=ieee-alphabetic]{biblatex}

\usepackage{titlesec}

\newlength\droppp
\newcommand*{\titleEI}{\begingroup
\droppp = 0.3\textheight
\vspace*{\droppp}
\raggedright
{\LARGE \emph{{\huge S}ome aspects of}  \\
\hspace*{20pt} \textbf{\huge {T}opological dynamics of Polish \\ \hspace*{30pt}groups} \par}
\vspace{1.5\baselineskip}
{\large \hspace*{25pt} \emph{{\LARGE W}ith an introduction to descriptive set theory} \\[.1\baselineskip]}
\vspace{3\baselineskip}
\hspace*{15pt}{ \large {\textit{\footnotesize BY}} \Large{Julien Melleray}}\par
\vfill
\endgroup}

\hypersetup{citecolor={navyblue}, urlcolor={navyblue}}

\DeclareMathOperator{\ucb}{UC_b}
\DeclareMathOperator{\cli}{C_L}

\newcommand{\cl}[2][]{\overline{#2}^{#1}}
\newcommand{\sub}{\subseteq}
\newcommand{\ana}{\mathbf{\Sigma}_{1}^{1}}
\newcommand{\coana}{\mathbf{\Pi}_{1}^{1}}
\newcommand{\sg}[1]{\underset{#1}{+}}

\makeindex
\begin{document}

\frontmatter
\thispagestyle{empty}

\title{Some aspects of topological dynamics of Polish groups \\ \begin{small} With an introduction to descriptive set theory  \end{small}}
\author{Julien Melleray}
\date{}

\begin{titlingpage}

\titleEI

\end{titlingpage}

\newpage
\thispagestyle{empty}
\setlength{\epigraphwidth}{0.45\textwidth}
\epigraph{\flushright \bf \Large{Foreword}}{} 

These notes grew out of a semester course given at the university of Turin during the spring of 2024 (the content of that course corresponds roughly to Chapters 1--8). They are aimed at graduate (or advanced undergraduate) students and researchers who want an introduction to some topics around the theory of \emph{Polish groups}, i.e., completely metrizable second-countable topological groups. This class of groups contains all second-countable locally compact Hausdorff topological groups; Polish groups are ubiquitous in functional analysis, ergodic theory, logic.... The primary focus of the first part of the book is an introduction to some aspects of topological dynamics of Polish group actions.  We are concerned with continuous actions of Polish groups by homeomorphisms on Hausdorff compact spaces, and develop some of the general theory. One of our main objectives in that part of the book is to establish a theorem due to Kechris, Pestov and Todor\v{c}evic, which characterizes automorphism groups of countable structures which are \emph{extremely amenable}, which means that every con\-ti\-nu\-ous action of such a group on a compact Hausdorff space has a fixed point. This theorem establishes a connection between topological dynamics and combinatorics and has proved to be extremely influential.

We do not take the shortest route to this result, instead taking time to discuss some of the background and present various methods, in particular based on Baire-category-theoretic notions, that are useful when working with Polish groups. In the first part of the book we sometimes need standard results from \emph{descriptive set theory}, which we use as black boxes; in the second part all those results (and more) are proved.

Descriptive set theory is, in a nutshell, the study of the properties of ``definable'' subsets of Polish spaces; this classical theory goes back more than a century and many topics are not covered here. I chose which topics to present based on personal taste and ease of presentation, with an eye on applicability in various areas of mathematics. The culminating point of the second part of the book is a proof, due to B. Miller, of a fundamental dichotomy theorem of Kechris, Solecki and Todor\v{c}evic; this dichotomy theorem has striking structural consequences, some of which we briefly discuss. My objective in this part is twofold: developing some useful tools, and hopefully motivating some readers to dig deeper into descriptive set theory and its recent developments. Several classical, important and beautiful parts of descriptive set theory are not covered at all in these notes; I hope the reader will want to know more and use this book as a first step towards more thorough presentations of the topic such as \cite{Kechris1995}.

This text uses some set-theoretic notations (notably, $\omega$ for the set of non negative integers); but there is little actual set-theoretic background, and I give a quick presentation of that background (as well as the axiomatic framework which we are using; notably, we freely use the axiom of choice) in an appendix. I tried to make this text accessible to readers without a background in logic, set theory or dynamics - some familiarity with point set topology (including the Baire category theorem) as well as group actions is however required.

The text is sprinkled with exercises of varying difficulty (some of them should be very easy, while others seem quite difficult to me); solutions to these exercises are given at the end of the text. Sometimes the reader is expected to fill in some of the details of the argument. It is my belief that at least attempting to solve many of these exercises while going through the text will help the reader get a good grasp on the material. 

Each chapter ends with a few comments, mostly providing references for the material presented and proposing some further reading.

I would like to thank the students who attended my lectures in Turin and worked on many of the exercises; in particular Paolo Boldrini, Fabio Carrozzo and Laura Paglietta deserve a mention here. I am also grateful to Raphaël Carroy for organizing my visit, to Luca Motto Ros for graciously sharing his office with me, and to the whole logic group for their welcome. This visit was made possible by the Unito Visiting Professor Program 2023/2024. 

Bruno Duchesne suggested turning my lecture notes into a publishable text; I am not sure I should thank him for all the additional work that this entailed, but will do so nevertheless. I am beholden to Sophie Grivaux for her encouragement and her very efficient handling of the manuscript.

I owe a debt of gratitude to Todor Tsankov, who explained to me much of what I know about topological dynamics of Polish groups. I am also indebted to an anonymous reviewer, who made many useful suggestions and caught an erroneous argument of mine; and to Andy Zucker for helping me understand why that argument was in fact incorrect. François Le Maître was of great help when I was working on the final version of this manuscript; for this and our many discussions over the years I thank him. Finally, I am much obliged to Johan Koschall and Flore Le Roux for their comments and feedback. 

There is unfortunately little doubt that some errors are still lurking - finding and correcting them is left as an exercise for the reader.

\newpage 

\tableofcontents*

\mainmatter
\epigraphhead{}
\clearpage
\epigraphhead{}
\pagestyle{plain}
\setlength{\epigraphwidth}{0.45\textwidth}

\part{Topological dynamics of Polish groups}

\clearpage


\chapter{Polish groups and their actions}

In the first part of this book, we are mainly concerned with actions of Polish groups (defined below) on compact Hausdorff spaces. In this chapter we cover some basics of Polish group theory; we use some descriptive-set-theoretic results which are only proved in the second part of the book.

We make the following convention, which will simplify the presentation a little : in this text, unless explicitly stated otherwise, {\bf we only consider Hausdorff topological spaces}.

\begin{defin}
A \emph{topological group}\index{topological group} is a group $(G,\cdot)$ endowed with a topology $\tau$ for which the group operations $g \mapsto g^{-1}$ and $(g,h) \mapsto g \cdot h$ are continuous.
\end{defin}

This could be stated more concisely by requiring $(g,h) \mapsto g \cdot h^{-1}$ to be continuous; in what follows we will simply write $gh$ instead of $g \cdot h$ since the group law should always be clear from the context. Similarly, when working with actions we will either write $g \cdot x$ or $gx$ if there is no risk of ambiguity. We denote the neutral element of $G$ either by $1$ or $1_G$.

\begin{defin}
Let $X$ be a topological space. We say that $(X,\tau)$ is \emph{Polish}\index{Polish space} if the topology of $X$ is separable and is induced by a complete metric.
\end{defin}
Note that separability of $X$ in this context is equivalent to $X$ being second-countable, i.e., that there exists a countable basis for the topology of $X$ (since for metrizable spaces, this condition is equivalent to separability).

Polish spaces (and their subsets)  form the backdrop of much of classical descriptive set theory, which is covered in more detail in the second part of the book. In this first part, concerned with dynamics, Polishness is a condition that we require of the topology of the acting group; completeness enables us to use the tools associated with the Baire category theorem.

\begin{defin}
A \emph{Polish group}\index{Polish group} is a topological group whose underlying topology is Polish.
\end{defin}

\begin{example}
Consider the group $\Sinf$ of all permutations of $\omega$ (the group law being given by composition of maps). See it as a subset of the Baire space $\omega^\omega$, and endow it with the induced topology. Explicitly, a basis of neighborhoods of the identity is given by the following clopen subsets (subgroups, actually):
\[ U_F = \lset \sigma \in \Sinf : \forall i \in F \ \sigma(i)=i \rset \]
where $F$ ranges over all finite subsets of $\omega$.

A compatible metric for this topology is given by 
\[d(\sigma,\tau) = \inf \lset 2^{-n} : \forall i < n \ \sigma(i)= \tau(i) \rset .\]
\end{example}

\begin{exo}\label{exo1}$ $

\begin{enumerate}
\item Show that the group operations on $\Sinf$ are continuous.
\item Use the sequence $(\sigma_i)_{i < \omega}$ defined by $\sigma_i(n)= n+1$ for $n \le i$, $\sigma_i(i+1)=0$, $\sigma_i(n)=n$ for $n >i+1$ to show that $d$ is not complete.
\item Show however that the metric $\rho$ defined by $\rho(\sigma,\tau)= d(\sigma,\tau)+ d(\sigma^{-1},\tau^{-1})$ is complete and conclude that $\Sinf$ is a Polish group.
\item Prove that $\Sinf$ is a $G_\delta$ subset of $\omega^\omega$ to obtain another proof that $\Sinf$ is Polish.
\end{enumerate}
\end{exo}

To solve the last question of the previous exercise, one needs to know that a subset $X$ of a completely metrizable topological space $Y$ is completely metrizable for the induced topology iff $X$ is a $G_\delta$ subset of $Y$. This is an important fact, which we will use repeatedly. This result, like the other results from descriptive set theory that are used without proof in the first part of this book, is proved in the second part (for this particular theorem, see Corollary \ref{coro:G_delta=Polish}). Most, if not all, of these descriptive-set-theoretic facts only come up in the first chapter (though their consequences are used later on).

\begin{exo}\label{exo2}
Prove that a closed subgroup $G$ of $\Sinf$ is compact if, and only if, the $G$-orbit of every element of $\omega$ is finite. 
\end{exo}

\begin{rem}
We see that the topology of $\Sinf$ is induced by a metric which is left-invariant but not complete (the metric $d$ above) as well as a metric $\rho$ which is complete but neither left- nor right-invariant. We will see shortly that this is part of a broader phenomenon.
\end{rem}

\begin{exo}\label{exo3}
Prove that a (at most) countable product of Polish groups, endowed with the product topology, is a Polish group.
\end{exo}

Clearly a closed subgroup of a Polish group is itself a Polish group; the class of all closed subgroups of $\Sinf$ is a rich topic of study, with a strong interaction with model theory.

\begin{exo}\label{exo4}
Let $G$ be a Polish group. Show that $G$ is isomorphic, as a topological group, to a closed subgroup of $\Sinf$ iff $1_G$ has a basis of neighborhoods consisting of open subgroups. Such Polish groups are called \emph{nonarchimedean}\index{nonarchimedean Polish group}.
\end{exo}

Nonarchimedean Polish groups are those that we will be most interested in in these notes; let us discuss briefly some other examples.

\begin{prop}
Let $(X,d)$ be a Polish metric space, and $G$ be its isometry group. Then $G$, endowed with the pointwise convergence topology, is a Polish group.
\end{prop}

\begin{proof}
As in the case of $\Sinf$ (which is the isometry group of $\omega$ endowed with the discrete metric) one can either argue by defining a complete metric on $G$, or by proving that it is a $G_\delta$ subset of a Polish space. We will use the second approach here (of course in both approaches one also needs to prove that group operations are continuous).

Fix a countable dense subset $(x_i)_{i < \omega}$ of $X$, then define
$\Phi \colon G \to X^\omega$ by $\Phi(g)(i)=g(x_i)$. 

By definition of the pointwise convergence topology (which is the topology on $G$ induced by the product topology on $X^X$) the map $\Phi$ is continuous. It is also injective: if $f \ne g \in G$ then there exists $x \in X$ and $\varepsilon >0$ such that $d(f(x),g(x)) > 2 \varepsilon$, and then the triangle inequality implies that $f(x_i) \ne g(x_i)$ for any $i$ satisfying $d(x_i,x)< \varepsilon$, whence $\Phi(f) \ne \Phi(g)$.

Then we note that $\Phi$ is a homeomorphism onto its image; let $U$ be open in $G$, without loss of generality we may assume that 
\[U=\lset g \in G : \forall a \in A  \ d(g(a),f(a))< \varepsilon_a \rset\] where $f \in G$, $A \subset X$ is finite and each $\varepsilon_a$ is $>0$. 

We have to show that $\Phi(U)$ is open in $\Phi(G)$, to do that we first fix $g \in U$. Then let $\varepsilon >0$ be such that $d(g(a),f(a)) + 3 \varepsilon < \varepsilon_a$ for all $a \in A$. For each $a \in A$ pick $i_a$ such that $d(x_{i_a},a) < \varepsilon$. Any $h \in G$ such that $d(g(x_{i_a}),h(x_{i_a})) < \varepsilon$ for all $a \in A$ satisfies $d(h(a),f(a)) < \varepsilon_a$ for all $a \in A$, and this gives us an open neighborhood of $\Phi(g)$ contained in $\Phi(U)$, as desired.

What we have proved so far is that the pointwise topology on $G$ is completely understood by looking only at countably many coordinates, hence is metrizable (more formally, $\Phi(G)$ is metrizable since it is a subspace of the metrizable topological space $X^\omega$). Recall that convergence of sequences for a product topology is easily described: $(g_n)_{n < \omega}$ converges to $g$ if, and only if, $(g_n(x))_{n< \omega}$ converges to $g(x)$ for all $x \in X$. 

Let us now check that the group operations are continuous. Assume that $(g_n)_{n < \omega}$ converges to $g$ in $G$. Fix $a \in X$ then let $b=g^{-1}(a)$. Fix $\varepsilon >0$. We have that $a_n=g_n(b)$ converges to $g(b)=a$, so for $n$ large enough $d(a_n,a) < \varepsilon$. Using that $g_n^{-1}$ is an isometry we obtain that for $n$ large enough $d(g_n^{-1}(a), b) < \varepsilon$, so $(g_n^{-1}(a))_{n < \omega}$ converges to $b$, and this proves that $g \mapsto g^{-1}$ is continuous. Using the fact that 
\begin{align*}
d(g_nh_n(x),gh(x)) &\le d(g_nh_n(x),g_nh(x))+ d(g_nh(x),gh(x)) \\
&= d(h_n(x),h(x))+ d(g_nh(x),gh(x)),
\end{align*}
continuity of $(g,h) \mapsto gh$ is easily proved and we leave the details to the reader.

Finally, we note that $\Phi(G)$ is a $G_\delta$ subset of $X^\omega$. Indeed, $f \in X^\omega$ belongs to $\Phi(G)$ if, and only if, it satisfies the following two conditions:
\begin{itemize}
\item For all $i,j \in \omega$  $d(f(i),f(j)) = d(x_i,x_j)$ (by completeness, a distance-preserving map defined on a dense subset extends to a distance-preserving map defined on the whole space) 
\item $\lset f(i) : i \in \omega\rset$ is dense in $X$ (since $X$ is complete, the image of a distance-preserving map is closed, so the map is surjective as soon as its image is dense).
\end{itemize}
This yields the following description of $\Phi(G)$:
\[ \Phi(G)= \bigcap_{i,j < \omega} \lset f : d(f(i),f(j)) = d(x_i,x_j) \rset \cap \bigcap_{i < \omega} \bigcap_{\varepsilon \in \Q^+} \bigcup_{j < \omega} \lset f : d(f(j),x_i) < \varepsilon \rset \] 

Since closed subsets of a metrizable space are $G_\delta$, and a countable intersection of $G_\delta$ subsets is again $G_\delta$, we conclude as expected that $\Phi(G)$ is $G_\delta$ in $X^\omega$.
\end{proof}

We already mentioned that closed subgroups of Polish groups are themselves Polish groups when endowed with the induced topology; actually, the converse is also true.

\begin{thm}
Let $G$ be a Polish group, and $H$ be a Polish subgroup of $G$, i.e.~a subgroup of $G$ which is a Polish group when endowed with the induced topology. Then $H$ is a closed subgroup of $G$.
\end{thm}

\begin{proof}
Since $H$ is Polish, it is a $G_\delta$ subset of $G$. For any $g \in \overline{H}$, $H$ and $gH$ are then dense $G_\delta$ subsets of the Polish space $\overline{H}$. The Baire category theorem then implies that $H \cap gH \ne \emptyset$, whence $g \in H$ and it follows that $\overline{H}=H$.
\end{proof}

\begin{defin}\label{def:Baire_measurable}
Let $X$ be a Polish space. We say that $A \subseteq X$ is \emph{Baire measurable}\index{Baire measurable subset} if there exist an open subset $O$ and a meager subset $M$ of $X$ such that $A=O \mathrel{\Delta} M$.
\end{defin}

We recall that Baire measurable subsets of $X$ form a $\sigma$-algebra; this $\sigma$-algebra is the smallest which contains open subsets as well as meager sets, and it contains all Borel subsets of $X$ (as well as the analytic and coanalytic subsets). For more details on the corresponding descriptive set-theoretic background, see the second part of these notes and references therein.

In contexts where there is no quasi-invariant $\sigma$-finite measure (as opposed to the case of locally compact groups which come endowed with the Haar measure) this is a very useful $\sigma$-algebra. The meager sets provide us a well-behaved notion of smallness, and the Kuratowski--Ulam theorem (which we discuss in some detail later) is a suitable analogue of the Fubini theorem, though of course the Baire-category-theoretic notions are much less quantitative than their measured counterparts.

We recall that a Borel measure on a Polish group $G$ is \emph{left quasi-invariant} if for any $g$ and any Borel $A$ we have $\mu(A)=0 \Leftrightarrow \mu(gA)=0$. The next result shows that, within the class of Polish groups, existence of a nontrivial left-quasi invariant measure is already enough to imply that $G$ is locally compact; so one really needs to accept that there will in general not be a single measure providing a well-behaved notion of smallness of subsets of $G$.

\begin{thm}[Weil]
Let $G$ be a Polish group, and $\mu$ be a nontrivial Borel $\sigma$-finite measure on $G$ which is left quasi-invariant, i.e.,such that for any Borel subset $A$ such that $\mu(A)=0$ one also has $\mu(gA)=0$ for all $g \in G$. Then $G$ is $\sigma$-compact (hence also locally compact since it is Polish).
\end{thm}

\begin{proof}
Let $A$ be a Borel subset of $G$ such that $0 < \mu(A) < + \infty$. There exists a compact subset $K \subseteq A$ such that $0 < \mu(K)$ (this follows from instance from the fact that one can refine the topology of $A$, without changing its Borel structure, to turn $A$ into a Polish space, combined with the inner regularity of finite Borel measures on Polish spaces; for this last property see Theorem 17.11 of \cite{Kechris1995}). Let $H$ be the subgroup of $G$ which is generated by $K$; since $H= \bigcup_{n < \omega} (K \cup K^{-1})^n$ we have that $H$ is $\sigma$-compact.

If $H$ has uncountable index in $G$, there is an uncountable $F \subset G$ such that $aH \cap bH = \emptyset$ for each $a \ne b$ in $F$, hence also $aK \cap bK = \emptyset$ for every $a \ne b \in F$. Note that for each $a \in F$ we have $\mu(aK)>0$.

Since $\mu$ is assumed to be $\sigma$-finite, there exists a sequence $(B_n)_{n< \omega}$ of subsets of $G$ such that each $\mu(B_n)$ is finite and $G= \bigcup_{n < \omega} B_n$. 

For each $n$, the family $(\mu(B_n \cap aK))_{a \in F}$ is summable, whence for all $n$ there are only at most countably many $a \in F$ such that $\mu(B_n \cap aK) >0$. But then there are only at most countably many $a \in A$ such that $\mu(aK) >0$, a contradiction.

Thus $H$ has countable index in $G$; since $H$ is $\sigma$-compact it follows that $G$ is also $\sigma$-compact, hence locally compact (applying the Baire category theorem shows that some compact subset of $G$ must have nonempty interior, so every element of $G$ has a compact neighborhood).
\end{proof}

\begin{defin}
Let $X$ be a Polish space, and $A$ be a subset of $X$. We denote\index{$U(A)$}
\[U(A) = \bigcup \lset O \text{ open in } X : O \setminus A \text{ is meager} \rset . \]
\end{defin}

Since $X$ is second-countable, there exists a countable family $(O_n)_{n< \omega}$ of open subsets of $X$ such that $O \setminus A$ is meager for all $n$, and $U(A)= \bigcup_n O_n$ (see Exercise \ref{exo5}). Thus $U(A) \setminus A$ is meager, without any definability assumption on $A$.

\begin{thm}
Let $X$ be a Polish space and $A$ a subset of $X$. Then $A$ is Baire measurable iff $U(A) \mathrel{\Delta} A$ is meager.
\end{thm}

\begin{proof}
Clearly if $U(A) \mathrel{\Delta} A$ is meager then $A$ is Baire measurable; conversely, assume that $A$ is Baire measurable and let $O$ be open such that $A \mathrel{\Delta} O$ is meager. By definition $O$ is contained in $U(A)$, whence $A \setminus U(A)$ is meager, so $A \mathrel{\Delta} U(A)$ is meager.
\end{proof}

The proof above applies only to second-countable spaces; one can however generalize this statement to any topological space, see e.g. Theorem 8.29 of \cite{Kechris1995}.

Continuity of group operations implies that, for a Polish group $G$, a subset $A$ of $G$ and $g \in G$ we have $U(A^{-1})=U(A)^{-1}$ as well as $U(gA)=gU(A)$.

\begin{lem}[Pettis]\index{Pettis lemma}
Let $G$ be a Polish group, and $A,B$ be two subsets of $G$.

 Then $U(A) U(B) \subseteq AB$.
\end{lem}

\begin{proof}
Let $g$ belong to $U(A)U(B)$. Then $gU(B)^{-1} \cap U(A) \ne \emptyset$, i.e.,$U(gB^{-1}) \cap U(A) \ne \emptyset$. Let $V$ denote the nonempty open subset $U(gB^{-1}) \cap U(A)$.

Since $A$ is comeager in $U(A)$, and $V$ is an open subset of $U(A)$, $A$ is also comeager in $V$; similarly $gB^{-1}$ is comeager in $V$. Applying the Baire category theorem in $V$, which is a Polish space for the induced topology, we obtain that $A \cap gB^{-1} \ne \emptyset$, equivalently $g \in AB$. This concludes the proof.
\end{proof}

\begin{coro}
Let $G$ be a Polish group and $A$ a Baire measurable, non-meager subset of $G$. Then $1$ belongs to the interior of $A A^{-1}$.
\end{coro}

\begin{proof}
This is immediate: since $U(A)$ is nonempty, $1 \in U(A)U(A)^{-1}=U(A)U(A^{-1})$. Pettis' lemma then tells us that $U(A)U(A)^{-1}$ is a neighborhood of $1$ contained in $A A^{-1}$.
\end{proof}

\begin{thm}[Banach]\label{thm:Banach}
Let $G,H$ be Polish groups and $\varphi \colon G \to H$ be a Baire measurable group homomorphism. Then $\varphi$ is continuous.
\end{thm}

Note that this implies that every Borel group homomorphism between Polish groups is continuous, since every Borel map is Baire measurable.

\begin{proof}
It is enough to prove that $\varphi$ is continuous at $1_G$. Let $V$ be an open neighborhood of $1_H$, and fix an open neighborhood $W$ of $1_H$ such that $W W^{-1} \subseteq V$ (continuity of the group operations, as well as the fact that $1_H 1_H=1_H$ give us the existence of $W$). 

Let $W' = \varphi(G) \cap W$, which is open in $\varphi(G)$.
We have $\varphi(G) = \bigcup_{h \in \varphi(G)} h W'$ whence, since the topology of $\varphi(G)$ admits a countable basis, there exists a sequence $(h_n)_{n< \omega}$ of elements of $\varphi(G)$ such that $\varphi(G)= \bigcup_{n < \omega} h_n W'$ (see if necessary the exercise right after this proof). 

Fix $(g_n) \in G$ such that $\varphi(g_n)=h_n$ for all $n$; we then have 
\[ G= \bigcup_{n < \omega} g_n \varphi^{-1}(W).\]
Since $\varphi$ is Baire measurable, $\varphi^{-1}(W)$ is Baire measurable; and it is not meager since countably many translates of it cover $G$. By Pettis' lemma, there is an open neighborhood $O$ of $1_G$ which is contained in $\varphi^{-1}(W) (\varphi^{-1}(W))^{-1}$. But then $\varphi(O)$ is contained in $W W^{-1} \subseteq V$. This proves that $\varphi$ is continuous at $1_G$, hence continuous outright.
\end{proof}

\begin{exo}\label{exo5}
Let $X$ be a topological space whose topology admits a countable basis, and let $(O_i)_{i \in I}$ be a family of open subsets of $X$. 
Prove that there is a (at most) countable subset $J$ of $I$ such that $\bigcup_{i \in I} O_i = \bigcup_{i \in J} O_i$.

A space $X$ such that every open cover admits a countable subcover is said to have the \emph{Lindelöff property}\index{Lindelöff property}.
\end{exo}

\begin{exo}\label{exo6}
Assume that $\varphi \colon \R \to \R$ is Baire measurable and satisfies $\varphi(x+y)=\varphi(x)+ \varphi(y)$ for all $x,y \in \R$. Prove that there exists $\alpha \in \R$ such that $\varphi(x)= \alpha x$ for all $x \in \R$.
\end{exo}

\begin{exo}\label{exo7}
\begin{enumerate}
\item Let $G$ be a Polish group and $H$ a subgroup of $G$ which is Baire measurable and non meager. Prove that $H$ is clopen in $G$ (i.e., both open and closed).
\item Let $G,H$ be two Polish groups and $\varphi \colon G \to H$ a Borel map which is an isomorphism of abstract  groups. Prove that $\varphi$ is a topological group isomorphism (i.e., prove that $\varphi$ is a homeomorphism).
\item Let $G$ be a group, and $\tau_1$, $\tau_2$ two Polish group topologies on $G$ such that $\tau_1 \subseteq \tau_2$. Prove that $\tau_1=\tau_2$.
\end{enumerate}
\end{exo}

To solve the second question of the above exercise, it is useful to know that the inverse of a Borel bijection between two Polish spaces is Borel (by the Lusin--Suslin theorem, see Theorem \ref{thm:Lusin--Suslin}).

\begin{exo}\label{exo7bis}
Let $\tau$ be the usual topology on $\Sinf$, and $\tau'$ be another topology such that $(G,\tau')$ is a Polish group.

\begin{enumerate}
\item Let $A \subset \omega$ be finite, and denote $U_A= \lset \sigma : \forall a \in A \  \sigma(a)=a \rset$.

Prove that $U_A$ is $\tau'$-Borel (first show that if $|A| \ge 3$ then $U_A$ is $\tau'$-closed)

\item Show that $\tau'= \tau$.
\end{enumerate}
\end{exo}

\begin{exo}\label{exo8}
Let $G,H$ be two Polish groups, and $\varphi \colon G \to H$ a continuous surjective homomorphism. Prove that $\varphi$ is an open map.
\end{exo}

Let us now come back to the question of existence of metrics inducing the topology of a given Polish group with additional properties such as left-invariance (or even invariance both on the left and on the right) or completeness.

\begin{thm}[Birkhoff--Kakutani]\index{Birkhoff--Kakutani theorem}
Let $G$ be a topological group with a countable basis of neighborhoods of $1$. Then there exists a left-invariant metric $d$ inducing the topology of $G$.
\end{thm}

Since uniform structures provide a natural framework to prove this result (in a slightly more general version), we postpone the proof of this fact to chapter \ref{uniform_structures_on_groups} where we discuss uniform structures on topological groups.

It follows in particular that any Polish group admits a left-invariant metric; however, it is often the case that there is no metric which is both complete and left-invariant. The group $\Sinf$ is an example, as follows from our earlier discussion combined with the following fact.

\begin{prop}
Let $G$ be a metrizable topological group, and $(g_n)_{n < \omega}$ be a sequence of elements of $G$. Then $(g_n)_{n < \omega}$ is a Cauchy sequence for \emph{some} compatible left-invariant metric on $G$ iff $(g_n)_{n< \omega}$ is a Cauchy sequence for \emph{any} compatible left-invariant metric on $G$.

It follows that $G$ admits a compatible left-invariant complete metric iff any left-invariant compatible metric on $G$ is complete.
\end{prop}

\begin{proof}
Fix some compatible left-invariant metric $d$ on $G$; a sequence $(g_n)_{n< \omega}$ is Cauchy for $d$ iff 
\[\forall \varepsilon > 0 \ \exists N < \omega \ \forall n,m \ge N \quad d(g_n,g_m) < \varepsilon .\]
Using left-invariance, this is equivalent to 
\[\forall \varepsilon > 0 \ \exists N < \omega \ \forall n,m \ge N \quad d(g_m^{-1}g_n,1) < \varepsilon .\]
Let $(U_i)_{i < \omega}$ be a basis of neighborhoods of $1$; the above property is equivalent to the following statement:
\[\forall i < \omega \ \exists N < \omega \ \forall n,m \ge N \quad g_m^{-1}g_n \in U_i .\]
The previous line only depends on the topology of $G$, and not on the choice of left-invariant metric inducing it, which proves our claim.
\end{proof}

So any two left-invariant metrics on a Polish group $G$ have the same Cauchy sequences (later on, we will prove that any two left-invariant metrics induce the same uniform structure on $G$); this suggests considering the completion of $(G,d)$ for some left-invariant metric $d$ on $G$. The group product extends to a continuous operation (this follows from the argument in the next proof), turning this object into a semigroup (well-worth studying in many cases!), however the inverse operation typically does not extend. It it interesting already to understand what happens in the case of $\Sinf$ and give an explicit description of this left completion of $\Sinf$ (see exercise \ref{ex:left_completion} in the next chapter for a more general statement).

\begin{thm}\label{t:extending_group_operations}
Let $G$ be a metrizable topological group, and $d$ a left-invariant metric inducing the topology of $G$. Define a new compatible metric $\rho$ on $G$ by setting $\rho(g,h)= d(g,h)+ d(g^{-1},h^{-1})$. 
Denote by $(\hat{G},\rho)$ the completion of $(G,\rho)$.

The group operations on $G$ extend to continuous operations on $\hat G$, and $(\hat{G},\rho)$ is a topological group.
\end{thm}

\begin{proof}
The map $g \mapsto g^{-1}$ is a distance-preserving map from $(G,\rho)$ to itself, thus it extends to a distance-preserving map from $(\hat{G},\rho)$ to itself. Its image is dense since it contains $G$, and is closed since $\hat{G}$ is complete, so the inverse map extends to an isometry of $\hat{G}$. 

Let us prove that $(g,h) \mapsto gh$ also extends continuously; for this it is enough to prove that if $(g_n)_{n< \omega}$ and $(h_n)_{n<\omega}$ are Cauchy sequences in $(G,\rho)$ then $(g_n h_n)_{n< \omega}$ is also Cauchy in $(G,\rho)$, and actually because of the definition of $\rho$ it is enough to prove that $(g_n h_n)_{n < \omega}$ is Cauchy in $(G,d)$ (since one can then apply this fact to $h_n^{-1} g_n^{-1}$).

We prove the slightly stronger fact that $(g_n h_n)_{n< \omega}$ is $d$-Cauchy as soon as $(g_n)_{n< \omega}$ and $(h_n)_{n < \omega}$ are both $d$-Cauchy\footnote{Note however that in general $(g,h) \mapsto gh$ is not $d$-uniformly continuous!}. Note that for all $n,m,p < \omega$ we have
\begin{eqnarray*}
d(g_n h_n, g_m h_m) & \le & d(g_nh_n,g_n h_p)+ d(g_nh_p,g_m h_p) + d(g_mh_p,g_m h_m) \\
					& \le & d(h_n,h_p) + d(g_n h_p, g_m h_p) + d(h_p,h_m).
\end{eqnarray*}

Fix $\varepsilon >0$, then $p < \omega$ such that for any  $n \ge p$ one has $d(h_p,h_n) \le \varepsilon$.

The map $g \mapsto g h_p$ is continuous at $1$, hence there exists $\delta >0$ such that 
\[ d(g,1) \le \delta \Rightarrow d(gh_p,h_p) \le \varepsilon .\]
Since $(g_n)_{n \in \omega}$ is $d$-Cauchy there exists $N \ge p$ such that for any $n,m \ge N$ we have $d(g_n,g_m)=d(g_m^{-1}g_n,1) \le \delta$. 

Then $d(g_n h_p,g_m h_p) = d{(g_m}^{-1}g_n h_p,h_p) \le \varepsilon$, and it follows that for any $n,m \ge N$ we have $d(g_n h_n, g_m h_m) \le 3 \varepsilon$.

So the map $(g,h) \mapsto gh$ extends continuously to $\hat{G} \times \hat{G}$; associativity of this binary operation on $G$, allied with continuity and density, gives us associativity of the operation on $\hat{G}$. The equality $1g=g$ for all $g \in G$ also extends to $\hat{G}$; similarly the equality $g g^{-1}= 1=g^{-1}g$ extends to $\hat{G}$, proving that every element of $\hat{G}$ has an inverse.

Thus $(\hat{G},\rho)$ is a topological group.
\end{proof}

\begin{coro}
Let $G$ be a Polish group, and $d$ be a compatible left-invariant metric on $G$. Then the metric $\rho$ defined by $\rho(g,h)= d(g,h)+ d(g^{-1},h^{-1})$ is a compatible complete metric on $G$.
\end{coro}

\begin{proof}
It is clear that $\rho$ is a compatible metric on $G$ since the inverse map is continuous. But then $G$ is a Polish subgroup of $(\hat{G},\rho)$, hence $G$ is closed in $\hat{G}$, which is only possible if $G= \hat{G}$.
\end{proof}

Note that, by a theorem of Struble, any locally compact Polish group admits a compatible metric which is both left-invariant and proper, i.e., such that closed balls are compact; such a metric is automatically complete. So the lack of compatible complete left-invariant  metric is again a phenomenon which happens only in ``large'' Polish groups (see Exercise \ref{exo:Struble}).

Now we want to understand how to form quotients of Polish groups. For that, we recall the following classical result, variously attributed to Sierpinski or Hausdorff.

\begin{thm}
Let $X$ be a Polish space, $Y$ a metrizable space and $f \colon X \to Y$ a continuous, surjective and open map. Then $Y$ is Polish.
\end{thm}

Several constructions in the rest of this chapter would be more naturally written down using uniform structures; we advise the reader to revisit them once we have covered that topic later on. We still work trough them now using metrics as a way to prepare ourselves for the use of uniformities later.

\begin{defin}\label{def:quotient_metric}
Let $G$ be a Polish group, and let $\rho$ be a right-invariant compatible metric on $G$. Let $H$ be a closed subgroup of $G$. We endow $G \lcoset H$ with the metric
\[d(fH,gH) = \inf_{h_1,h_2 \in H} \rho(fh_1,gh_2) = \inf_{h \in H} \rho(fh,g) \]
and denote by $\pi \colon G \to G/H$ the natural surjection $g \mapsto gH$.
\end{defin}

\begin{exo}\label{exo9}
Show that $d$ is indeed a metric, and that $d$ induces the quotient topology on $G/H$.

(since $\pi$ is clearly continuous, the last statement amounts to claiming that $\pi$ is open).
\end{exo}

Then $(G/H,d)$ is metrizable and is a continuous open image of a Polish space, thus it is Polish in its own right.

\begin{thm}
Let $G$ be a Polish group, and $H$ a closed normal subgroup of $G$. Then $G/H$, endowed with the quotient topology, is a Polish group.
\end{thm}

\begin{proof}
We have to prove that the group operations are continuous on $G/H$. Let $(f_n)_{n< \omega}$ and $(g_n)_{n< \omega}$ be such that $f_n H \to fH$ and $g_n H \to gH$ in $(G/H,d)$. 

By definition of $d$, we can find sequences $(h_n)_{n < \omega}$ and $(\tilde h_n)_{n < \omega}$ in $H$  such that $f_n h_n$ converges to $f$ and $g_n \tilde h_n$ converges to $g$ in $G$.

Then $(f_n h_n)(g_n \tilde h_n)^{-1}$ converges to $f g^{-1}$; but we have
\[(f_n h_n)(g_n \tilde h_n)^{-1}= f_n h_n \tilde h_n^{-1} g_n^{-1} = f_n g_n ^{-1} (g_n h_n \tilde h_n^{-1} g_n^{-1}). \]
Since $g_n h_n \tilde h_n^{-1} g_n^{-1}$ belongs to $H$ for all $n < \omega$, we conclude that $f_n g_n^{-1} H =(f_n H)(g_n H)^{-1} $ converges to $f g^{-1} H$. 

This establishes continuity of $(fH,gH) \mapsto (fH)(gH)^{-1}$, and we are done.
\end{proof}

\begin{coro}
Let $G,H$ be two Polish groups and $\varphi \colon G \to H$ a surjective homomorphism. Then $\varphi$ induces an isomorphism of topological groups from $G/\ker{\varphi}$ onto $H$.
\end{coro}

\begin{proof}
By definition of the quotient topology, $\varphi$ induces a continuous injective morphism $\tilde \varphi$ from $G/\ker{\varphi}$ to $H$, which is onto since $\varphi$ is onto. Thus $\tilde \varphi$ is an isomorphism of abstract groups between two Polish groups which is continuous, and we saw in an earlier exercise that this implies that $\tilde \varphi$ is an isomorphism of topological groups.
\end{proof}

Of course, now that we have singled out a class of groups we are interested in, we want to make them act on structures; in these notes the main focus will be continuous actions on compact spaces but many other examples are of interest, such as the diagonal conjugation action of $G$ on $G^n$ given by $g \mapsto (g_1,\ldots,g_n)=(gg_1g^{-1},\ldots,gg_ng^{-1})$, unitary representations, measure-preserving actions, actions by permutations on countable sets... 

\begin{exo}\label{exo10}
Let $X$ be a Polish space, and $G$ be a group of homeomorphisms of $X$. Prove that the following conditions are equivalent:
\begin{enumerate}[(i)]
\item There exists $x \in X$ with a dense orbit.
\item For any two nonempty open subsets $U$, $V$ of $X$ there exists $g \in G$ such that $gU \cap V \ne \emptyset$.
\end{enumerate}
If these two equivalent conditions are satisfied, one says that $G$ acts \emph{topologically transitively}\index{topologically transitive action}.
\end{exo}

\begin{exo}[The first $0-1$ topological law]\label{exo11}
Let $X$ be a Polish space, and $G$ be a group of homeomorphisms acting on $X$ topologically transitively. Prove that any subset of $X$ which is both Baire measurable and $G$-invariant is either meager or comeager.
\end{exo}

In the previous exercise we did not even require the action to be jointly continuous (only that each $x \mapsto gx$ is continuous) but that is what we will require most of the time.

\begin{defin}
Let $G$ be a topological group acting on a topological space $X$. We say that the action is \emph{continuous}\index{continuous group action} if $(g,x) \mapsto gx$ is continuous.
\end{defin}

Note that, if a Polish group $G$ acts on a Polish space $X$ continuously and topologically transitively, then each $G$-orbit is Baire measurable (it is analytic since it is the image of $G$ under the continuous map $g \mapsto g \cdot x$; actually it is even Borel, see below) thus each orbit is either meager or comeager. Since distinct orbits are disjoint, there can only exist at most one comeager orbit.

We now turn to discussing some related questions; even though we will not need it, we note the following fact.

\begin{thm}[Miller]
Let $X$ be a Polish space, and $G$ a Polish group acting in a Borel way on $X$ (i.e., $(g,x) \mapsto gx$ is Borel). Then for any $x$ in $X$ the stabilizer $G_x$ is a closed subgroup of $G$, whence the orbit $Gx$ is Borel.
\end{thm}

\begin{proof}
Begin by fixing $x \in X$. Once we have proved that $G_x$ is closed, we are done: indeed $G/G_x$ is a Polish space, and the map $g \mapsto gx$ induces an injective Borel map from $G/G_x$ to $X$ whose image is $Gx$, so $Gx$ is Borel as an injective Borel image of a Polish space (by the Lusin--Suslin theorem). 

Without loss of generality, we may assume that $G_x$ is dense in $G$; since $G_x$ is Borel, hence Baire measurable, Pettis' lemma ensures that $G_x=G$ as soon as $G_x$ is non meager. Let us assume for a contradiction that $G_x$ is meager; since $G_x$ is dense in $G$ an application of the $0-1$ topological law tells us that each Baire measurable $A \subseteq G$ such that $AG_x=A$ is either meager or comeager (the action of $G_x$ on $G$ by right translation is topologically transitive).

Let us now fix a countable basis $(U_n)_{n< \omega}$ for the topology of $X$, and consider 
\[A_n = \lset g \in G : gx \in U_n \rset . \]
Each $A_n$ is Borel, and either meager or comeager by the remark above since $A_n G_x = A_n$. Further, for each $g \in G$ we have 
\[g G_x = \lset h : h x = g x \rset =  \bigcap_{\lset n : g \in A_n\rset } A_n .\]
Since $gG_x$ is meager, it follows that for each $g$ there exists some meager $A_n$ containing $g$. But then $G$ is contained in a countable union of meager sets, contradicting the Baire category theorem.
\end{proof}

Let us now give a criterion for the existence of a comeager orbit.

\begin{lem}[Rosendal (\cite{Rosendal11}, Lemma 9)]\label{l:Rosendal}\index{Rosendal's criterion}
Assume that $G$ is a Polish group acting continuously and topologically transitively on a Polish space $X$. Then the following conditions are equivalent:
\begin{enumerate}[(i)]
\item There exists a comeager orbit.
\item For any nonempty open subset $V$ of $G$, the set $\lset x \in X : Vx \text{ is somewhere dense} \rset$ is dense in $X$.
\item For any open $V \ni 1$ in $G$ and any nonempty open subset $U$ of $X$ there exists a nonempty open $U' \subseteq U$ such that for every nonempty open $W_1$, $W_2 \subseteq U'$ one has $VW_1 \cap W_2 \ne \emptyset$.
\end{enumerate}
\end{lem}

\begin{proof}
$(i) \Rightarrow (ii)$. Assume that $Gx$ is comeager, and let $V \subseteq G$ be nonempty open. Then there exists $(g_n)_{n< \omega}$ in $G$ such that $G = \bigcup_n g_n V$, so $Gx = \bigcup_n g_n Vx$, whence $Vx$ is not meager and thus somewhere dense. Since this property holds for every $y \in Gx$ and $Gx$ is dense, we are done.

\medskip
$(ii) \Rightarrow (iii)$. Assume that the second condition above holds, and fix $1 \in V \subseteq G$ open and $U \subseteq X$ nonempty open. Using continuity of the action, we may find an open neighborhood $\tilde V$ of $1$ contained in $V$ and a nonempty $\tilde U \subseteq U$ such that $\tilde V \tilde U \subseteq U$.

\medskip
Pick some $1 \in V_1 \subseteq \tilde V$ symmetric open such that $V_1 V_1 \subseteq \tilde V$. There exists $x \in \tilde U$ such that $\overline{V_1 x}$ has nonempty interior, and $V_1x \subseteq U$ so its closure contains some nonempty open $U' \subseteq U$. Pick $W_1$, $W_2$ nonempty open and contained in $U'$. There exists $g, h \in V_1$ such that $gx \in W_1$ and $h x \in W_2$, whence $hg^{-1} W_1 \cap W_2 \ne \emptyset$. Since $h g^{-1} \in V_1 V_1^{-1} = V_1 V_1 \subseteq \tilde V$, we have proved the second implication.

$(iii) \Rightarrow (i)$. Assume $(i)$ is false, i.e., there is no comeager orbit. Since the action of $G$ is topologically transitive, all orbits are meager. For every $x \in X$, there exists a countable family $(F_n)_{n<\omega}$ of closed subsets with empty interior such that $Gx \subseteq \bigcup_n F_n$. For some $n$ the set $\{g : g x \in F_n\}$ must have nonempty interior in $G$ by Baire's theorem, so there exists some nonempty open subset $O$ of $G$ such that $Ox$ is nowhere dense.

Translating $O$ if necessary, we obtain that for all $x$ there exists an open $V \ni 1$ such that $Vx$ is nowhere dense. We may restrict $V$ to range over some fixed countable basis of neighborhoods $\mathcal V$ of $1$, and we have 
obtained
\[X = \bigcup_{V \in \mathcal V} \lset x : Vx \text{ is nowhere dense}  \rset .\]
Applying Baire's theorem again, there exists $V$ in $\mathcal V$ such that $\lset x : Vx \text{ is nowhere dense} \rset$ is not meager; this set is Borel so it must be comeager in some nonempty open subset $U$ of $X$. 

Assume that the third condition above holds, and apply this to $V$ and $U$ to find a nonempty $U'$ witnessing that. We have that $\lset x \in U' : Vx \cap W \ne \emptyset \rset$ is dense open for any nonempty open $W$ in $U'$, whence for a generic element of $U'$ the closure of $Vx$ contains $U'$. But for a generic element of $U$ (hence also of $U'$) the closure of $Vx$ is nowhere dense. This yields the desired contradiction.
\end{proof}

Rosendal's criterion enables us to detect whether there exists a comeager orbit without knowing a priori which $x$ is such that $Gx$ is comeager. If we want to understand whether a given orbit $Gx$ is comeager, or use to our advantage the fact that it is comeager, then the next theorem is very useful.

\begin{thm}[Effros]\index{Effros theorem}
Let $G$ be a Polish group acting continuously on a Polish space $X$. For every $x \in X$ the following conditions are equivalent:
\begin{enumerate}[(i)]
\item The map $g \mapsto gx$, from $G$ to $Gx$, is open.
\item $Gx$ is a $G_\delta$ subset of $X$.
\item $Gx$ is non meager in its relative topology.
\item $Gx$ is comeager in $\overline{Gx}$.
\end{enumerate}
\end{thm}

\begin{proof}
$(i) \Rightarrow (ii)$. If $g \mapsto gx$ is open then $Gx$ is a metrizable, continuous open image of a Polish space hence it is Polish, thus $G_\delta$ in $X$. 

\medskip $(ii) \Rightarrow (iii)$. Similarly, if $Gx$ is $G_\delta$ then it is itself Polish, thus non meager in its relative topology.

\medskip $(iii) \Rightarrow (iv)$.
If $Gx$ is non meager in its relative topology then it is non meager in $\overline{Gx}$, so it is comeager since the left translation action of $G$ on $\overline{Gx}$ is topologically transitive.

\medskip $(iv) \Rightarrow (i)$.
Assume that $Gx$ is comeager in $\overline{Gx}=:Y$. 

Denote by $G_x$ the stabilizer of $x$. The orbit map $\varphi \colon G \to Gx$ is continuous, and induces an injective continuous map $\tilde \varphi$ from the Polish space $G/G_x$ to $Gx$. Thus $Gx$ is Borel in $Y$ (we already knew that) and the map $\psi \colon gx \mapsto gG_x$ is Borel. We want to prove that $\psi$ is actually continuous (this immediately implies the desired result since the quotient map from $G$ to $G/G_x$ is open by definition of the quotient topology).

Since $\psi$ is Borel, there is a dense $G_\delta$ subset $\Omega$ of $Y$ such that $\psi$ is continuous on $\Omega \cap Gx$ (extend $\psi$ to be constant on $Y \setminus Gx$, then use the general fact that for any Borel map $f$ between Polish spaces there exists a dense $G_\delta$ set on which the restriction of $f$ is continuous; see Theorem \ref{thm:Baire_measurable_continuous_on_dense_Gdelta} for a more general statement).

Using the notation $\forall^* x \ P(x)$ to mean that $\lset x : P(x) \text{ is true} \rset$ is comeager, we have
\[\forall g \in G \, \forall^* y \in Y \ gy \in \Omega . \]
The Kuratowski--Ulam theorem (which we discuss after this proof) then gives
\[\forall^* y \in Y \, \forall^* g \in G \ gy \in \Omega .\]
In other words, the set $\Sigma = \lset y \in Y : \forall^* g \in G \ gy \in \Omega \rset$ is comeager in $Y$, and this set is $G$-invariant by definition. Since $Gx$ is non meager it must intersect $\Sigma$, thus be contained in $\Sigma$ by $G$-invariance of $\Sigma$.

Now, let $(y_n)_{n< \omega}$ be a sequence of elements of $Y$ which converges to some $y \in Y$. For each $n$ the set $\{g : g y_n \in \Omega\}$ is comeager in $G$, as is $\{g : g y \in \Omega\}$. Hence there exists $g \in G$ such that $gy_n \in \Omega$ for all $n$ and $gy \in \Omega$. Since $\psi$ is continuous on $\Omega \cap Y$, we conclude that $\psi(gy_n)$ converges to $\psi(gy)$ in $G/G_x$. Since $\psi$ is $G$-equivariant we conclude as desired that $\psi(y_n)$ converges to $\psi(y)$.
\end{proof}

To end this chapter, we pursue the analogy between Baire category and measure and prove a useful generalization of the Kuratowski--Ulam theorem. First we discuss a notion which provides an analogue of measure-preserving maps.

\begin{defin}
Let $X$, $Y$ be Polish spaces and $f \colon X \to Y$ a continuous map. We say that $x \in X$ is \emph{locally dense}\index{locally dense point} for $f$ if for every neighborhood $U$ of $x$ the set $\overline{f(U)}$ is a neighborhood of $f(x)$.
\end{defin}

\begin{exo}[King]\label{exo12}
Let $X$, $Y$ be Polish spaces and $f \colon X \to Y$ a continuous map.

\begin{enumerate}
\item For each $r >0$ set $U_r= \lset x \in X : \exists \varepsilon < r \  f(x) \in \mathrm{Int}\overline{f(B(x,\varepsilon))}\rset$
where $B(x,\varepsilon)$ is the open ball of radius $\varepsilon >0$ for some (fixed) compatible metric on $X$.

Show that each $U_r$ is open.
\item Prove that $\bigcap_{r>0} U_r$ coincides with the set of points of local density for $f$. Hence this set is $G_\delta$ in $X$.
\end{enumerate}
\end{exo}

\begin{lem}[``Dougherty's lemma'']\index{Dougherty's lemma}
Assume $X,Y$ are Polish spaces, $f \colon X \to Y$ is continuous and the set of points which are locally dense for $f$ is dense in $X$. Then $f(X)$ is not meager.
\end{lem}

\begin{proof}
Assume for a contradiction that points of local density are dense and $f(X) \subseteq \bigcup_n F_n$, with each $F_n$ a closed subset of $Y$ with empty interior. Since each $f^{-1}(F_n)$ is closed and those sets cover $X$, the Baire category theorem assures us that for some $n$ there is a nonempty open $U$ contained in $f^{-1}(F_n)$. But then $U$ must contain a point of local density for $f$, so $\overline{f(U)} \subseteq F_n$ is a neighborhood of $f(x)$, contradicting the fact that $F_n$ has empty interior.
\end{proof}

This notion is often used to prove that maps do not have a meager image; of course one does not need points of local density to be dense for the image to be non meager, but the next lemma proves that those points do have to exist.

\begin{lem}\label{image_non_loc_dense_is_meager}
Let $f \colon X \to Y$ be a continuous map between Polish spaces, and let $A \subseteq X$ be the set of points which are not locally dense for $f$. Then $f(A)$ is meager.
\end{lem}

\begin{proof}
Fix a basis $(U_n)_{n< \omega}$ for the topology of $Y$. For every $x \in A$ there exists $n$ such that $x \in U_n$ and $\overline{f(U_n)}$ is not a neighborhood of $f(x)$. Hence $f(x) \in \overline{f(U_n)} \setminus \mathrm{Int}(\overline{f(U_n)})$. Each $F_n=\overline{f(U_n)} \setminus \mathrm{Int}(\overline{f(U_n)})$ is closed, has empty interior, and we have shown that $f(A)$ is contained in $\bigcup_n F_n$.
\end{proof}

\begin{exo}\label{exo13}
Let $X,Y$ be two Polish spaces, and $f \colon X \to Y$ be continuous. Show that the following conditions are equivalent:
\begin{enumerate}[(i)]
\item\label{exo13_i} For every meager $A \subseteq Y$, $f^{-1}(A)$ is meager.
\item\label{exo13_ii} For every comeager $A \subseteq Y$, $f^{-1}(A)$ is comeager.
\item\label{exo13_iii} For every dense open $A \subseteq Y$, $f^{-1}(A)$ is dense.
\item\label{exo13_iv} For every nonempty open $U \subseteq X$, $f(U)$ is not meager.
\item\label{exo13_v} For every nonempty open $U \subseteq X$, $f(U)$ is somewhere dense.
\end{enumerate}
\end{exo}

As an example, every continuous open map satisfies the previous conditions.

\begin{defin}
Let $X, Y$ be Polish spaces and $f \colon X \to Y$ be a continuous map. We say that $f$ is \emph{category-preserving}\index{category-preserving map} if $f$ satisfies one of the equivalent conditions appearing in the statement of the previous exercise.
\end{defin}

\begin{prop}
Assume that $X,Y$ are Polish spaces and $f \colon X \to Y$ is continuous. Then $f$ is category-preserving iff the set of points which are locally dense for $f$ is dense in $X$.
\end{prop}

\begin{proof}
Assume that points of local density are dense. Then any nonempty open $U$ contains a point of local density, so $f(U)$ is somewhere dense. Hence $f$ is category-preserving.

Conversely, assume that $U \subseteq X$ is nonempty, open, and does not contain any point of local density for $f$. Then by Lemma \ref{image_non_loc_dense_is_meager} $f(U)$ is meager, so $f$ is not category-preserving.
\end{proof}

\begin{lem}
Let $X$, $Y$ be Polish spaces and $f \colon X \to Y$ be a continuous map. Then there exists a dense $G_\delta$ subset $A$ of $Y$ such that $f \colon f^{-1}(A) \to A$ is open.
\end{lem}
Of course it could happen that $f^{-1}(A)$ is empty.

\begin{proof}
Fix a countable basis $(U_n)_{n< \omega}$ for the topology of $X$. Each $f(U_n)$ is analytic, so we may find an open $O_n$ and a meager $M_n$ in $Y$ such that $f(U_n)=O_n \Delta M_n$.

Let $B= Y \setminus \bigcup_n M_n$. Since $B$ is comeager, it contains a dense $G_\delta$ subset $A$. By construction, for every $n$ we have  $f(U_n \cap f^{-1}(A)) = O_n \cap A$, so each $f(U_n \cap f^{-1}(A))$ is open in $A$. This implies that for every open $U \subseteq X$ the set $f(U \cap f^{-1}(A))$ is open in $A$.
\end{proof}

We conclude this chapter by proving a generalization of the Kuratowski--Ulam theorem, which corresponds to the case where $f$ below is the projection map from $X=X_1 \times X_2$ to $Y=X_1$; note that projection maps are continuous and open by definition of the product topology, hence category-preserving. 

The Kuratowski--Ulam theorem\index{Kuratowski--Ulam theorem} is the statement that, for a Baire measurable $\Omega \subseteq X_1 \times X_2$ one has 
\[(\forall^* (x_1, x_2) \ (x_1,x_2) \in  \Omega ) \Leftrightarrow ( \forall^* x_1 \, \forall^* x_2 \ (x_1,x_2) \in \Omega).\]
This is the analogue, for Baire category, of the Fubini theorem.

\begin{thm}\label{t:splitting_category}
Let $X,Y$ be Polish spaces and $f \colon X \to Y$ be continuous and category-preserving. Assume that $\Omega \subseteq X$ is Baire measurable. The following assertions are equivalent:
\begin{enumerate}[(i)]
\item\label{cat_pres_1} $\Omega$ is comeager in $X$.
\item\label{cat_pres_2} $\displaystyle \lset y : \Omega \cap f^{-1}(\lset y\rset) \textrm{\emph{ is comeager in} } f^{-1}(\lset y \rset ) \rset$ is comeager in $Y$.
\end{enumerate}
Using category quantifiers:
\[(\forall^* x \in X \ \Omega(x)) \Leftrightarrow (\forall^*y \in Y \ \forall^* x \in f^{-1}(\lset y \rset )\  \Omega(x)  ) .\] 
\end{thm}

\begin{proof}
We give the proof for $f$ open, and leave the general case as an exercise. 

We begin by proving \ref{cat_pres_1} $\Rightarrow$ \ref{cat_pres_2}. By Baire's theorem, it is enough to prove that implication for $\Omega$ dense and open in $X$.

So we assume that $\Omega$ is dense open in $X$. Since each $\Omega \cap f^{-1}(\lset y \rset)$ is then open in $f^{-1}(\lset y \rset)$, it is enough to prove that 
\[\forall^* y \in Y \  \Omega \cap f^{-1}(\lset y\rset) \textrm{ is dense in } f^{-1}(\lset y \rset ) .\]
Fix a countable basis $(U_n)_{n< \omega}$ for the topology of $X$. The previous condition is equivalent to 
\[\forall^*y \in Y  \ \forall n \ (f^{-1}(\lset y \rset )\cap U_n \ne \emptyset) \Rightarrow (\Omega \cap f^{-1}(\lset y \rset )\cap U_n \ne \emptyset ).\]
Applying the Baire category theorem, this in turn amounts to 
\[ \forall n \ \forall^*y \in Y \ (f^{-1}(\lset y \rset )\cap U_n \ne \emptyset) \Rightarrow (\Omega \cap f^{-1}(\lset y \rset )\cap U_n \ne \emptyset ) .\]

Fix $n< \omega$, and denote 
\[A_n = \lset y : (f^{-1}(\lset y \rset )\cap U_n \ne \emptyset) \Rightarrow (\Omega \cap f^{-1}(\lset y \rset )\cap U_n) \ne \emptyset \rset  .\]
Since $f$ is open, $A_n$ is the union of an open set and a closed set, so it is $G_\delta$ (open subsets as well as closed subsets are $G_\delta$, and the union of two $G_\delta$ subsets is $G_\delta$). Our aim is to prove that it is comeager, so our job amounts to proving that each $A_n$ is dense. 

Pick a nonempty open $V$ in $Y$. If $V \not \subseteq f(U_n)$ then it is immediate that $V$ meets $A_n$ (any element not in $f(U_n)$ belongs to $A_n$) so we may assume $V \subseteq f(U_n)$.

Then $f^{-1}(V) \cap U_n$ is non-empty open, so $f^{-1}(V)\cap U_n \cap \Omega$ is nonempty. Pick some element $x \in f^{-1}(V)\cap U_n \cap \Omega$ and let $y=f(x)$. Then $y \in A_n \cap V$, proving that $A_n$ is dense.

\medskip To prove \ref{cat_pres_2} $\Rightarrow$ \ref{cat_pres_1}, assume for a contradiction that $\Omega$ is Baire measurable, satisfies the condition of \ref{cat_pres_2} but is not comeager. Then there exists a nonempty open $O$ in $X$ such that $\Omega \cap O$ is meager. Applying \ref{cat_pres_1} $\Rightarrow$ \ref{cat_pres_2} in the Polish space $O$, we obtain 
\[\forall^* y \in f(O) \quad \Omega \cap O \cap f^{-1}(\lset y \rset) \textrm{ is meager in } f^{-1}(\lset y \rset) .\] 
But \ref{cat_pres_2} and the fact that $f(O)$ is open imply that 
\[\forall^*y \in f(O) \quad  \Omega \cap f^{-1}(\lset y \rset) \textrm{ is comeager in } f^{-1}(\lset y \rset).\]
So there exists $y \in f(O)$ such that $\Omega \cap f^{-1}(\lset y \rset)$ is both meager and comeager in $O \cap f^{-1}(\lset y \rset)$, a contradiction.
\end{proof}

\begin{exo}\label{exo14}
Complete the proof of Theorem \ref{t:splitting_category}.

(Hint: combine the result for open $f$ with the existence of a dense $G_\delta$ subset $A$ of $Y$ such that $f \colon f^{-1}(A) \to A$ is open)
\end{exo}

\emph{Comments.} The terminology ``category-preserving map'' comes from \cite{MellerayTsankov}; these maps are further studied in \cite{Melleray14}. While the authors were not aware of the fact, this notion had already been considered long before. The earliest reference I know of is \cite{Herrlich}, where these maps are called ``demi-open'' (see the remark below Proposition 4 there); this paper is contemporary with the monograph \cite{Mioduszewski} where the term ``skeletal'' is used. One can also find the terminology ``pseudo-open'', ``almost open'' or ``weakly open'' in more recent works, though these terms are also used for other, related notions.

The fact that points of local density of a continuous function always form a $G_\delta$ subset comes from \cite{King2000}. 

Following work of Hodges, Hodkinson, Lascar and Shelah \cite{HHLK}, and later of Kechris and Rosendal \cite{KechrisRosendal}, a lot of attention has been devoted to the problem of determining whether some Polish groups of interest admit dense, or comeager, conjugacy classes (or, more generally, whether certain Polish group actions have comeager conjugacy classes). We refer the reader to \Cite{GlasnerWeiss2008} for more on this topic.

\chapter{Fraïssé limits and their automorphism groups}\label{chapter:Fraïssé}
Automorphism groups of countable first-order structures form a rich class of Polish groups to study, which is our main focus in these notes.
We begin by briefly reviewing some model-theoretic vocabulary.

\begin{defin}
A \emph{language}\index{language} consists of the following data:
\begin{itemize}
\item A set of constant symbols $(c_i)_{i \in I}$.
\item A set of relational symbols $(R_j)_{j \in J}$ of arity $n_j \in \omega \setminus \{0\}$.
\item A set of function symbols $(f_k)_{k \in K}$ of arity $m_k \in \omega \setminus \{0\}$.
\end{itemize}
The language is said to be \emph{countable} if $I$, $J$ and $K$ are each (at most) countable.
\end{defin}

\begin{defin}
Given some language $\mcL$ as above, a $\mcL$-structure\index{structure} $\bM$ is a set $M$ with:
\begin{itemize}
\item For each $i \in I$ some element $c_i^\bM$ of $M$.
\item For each $j \in J$ some subset $R_j^\bM$ of $M^{n_j}$.
\item For each $k \in K$ some function $f_k^\bM$ from $M^{m_k}$ to $M$.
\end{itemize}
\end{defin}

Since there should be little risk of confusion we will usually omit the superscript ${}^\bM$ in our notations, and for instance both use $c_i$ for the constant symbol of the language and its interpretation in the structure we are working with.

We always denote by $M$ the underlying set of a structure $\bM$; we say that $\bM$ is \emph{countable} if $M$ is (at most) countable. We will only work with countable languages and structures. We always assume that our languages contain a special binary symbol $=$, which is always interpreted by the equality on $M$; most of the time we will not bother mentioning this symbol but it is always there (for instance in the examples below).

\begin{example}
\begin{itemize}
\item We may view pure sets as structures in the language containing only $=$.
\item Using an additional binary relational symbol, we can consider ordered sets, graphs...
\item Using a binary functional symbol $\cdot$ as well as a constant symbol $e$, we may consider the class of groups (with $e$ being interpreted by the neutral element)
\item The language $(0,1,+,\cdot)$ is well-suited to study rings and fields (and for fields one might be tempted to add a unary functional symbol for multiplicative inverse).
\item The language $(0,1,\wedge, \vee, {}^c)$ is used to study Boolean algebras.
\end{itemize}
\end{example}

We have a natural notion of \emph{substructure}\index{substructure} of a structure $\bM$: a subset $N$ containing each $c_i$ and such that for every $k \in K$ and every $\bar x=(x_1,\ldots,x_{m_k})  \in N^{m_k}$ we have $f_k(x_1,\ldots,x_{m_k}) \in N$. Then one turns $N$ into an $\mcL$-structure $\bN$ by restricting the relations and functions of $\mcL$ to $N$.

Note that the language employed has an influence on the notion of substructure (with our choices of language above to talk about groups, a substructure of a group is not necessarily a subgroup; to address this we might want to add a symbol for the inverse map in our language)

Given an $\mcL$-structure $\bM$ and some subset $A$ of $M$, the structure $\langle A \rangle$ \emph{generated by} $A$ is the smallest substructure of $\bM$ containing $A$. 

We are not going to do any model theory with our structures; our concern will be with their automorphism groups (we should note here that sometimes different structures may induce the same automorphism group).

Given $\varphi \colon M \to N$ and $\bar x= (x_1,\ldots,x_p) \in M^p$ we denote $\varphi(\bar x)= (\varphi(x_1),\ldots,\varphi(x_p))$.

\begin{defin}
Let $\bM$, $\bN$ be two $\mcL$-structures. A map $\varphi \colon M \to N$ is a \emph{homomorphism}\index{homomorphism between first-order structures} if:
\begin{itemize}
\item $\forall i \in I \ f(c_i^\bM)= c_i^\bN$.
\item $\forall j \in J \ \forall \bar x \in M^{n_j} \quad \bar x \in R_j^\bM \Rightarrow \varphi(\bar x) \in R_j^\bN$ (we say that $\varphi$ is an \emph{embedding}\index{embeddings of first order structures} if this implication is an equivalence; embeddings are injective since $x=y \Leftrightarrow \varphi(x)=\varphi(y)$).
\item $\forall k \in K \ \forall \bar x \in M^{m_k} \quad f^\bN(\varphi(\bar x))= \varphi(f^\bM(\bar x))$.
\end{itemize}
\end{defin}

\begin{defin}
Let $\bM$, $\bN$ be $\mcL$-structures. An \emph{isomorphism}\index{isomorphism of first-order structures} is a surjective embedding $f \colon \bM \to \bN$ (and then $f^{-1}$ is an isomorphism from $\bN$ to $\bM$).

An \emph{automorphism}\index{automorphism} of $\bM$ is an isomorphism from $\bM$ to $\bM$.
\end{defin}
The automorphisms of $\bM$ form a group which we denote $\Aut(\bM)$. 

The universe of $\bM$ will almost always be countably infinite, and then we assume it is $\omega$ and view $\Aut(\bM)$ as a subgroup of $\Sinf$.

\begin{prop}\label{p:auts_countable_structure}
Let $\bM$ be a $\mcL$-structure with universe $\omega$. Then $\Aut(\bM)$ is a closed subgroup of $\Sinf$.
Conversely, every closed subgroup of $\Sinf$ is of the form $\Aut(\bM)$ for some structure with universe $\omega$ in some countable relational language.
\end{prop}

\begin{proof}
Assume that $g \in \Sinf \setminus \Aut(\bM)$. Then:
\begin{itemize}
\item Either $g(c_i) \neq c_i$ for some $i \in I$, and then any $h$ such that $h(c_i)= g(c_i)$ does not belong to $\Aut(\bM)$;
\item Or for some $j \in J$ and some $\bar x \in \omega^{n_j}$ we have $\left( \bar x \in R_j \right) \not \Leftrightarrow (g(\bar x) \in R_j)$, and then any $h$ such that $h(\bar x)= g(\bar x)$ does not belong to $\Aut(\bM)$;
\item Or for some $k \in K$ and some $\bar x \in \omega^{m_k}$ we have $f_k(g(\bar x)) \ne g(f_k(\bar x))$ and then any $h$ such that $h(\bar x)= g(\bar x)$ does not belong to $\Aut(\bM)$.
\end{itemize}
In each of the cases above, we found an open subset of $\Sinf$ which contains $g$ and has an empty intersection with $\Aut(\bM)$. This proves that $\Aut(\bM)$ is closed.

Next, let $G$ be a closed subgroup of $\Sinf$. For every $k$, consider the \emph{diagonal action}\index{diagonal action} $G \actson \omega^k$, defined by 
\[g \cdot (n_1,\ldots,n_k)= (g(n_1),\ldots,g(n_k)) , \]
and let $J_k$ denote the set of orbits for this action. We may then form a countable relational language, with a $k$-ary relational symbol $R_O$ for each $O \in J_k$, and consider the $\mcL$-structure$\bM$ with universe $\omega$ where each $R_O$ is interpreted by $O$ (namely $\bar x \in R_O^{\bM} \Leftrightarrow \bar x \in O$).

By definition, we have $G \le \Aut(\bM)$. To see the converse, pick some $h \in \Aut(\bM)$ and some $\bar x \in \omega^k$. Then $h(\bar x) \in G \bar x$ (since the orbit of $x$ is named in the structure $\bM$), i.e., there exists $g \in G$ such that $g(\bar x)= h(\bar x)$. This proves that $h \in \overline{G}$, which is equal to $G$ since $G$ is a Polish (hence closed) subgroup of $\Sinf$.
\end{proof}

A particularly interesting case, which will come up later in these notes, is the case where each action $\Aut(\bM) \actson \omega^k$ only has finitely many orbits (so in the language above there are only finitely many relational symbols of any fixed arity). These groups are called \emph{oligomorphic}\index{oligomorphic group} and correspond to automorphism groups of $\aleph_0$-categorical first-order structures ((i.e, structures whose isomorphism type is completely determined by their theory; we refer to \cite{TentZiegler} for model-theoretic definitions and theorems, which will not be needed in these notes).

\begin{defin}
A $\mcL$-structure $\bM$ is \emph{ultrahomogeneous}\index{ultrahomogeneous structure} if it satisfies the following condition: for any two finitely generated substructures $\bN_1$, $\bN_2$ of $\bM$, any isomorphism $\varphi \colon \bN_1 \to \bN_2$ extends to an automorphism of $\bM$.
\end{defin}

This is a very strong condition, often only considered in the case of relational structures where ``finitely generated'' above is equivalent to ``finite''.

\begin{exo}\label{exo15}
Prove that:
\begin{enumerate}
\item Any pure set is ultrahomogeneous.
\item The set of all rational numbers, seen as an ordered set with its usual order, is ultrahomogeneous.
\item Any vector space over an at most countable field is ultrahomogeneous.
\item An infinite countable, atomless Boolean algebra is ultrahomogeneous.
\end{enumerate}
(In each of the above cases, begin by describing precisely the language you are using)

For the second and the fourth cases above, use a \emph{back-and-forth} argument, i.e., build the desired automorphism as a union of isomorphisms $f_n \colon A_n \to B_n$ between finite substructures, using the even steps to ensure that $\bigcup A_n$ exhausts the structure, and the odd steps to ensure that $\bigcup B_n$ is also exhaustive.
\end{exo}

\begin{exo}\label{exo16}
Prove that every closed subgroup of $\Sinf$ is the automorphism group of some ultrahomogeneous structure on $\omega$ in a countable relational language.

(It is enough to prove that the structure built in the proof of \ref{p:auts_countable_structure} is ultrahomogeneous!)
\end{exo}

\begin{exo}\label{ex:left_completion}
Let $\bM$ be a countable ultrahomogeneous structure, and $G=\Aut(\bM)$. Prove that the left completion of $G$ is naturally identified with the set of embeddings of $\bM$ into itself.
\end{exo}

For instance, it follows from the result of the previous exercise that the left completion of $\Sinf$ is the set of all injective maps from $\omega$ to $\omega$.

\begin{defin}
Let $\bM$ be a $\mcL$-structure. The \emph{age} of $\bM$\index{age of a structure}, denoted by $\age(\bM)$, is the class of all finitely generated $\mcL$-structures which are isomorphic to a substructure of $\bM$.
\end{defin}

First, some easy observations: fix a $\mcL$-structure $\bM$ and let $\mcK= \age(\bM)$.
Then:
\begin{itemize}
\item For any finitely generated $\mcL$-structures $\bA$, $\bB$, if $\bA$ embeds into $\bB$ and $\bB \in \mcK$ then $\bA$ also belongs to $\mcK$. We say that $\mcK$ is \emph{hereditary}\index{hereditary class of structures}.
\end{itemize}

Here, one needs to pay attention to the fact that a substructure of a finitely generated structure need not be finitely generated; for instance the free group on two generators admits a subgroup which is free on countably many generators. Of course this phenomenon cannot occur with relational structures. This is why we specified that $\bA$ is finitely generated.
\begin{itemize}
\item For any $\bA$, $\bB \in \mcK$ there exists $\bC \in \mcK$ such that both $\bA$ and $\bB$ embed in $\bC$. We say that $\mcK$ has the \emph{joint embedding property}\index{joint embedding property}.
\end{itemize}

\begin{defin} Let $\mcK$ be a class of finitely generated $\mcL$-structures. We say that $\mcK$ has the \emph{amalgamation property}\index{amalgamation property} if for any $\bA$, $\bB$, $\bC \in \mcK$ and any embeddings $\alpha \colon \bA \to \bB$, $\beta \colon \bA \to \bC$, there exist $\bD \in \mcK$ and embeddings $i \colon \bB \to \bD$ and $j \colon \bC \to \bD$ such that the following diagram commutes:
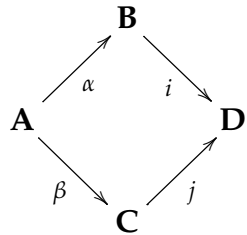
\[
\xymatrix{
& \bB \ar[dr]_i& \\
\bA \ar[ur]_\alpha \ar[dr]_\beta & & \bD \\
& \bC \ar[ru]_j & \\		
}
\]
\end{defin}

\begin{exo}\label{exo18}
Prove that the following classes of structures have the amalgamation property.
\begin{enumerate}
\item The class of all finite graphs.
\item The class of all finite linear orders.
\item The class of all finite groups.
\end{enumerate}
For the first two, the language has one binary relational symbol (besides equality); for the third one there is a binary symbol for the group operation as well as a constant symbol for the neutral element. At this stage in the notes, the third example requires some group-theoretic knowledge (e.g., that an amalgamated free product of finite groups is residually finite), an alternative argument will come up shortly (see Exercise \ref{exo:Hall}).
\end{exo}

\begin{thm}[Fraïssé]\label{t:Fraïssé}
Let $\bM$ be a ultrahomogeneous structure. Then $\age(\bM)$ has the amalgamation property.
\end{thm}

\begin{proof}
Fix $\alpha \colon \bA \to \bB$ and $\beta \colon \bA \to \bC$ as in the definition. We may assume that $\bA, \bB, \bC$ are substructures of $M$. Then $\beta$ is a partial isomorphism of $\bM$ with domain $A$ and image $\beta(A) \subseteq C$. 

By definition of ultrahomogeneity, there exists $g \in \Aut(\bM)$ such that $g(\alpha(a))=\beta(a)$ for all $a \in A$. Let $\bD$ be the substructure of $\bM$ generated by $\bB$ and $g^{-1}(\bC)$. Then $\bD$ is finitely generated; letting $i(b)=b$ for all $b \in B$ and $j(c)=g^{-1}(c)$ for all $c \in \bC$ gives us the desired maps $i,j$.
\end{proof}

\begin{defin}
We say that a structure $\bM$ has the \emph{extension property}\index{extension property} if for any finitely generated substructure $\bA$ of $\bM$ and any embedding $\alpha \colon \bA \to \bB$ with $\bB \in \age(\bM)$ there exists an embedding $\varphi \colon \bB \to \bM$ such that $\varphi( \alpha(a))= a$ for all $a \in A$.
\end{defin}

In words: any abstract extension of a copy of $\bA$ which is contained in $\age(\bM)$ can be realized inside $\bM$ by a substructure which contains $\bA$. 

\begin{exo}\label{exo19}
Prove that ultrahomogeneous structures have the extension property (use the same argument as in the proof of Theorem \ref{t:Fraïssé}).
\end{exo}

\begin{thm}
Let $\bM$, $\bN$ be two countable ultrahomogeneous $\mcL$-structures.

Assume that $\age(\bM)=\age(\bN)$, let $\bA$ be a finitely generated substructure of $\bM$ and $\alpha \colon \bA \to \bN$ an embedding. Then there exists an isomorphism $g \colon \bM \to \bN$ such that $g_{|\bA}=\alpha$.
\end{thm}

In particular, any two ultrahomogeneous structures with the same age are isomorphic.

\begin{proof}
Fix enumerations $(m_k)_{k< \omega}$, $(n_k)_{k< \omega}$ of $M$, $N$ respectively. We claim that, using the extension property, it is possible to build inductively an increasing sequence of finitely generated substructures $\bA_k$, $\bB_k$ of $\bM$, $\bN$ and isomorphisms $\alpha_k \colon \bA_k \to \bB_k$ with the following properties:
\begin{itemize}
\item $\bA_0=\bA$, $\bB_0=\alpha(\bA)$, $\alpha_0=\alpha$. 
\item For each $k$, $\alpha_{k+1}$ extends $\alpha_k$.
\item For all $k$, $m_k \in \bA_{2k+1}$.
\item For all $k$, $n_k \in \bB_{2k+2}$.
\end{itemize}
Assume for now that this is possible. Then $g=\bigcup \alpha_k$ is an isomorphism from $\bigcup_k \bA_k= \bM$ to $\bigcup \bB_k=\bN$ which extends $\alpha$ (the third condition above is there to ensure that the domain of $g$ is $M$, and the fourth one guarantees that the image of $g$ is $N$).

To see why this is indeed possible, assume we have built $\bA_k$, $\bB_k$ and $\alpha_k$ up to some rank $p$. Assume also that $p=2q$ is even (the odd case is similar). If $m_q \in A_p$ then we have nothing to do and simply set $\bA_{p+1}= \bA_p$, $\bB_{p+1}= \bB_p$ and $\alpha_{p+1}=\alpha_p$.

If $m_q \not \in A_p$, then we set $\bA_{p+1}= \langle A_p, m_q \rangle$ and use the extension property of $\bN$ and the fact that $\bA_{p+1} \in \age(\bM)=\age(\bN)$ to find $n \in N$ such that $\alpha_p$ extends to an isomorphism $\alpha_{p+1}$ from $\langle A_p, m_q \rangle$ to $\langle B_p, n \rangle$ such that $\alpha_{p+1}(m_q)=n$. We let $\bB_{p+1}= \langle B_p, n \rangle$ and move on to the next step.
\end{proof}

Given that we only used the extension property of $\bM$, $\bN$ above, the same back-and-forth argument implies that any countable structure with the extension property is ultrahomogeneous (so the extension property is equivalent to ultrahomogeneity).

\begin{defin}
A class $\mcK$ of finitely generated $\mcL$-structures is said to be a \emph{Fraïssé class}\index{Fraïssé class} if it satisfies the following conditions:
\begin{itemize}
\item $\mcK$ contains at most countably many isomorphism types.
\item $\mcK$ has the joint embedding property.
\item $\mcK$ is hereditary.
\item $\mcK$ has the amalgamation property.
\end{itemize}
\end{defin}

Note that here we allow the degenerate case where every structure with age contained in $\mcK$ is finitely generated (in particular, $\mcK$ could be a finite set of relational structures and have a largest element); there is no reason to exclude this possibility it at this point.

\begin{exo}\label{exo20}
Let $\mcL=(R_q)_{q \in \Q^+}$ be the countable language with a binary relational symbol $R_q$ for each positive rational $q$. Any metric space whose metric takes values in $\Q$ can be seen as an $\mcL$-structure, by setting $(x,y) \in R_q \Leftrightarrow d(x,y)=q$. 

Show that the class of finite $\Q$-valued metric spaces, seen as $\mcL$-structures as above, is a Fraïssé class.
\end{exo}

We already saw that the age of any ultrahomogeneous structure is a Fraïssé class. We turn to establishing a converse to that statement.

\begin{thm}[Fraïssé]
Let $\mcK$ be a Fraïssé class in a countable language $\mcL$. There there exists a unique (up to isomorphism) $\mcL$-structure $\bF_{\mcK}$ which is ultrahomogeneous, (at most) countable, and such that $\age(\bF_{\mcK})=\mcK$.

This structure is called the \emph{Fraïssé limit}\index{Fraïssé limit} of $\mcK$.
\end{thm}

Note that uniqueness (up to isomorphism) of a ultrahomogeneous structure with a given age has already been established.

The previous theorem is really only interesting if there are non-finitely generated structures whose age is contained in $\mcK$. 

\begin{proof}
We sketch the argument; see for instance Theorem 4.4.4 of \cite{TentZiegler} for a somewhat more detailed write-up.

We claim that one can build an increasing sequence of structures $\bF_i \in \mcK$ in such a way that
if $\bA \le \bB \in \mcK$ and $f \colon \bA \to \bF_i$ is an embedding, then there exists $j$ and an embedding $g \colon \bB \to \bF_j$ which extends $f$. 

Assume for the moment that such a construction is possible, and let $\bF= \bigcup_i \bF_i$. Then $\age(\bF) \subseteq \mcK$. For any $\bA \in \mcK$ there exists (by the joint embedding property) some $\bB \in \mcK$ such that both $\bF_0$ and $\bA$ embed in $\bB$. Then the identity map from $\bF_0$ to itself extends to an embedding of $\bB$ in some $\bF_j$; hence $\bF_j$ contains a copy of $\bA$ and $\age(\bF)=\mcK$. Then the above condition precisely says that $\bF$ satisfies the extension property, and we are done.

Note that there are (up to isomorphism) only countably many quadruples $(\bA,\bB,\bC,f)$ with $\bA,\bB, \bC \in \mcK$, $\bA \le \bB$ and $f\colon \bA \to \bC$ an embedding. So in order to perform the desired construction we simply need to be able, given such a quadruple, to produce $\bD$ containing $\bC$ such that $f$ extends to an embedding of $\bB$ into $\bD$; and this is precisely what is provided by the following amalgamation diagram,where $i$ is the inclusion map from $A$ to $B$:

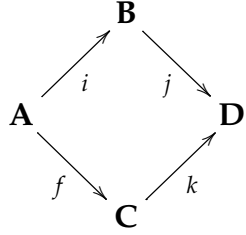
\[
\xymatrix{
& \bB \ar[dr]_j& \\
\bA \ar[ur]_i \ar[dr]_f & & \bD \\
& \bC \ar[ru]_k & \\		
}
\]
Indeed, identifying $\bC$ with $k(\bC)$, $j$ is the desired embedding of $\bB$ into $\bD$.
\end{proof}

\begin{exo}\label{exo21}
\begin{enumerate}
\item Let $\mcK$ be the class of all finite graphs (i.e., sets endowed with a symmetric, irreflexive binary relation) and let $\bR$ be its Fraïssé limit. Show that $\bR$ is (up to isomorphism) the unique countable graph such that for any finite $A,B \subset R$ there exists an element $x \in R$ which is adjacent to every element of $A$ and to no element of $B$.
\item We build a graph with vertex set $\omega$ as follows: for each $i<j \in \omega$ we put an edge between $i$ and $j$ with probability $1/2$ (all those choices being made independently). Prove that almost surely the resulting graph is isomorphic to $\bR$ (this explains why $\bR$ is often called the \emph{random graph}\index{random graph}). 
\end{enumerate}
\end{exo}

One could prove, similarly to what is done in Exercise \ref{exo20}, that the class of finite metric spaces whose metric takes its values in $\{0,\ldots,n\}$ for some $1 \le n \le \omega$ is a Fraïssé class. For $n=1$, the corresponding Fraïssé limit is a pure set; for $n=2$ it is naturally identified with the random graph : say that there is an edge between $x$ and $y$ iff $d(x,y)=1$. When the set of allowed distances is equal to $\Q$, the corresponding Fraïssé limit is called the \emph{rational Urysohn space}.

\begin{exo}\label{exo22}
A \emph{tournament}\index{tournament} is an oriented graph (asymmetric, irreflexive binary relation) with the property that for any pair of distinct points $x,y$ there is an edge from $x$ to $y$ or an edge from $y$ to $x$ (equivalently, a complete graph with an orientation of its edges). Prove that the class of finite tournaments is a Fraïssé class.
\end{exo}

\begin{exo}\label{exo:Hall}
We build an increasing sequence of groups $(G_i)_{i< \omega}$ as follows: let $G_0$ be the permutation group on $3$ elements. Assuming $G_i$ has been built, embed $G_i$ in the permutation group $\mathfrak S(G_i)$ via the left translation action of $G_i$ on itself, and set $G_{i+1}=\mathfrak S(G_i)$. Then define $\bH = \bigcup_i G_i$.

\begin{enumerate}
\item Prove that $\bH$ is locally finite (i.e., any finitely generated subgroup is finite), ultrahomogeneous, and that $\age(\bH)$ is the class of all finite groups. The group $\bH$ is called \emph{Hall's universal locally finite group}\index{Hall's universal locally finite group}.
\item Prove that the class of finite groups satisfies the amalgamation property.
\item Prove that any two elements of $\bH$ with the same order are conjugate in $\bH$.
\item Prove that $\bH$ is simple.
\end{enumerate}
\end{exo}

\begin{exo}\label{exo24}
We already saw that the class of finite linear orders is a Fraïssé class. Show that its limit $\bM=(M,<)$ is the unique countable linear ordering which is 
\begin{itemize}
\item \emph{dense}, i.e., for any $x<y \in M$ there exists $z \in M$ such that $x<z<y$.
\item \emph{without endpoints}, i.e., $M$ has neither a maximum nor a minimum.
\end{itemize}
Conclude that $\Q$, with its usual ordering, is the Fraïssé limit of finite linear orders, and that it is up to isomorphism the unique countable linear order which is dense and without endpoints (Cantor's theorem).
\end{exo}

It is interesting to analyze the interplay between combinatorial properties of $\mcK$ and properties of the automorphism group of its Fraïssé limit. Let us give an example to conclude this chapter.

\begin{defin}
Let $\mcK$ be a Fraïssé class in a relational language $\mcL$. We say that $\mcK$ has the \emph{extension property for partial automorphisms}\index{extension property for partial automorphisms (EPPA)}, abbreviated EPPA, if for any $\bA \in \mcK$ there exists $\bB \in \mcK$ such that $\bA$ embeds in $\bB$ and every partial automorphism of $\bA$ extends to an isomorphism of $\bB$.
\end{defin}

Hrushovski \cite{Hrushovski1992} proved that the class of finite graphs has the EPPA; since then many other examples have been obtained, in particular following work of Herwig--Lascar (for instance Solecki proved that this property holds for the class of finite metric spaces with rational distances). This property is however in general very difficult to establish.

\begin{thm}[Kechris--Rosendal (\cite{KechrisRosendal}, Proposition~6.4)] Let $\mcK$ be a Fraïssé class in a relational language (thus all elements of $\mcK$ are finite), $\bF$ its Fraïssé limit, and $G=\Aut(\bF)$. Then the following are equivalent:
\begin{itemize}
\item $\mcK$ has the extension property for partial automorphisms.
\item There exists an increasing sequence of compact subgroups $G_n \le \Aut(\bF)$ whose union is dense in $\Aut(\bF)$.
\end{itemize}
\end{thm}

As usual, this result is only interesting if $\mcK$ contains structures of arbitrarily large finite cardinality, i.e., if the underlying set of $\bF$ is infinite (otherwise $\bF$ is a finite structure, and both properties above trivially hold).

\begin{proof}
We assume that $\bF$ is infinite (and countable, of course). We assume that $\mcK$ has the EPPA and denote $G=\Aut(\bF)$. Given $n$, consider the set 
\[\Omega_n =\lset(g_1,\ldots,g_n) \in G^n : \forall x \in F \ \langle g_1,\ldots,g_n \rangle x\text{ is finite} \rset .\]
Given some finite $A \subset F$, the set $\lset \bar g \in G^n : \langle \bar g \rangle A =A \rset$ is open, since if $\bar g$ is in that set and $\bar h$ coincides with $\bar g$ on $A$ then $\bar h$ also belongs to that set. So saying that a given $x$ has a finite orbit under $\langle \bar g \rangle$ is an open condition. Thus $\Omega_n$ is a countable intersection of open sets, hence is $G_\delta$. 

The EPPA guarantees that $\Omega_n$ is dense in $G^n$; thus a generic element $\bar g$ of $G^n$ belongs to $\Omega_n$, equivalently the closed subgroup generated by $\bar g$ is compact (recall that a closed subgroup $G \le \Sinf$ is compact iff every point has a finite $G$-orbit).

It follows from this, and an application of the Kuratowski--Ulam theorem, that for a generic sequence $(g_n)_{n< \omega} \in G^\omega$, we have for all $n$ that $\overline{\langle g_1,\ldots,g_n \rangle}$ is compact.

In any Polish group, the set of all $(g_n) \in G^\omega$ which generate a dense subgroup is dense $G_\delta$; applying the Baire category theorem, we conclude that a generic sequence in $G^\omega$ induces an increasing sequence of compact subgroups with a dense union.

\medskip
Conversely, assume that $(G_i)_{i< \omega}$ is an increasing sequence of compact subgroups with a dense union, and let $\bA \in \mcK$. We may view $\bA$ as a (finite) substructure of $\bF$. Let $f_1,\ldots,f_k$ enumerate all partial automorphisms of $\bA$; by ultrahomogeneity of $\bF$ we can extend each $f_j$ to some $\tilde f_j \in \Aut(\bF)$. By density of $\bigcup_i G_i$, for all $j$ there exists some $i_j$ and $g_j \in G_{i_j}$ which coincides with $\tilde f_j$ on $A$, hence coincides with $f_j$ on its domain. Let $n = \max\lset i_j : j \in \lset 1,\ldots k \rset \rset$. Then $g_j \in G_n$ for each $n$; and each $g_j$ induces an automorphism of the finite structure $G_n \bA$ (since $G_n$ is compact, every element of $\bF$ has a finite $G_n$-orbit). So every partial automorphism of $\bA$ extends to an automorphism of $G_n \bA$, and we are done.
\end{proof}

\emph{Comments.} For more on Fraïssé classes the reader can consult courses on model theory such as \cite{TentZiegler} or \cite{Hodges1993a}. The survey \cite{MacPherson} is also a good introductory source for some of the material mentioned here, and much more. The connection between automorphism groups of ultrahomogeneous structures and closed subgroups of $\Sinf$ comes from \cite{Becker1996}, which contains a lot of interesting material on Polish group actions. Regarding the EPPA, much work has built on \cite{Herwig2000}, including very recent papers like \cite{EvansHubickaNesetril} and \cite{HubickaKonecnyNesetril}. The fact that the class of finite rational metric spaces has the EPPA comes from \cite{Solecki2005b}. 
Whether the class of finite tournaments satisfies the EPPA is a notorious (and notoriously difficult) open problem.

\chapter{Uniform spaces}
To study topological groups, even metrizable, it is very useful to understand the notion of a uniform structure, which we study in some detail now. Later we will consider various uniformities on a given Polish group; in particular, the algebra of all continuous, bounded functions for the so-called right uniformity will play a major part. Before delving into that aspect of the theory, we need to set up the basic concepts.

\begin{defin}
Let $X$ be a set; denote $\Delta_X=\lset (x,x) : x \in X \rset$. A \emph{uniform structure}\index{uniform structure}, also called a \emph{uniformity}, on $X$, is a set $\mcU$ of subsets of $X \times X$ such that:
\begin{itemize}
\item For all $U \in \mcU$ one has $\Delta_X \subseteq U$.
\item For all $U \in \mcU$, $U^{-1} \in \mcU$ (where $U^{-1}= \lset (x,y) : (y,x) \in U \rset$).
\item For all $U,V \in \mcU$, $U \cap V \in \mcU$.
\item For all $U$, $V$ $(U \in \mcU \text{ and } U \subseteq V ) \Rightarrow V \in \mcU$.
\item For all $U \in \mcU$ there exists $V \in \mcU$ such that $V \circ V \subseteq U$.

(where $V \circ V= \lset (x,y) : \exists z \ (x,z) \in V \text{ and } (z,y) \in V \rset$)
\end{itemize}
\end{defin}
Elements of a uniform structure $\mcU$ are called \emph{entourages}\index{entourage}; we think of them as being neighborhoods of $\Delta_X$. A fundamental example is given by metric spaces: given a metric space $(X,d)$ one can consider the uniformity $\mcU_d$ whose entourages are the subsets of $X \times X$ containing a set of the form $\lset (x,y) : d(x,y)< r \rset$ for some $r>0$.

\begin{exo}\label{exo25}
Show that $\mcU_d$ is indeed a uniformity. Prove that two metrics $d_1$, $d_2$ on $X$ are uniformly equivalent iff $\mcU_{d_1}= \mcU_{d_2}$.

(Recall that two metrics on $X$ are uniformly equivalent iff the identity map is uniformly continuous in both directions)
\end{exo}

\begin{defin}
Given a uniformity $\mcU$ on a set $X$, we endow $X$ with a topology by stating that $V \subseteq X$ is open iff
\[\forall x \in V \ \exists U \in \mcU \ \lset y : (x,y) \in U \rset \subseteq V .\]
\end{defin}
We will denote $U[x]=\lset y : (x,y) \in U \rset$.

\begin{exo}\label{exo26}
Let $d$ be a metric on $X$, and $\mcU_d$ the metric uniformity. Show that the topology induced by $\mcU_d$ is the same as the topology induced by $d$.
\end{exo}

\begin{exo}\label{exo27}
Let $(X,\mcU)$ be a uniform space, and $x \in X$. Show that the neighborhoods of $x$ for the topology induced by $\mcU$ are exactly the sets of the form $U[x]$ for some $U \in \mcU$.
\end{exo}

\begin{defin}
Let $(X,\mcU)$ and $(Y,\mcV)$ be two uniform spaces. A map $f \colon X \to Y$ is \emph{uniformly continuous} if 
\[\forall V \in \mcV \ \exists U \in \mcU \ \forall x,y \in X \quad (x,y) \in U \Rightarrow (f(x),f(y)) \in V . \]
\end{defin}

We say that a uniform structure is \emph{metrizable}\index{metrizable uniform structure} if there is a metric which induces it. Not all uniform structures are metrizable (though they are always induced by a family of pseudometrics); even in metrizable settings, it is sometimes the case that uniform structures are more natural objects that metrics (i.e., there are some natural choices of uniform structure, but no canonical choice of metric; we already saw an example of this phenomenon when discussing left-invariant metrics on Polish groups).

\begin{defin}
A uniform structure is Hausdorff if the topology it induces is Hausdorff.
\end{defin}

Since our convention is that all the topological spaces we consider are Hausdorff, we need to clarify which uniform structures we are working with. That turns out to be fairly straightforward.

\begin{prop}
$(X,\mcU)$ is Hausdorff iff $\displaystyle \bigcap_{U \in \mcU} U = \Delta_X$.
\end{prop}

\begin{proof}
Assume $(X,\mcU)$ is Hausdorff. Let $x \ne y$; there exists $U \in \mcU$ such that $y \not \in U[x]$, equivalently $(x,y) \not \in U$ so $\bigcap_{U \in \mcU} U \subseteq \Delta_X$. The other inclusion is an immediate consequence of the definition of a uniform structure.

Conversely, assume that $\bigcap_{U \in \mcU} U = \Delta_X$ and let $x \ne y \in X$. There exists $U \in \mcU$ such that $(x,y) \not \in U$, and $V \in \mcU$ such that $V \circ V^{-1} \subseteq U$. Then $V[x]$ and $V[y]$ are neighborhoods of $x,y$ respectively; if $z \in V[x] \cap V[y]$ then $(x,z) \in V$ and $(z,y) \in V^{-1}$, whence $(x,y) \in U$, a contradiction. So $V[x]$ and $V[y]$ are disjoint.
\end{proof}

The above proof shows that a uniform space is Hausdorff as soon as for any two points $x \ne y$ there exists a neighborhood of $x$ which does not contain $y$, a property which is in general weaker than being Hausdorff. From now on we only work with Hausdorff uniform spaces.

If the topology on a topological space $X$ is induced by a uniformity, one says that it is \emph{uniformizable}; not every topology is uniformizable. We mention without proof the fact that a topological space is uniformizable iff it is completely regular, i.e., points can be separated from closed sets by continuous functions.

\begin{exo}\label{exo28}
Let $(X,\mcU)$ be a uniform space. Prove that for any $U \in \mcU$ the interior of $U$ (for the product topology on $X \times X$ induced by $\mcU$) belongs to $\mcU$.
\end{exo}

\begin{defin}
Let $(X,\mcU)$ be a uniform space. A \emph{fundamental system}\index{fundamental system} is a family $\mcE$ of entourages such that any element of $\mcU$ contains an element of $\mcE$.
\end{defin}

Note that every uniform structure admits a fundamental system consisting of open entourages (start from any fundamental system then consider the interiors of its elements).

\begin{thm}\label{t:metrizable_uniformity}
Let $(X,\mcU)$ be a uniform space. Then $\mcU$ is metrizable iff it admits a countable fundamental system.
\end{thm}

\begin{proof}
If $d$ induces $\mcU$ then the family of entourages $\lset (x,y) : d(x,y) < \varepsilon\rset$, where $\varepsilon$ ranges over positive rational numbers, is a countable fundamental system. This proves one implication.

Conversely, assume that $\mcU$ admits a countable fundamental system. One can then produce another countable fundamental system $\lset U_n \rset_{n< \omega}$ with the following properties:
\begin{itemize}
\item $U_0= X \times X$.
\item For all $n$ one has $U_n^{-1}=U_n$.
\item For all $n$ one has $U_{n+1} \circ U_{n+1} \circ U_{n+1} \subseteq U_n$.
\end{itemize}
We may then define, for $x,y \in X$, $\rho(x,y) =\inf \lset 2^{-n} : (x,y) \in U_n \rset$.

Then $\rho$ is symmetric, and since $\mcU$ is Hausdorff we have $\rho(x,y)=0$ iff $x=y$ (since $(U_n)_n$ is a fundamental system we have $\bigcap_n U_n = \Delta_X$). However $\rho$ need not satisfy the triangle inequality. Fortunately, there is a way to produce a metric from a symmetric weight function such as $\rho$: we let 
\[d(x,y)= \inf \lset \sum_{i=0}^n \rho(x_i,x_{i+1}) : x_0=x, \  x_{n+1}=y \rset .\]
Clearly $d$ is now a pseudometric and $d \le \rho$. We claim that $2 d \ge \rho$ (which implies that $d$ is actually a metric); accept this for the moment. 

Choose $\varepsilon$ such that $0< \varepsilon \le \frac{1}{2}$, and let $i \ge 1$ be such that $2^{-i-1} \le \varepsilon \le 2^{-i}$.  

If $(x,y) \in U_{i+1}$ then $d(x,y) \le \rho(x,y) \le 2^{-i-1} \le \varepsilon$. Thus $\lset (x,y) : d(x,y) \le \varepsilon \rset$ contains $U_{i+1}$, hence belongs to $\mcU$. 

Conversely, if $d(x,y) \le \varepsilon$ then $\rho(x,y) \le 2 d(x,y) \le 2^{-i+1}$, whence $(x,y) \in U_{i-1}$ so $U_{i-1}$ belongs to the uniformity generated by $d$. It follows that $d$ induces $\mcU$.

We still have to prove that $2d \ge \rho$. For this, we prove by induction on $n \ge 1$ that for any $x_0,\ldots,x_{n+1} \in X$ one has 
\[ \sum_{i=0}^n \rho(x_i,x_{i+1}) \ge \frac{\rho(x_0,x_{n+1})}{2} .\]
Note that the fact that $U_{n+1} \circ U_{n+1} \circ U_{n+1} \subseteq U_n$ for all $n$ implies that 
\begin{equation}\label{inductive_metrizable_unifrmity}\tag{$*$}
\forall x_1,\, x_2, \, x_3 \ \forall \varepsilon >0 \  (\rho(x_0,x_1) \le \varepsilon \textrm{ and } \rho(x_1,x_2) \le \varepsilon \textrm{ and } \rho (x_2,x_3) \le \varepsilon ) \Rightarrow \rho(x_0,x_3) \le 2 \varepsilon .
\end{equation}
This takes care of the cases $n=0,1,2$ in our induction. Assume we have established our property for all $m \le n-1$ for some integer $n \ge 3$. 

Fix $x_0,\ldots,x_{n+1}$ and let $\displaystyle r= \sum_{i=0}^n \rho(x_i,x_{i+1})$.
If $\displaystyle \sum_{i=0}^{n-1} \rho(x_i,x_{i+1}) \le \frac{r}{2}$ then (by induction) we have $\rho(x_0,x_n) \le r$, whence we conclude by \eqref{inductive_metrizable_unifrmity} that $\rho(x_0,x_{n+1}) \le 2r$. Similarly, we obtain the desired conclusion if $\displaystyle \sum_{i=1}^{n} \rho(x_i,x_{i+1}) \le \frac{r}{2}$.

Otherwise, there exists some $i \in \lset 1,\ldots n-1 \rset$ such that 
\[\sum_{j=0}^{i-1} \rho(x_j,x_{j+1}) \le \frac{r}{2} \text{ and } \sum_{j=i+1}^n \rho(x_j,x_{j+1}) \le \frac{r}{2} .\] (choose the largest $i$ for which $\sum_{j=0}^{i-1} \rho(x_j,x_{j+1}) \le \frac{r}{2}$). 

Again using the inductive assumption, we have $\rho(x_0,x_i) \le r$, $\rho(x_{i+1},x_{n+1}) \le r$ and of course we also have $\rho(x_i,x_{i+1}) \le r$. Hence \eqref{inductive_metrizable_unifrmity} gives $\rho(x_0,x_{n+1}) \le 2 r$.
\end{proof}

\begin{prop}\label{p:description_uniform_structure_compact_space}
Let $X$ be a compact space. Then the set of all neighborhoods of $\Delta_X$ (in $X \times X$ endowed with the product topology) is a uniform structure on $X$.
\end{prop}

\begin{proof}
Reviewing the definition of a uniform structure, we see that the only non immediate fact that we have to prove is that if $U$ is a neighborhood of $\Delta_X$ then there exists a neighborhood $V$ of $\Delta_X$ such that $V \circ V \subseteq U$. 

We will make use of the fact that compact Hausdorff spaces are \emph{completely regular}, i.e., for every $x \in F$ and any closed $F \subseteq X$ such that $x \not \in F$, there exists a continuous function $f \colon X \to \R$ such that $f(x)=0$ and $f$ is constant equal to $1$ on $F$.
We claim that every neighborhood of $\Delta_X$ contains a set of the form
$\lset (x,y) : \forall i \in I \ |f_i(x)-f_i(y)| < \varepsilon_i \rset$ with $f_i \in C(X,\R)$  and $ \varepsilon_i >0 $
for some finite index set $I$ (also, note that each of these sets is an open neighborhood of $\Delta_X$).
Once this is proved, the desired result follows easily. 

We check the claim. Let $A$ be a neighborhood of $\Delta_X$; by compactness, there is an open covering $(O_i)_{1 \le i \le n}$ of $X$ such that $\bigcup_{i=1}^n O_i \times O_i \subseteq A$. For each $x \in X$, find $i$ such that $x \in O_i$ and then a continuous map $f_x \colon X \to \R$ such that $f_x(x)=0$ and $f_x(y)=1$ for every $y \not \in O_i$. Then let $V_x= \lset y : f_x(y) < \frac{1}{2} \rset$.

By compactness again, there exist $x_1,\ldots,x_p$ such that $V_{x_1}, \ldots, V_{x_p}$ cover $X$. For $j \in \{1,\ldots,p\}$ denote $g_j=f_{x_j}$. Assume that $(x,y)$ is such that $|g_j(x)-g_j(y)| < \frac{1}{2}$ for all $j$. Then there exists $j \in \{1,\ldots,p\}$ such that $x \in V_{x_j}$. In turn, there exists $k \in \lset 1,\ldots,n \rset$ such that $x_j \in O_k$; we then have $f_{x_j}(x) < \frac{1}{2}$ and $f_{x_j}(y) < 1$, so $(x,y) \in O_k \times O_k \subseteq A$.
\end{proof}

\begin{thm}\label{t:uniform_continuity_compact}
Let $(X,\mcU)$ be a compact uniform space, $(Y,\mcV)$ be a uniform space, and $f \colon X \to Y$ a continuous map. Then $f$ is uniformly continuous.
\end{thm}

\begin{proof}
Fix $V \in \mcV$, and choose a symmetric $V'$ such that $V' \circ V' \subseteq V$. Continuity of $f$ implies that for every $x \in X$ there exists $U_x \in \mcU$ such that $(x,x') \in U_x \Rightarrow (f(x),f(x')) \in V'$. Find some symmetric $U_x' \in \mcU$ such that $U_x' \circ U_x' \subseteq U_x$.  

Recall that $U_x'[x]$ is a neighborhood of $x$. Hence by compactness there exist $x_1,\ldots,x_n \in X$ such that $X= \bigcup_{i=1}^n U'_{x_i}[x_i]$; let $U = \bigcap_{i=1}^n U_{x_i}'$. By definition of a uniform structure, $U$ belongs to $\mcU$.

Now, assume that $(x,y) \in U$. Then for some $i$ we have $x \in U_{x_i}'[x_i]$. Also $(x,y) \in U'_{x_i}$ since $U \subset U'_{x_i}$. Thus we have that $(x_i,x) \in U_{x_i}'$ and $(x_i,y) \in U'_{x_i} \circ U'_{x_i} \subseteq U_{x_i}$.

It follows that $(f(x_i), f(x)) \in V'$ and $(f(x_i),f(y)) \in V'$, so $(f(x),f(y)) \in V$.
\end{proof}

\begin{exo}\label{exo29}
Prove that every compact space admits a unique compatible uniformity.
\end{exo}

The following fact will be handy when we turn to topological dynamics and need to check uniform continuity of certain maps defined on dense subsets of compact spaces in order to extend them.

\begin{prop}\label{p:carac_unif_continuity_compact}
Let $(X,\mcU)$ be a uniform space, and $Y$ a compact topological space. 

Then $f \colon (X,\mcU) \to Y$ is uniformly continuous if, and only if, $g \circ f$ is uniformly continuous for every continuous function $g \colon Y \to \R$.
\end{prop}

\begin{proof}
Note that we did not bother mentioning which uniform structure we put on $Y$, since there is exactly one. One direction in the above equivalence is obvious, since every continuous function from $Y$ to $\R$ is uniformly continuous and a composition of two uniformly continuous functions is uniformly continuous.

Conversely, assume that $f$ satisfies the above condition, and let $V$ be an entourage for the uniformity on $Y$. We know that $V$ is a neighborhood of the diagonal, and that reducing $V$ if necessary we may assume that 
\[V = \lset (y_1,y_2) \in Y : \forall i \in I \ |g_i(y_1)-g_i(y_2)| < 1 \rset \]
for some finite index set $I$ and some continuous functions $g_i \colon Y \to \R$. By assumption, each $g_i \circ f$ is uniformly continuous on $(X,\mcU)$, whence for each $i \in I$ there exists $U_i \in \mcU$ such that for all $(x_1,x_2) \in U_i$ one has $|g_i \circ f(x_1) - g_i \circ f(x_2)| < 1$.

Then $U= \bigcap_i U_i$ belongs to $\mcU$ and for all $(x_1,x_2) \in \mcU$ we have $(f(x_1), f(x_2)) \in V$.
\end{proof}

We now need to present a notion of completion of a (Hausdorff, as always) uniform space; for this we use filters. We first recall the basic definitions.

\begin{defin}
Let $X$ be a set. A subset $\mcF$ of $\mcP(X)$ is a \emph{filter}\index{filter} if :
\begin{itemize}
\item $X \in\mcF$ and $\emptyset \not \in \mcF$.
\item If $F_1$, $F_2 \in \mcF$ then $F_1 \cap F_2$ belongs to $\mcF$.
\item If $F \in \mcF$ and $F \subseteq A$ then $A \in \mcF$.
\end{itemize}
An \emph{ultrafilter}\index{ultrafilter} is a filter which is maximal (among filters) for inclusion.
\end{defin}

Intuitively, a filter provides a notion of ``large subset'': the whole set is large, the empty set is not large, an intersection of two large sets is still large, and a subset which contains a large subset is large itself.

\begin{exo}\label{exo30}
Let $\mcF$ be a filter on $X$. Then $\mcF$ is an ultrafilter iff for any $A \subseteq X$ one has either $A \in \mcF$ or $X \setminus A \in \mcF$.
\end{exo}

A filter on $X$ is \emph{principal} if there exists some $A \subseteq X$ such that $\mcF= \lset F : A \subseteq F \rset$.

More interesting is the \emph{Fréchet filter}\index{Fréchet filter}, which is the set of all subsets of $X$ with finite complement (on an infinite set $X$). Any ultrafilter is either principal (and contains a singleton) or contains the Fréchet filter. Principal filters and the Fréchet filter may both be seen as particular cases of neighborhood filters, which we define now.

\begin{defin}
Let $X$ be a topological space and $x \in X$. The \emph{neighborhood filter}\index{neighborhood filter} $\mcV_x$ is the family of all neighborhoods of $x$.
\end{defin}

\begin{defin}
Let $X$ be a topological space, and $\mcF$ be a filter on $X$. We say that $\mcF$ \emph{converges to $x$}\index{convergent filter} if $\mcF$ contains the neighborhood filter $\mcV_x$ of $x$.
\end{defin}

\begin{defin}
Let $X$, $Y$ be two sets, $\mcF$ a filter on $X$ and $g \colon X \to Y$. The \emph{image filter} $g(\mcF)$ is $\{A \subseteq Y : g^{-1}(A) \in \mcF\}$.
\end{defin}

\begin{exo}\label{exo31}
Assume that $X$ is a (as always, Hausdorff) topological space and $\mcF$ is a filter on $X$ which converges to both $x$ and $y$. Show that $x=y$.
\end{exo}

\begin{exo}\label{exo32}
Let $X,Y$ be two topological spaces, and $f \colon X \to Y$ be a function. Show that $f$ is continuous at $x \in X$ iff the image of any filter that converges to $x$ is a filter that converges to $f(x)$.
\end{exo}

\begin{exo}\label{exo33}
Let $X$ be a topological space, and $(x_n)_{n < \omega} \in X^\omega$. Let $\varphi \colon \omega \to X$ be defined by $\varphi(n)=x_n$. Let $\mcF$ denote the image under $\varphi$ of the Fréchet filter on $\omega$.

Show that $(x_n)_{n \in \omega}$ converges to $x$ iff $\mcF$ converges to $x$.
\end{exo}

So convergence of filters generalizes convergence of sequences, and enables us to capture the topology of spaces which are not first countable. One also often sees generalized sequences or nets used for the same purpose; filters have the advantage of being very amenable to set-theoretic manipulations, and probably the disadvantage of being less intuitive.

\begin{prop}
Let $X$ be a compact topological space. Then any ultrafilter on $X$ is convergent.
\end{prop}

\begin{proof}
Let $\mcF$ be an ultrafilter on $X$, and set $\mcA=\lset \overline{F} : F \in \mcF \rset$. Then any finite intersection of elements of $\mcA$ is nonempty, hence by compactness there exists $x_0 \in X$ such that $x_0 \in \overline F$ for any $F \in \mcF$.

Let $V$ be a neighborhood of $x_0$; if $V \not \in \mcF$ then $X \setminus V \in \mcF$, and $x_0 \not \in \overline{X \setminus V}$, a contradiction. So $V \in \mcF$, proving that $\mcF$ converges to $x_0$.
\end{proof}

\begin{exo}\label{exo34}
Prove that the previous property characterizes compact spaces, i.e., that a topological space is compact iff any ultrafilter on $X$ is convergent. Use this to prove Tychonoff's theorem: any product of compact spaces is compact (first characterize convergence of a filter on a product space by the convergence of each of its projections).
\end{exo}

\begin{defin}
Let $(X,\mcU)$ be a uniform space, and $\mcF$ be a filter on $X$. We say that $\mcF$ is a \emph{Cauchy filter}\index{Cauchy filter} if :
\[\forall U \in \mcU \ \exists F \in \mcF \quad F \times F \subseteq U . \]
\end{defin}

\begin{exo}\label{exo35}
Let $(X,d)$ be a metric space, and $(x_n)_{n < \omega} \in X^\omega$. Let $\varphi \colon \omega \to X$ be defined by $\varphi(n)=x_n$. Let $\mcF$ denote the image under $\varphi$ of the Fréchet filter on $\omega$.

Show that $(x_n)_{n < \omega}$ is a Cauchy sequence iff $\mcF$ is a Cauchy filter.
\end{exo}
\begin{defin}
A uniform space $(X,\mcU)$ is \emph{complete}\index{complete uniform space} iff any Cauchy filter on $X$ is convergent.
\end{defin}

\begin{exo}\label{exo36}
Let $(X,\mcU)$ be a uniform space.
\begin{enumerate}
\item Prove that any convergent filter is Cauchy (this amounts to proving that any neighborhood filter is a Cauchy filter).
\item Assume that $(X,\mcU)$ is metrizable by some metric $d$. Prove that $(X,\mcU)$ is complete iff $(X,d)$ is complete. 
\end{enumerate}
\end{exo}

\begin{defin}
Let $(X,\mcU)$ be a uniform space, and $\mcF$ be a filter on $X$. We say that $x$ is an \emph{adherent point} of $\mcF$ if for any neighborhood $V$ of $x$ and any $F \in \mcF$ we have $V \cap F \ne \emptyset$.
\end{defin}

\begin{exo}\label{exo37}
Show that a Cauchy filter with an adherent point is convergent (to that point).
\end{exo}

It follows from this that any compact uniform space is complete (if $\mcF \subseteq \mcU$ are two filters and $\mcU$ converges to $x$ then $x$ is an adherent point of $\mcF$, so on a compact space every filter admits an adherent point).

\begin{exo}\label{exo38}
We say that a uniform space $(X,\mcU)$ is \emph{totally bounded} if for any $U \in \mcU$ there exists a finite $A \subset X$ such that $X= U[A]$ ($=\bigcup_{x \in A} U[x]$).

Show that $(X,\mcU)$ is compact iff it is both complete and totally bounded.
\end{exo}

We now discuss how to build the completion of a uniform space. I chose to include the details in these notes for completeness (no pun intended) but they may be skipped without hindering comprehension of the next chapters (one then just needs to accept as a black box that one can complete Hausdorff uniform spaces in essentially the same way that one completes metric spaces).

\begin{lem}\label{l:Cauchy_filter_on_a_subspace}
Let $(X,\mcU)$ be a uniform space. Assume that $Y$ is a dense subset of $X$ such that for every Cauchy filter $\mcF$ on $Y$ the filter on $X$ generated by $\mcF$ converges in $X$. Then $(X,\mcU)$ is complete.
\end{lem}

The filter on $X$ generated by $\mcF$ mentioned in the statement above is the set of subsets of $X$ which contain an element of $\mcF$ (equivalently, the image filter of $\mcF$ under the inclusion map from $Y$ to $X$).

\begin{proof}
Let $\mcF$ be a Cauchy filter on $X$. Let $\mcF^\mcU$ be the filter generated by all the subsets $U[F]$ ($=\bigcup_{x \in F} U[x]$) for some $U \in \mcU$ and $F \in \mcF$.

Then $\mcF^\mcU$ is contained in $\mcF$, and it is a Cauchy filter: fix $U \in \mcU$ and choose a symmetric $V \in \mcU$ such that $V \circ V \circ V \subseteq U$. Since $\mcF$ is Cauchy, there is $F \in \mcF$ such that $F \times F \subseteq V$. Then $V[F] \times V[F] \subseteq U$: if $(x,y) \in V[F] \times V[F]$ there exists $f_1,f_2 \in F$ such that $(x,f_1) \in V$ and $(y,f_2) \in V$, and since $(f_1,f_2) \in V$ we get $(x,y) \in V \circ V \circ V \subseteq U$.

For any element $A$ of $\mcF^\mcU$ we have $A \cap Y \ne \emptyset$ since $Y$ is dense in $X$, so $\mcF^\mcU$ induces a Cauchy filter on $Y$ (for the uniformity on $Y$ induced by $\mcU$). By assumption, we obtain that there exists $x \in X$ such that for every neighborhood $N$ of $x$, every $F \in \mcF$ and every $U \in \mcU$ one has $U[F] \cap N \ne \emptyset$. 

Let $N=V[x]$ be a neighborhood of $x$ with $V \in \mcU$, and let $U \in \mcU$ be symmetric and such that $U \circ U \subseteq V$. Given $F \in \mcF$, we can pick $z \in U[F] \cap U[x]$. Then $(x,z) \in U$ and there exists $f \in F$ such that $(z,f) \in U$, hence $(x,f) \in V$, so $F \cap V[x] \ne \emptyset$.

Finally, $x$ is an adherent point of the Cauchy filter $\mcF$, hence $\mcF$ converges to $x$ and $(X,\mcU)$ is complete.
\end{proof}

Given a uniform space $(X,\mcU)$, we let $\yhwidehat{X}$ denote the set of all Cauchy filters which are minimal for inclusion (among Cauchy filters).

\begin{lem}
For any $x \in X$, the neighborhood filter $\mcV_{x}$ is a minimal Cauchy filter.
\end{lem}

\begin{proof}
We already know that $\mcV_x$ is a Cauchy filter. Assume that $\mcF$ is a Cauchy filter contained in $\mcV_x$. Let $O$ be a neighborhood of $x$, which is of the form $U[x]$ for some $U \in \mcU$.

Since $\mcF$ is Cauchy, there exists $F \in \mcF$ such that $F \times F \subseteq U$; since $\mcF$ is contained in $\mcV_x$ we know that $F$ is a neighborhood of $x$, in particular $x \in F$. Thus $\{x\} \times F \subseteq U$, whence $F \subseteq U[x]=O$. It follows that $O \in \mcF$ and we are done. 
\end{proof}

Note that a minimal Cauchy filter which converges to some $x$ must coincide with the neighborhood filter $\mcV_x$; hence any non convergent minimal Cauchy filter witnesses that our uniform space is not complete. 

\begin{lem}
Let $\mcF$ be a Cauchy filter on $(X,\mcU)$. Then $\mcF$ contains a unique minimal Cauchy filter.
\end{lem}

\begin{proof}
Given $U \in \mcU$ and $F \in \mcF$, consider 
\[U[F]= \{y \in X \colon \exists x \in F \ (x,y) \in U \}.\]
It follows from the axioms defining uniformities and filters that any finite intersection of sets of this form is nonempty, so they generate a filter which we call $\mcG$.

To see that $\mcG$ is a Cauchy filter, fix $U \in \mcU$. Choose a symmetric $V \in \mcU$ such that $V \circ V \circ V \subseteq U$, and $F \in \mcF$ such that $F \times F \subseteq V$. Then $V[F] \times V[F] \subseteq U$; since $V[F] \in \mcG$ we have shown that $\mcG$ is a Cauchy filter.

Assume that $\mcG'$ is another Cauchy filter contained in $\mcF$. Fix $U \in \mcU$ and $F \in \mcF$. There exists $A \in \mcG'$ such that $A \times A \subseteq U$; since $A \in \mcG' \subseteq \mcF$ we have that $F'= A \cap F$ belongs to $\mcF$. In particular $A \cap F \ne \emptyset$, so $A \subseteq U[F]$. It follows that $U[F]$ belongs to $\mcG'$, hence $\mcG$ is contained in $\mcG'$. This proves both that $\mcG$ is minimal and that it is the unique minimal Cauchy filter contained in $\mcF$.
\end{proof}

It is then natural to consider the map $i \colon X \to \yhwidehat{X}$ which maps $x \in X$ to $\mcV_x$; this map is injective because of our standing assumption that all spaces are Hausdorff.

We now want to endow $\yhwidehat{X}$ with a uniform structure $\yhwidehat{\mcU}$ such that $i(X)$ is dense in $\yhwidehat{X}$ and $(\yhwidehat{X},\yhwidehat{\mcU})$ is complete. For some fixed entourage $U \in \mcU$, we let
\[C(U) = \lset (\mcF,\mcG) \in (\yhwidehat{X} \times \yhwidehat{X}): \exists A \in \mcF \cap \mcG \ A \times A \subseteq U \rset .\]
By definition of a Cauchy filter, $C(U)$ contains $\Delta_{\yhwidehat{X}}$ for every $U \in \mcU$. We want to prove that the sets $C(U)$ generate a uniform structure $\yhwidehat{\mcU}$ on $\yhwidehat{X}$; we let $\yhwidehat{\mcU}$ be the set of all subsets of $\yhwidehat{X} \times \yhwidehat{X}$ which contain some $C(U)$. The only thing that really requires checking is the existence of ``square roots'' in $\yhwidehat{\mcU}$, i.e., the last axiom in the definition of a uniformity. The next lemma shows that this property holds for $\yhwidehat{\mcU}$.

\begin{lem}
Assume that $V \circ V \subseteq U$. Then $C(V)\circ C(V) \subseteq C(U)$.
\end{lem}

\begin{proof}
Assume that $(\mcF,\mcG) \in C(V)$ and $(\mcG,\mcH) \in C(V)$. Then we can find $A \in \mcF \cap \mcG$ such that $A \times A \subseteq V$, and $B \in \mcG \cap \mcH$ such that $B \times B \subseteq V$. Let $D=A \cup B$; it belongs to $\mcF \cap \mcH$. Now, let $(x,y)$ belong to $D \times D$. There are four possibilities to consider:
\begin{enumerate}[(i)]
\item $x \in A$ and $y \in A$. Then $(x,y) \in A \times A \subseteq V \subseteq U$. 
\item $x \in B$ and $y \in B$. Then $(x,y) \in B \times B \subseteq V \subseteq U$.
\item $x \in A$ and $y \in B$. Since both $A$ and $B$ belong to $\mcG$, $A \cap B$ is nonempty so we can pick $z \in A \cap B$. Then $(x,z) \in A \times A$ and $(z,y) \in B \times B$ so $(x,z) \in V$ and $(z,y) \in V$ whence $(x,y) \in V \circ V \subseteq U$.
\item The case $x \in B$ and $y \in A$ is dealt with in the same way. \qedhere
\end{enumerate}
\end{proof}

\begin{lem}
$i(X)$ is dense in $\yhwidehat{X}$.
\end{lem}

\begin{proof}
Let $\mcF \in \yhwidehat{X}$ be a minimal Cauchy filter, and $U \in \mcU$. Choose a symmetric $V \in \mcU$ such that $V \circ V \subseteq U$.

Since $\mcF$ is Cauchy, there exists $F \in \mcF$ such that $F \times F \subseteq V$. In particular for all $x \in F$ we have $F \subseteq V[x]$ so $V[x] \in \mcF$. Since $V[x]$ is both a neighborhood of $x$ and an element of $\mcF$, and $V[x] \times V[x] \subseteq U$, we conclude that $(\mcF,\mcV_x) \in C(U)$. Hence $\mcV_x = i(x)$ belongs to the neighborhood $C(U)[\mcF]$. 
\end{proof}

\begin{lem}\label{l:uniform_isomorphism}
The map $i$ is uniformly continuous. Conversely, for any $U \in \mcU$ and $x,y \in X$ if $(\mcV_x,\mcV_y) \in C(U)$ then $(x,y) \in U$.
\end{lem}

It follows from this that $i$ is a uniform isomorphism from $(X,\mcU)$ to $i(X)$ endowed with the uniformity induced by $\yhwidehat{\mcU}$.

\begin{proof}
Fix $U \in \mcU$, and let $V \in \mcU$ be open, symmetric and such that $V \circ V \subseteq U$. Assume that $(x,y) \in V$. Then $V[y]$ is a neighborhood of $x$ and $y$; and $V[y] \times V[y] \subseteq U$. So $(\mcV_x,\mcV_y) \in C(U)$, and this proves that $i$ is uniformly continuous.

The second assertion follows immediately from the definition: if $(\mcV_x,\mcV_y) \in C(U)$ there exists a neighborhood $A$ of both $x$ and $y$ such that $A \times A \subseteq U$, in particular $(x,y) \in U$.
\end{proof}

\begin{lem}
The uniform space $(\yhwidehat{X},\yhwidehat{\mcU})$ is Hausdorff.
\end{lem}

\begin{proof}
Assume that $(\mcF,\mcG)\in \bigcap_{U \in \mcU} C(U)$; we have to prove that $\mcF=\mcG$.

Let $\mcH$ be the filter generated by sets of the form $F \cup G$ for $F \in \mcF$ and $G \in \mcG$. Clearly $\mcH$ is a filter which is contained in both $\mcF$ and $\mcG$; if we prove that $\mcH$ is a Cauchy filter then by minimality of $\mcF$ and $\mcG$ we will conclude that $\mcF=\mcH=\mcG$.

Choose $U \in \mcU$. There exists $F \in \mcF$ such that $F \times F \subseteq U$, as well as $G \in \mcG$ such that $G \times G \subseteq U$; by assumption on $(\mcF,\mcG)$ there also exists $A \in \mcF \cap \mcG$ such that $A \times A \subseteq U$. Replacing $F$ by $F \cap A$, $G$ by $G \cap A$ we obtain that there exist $F \in \mcF$, $G \in \mcG$ such that $F \times F \subseteq U$, $G \times G \subseteq U$, $F \times G \subseteq U$ and $G \times F \subseteq U$.

It follows that $(F \cup G) \times (F \cup G) \subseteq U$, whence $\mcH$ is a Cauchy filter.
\end{proof}

\begin{lem}
The uniform space $(\yhwidehat{X},\yhwidehat{\mcU})$ is complete. It is called the \emph{completion}\index{completion of a uniform space} of $(X,\mcU)$.
\end{lem}

\begin{proof}
It is enough to prove that for any Cauchy filter $\mcF$ on $i(X)$ the filter $\yhwidehat{\mcF}$ that it generates on $\yhwidehat{X}$ is convergent in $\yhwidehat{X}$ (see Lemma \ref{l:Cauchy_filter_on_a_subspace}, whose notations we reuse here). Let $\mcG=i^{-1}(\mcF)$, which is a Cauchy filter on $X$ by Lemma \ref{l:uniform_isomorphism}. Then $\mcG^\mcU$ is a Cauchy filter on $X$, contained in $\mcG$ (see \ref{l:Cauchy_filter_on_a_subspace}).

Let us show that $\mcG^\mcU$ is a minimal Cauchy filter: assume that $\mcH \subseteq \mcG^\mcU$ is also a Cauchy filter, then fix $G \in \mcG$ and $U \in \mcU$. There exists $H \in \mcH$ such that $H \times H \subseteq U$; replacing $H$ by $G \cap H$ if necessary we may assume that $H \subseteq G$. Since $H \times H \subseteq U$ we have $H \subseteq U[H] \subseteq U[G]$, hence $U[G] \in \mcH$ and we obtain $\mcG^\mcU \subseteq \mcH$.

It follows as promised that $\mcG^\mcU$ is a minimal Cauchy filter (even, the unique minimal Cauchy filter contained in $\mcG$), and our last step is to prove that $\yhwidehat{\mcF}$ converges to it in $\yhwidehat{X}$, i.e., that the filter of neighborhoods of $\mcG^\mcU$ is contained in $\yhwidehat{\mcF}$. 

By definition of $\yhwidehat{\mcF}$, we need to show that for every $U \in \mcU$ we have $C(U)[\mcG^\mcU]  \cap i(X) \in \mcF$, equivalently that $\Sigma:=\lset x \in X : \mcV_x \in C(U)[\mcG^\mcU] \rset \in \mcG$. 

Explicitly, we have
\[\Sigma = \lset x \in X : \exists A \in \mcV_x \cap \mcG^\mcU \ A \times A \subseteq U \rset .\]

Pick $U \in \mcU$ and let $V \in \mcU$ be symmetric and such that $V \circ V \circ V \subseteq U$. There exists $G \in \mcG$ such that $G \times G \subseteq V$, hence $V[G] \times V[G] \subseteq V \circ V \circ V \subseteq U$.

For every $x \in G$, $V[G]$ is a neighborhood of $x$, it belongs to $\mcG^\mcU$ and $V[G] \times V[G] \subseteq U$; this proves that $G \subseteq \Sigma$, hence $\Sigma \in \mcG$ and we are done.
\end{proof}

\begin{exo}\label{exo39}
Let $(X,\mcU)$ be a uniform space, and $(\yhwidehat{X}, \yhwidehat{\mcU})$ be its completion. 

Show that for any complete (and, as usual, Hausdorff) uniform space $(Y,\mcV)$ and any uniformly continuous $f \colon (X,\mcU) \to (Y,\mcV)$ there exists a uniformly continuous $\hat{f} \colon \yhwidehat{X} \to Y$ such that $\hat{f} \circ i = f$.

Explain in what sense this property characterizes $(\yhwidehat{X},\yhwidehat{\mcU})$.
\end{exo}

\begin{exo}\label{exo40}
Let $(X,\mcU)$ be a uniform space, and $Y$ be a dense subspace of $X$. Denote by $\mcU_Y$ the uniformity on $Y$ induced by $\mcU$. Prove that $(\yhwidehat{Y},\yhwidehat{\mcU_Y})=(\yhwidehat{X},\yhwidehat{\mcU})$.
\end{exo}

\begin{exo}\label{exo41}
Let $(X,\mcU)$ be a uniform space. Show that $(X,\mcU)$ is totally bounded iff $(\yhwidehat{X},\yhwidehat{\mcU})$ is compact.
\end{exo}

\emph{Comments.} Uniform structures were introduced by Weil as a means to study topological groups (though similar ideas had occurred earlier). There are many sources containing more material than is presented here, for instance \cite{Engelking}. The books \cite{Tkachenko} and\cite{RoelckeDierolf} are specifically concerned with topological groups. Our presentation is influenced by the seminal paper \cite{Samuel48}.
\chapter{Uniform structures on Polish groups}\label{uniform_structures_on_groups}
As promised, we can now turn to a discussion of uniform structures on topological groups. It turns out that there are four natural uniformities to consider; we introduce them now.

\begin{defin}
Let $G$ be a topological group. We consider the following four uniformities on $G$:
\begin{enumerate}[(i)]
\item The \emph{left uniformity}\index{left uniformity on a topological group} $\mcU_l$ is generated by the entourages $\lset (g,h) : g^{-1}h \in U \rset$, where $U$ ranges over all neighborhoods of $1$.
\item The \emph{right uniformity}\index{right uniformity on a topological group} $\mcU_r$ is generated by the entourages $\lset (g,h) : gh^{-1} \in U \rset$, where $U$ ranges over all neighborhoods of $1$.
\item The \emph{upper uniformity}\index{upper uniformity on a topological group} $\mcU_+$ is the coarsest uniformity refining both the left and right uniformities; the sets 
$\lset (g,h) : gh^{-1} \in U \text{ and } g^{-1}h \in U \rset$, where $U$ ranges over all neighborhoods of $1$, form a fundamental system of entourages for $\mcU_+$.
\item The \emph{lower uniformity} or \emph{Roelcke uniformity}\index{Roelcke uniformity on a topological group} $\mcU_{Roelcke}$ is the finest uniformity coarser than both the left and right uniformities; a fundamental system  for $\mcU_{Roelcke}$ is given by sets of the form 
$\lset (g,h) : h \in U g U\rset$, where $U$ ranges over all neighborhoods of $1$.
\end{enumerate}
\end{defin}

The convention we are using for ``left'' and ``right'' uniformity is sometimes switched in the literature; our choice of terminology is motivated by the fact that left translations map every element of the fundamental system given above to itself, and similarly for right translations and the right uniformity. But one needs to pay attention to the fact that $(g,h)$ being close for (say) the right uniformity means that one can go from $g$ to $h$ by multiplying \emph{on the left} by an element of $G$ which is close to $1$. As we will see below, the left uniformity on a metrizable group is induced by a left-invariant metric, and similarly for the right uniformity.

Note that in general $\mcU_l$ and $\mcU_r$ do not coincide and then all four uniformities above are distinct. Clearly left translations are uniform isomorphisms of $(G,\mcU_l)$, right translations are uniform isomorphisms of $(G,\mcU_r)$ and the inverse map is a uniform isomorphism of $(G,\mcU_+)$. Given that $\mcU_r$ and $\mcU_l$ do not coincide in general, the result of the next exercise is maybe a little surprising.

\begin{exo}\label{exo:right_translations_left_UC}
Show that for each $g$ the map $h \mapsto hg$ is a uniform isomorphism of $(G,\mcU_l)$, and the map $h \mapsto gh$ is a uniform isomorphism of $(G,\mcU_r)$. 

Prove that each left translation, as well as each right translation, is a uniform isomorphism for $\mcU_+$ and $\mcU_{Roelcke}$. 
\end{exo}

\begin{prop}
Let $(G,\tau)$ be a topological group. Each of the four uniformities  above induces the topology of $G$.
\end{prop}

\begin{proof}
Clearly, $\tau_{Roelcke} \subseteq \tau_l, \, \tau_r \subseteq \tau_+$. 

Assume that $A \subseteq G$ is $\tau_+$-open and let $g \in A$. Then there exists a $\tau$-neighborhood $U$ of $1$ such that $A$ contains $\lset h : g^{-1}h \in U \text{ and } gh^{-1} \in U \rset = gU \cap U^{-1}g$, a $\tau$-neighborhood of $g$. Hence $\tau_+ \subseteq \tau$.

Conversely, let $A \subseteq G$ be $\tau$-open and let $g \in A$. The map $\varphi \colon (f_1,f_2) \mapsto f_1 g f_2$ is $\tau$-continuous, and $\varphi(1,1)=g$. Hence there exists a neighborhood $U$ of $1$ such that $\varphi(U \times U) \subseteq A$. So $UgU \subseteq A$, and $UgU$ is a $\tau_{Roelcke}$-neighborhood of $g$. This proves that $\tau \subseteq \tau_{Roelcke}$, which concludes the proof.
\end{proof}

\begin{thm}[Birkhoff--Kakutani]\index{Birkhoff--Kakutani theorem}
Let $G$ be a topological group such that $1$ has a countable basis of neighborhoods. Then $G$ admits a compatible left-invariant metric.
\end{thm}

\begin{proof}
All four uniformities defined above are Hausdorff and admit a countable fundamental system of entourages as soon as $G$ is first-countable. So Theorem \ref{t:metrizable_uniformity} implies that all four uniformities are metrizable under our assumptions on $G$.

Now, let $d$ be a bounded metric inducing $\mcU_l$. While there is no reason for $d$ to be left-invariant, we can consider a new metric $\rho$ on $G$ defined by $\rho(g,h)= \sup \lset d(kg ,kh) : k \in G \rset$. This metric is left-invariant and $\mathrm{id} \colon (G,\rho) \to (G,d)$ is $1$-Lipschitz hence $\mathrm{id} \colon (G,\rho) \to (G,\mcU_l)$ is uniformly continuous.

To see the converse, fix $\varepsilon >0$. Since $d$ induces $\mcU_l$, there exists a neighborhood $U$ of $1$ such that for all $h$ one has $g^{-1}h \in U \Rightarrow d(g,h) \le \varepsilon$. But then whenever $g^{-1}h \in U$ we have $\rho(g,h) \le \varepsilon$, and this proves that $\mathrm{id} \colon (G,\mcU_l) \to (G,\rho)$ is uniformly continuous.
\end{proof}

\begin{exo}\label{exo43}
Let $G$ be a metrizable topological group and $d$ be a left-invariant compatible metric. Prove that the uniformity induced by $d$ coincides with $\mcU_l$.

Use this to recover the result that two left-invariant metrics inducing the same topology have the same Cauchy sequences.
\end{exo}

\begin{exo}\label{exo:Struble}
Let $G$ be a locally compact, first-countable topological group. 
\begin{enumerate}
\item Prove that for any left-invariant metric $d$ on $G$ there is $r>0$ such that closed balls of radius less than $r$ are compact. 
Use this to prove that any left-invariant metric on $G$ is complete.
\item (Struble)\index{Struble's theorem} Assume further that $G$ is second-countable. Show that $G$ admits a compatible, left-invariant metric in which every closed ball is compact. 

To establish this result, first show that there exists a sequence $(V_n)_{n \in \Z}$ of identity neighborhoods such that $\bigcap_{n \in \Z} V_n = \{1_G\}$, $\bigcup_{n \in \Z} V_n= G$, and $V_n^3 \subseteq V_{n+1}$ for each $n \in \Z$. Then apply the argument used in the proof of Theorem \ref{t:metrizable_uniformity}.
\end{enumerate}
\end{exo}

\begin{exo}\label{exo45}
Let $G$ be a first-countable topological group. Prove that $G$ admits a compatible bi-invariant metric if, and only if, $1$ admits a basis of neighborhoods which are conjugacy invariant (i.e., $gUg^{-1}=U$ for each $g \in G$ and each $U$ in the basis).
\end{exo}

\begin{exo}\label{exo46}
Let $G$ be a first-countable topological group. Prove that $(g,h) \mapsto gh$ is left-uniformly continuous  (i.e., uniformly continuous from $(G,\mcU_l)\times(G,\mcU_l)$ to $(G,\mcU_l)$) iff it is right-uniformly continuous iff $G$ admits a compatible bi-invariant metric.

Contrast this with the result of Exercise \ref{exo:right_translations_left_UC}.
\end{exo}

\begin{exo}\label{exo47}
Let $G$ be a Polish group. Prove that $(G,\mcU_+)$ is complete (one says that $G$ is \emph{Raikov-complete})\index{Raikov-complete topological group}.

Prove that if $\mcU_l$ and $\mcU_r$ coincide then they are both complete (note, however, that this is only an implication and not an equivalence).
\end{exo}

\begin{prop}[Solecki]\label{p:Solecki}
Let $G$ be a Polish group such that $(G,\mcU_l)$ is precompact (i.e., its completion is compact). Then $G$ is compact.
\end{prop}

\begin{proof}
Here we may simply work with sequences since $\mcU_l$ is metrizable (and the completion is obtained by taking the metric completion of $(G,d_l)$ for some compatible left-invariant metric $d_l$ on $G$). Let $(g_n)_{n < \omega}$ be a sequence of elements of $G$.

Since $(G,\mcU_l)$ is precompact, $(g_n)_{n< \omega}$ admits a subsequence $(g_{\varphi(n)})_{n < \omega}$ which is Cauchy in $\mcU_l$; applying this to $(g_{\varphi(n)}^{-1})$, we obtain that $(g_n)_{n < \omega}$ admits a subsequence $(g_{\psi(n)})_{n <\omega}$ which is Cauchy both for $\mcU_l$ and for $\mcU_r$, hence Cauchy for $\mcU_+$. Since $\mcU_+$ is complete, we conclude that $(g_{\psi(n)})_{n < \omega}$ is convergent. Hence $G$ is compact.
\end{proof}

The situation is completely different for the Roelcke uniformity.

\begin{defin}
A topological group is \emph{Roelcke precompact}\index{Roelcke precompact topological group} if the completion of $(G,\mcU_{Roelcke})$ is compact. 
\end{defin}

Note that this amounts to saying that $(G,\mcU_{Roelcke})$ is totally bounded, i.e., that for every $V \in \mcU_{Roelcke}$ there exists a finite subset $F \subseteq G$ such that $G=V[F]$. Unpacking this a little further, we obtain that $G$ is Roelcke precompact iff for every nonempty open $U$ there exists a finite $F$ such that $UFU=G$.

\begin{prop}
Let $G$ be a topological group which is both locally compact and Roelcke precompact. Then $G$ is compact.
\end{prop}

\begin{proof}
Since $G$ is locally compact, there exists a compact subset $K$ of $G$ with a nonempty interior. By Roelcke precompactness $G=KFK$ for some finite set $F \subseteq G$, and $KFK$ is compact since $(f,g,h) \mapsto fgh$ is continuous.
\end{proof}

However, many large Polish groups of interest are Roelcke precompact, and there is a rich connection with model theory.

\begin{thm}\label{thm:oligomorphic_is_Roelcke}
Every oligomorphic subgroup of $\Sinf$ is Roelcke precompact.
\end{thm}

Recall that $G \le \Sinf$ is oligomorphic iff there are finitely many orbits for the action $G \actson \omega^k$ for each $k$; oligomorphic closed subgroups of $\Sinf$ are exactly the automorphism groups of $\aleph_0$-categorical countable structures (in particular, automorphism groups of Fraïssé limits in a finite relational language are oligomorphic).

\begin{proof}
Let $U$ be a neighborhood of $1$; we may assume that $U=\lset g \in G : \forall i \le n \ g(i)=i\rset$ is a clopen subgroup. We need to find a finite subset $F$ of $G$ such that $G=UFU$, i.e., prove that there are only finitely many disjoint double cosets $UgU$.

Double cosets $UgU$ are in bijection with orbits for the diagonal action $G \actson G/U \times G/U$. Indeed, any orbit for this action contains an element of the form $(U,gU)$ (by translating on the first coordinate); and $(U,gU)$, $(U,hU)$ belong to the same orbit iff there exists $k$ such that $kU=U$ and $kgU=hU$, i.e., iff there exists $k \in U$ such that $kgU=hU$, which is equivalent to $g \in UhU$. So assigning to the orbit of $(U,gU)$ the double coset $UgU$ gives the desired bijection.

Since $U$ is the pointwise stabilizer of $\lset 0,\ldots,n \rset$, $G \actson G/U$ can be seen as the action of $G$ on a subset of $\omega^{n+1}$ (the orbit of $(0,\ldots,n)$); hence $G \actson G/U \times G/U$ is a subaction of $G \actson \omega^{n+1} \times \omega^{n+1}$, which has only finitely many orbits since $G$ is oligomorphic. So there are finitely many double cosets $UgU$ and we are done.
\end{proof}

For instance, $\Sinf$ itself is Roelcke precompact, as are the automorphism groups of the random graph, of $(\Q,<)$, and of any Fraïssé limit in a finite relational language. This gives many examples of non compact, Roelcke precompact Polish groups (there are also some interesting connected examples, and a strong connection with continuous logic, though we will not develop that here).

\begin{exo} \label{exo48bis}
We recall the definition of an inverse limit of topological groups: let $(I,\le)$ be an ordered set in which every finite subset has an upper bound, $(G_i)_{i \in I}$ be a family of topological groups and assume that for each $i \le j$ we are given a continuous homomorphism $\pi_{i,j} \colon G_j \to G_i$ such that $\pi_{i,j} \circ \pi_{j,k}= \pi_{i,k}$ whenever $i \le j \le k$. Assume further that $\pi_{i,i}=\mathrm{id}_{G_i}$ for all $i$. We call the maps $\pi_{i,j}$ the \emph{transition maps} and say that $((G_i)_{i \in I}, (\pi_{i,j})_{i \le j \in I})$ is an \emph{inverse system}.

Then the associated inverse limit is 
\[ \lset x \in \prod_{i \in I} G_i : \forall i \le j \ x(i)= \pi_{i,j}(x(j) \rset ,\]
which is a closed subset of $\prod_{i \in I} G_i$.

We say that the inverse system is \emph{surjective} if each projection map $\pi_i \colon G \to G_i$ is surjective.

Prove that the inverse limit of a surjective inverse system of Roelcke precompact groups is Roelcke precompact (in particular, the product of a family of Roelcke precompact groups is Roelcke precompact).

\end{exo}

\begin{exo}[Tsankov (\cite{Tsankov12}, Theorem 2.1)]\label{exo48} Prove that the following conditions are equivalent, for a subgroup $G \le \Sinf$:
\begin{enumerate}[(i)]
\item $G$ is Roelcke precompact.
\item For every continuous action $G \actson X$ of $G$ on a countable, discrete set $X$ with finitely many orbits, the diagonal action $G \actson X^n$ has finitely many orbits for all $n$.
\item $G$ is isomorphic (as a topological group) to the inverse limit of a surjective inverse system of oligomorphic subgroups of $\Sinf$. 
\end{enumerate}
\end{exo}

\begin{exo}\label{exo49}
Prove that a Roelcke precompact subgroup of $\Sinf$ has only (at most) countably many open subgroups.
\end{exo}

\begin{thm}[Tsankov] Assume that $G$ is a Roelcke precompact topological group which acts isometrically on a metric space $(X,d)$ so that $g \mapsto g \cdot x$ is continuous for each $x \in X$. Then every $G$-orbit is bounded.

In particular, a left-invariant continuous pseudometric on any Roelcke precompact topological group is always bounded (this applies to every left-invariant continuous metric on a Roelcke precompact Polish group).
\end{thm}

\begin{proof}
Fix $x_0 \in X$ and let $\varphi(g)=d(x_0,g \cdot x_0)$. We claim that $\varphi$ is Roelcke-uniformly continuous (that is, both left- and right-uniformly continuous). 

Granting that, there exists $\hat \varphi \colon \yhwidehat{(G,\mcU_{Roelcke})} \to \R$ uniformly continuous such that $\hat \varphi \circ i= \varphi$ (where $i$ denotes the map $x \mapsto \mcV_x$ from $(G,\mcU_{Roelcke})$ to its completion). Since $\yhwidehat{(G,\mcU_{Roelcke})}$ is compact, the image of $\varphi$ is then contained in a compact subset of $\R$, hence is bounded.

Now we establish our claim. Let $f,g,h$ belong to $G$. We have
\begin{eqnarray*}
d(x_0,fgh \cdot x_0) & = & d(f^{-1} \cdot x_0, gh \cdot x_0) \\
					 & \le & d(f^{-1}\cdot x_0, x_0) + d(x_0,g\cdot x_0) + d(g\cdot x_0,gh \cdot x_0) \\
					 & \le & d(x_0, f \cdot x_0) + d(x_0,g \cdot x_0) + d(x_0, h \cdot x_0)
\end{eqnarray*}
We conclude that for all $f$, $g$, $h$ one has
\[d(x_0,fgh \cdot x_0) - d(x_0,g \cdot x_0) \le d(x_0, f \cdot x_0) + d(x_0,h \cdot x_0)\]
Applying this inequality to $f^{-1}$, $fgh$, $h^{-1}$ we obtain 
\[d(x_0,g \cdot x_0) - d(x_0,fgh \cdot x_0) \le d(x_0,f^{-1} \cdot x_0)+ d(x_0,h^{-1}\cdot x_0)= d(x_0,f \cdot x_0)+ d(x_0,h \cdot x_0)\]
Finally we obtain that for every $f$, $g$, $h$ 
\[|d(x_0,fgh \cdot x_0) - d(x_0,g \cdot x_0)| \le d(x_0, f \cdot x_0) + d(x_0,h \cdot x_0)\]

Fix $\varepsilon >0$. Since $g \mapsto g \cdot x_0$ is continuous, there exists a neighborhood $U$ of $1$ such that $d(u \cdot x_0,x_0) \le \varepsilon$ for all $u \in U$. This yields that for all $(g,h)$ such that $h \in UgU$ we have $|\varphi(g) - \varphi(h)| \le 2 \varepsilon$, proving that $\varphi$ is Roelcke-uniformly continuous as claimed.
\end{proof}

\emph{Comments.} The interest in Roelcke precompact Polish groups stems in large part from work of Uspenskij \cite{Uspenskij} and then Rosendal \cite{Rosendal09}, Tsankov \cite{Tsankov12} and Ben Yaacov--Tsankov \cite{BenYaacovTsankov}. Recently this notion has been investigated quite thoroughly, and some weaker notions (such as locally Roelcke precompact groups, see \cite{Zielinski21}) have also been studied.  See Rosendal's book \cite{Rosendal22} for more information and connections with large scale geometry of Polish groups.

\chapter{Compactifications}
We can now make a link between compactifications of a given topological space $X$, and subalgebras of the algebra of all bounded continuous functions on $X$. This duality is key in our understanding of topological dynamics.

\begin{defin}
We let $\beta \omega$\index{$\beta \omega$} be the space of all ultrafilters on $\omega$, which we endow with the topology whose basic open sets are of the form $[A] =\lset p \in \beta \omega : A \in p \rset$, with $A$ a subset of $\omega$.
\end{defin}

Equivalently, we endow the space of ultrafilters on $\omega$ with the topology induced by the product topology on $2^{\mcP(\omega)}$, identifying each ultrafilter with its characteristic function. Indeed, a basic open set for this topology is of the form \[ U_{F_1,F_2} =\lset p : \forall A \in F_1 \ p(A)= 1 \text{ and } \forall B \in F_2 \ p(B)=0 \rset \]
where $F_1$, $F_2$ are finite subsets of $\mcP(\omega)$. A filter contains every element of $F_1$ iff it contains $B=\bigcap_{A \in F_1} A$; and an ultrafilter contains no element of $F_2$ iff it contains $C=\bigcap_{A \in F_2} (\omega \setminus A)$. Letting $A=B \cap C$, we see that $\beta \omega \cap U_{F_1,F_2} = [A]$.

\begin{prop}\label{p:extremally_disconnected}
\begin{enumerate}[(i)]
\item The space $\beta\omega$ is compact Hausdorff. 
\item The set of principal ultrafilters is dense in $\beta \omega$.

\item The closure of any open subset of $\beta \omega$ is open (one says that $\beta \omega$ is \emph{extremally disconnected}\index{extremally disconnected topological space}). Clopen subsets of $\beta \omega$ are of the form $[A]$ for some $A \subseteq \omega$.
\item For any two disjoint open subsets $U$, $V$ of $\beta \omega$ one has $\overline{ U }\cap \overline{ V }= \emptyset$.
\item For any convergent sequence $(p_n)_{n < \omega}$ in $\beta \omega$, there exists some $N$ such that $p_n=p_N$ for each $n>N$.
\end{enumerate}
\end{prop}

Below we often identify an element of $\omega$ with the corresponding principal ultrafilter, thus viewing $\omega$ as a dense subset of $\beta \omega$.

\begin{proof}
(i) Each of the following sets is a clopen subset of $2^{\mcP(\omega)}$ (below $A$, $B$ are subsets of $\omega$):
\begin{itemize}
\item $\Sigma_1 = \lset p : p(\omega)=1 \text{ and } p(\emptyset)=0\rset$.
\item $\Sigma_2(A,B) = \lset p : (p(A)= 1 \text{ and } p(B)=1) \Rightarrow p(A \cap B)=1 \rset$.
\item $\Sigma_3(A,B) = \lset p : (p(A)=1 \text{ and }A \subseteq B) \Rightarrow p(B)=1 \rset .$
\item $\Sigma_4(A) = \lset p : p(A)= 1 \text{ or } p(\omega \setminus A)=1 \rset .$
\end{itemize}
And we have 
\[ \beta \omega = \Sigma_1 \cap \bigcap_{A,B} \Sigma_2(A,B) \cap \bigcap_{A,B} \Sigma_3(A,B) \cap \bigcap_A \Sigma_4(A) .\]
This proves that $\beta \omega$ is closed in $2^{\mcP(\omega)}$, hence compact Hausdorff.

\medskip (ii) Let $A$ be a nonempty subset of $\omega$. For any $n \in A$ the principal ultrafilter associated to $n$ belongs to $[A]$.

\medskip (iii) Let $U$ be open in $\beta \omega$, and set $A=U \cap \omega$. Since $U$ is open and $\omega$ is dense we have $\overline{U} = \overline{A}$. We claim that $\overline A = [A]$, which is clopen (it is open and its complement is $[X \setminus A]$, which is open). To see this, note first that $A \subset [A]$ by definition. Conversely, let $p \in \beta \omega$ contain $A$, and let $[B]$ be a basic open neighborhood of $p$. Then both $A$ and $B$ belong to $p$, so $A \cap B \in p$ and is in particular nonempty. For any $n \in A \cap B$ we have $n \in A \cap [B]$, whence $A \cap [B] \ne \emptyset$. It follows that $p \in \overline{A}$. So $\overline{A}= [A]$ as promised.

\medskip (iv) Assume $U,V$ are disjoint open. Then $U \cap \omega = A$ and $V \cap \omega = B$ are disjoint nonempty subsets of $\omega$, hence $[A] \cap [B] = \emptyset$. Given what we just proved above, $\overline{U} \cap \overline{V} = [A] \cap [B] = \emptyset$.

\medskip (v) Finally, let $(p_n)_{n< \omega}$ be a convergent sequence. If it is not eventually constant, we may up to some extraction assume that it is an injective sequence which converges to some $p \not \in \lset p_n: n < \omega\rset$. Since $\beta \omega$ is Hausdorff, we may then inductively build a sequence of open subsets $(O_n)_{n< \omega}$ of $\beta \omega$ with $p_n \in O_n$ for all $n$ and $O_n \cap O_m = \emptyset$ for all $n \ne m$. Let $U= \bigcup_n O_{2n}$ and $V=\bigcup_n O_{2n+1}$. Then $U$ and $V$ are disjoint open but $p \in \overline{U} \cap \overline{V}$, a contradiction.
\end{proof}

\begin{prop}
Let $X$ be a compact topological space, and $f \colon \omega \to X$ a function. Then there exists a continuous $\hat{f} \colon \beta \omega \to X$ which extends $f$. 
\end{prop}

\begin{proof}
Given $p \in \beta \omega$, $f(p)$ is an ultrafilter on $X$, hence is convergent since $X$ is compact. We may then set $\hat{f}(p)= \lim f(p)$. Clearly $\hat{f}$ extends $f$, and it remains to prove that $\hat{f}$ is continuous.

Fix $p \in \beta \omega$, let $x=\hat{f}(p)$, and let $V$ be a neighborhood of $x$. Since $X$ is compact it is regular, i.e., there exists a neighborhood $W$ of $x$ such that $\overline{W} \subseteq V$.

By definition of convergence of filters, $W$ belongs to $f(p)$, equivalently $\lset n : f(n) \in W \rset$ belongs to $p$. Let $A=\lset n : f(n) \in W \rset$; $A \in p$ so $p \in [A]$. For every $q \in [A]$ we have $W \in f(q)$, so $\lim f(q) \in \overline{W}$ hence $\hat f(q) \in V$ for all $q \in [A]$, proving that $\hat f$ is continuous at $p$.
\end{proof}

\begin{exo}\label{exo50}
Explain in what sense the above proposition characterizes $\beta \omega$ among all compact spaces in which $\omega$ densely embeds.
\end{exo}

The space $\beta \omega$ is called the \emph{Stone-\v{C}ech} compactification \index{Stone-\v{C}ech compactification of $\omega$} of $\omega$. As we saw, it is a compact space to which $\omega$ maps densely (and injectively) so that all bounded functions from $\omega$ to $\R$ (for instance) extend continuously. This leads us to the general notion of a \emph{compactification}.

\begin{defin}
Let $X$ be a topological space. A \emph{compactification}\index{compactification of a topological space} is a pair $(Y,\varphi)$, where $Y$ is compact and $\varphi \colon X \to Y$ is continuous and has a dense image.

We say that $(Y,\varphi)$ is \emph{proper}\index{proper compactification} if $\varphi$ is a homeomorphism onto $\varphi(X)$.
\end{defin}

Note that we do not even ask that $\varphi$ is injective, for instance the constant map from $X$ to $\lset 0 \rset$ is a compactification in the above sense (beware that in the literature terminologies regarding compactifications vary!). 

Given a topological space $X$, we denote by $C_b(X)$\index{$C_b(X)$} the algebra of all bounded continuous functions from $X$ to $\R$ (and simply denote it by $C(X)$ when $X$ is compact). Endowed with the sup norm $\| \cdot \|_\infty$ this is a unital, commutative Banach algebra.

\begin{prop}
Let $(Y,\varphi)$ be a compactification of a topological space $X$. Then $\lset f \circ \varphi : f \in C(Y) \rset$ is a closed unital subalgebra of $C_b(X)$.
\end{prop}

\begin{proof}
It is immediate that $\lset f \circ \varphi : f \in C(Y) \rset$ is a subalgebra which contains $1$; since $\varphi(X)$ is dense the map $\Phi \colon f \mapsto f \circ \varphi$ is norm preserving from $(C(Y),\|\cdot\|_\infty)$ to $(C_b(X),\|\cdot\|_\infty)$. Since $C(Y)$ is complete its image under $\Phi$ also is, hence it is closed in $C_b(X)$.
\end{proof}

\begin{thm}[Gelfand--Naimark]
Let $X$ be a topological space, and $A \subseteq C_b(X)$ be a closed subalgebra containing the constant functions.

Then there exists a compactification $(Y,\varphi)$ such that $A=\{f \circ \varphi : f \in C(Y) \}$.

This compactification is proper iff $A$ separates points and closed sets, i.e., iff for every $x \in X$ and every closed $F \not \ni x$ there exists $a \in A$ such that $a(x) \not \in \overline{a(F)}$.
\end{thm}

Note that in the theorem above $A$ must be isomorphic (as a Banach algebra) with $C(Y)$. 

\begin{proof}
Fix $A$ as in the statement of the theorem. For every $a \in A$, denote $I_a = \overline{a(X)}$, which is a compact subset of $\R$. Let $Z= \prod_{a \in A} I_a$ and define $\varphi \colon X \to Z$ so that $\varphi(x)=(a(x))_{a \in A}$. Let $Y= \overline{\varphi(X)}$, which is compact. Clearly $(Y,\varphi)$ is a compactification of $X$.

Given $a \in A$, let $f_a \colon Y \to \R$ be defined by $f_a(y)=y(a)$. It is continuous and by definition for every $x \in X$ we have $f_a(\varphi(x))= a(x)$, in other words $a= f_a \circ \varphi$.

Note that $a \mapsto f_a$ is a homomorphism of Banach algebras from $A$ to $C(Y)$. Indeed, for every $a_1,a_2, a_3 \in A$, every $\lambda \in \R$ and every $x \in X$ we have 
\[f_{(\lambda a_1+a_2)a_3}(\varphi(x)) = ((\lambda a_1 + a_2)a_3)(x) = ((\lambda f_{a_1} + f_{a_2})f_{a_3})( \varphi(x))\]
Since $\varphi(X)$ is dense in $Y$ we obtain $f_{(\lambda a_1+a_2)a_3}=(\lambda f_{a_1} + f_{a_2})f_{a_3}$. 

So $\yhwidehat{A}=\lset f_a : a \in A \rset$ is a subalgebra of $C(Y)$, and for every $a \in A$ we have (again by density of $\varphi(X)$ in $Y$)
\[\|f_a\|_\infty = \sup\lset |f_a (\varphi(x))| : x \in X \rset = \sup \lset  |a(x)| : x \in X \rset =\|a\|_\infty . \]

It follows that $\yhwidehat{A}$ is closed in $C(Y)$ (it is isometric to $A$, which is complete). By definition, $\yhwidehat{A}$ separates points: if $y_1 \ne y_2 \in Y$ then for some $a$ we have $y_1(a) \ne y_2(a)$, i.e.,$f_a(y_1) \ne f_a(y_2)$. By the Stone--Weierstrass theorem, it follows that $\yhwidehat{A} = C(Y)$. Finally, 
\[ \lset f \circ \varphi : f \in C(Y) \rset = \lset f_a \circ \varphi : a \in A \rset =A .\]

This takes care of the first part of the statement. Now, assume that the compactification $(Y,\varphi)$ is proper and let $x \in X$, $F$ closed such that $x \not \in F$. Then $\overline{\varphi(F)}$ is a compact subset of $Y$ which does not contain $\varphi(x)$, so there exists $f \in C(Y)$ such that $f(\varphi(x)) \not \in \overline{f(\varphi(F))}$. Since $f \circ \varphi$ belongs to $A$, this proves that $A$ separates points and closed sets. 

To conclude the proof, assume that $A$ separates points and closed sets (whence $\varphi$ is injective), and let $U$ be an open subset of $X$ and $x \in U$. There exists $a \in A$ such that $a(x) \not \in \overline{a(X \setminus U)}$. 

For some $f \in C(Y)$ we have $a=f \circ \varphi$.
Then $V=\lset y : f(y) \not \in \overline{a(X \setminus U)} \rset$ is open in $Y$, contains $\varphi(x)$, and $\varphi(X) \cap V \subseteq \varphi(U)$. Hence $\varphi \colon X \to \varphi(X)$ is open.
\end{proof}

We have established that compactifications of a Hausdorff space $X$ correspond to closed unital subalgebras of $C_b(X)$. On the algebra side, we have a natural partial ordering given by inclusion (with $C_b(X)$ as its maximum). We now discuss what this partial order corresponds to for compactifications.

\begin{thm}
Let $X$ be a topological space, and $(Y_1,\varphi_1)$, $(Y_2,\varphi_2)$ be two compactifications of $X$. 
Let $A_1= \lset f \circ \varphi_1 \colon f \in C(Y_1) \rset$ and $A_2= \lset f \circ \varphi_2 : f \in C(Y_2) \rset$.

Then the following conditions are equivalent:
\begin{enumerate}[(i)]
\item $A_2 \subseteq A_1$.
\item There exists a continuous $\psi \colon Y_1 \to Y_2$ such that $\varphi_2 = \psi \circ \varphi_1$.
\end{enumerate}
\end{thm}

\begin{proof}
One implication is immediate: indeed, if there exists a continuous $\psi \colon Y_1 \to Y_2$ such that $\varphi_2 = \psi \circ \varphi_1$ then we have 
\[A_2= \lset f \circ \varphi_2 \colon f \in C(Y_2)  \rset =\lset (f \circ \psi) \circ \varphi_1 \colon f \in C(Y_2) \rset \subseteq \lset f \circ \varphi_1 \colon f \in C(Y_1)  \rset =A_1 .\]

Assume that $A_2 \subseteq A_1$. Note that if $x$, $x'$ are such that $\varphi_2(x) \ne \varphi_2(x')$ then there exists $a \in A_2$ so that $a(x) \ne a(x')$ (because $C_b(Y_2)$ separates points). Since $A_2 \subseteq A_1$ there exists $g \in C_b(Y_1)$ so that $g \circ \varphi_1 = a$, hence $\varphi_1(x) \ne \varphi_1(x')$.

This means that we can define $\psi \colon \varphi_1(X) \to \varphi_2(X)$ by setting $ \psi(\varphi_1(x))=\varphi_2(x)$. Certainly we guaranteed that $\psi \circ \varphi_1=\varphi_2$, but we do not know yet how to extend $\psi$ to $Y_1$, or whether it is continuous.

Both issues are taken care of at once if we prove that $\psi$ is uniformly continuous from $\varphi_1(X)$ to $Y_2$ (for the uniformities coming from the compact topologies of $Y_1$, $Y_2$) since $Y_2$ is complete for its unique compatible uniformity. By proposition \ref{p:carac_unif_continuity_compact} this is equivalent to proving that for every continuous $f \colon Y_2 \to \R$, $f\circ \psi$ is uniformly continuous on $\varphi_1(X)$.

By definition, we have for all $x$ that $f \circ \psi (\varphi_1(x))= f \circ \varphi_2(x)$. Thus $f \circ \psi \circ \varphi_1 \in A_2$, so the exists $g \in C(Y_1)$ such that $f \circ \psi \circ \varphi_1= g \circ \varphi_1$. Since $g$ is uniformly continuous on $Y_1$, $f \circ \psi$ is uniformly continuous on $\varphi_1(X)$. Hence $\psi$ is uniformly continuous, so it extends to a continuous map $\psi \colon \varphi_1(X) \to \varphi_2(X)$.
\end{proof}

In particular, two compactifications $(Y_1,\varphi_1)$ and $(Y_2,\varphi_2)$ which correspond to the same subalgebra of $C_b(X)$ are such that there exist a continuous $\psi \colon Y_1 \to Y_2$ such that $\psi(\varphi_1(x))= \varphi_2(x)$ for all $x \in X$, and a continuous $\tilde \psi \colon Y_2 \to Y_1$ such that $\tilde \psi(\varphi_2(x))= \varphi_1(x)$ for all $x$. By density it follows that $\tilde \psi = \psi^{-1}$, so $\psi$ is a homeomorphism from $Y_1 \to Y_2$ such that $\varphi_2 = \psi \circ \varphi_1$. When that situation occurs we say that the two compactifications are \emph{equivalent}\index{equivalent compactifications}.

Up to equivalence of compactifications, there is a unique compactification associated to any closed unital subalgebra $A$ of $C_b(X)$; for $A= C_b(X)$ we obtain the \emph{Stone--\v{C}ech compactification}\index{Stone--\v{C}ech compactification of a Hausdorff topological space} $(\beta_X,\beta)$ of $X$, which is characterized by the following universal property: for every continuous map $\varphi \colon X \to Y$ with $Y$ a compact space, there exists a continuous application $\psi \colon \beta_X \to Y$ such that $\varphi = \psi \circ \beta$.

In the particular case where $X=\omega$ with the discrete topology, we recover the space $\beta \omega$ which we already discussed, with the natural embedding of $\omega$ inside $\beta \omega$.

\begin{prop}
Let $(X,\mcU)$ be a uniform space. The set $A$ of all uniformly continuous bounded functions from $(X,\mcU)$ to $\R$ 
is a closed unital subalgebra of $C_b(X)$.
\end{prop}

\begin{proof}
Clearly $A$ contains the constant function $1$. To see that it is closed, one only needs to check that a uniform limit of uniformly continuous bounded functions is still uniformly continuous and bounded, and we leave that to the reader. It is also clear that $A$ is closed under linear combinations. Closure under products follows from a general result about uniform spaces, see exercise \ref{exo:product_unif_continuous} below.
\end{proof}

\begin{exo}\label{exo:product_unif_continuous}
Let $(X,\mcU)$ be a uniform space. Assume that $f,g \colon X \to \R$ are uniformly continuous and bounded. Prove that $fg$ is uniformly continuous.
\end{exo}

\begin{defin}
Let $(X,\mcU)$ be a uniform space. The \emph{Samuel compactification}\index{Samuel compactification} of $(X,\mcU)$ is the compactification associated to the algebra of all uniformly continuous bounded functions from $(X,\mcU)$ to $\R$.
\end{defin}

When $G$ is a topological group, it is tempting to consider the compactification $\beta G$ associated to $C_b(G)$, however there is in general no way to continuously extend the group operations to this compactification.
We now need to know which algebras are associated to continuous group actions.
Whenever $G$ acts on $X$ by homeomorphisms, it also acts on $C(X)$ via $(g \varphi)(x)= \varphi(g^{-1}x)$.

\begin{lem}\label{l:continuity_space_algebra}
Assume that $X$ is compact. Then the action $G \actson X$ is continuous iff $G \actson C(X)$ is continuous.
\end{lem}

\begin{proof}
Assume that $G$ acts on $X$ continuously. Certainly each $\varphi \mapsto g \varphi$ is continuous since it is isometric. So we need to prove continuity of the action at each $(1_G,\varphi)$, i.e., that for any $\varphi \in C(X)$ and $\varepsilon >0$ there exists a neighborhood $U$ of $1_G$ and $\delta >0$ such that
\[g \in U \text{ and } \|\psi - \varphi \| \le \delta \Rightarrow \|g \psi - \varphi \| \le \varepsilon . \]
Since $\|g \psi - g \varphi \| = \|\psi - \varphi\|$, this reduces to proving that for any $\varphi$ and $\varepsilon$ there exists a neighborhood $U$ of $1$ such that $\|g \varphi - \varphi\| \le \varepsilon$ for all $g \in U$.
Since $\varphi$ is uniformly continuous, there exists a neighborhood $V$ of $\Delta_X$ such that 
$|\varphi(x)- \varphi(y)| \le \varepsilon$ for all $(x,y) \in V$. 
By continuity of the action, for every $x \in X$ there exists an open $O_x \ni x$ and a symmetric neighborhood $U_x$ of $1_G$ such that $(y,gy) \in V$ for all $g \in U_x$ and all $y \in O_x$. Applying compactness we cover $X$ by $O_{x_1},\ldots,O_{x_n}$ and then $U=\bigcap_i U_{x_i}$ is a neighborhood of $1_G$ such that $(y,gy) \in V$ for all $y \in X$. It follows that for all $x \in X$ and all $g \in U$ we have $|\varphi(x) - \varphi(gx)| \le \varepsilon$.

Conversely, assume that $G \actson C(X)$ is continuous, and let $x \in X$ and $O \ni x$ be open in $X$. Since $X$ is compact, there exists a continuous $\varphi \colon X \to \R$ such that $\varphi(x)=0$ and $\{y \colon \varphi(y)< 1\} \subseteq O$. Let $V=\{y : \varphi(y) < \frac{1}{2} \}$, and $U$ a symmetric open neighborhood of $1_G$ such that $\|g \varphi - \varphi\| \le \frac{1}{2}$ for all $g \in U$. Then for every $g \in U$ and every $y \in V$ we have
\[\varphi(gy) \le \|g^{-1}\varphi - \varphi\| + \varphi(y) < 1 \]
hence $gy \in O$. This proves that $(g,x) \mapsto gx$ is continuous at $(1_G,x)$ and we are done.
\end{proof}

If $G$ acts continuously on a compact space $X$, and $x_0 \in X$, the pair $(\overline{G x_0},i_{x_0})$ defined by $i_{x_0}(g)= g x_0$ is a compactification of $G$. Then $\{g \mapsto f(gx_0) : f \in C(X)\}$ is a closed, unital subalgebra $A$ of $C_b(G)$, and by the previous proposition the natural action of $G$ on $A$ is continuous. In particular, for every $a \in A$ we have $\lim_{g \to 1} ga = a$ for every $a$ in $A$, which means that for every $\varepsilon >0$ there exists an open neighborhood $U$ of $1_G$ such that 
\[\forall g \in U \ \forall h \in G \ |a(g^{-1}h) - a(h)| \le \varepsilon .\]
This implies that $a$ is right-uniformly continuous (i.e., continuous for the right uniformity $\mcU_r$ on $G$). Indeed,   if $U$ is as above and $hg^{-1} \in U$ then $|a(g)-a(h)|= |a((hg^{-1})^{-1}h)- a(h)| \le \varepsilon$.

Let us sum up what we just observed.

\begin{prop}\label{p:right_uniform_continuity_compact_flow}
Let $G$ be a topological group acting continuously on a compact space $X$, and fix $x_0 \in X$.

Then the map $g \mapsto gx_0$ is uniformly continuous from $(G,\mcU_r)$ to $X$ (endowed with its unique uniform structure).   
\end{prop}

Note that this amounts to saying that for every $\varphi \in C(X)$ the map $g \mapsto \varphi(g x_0)$ is right-uniformly continuous, which is what we proved above.

\begin{defin} Let $G$ be a topological group. Denote by $\mathrm{RUC}_b(G)$ the $*$-algebra of all right-uniformly continuous bounded functions from $G$ to $\R$.
 \end{defin}

\begin{prop}
The following facts hold.
\begin{enumerate}[(i)]
\item\label{RUCb-2} $\mathrm{RUC}_b(G)$ is $G$-invariant for the action $G \actson C_b(G)$ (defined by $(g \varphi)(x)= \varphi(g^{-1}x)$).
\item\label{RUCb-3} The action of $G$ on $\mathrm{RUC}_b(G)$ is continuous.
\item\label{RUCb-4} If $A \subseteq C_b(G)$ is such that $\lim_{g \to 1} g \cdot f = f$ for every $f \in A$, then $A \subseteq \mathrm{RUC}_b(G)$.
\end{enumerate}
\end{prop}

\begin{proof}

\ref{RUCb-2} To check $G$-invariance, fix $g \in G$ and $\varphi \in \mathrm{RUC}_b(G)$. Given $\varepsilon >0$, there exists an open $U \ni 1$ such that $|\varphi(uh)-\varphi(h)| \le \varepsilon$ for all $u \in U$ and all $h \in G$.

This amounts to saying that $|(g\varphi)(guh)- (g\varphi)(gh)| \le \varepsilon$ for all $u \in U$ and $h \in G$, i.e., that $
|(g \varphi)(gug^{-1}k) - (g\varphi)(k)| \le \varepsilon$ for all $u \in U$ and all $k \in G$. Since $gUg^{-1}$ is an open neighborhood of $1$ we conclude that $g \varphi \in \mathrm{RUC}_b(G)$.

\medskip
\ref{RUCb-3} Fix $\varepsilon >0$ and $\varphi \in \mathrm{RUC}_b(G)$. We need to find $\delta >0$ and $U \ni 1_G$ open such that 
\[\forall \psi \ \forall g \quad \left(g \in U \text{ and } \|\varphi - \psi\| \le \delta \right) \Rightarrow \|g \psi  - \varphi\| \le \varepsilon .\] 
This is straightforward: we have for all $h \in G$
\[|g\psi(h) - \varphi(h)| \le \|\psi - \varphi\| + |\varphi(g^{-1}h) - \varphi(h)| .\]
Since $\varphi \in \mathrm{RUC}_b(G)$, there exists a neighborhood $U$ of $1_G$ such that for every $g \in U$ and every $h \in G$ we have $\|\varphi(g^{-1}h) - \varphi(h)| \le \frac{\varepsilon}{2}$. This $U$ along with $\delta=\frac{\varepsilon}{2}$ are what we were looking for.

\medskip
\ref{RUCb-4} is immediate from the definition of $\mathrm{RUC}_b(G)$ (and has already been noted above).
\end{proof}

\begin{defin}
Let $G$ be a topological group. We denote by $(S(G),i)$\index{$S(G)$} the Samuel compactification of $(G,\mcU_r)$.
\end{defin}

\begin{exo}\label{exo52}
Let $G$ be a Polish group. Prove that the compactification $(S(G),i)$ is proper.

(This fact actually holds for any Hausdorff topological group, which makes for a slightly more challenging exercise)
\end{exo}

Recall that we realized $S(G)$ ``concretely'' by considering 
\[Y= \prod_{\varphi \in \mathrm{RUC}_b(G)} \overline{\varphi(G)} \ , \quad \ i(g)(\varphi)=\varphi(g) \ , \quad S(G)=\overline{i(G)} .\]

Since $\mathrm{RUC}_b(G)$ is $G$-invariant, it follows that $G$ acts on $S(G)$ via $(g \cdot y)(\varphi)= y(g^{-1}\varphi)$. Note that each $y \mapsto g \cdot y$ is a homeomorphism.

\begin{exo}\label{exo53}
Prove that the above formula indeed defines an action of $G$ on $S(G)$, and that $i(gh)=g \cdot i(h)$ for all $g,h \in G$.
\end{exo}

\begin{prop}
The action $G \actson S(G)$ is continuous.
\end{prop}

\begin{proof}
Let $i \colon G \to S(G)$. The map $\varphi \mapsto \varphi \circ i$ induces an algebra isomorphism from $C(S(G))$ to $\mathrm{RUC}_b(G)$, and the previous exercise shows that $i$ is $G$-equivariant.

Since $G \actson \mathrm{RUC}_b(G)$ is continuous, we conclude that $G \actson C(S(G))$ is continuous ($i$ carries all the structure of one, including the action, onto the other) so $G \actson S(G)$ is continuous by Lemma \ref{l:continuity_space_algebra}.
\end{proof}

\begin{exo}\label{exo54}
Let $A \subseteq \mathrm{RUC}_b(G)$ be a $G$-invariant, closed, unital subalgebra. Let $X_A$ be the compactification associated to $A$. Prove that the left translation of $G$ on itself extends to a continuous action of $G$ on $X_A$.
\end{exo}

\begin{exo}\label{exo55}
Assume $G$ is discrete (countable if you wish) and let $\beta G$ be its Stone--\v{C}ech compactification. 

Prove that $G \actson S(G)$ is simply the translation action $G \actson \beta G$, where $A \in g \cdot p \Leftrightarrow g^{-1}A \in p$.
\end{exo}

\emph{Comments.} The book \cite{deVries} is an encyclopedic reference about topological dynamics, in particular Appendix D contains information about compactifications, Many books cover this material, for instance the very first chapter of Folland's book \cite{Folland} (though with a somewhat different formulation). Uspenskij's paper \cite{Uspenskij} gives good reasons to care about compactifications of topological groups.
The theory presented here fits into the much broader picture of operator algebras, see \cite{Blackadar06} for an in-depth introduction to this topic (though here we only worked with real Banach algebras because it simplifies the presentation a little).

\chapter{The greatest ambit and universal minimal flow}
In this chapter, we reap the rewards of our work so far, which enables us to prove the existence of two objects that capture many dynamical properties of a given topological group: its \emph{greatest ambit} and \emph{universal minimal flow}.

\begin{defin}
Let $G$ be a topological group. A \emph{$G$-flow}\index{flow} is a continuous action $G \actson X$ where $X$ is compact, nonempty and (as always) Hausdorff.

A $G$-\emph{ambit}\index{ambit} is a pair $(X,x_0)$ where $x_0 \in X$, $X=\overline{Gx_0}$ is compact and $G \actson X$ is continuous (i.e., a $G$-ambit is a $G$-flow where we named a point with a dense orbit).
\end{defin}

We identify $G$ with its image $i(G)$ in $S(G)$.

\begin{thm}
Let $(X,x_0)$ be a $G$-ambit. Then there exists a unique continuous map $\pi \colon S(G) \to X$ which is $G$-equivariant and such that $\pi(1_G)=x_0$. 
\end{thm}
Note that $\pi$ above is automatically surjective, since $S(G)$ is compact and $Gx_0$ is dense in $X$. We say that $(S(G),1_G)$ is the \emph{greatest ambit}\index{greatest ambit} of $G$.

\begin{proof}
Consider the map $\varphi \colon g \mapsto g \cdot x_0$. Then $(X,\varphi)$ is a compactification of $G$, and since the action $G \actson X$ is continuous it follows from Proposition \ref{p:right_uniform_continuity_compact_flow} that $\lset f \circ \varphi : f \in C(X) \rset \subseteq \mathrm{RUC}_b(G)$. This gives us the existence of a continuous $\pi \colon S(G) \to X$ such that $\pi \circ i = \varphi$.

We then have $\pi(1_G)=x_0$. Also, for every $g,h \in G$ we have $\pi(gh)= gh \cdot x_0 =  g \cdot \pi(h)$. By continuity of $G \actson S(G)$, density of $G$ on $S(G)$ and continuity of $\pi$ we conclude that for every $g \in G$ and every $p \in S(G)$ we have $\pi(g \cdot p)= g \cdot \pi(p)$. This proves that $\pi$ is $G$-equivariant.

Uniqueness of $\pi$ is immediate, since $G$ is dense in $S(G)$ and on $G$ we have $\pi(g)=g \cdot x_0$ by $G$-equivariance.
\end{proof}

\begin{defin}
Let $G$ be a topological group. A $G$-flow $X$ is \emph{minimal}\index{minimal flow} if every $G$-orbit is dense in $X$.
\end{defin}

Note that if $X$ is minimal then $(X,x_0)$ is an ambit for every $x_0 \in X$.

\begin{exo}\label{exo56}
\begin{enumerate}[(i)]
\item Prove that a $G$-flow $G \actson X$ is minimal iff the only closed $G$-invariant sets are $\emptyset$ and $X$ iff the only open $G$-invariant sets are $\emptyset$ and $X$ iff for every nonempty open $U \subseteq X$ there exists a finite $F \subseteq G$ such that $X=F \cdot U$.
\item Let $X$ be a $G$-flow. Prove that the set of all nonempty closed $G$-invariant subsets of $X$, ordered by $\supseteq$, is inductive. Apply Zorn's lemma to prove that there exists a nonempty closed $G$-invariant $Y \subseteq X$ such that $G \actson Y$ is minimal.
\end{enumerate}
\end{exo}

\begin{defin}
Let $G$ be a topological group. A minimal $G$-flow $X$ is \emph{universal}\index{universal minimal flow} if for any minimal $G$-flow $Y$ there exists a continuous equivariant map $\pi \colon X \to Y$.
\end{defin}

Note that, by minimality of $Y$, a map $\pi$ as above must be onto.

\begin{prop}
Any minimal subflow of $S(G)$ is universal.
\end{prop}

\begin{proof}
Let $Y$ be a minimal $G$-flow, and $y \in Y$. Let $M$ be a minimal subflow of $S(G)$. There exists a $G$-equivariant map $\pi \colon S(G) \to Y$. The restriction of $\pi$ to $M$ witnesses that $M$ is universal.
\end{proof}

To show the existence of a universal minimal flow, we could also have taken the product of a set of representatives of all possible  minimal $G$-flows, and taken a minimal component in that product. So there is nothing unexpected in the existence of a universal minimal flow; but it turns out to be unique up to isomorphism, which makes it an interesting object to study. To prove this we are going to introduce some additional structure on $S(G)$ (though again there are other possible arguments, see for instance \cite{Gutman_Li}).

\begin{defin}
Let $p,q$ belong to $S(G)$. 

There exists a unique continuous $G$-equivariant $\pi_q \colon S(G) \to S(G)$ such that $\pi_q(1)=q$. We define $p \cdot q = \pi_q(p)$\index{semigroup structure on the Samuel compactification} .
\end{defin}

The existence and uniqueness of $\pi_q$ come from considering $(\overline{Gq},q)$ as an ambit and applying the universal property of $S(G)$.

\begin{prop}
The map $(p,q) \mapsto p \cdot q$ is associative and extends the group action of $G$ on $S(G)$. Furthermore, $p \mapsto p \cdot q$ is continuous for all $q \in S(G)$.
\end{prop}
One says that $S(G)$ is a right topological semigroup (right translations are continuous). Beware however that in general left translations are not continuous!

\begin{proof}
Since $p \cdot q=\pi_q(p)$ and $\pi_q$ is continuous, the continuity of $p \mapsto p \cdot q$ for all $q \in S(G)$ is by definition.

Given $g \in G$ and $q \in S(G)$ we have $g \cdot q = \pi_q(g)= \pi_q(g1_G)= g \pi_q(1_G)=g q$.

To check associativity, let $p,q,r \in S(G)$. We have $\pi_r \circ \pi_q(1)= \pi_r(q) = q \cdot r$, and $\pi_r \circ \pi_q$ is $G$-equivariant. Hence $\pi_{q \cdot r}= \pi_r \circ \pi_q$, and it follows that 
\[p \cdot (q \cdot r)= \pi_{q \cdot r}(p)= \pi_r(\pi_q(p))= \pi_q(p) \cdot r=(p \cdot q) \cdot r .\qedhere \]
\end{proof}

\begin{defin}
We say that $I \subseteq S(G)$ is a \emph{left-ideal}\index{left-ideal in $S(G)$} if $I$ is nonempty, closed and $S(G)I \subseteq I$.

We say that a left-ideal is \emph{minimal} if it contains no proper left-ideal.
\end{defin}

Note that since $I$ is closed, $G$ is dense in $S(G)$ and the semigroup operation extends the action $G \actson S(G)$, $I$ is a left-ideal iff $GI =I$ iff $I$ is a subflow of $G \actson S(G)$. Similarly, minimal left-ideals correspond to minimal subflows of $G \actson S(G)$. So Zorn's lemma implies that minimal left-ideals exist in $S(G)$.

\begin{lem}[Ellis]\index{Ellis lemma} Let $I$ be a left-ideal in a compact right-topological semigroup. Then $I$ contains an \emph{idempotent}\index{idempotent element in $S(G)$}, i.e., there exists $e \in I$ such that $e \cdot e=e$.
\end{lem}
The previous lemma of course applies to $S(G)$.

\begin{proof}
Using Zorn's lemma, we find $S$ minimal among nonempty closed subsets of $I$ such that $S \cdot S \subseteq S$ (note that $I \cdot I \subseteq I$).

For any $x \in S$, we have $S \cdot x \subseteq S$; and $(S \cdot x) \cdot(S \cdot x) \subseteq S^3 \cdot x \subseteq S \cdot x$. By minimality of $S$, it follows that $S \cdot x= S$ for every $x \in S$.
Let $W_x = \lset y \in S : y \cdot x = x \rset$, which is nonempty since $S \cdot x =S$. Then $W_x$ is closed by continuity of $y \mapsto y \cdot x$. Further, by associativity we have $W_x \cdot W_x \subseteq W_x$; hence $W_x=S$ by minimality of $S$ again. Thus $y \cdot x = x$ for all $x \in S$, so $S=\lset e \rset$ with $e \cdot e=e$.
\end{proof}

\begin{lem}
Let $I$ be a minimal left-ideal in $S(G)$, and $e \in I$ be an idempotent. Then for every $p \in I$ we have that $p \cdot e=p$ ; for every $q \in M,$ $p \mapsto p \cdot q$ is a $G$-equivariant surjection of $I$ onto itself.
\end{lem}

\begin{proof}
Note that $I\cdot e$ is a left-ideal contained in $I$ since $I$ is a left ideal and $\cdot$ is associative. Hence $I \cdot e = I$. For every $ p \in I$ there exists some $q \in I$ such that $p= q\cdot e$, from which we obtain that $p \cdot e= q \cdot e^2= q \cdot e = p$. 

Similarly for every $q \in I$ we have $I \cdot q=I$, and $G$-equivariance of $p \mapsto p \cdot q$ is part of its definition.
\end{proof}

The uniqueness of a universal minimal flow is a consequence of the next result.

\begin{thm}[Ellis] Let $G$ be a topological group, and $M \subseteq S(G)$ be a minimal left-ideal. Then $M$ is \emph{coalescent}\index{coalescent flow}, i.e., every $G$-equivariant continuous $\pi \colon M \to M$ is bijective.
\end{thm}

The surjectivity of $\pi$ is immediate by minimality of $M$, so the interesting point here is that $\pi$ is injective.

\begin{proof}
Let $\pi \colon M \to M$ be $G$-equivariant and continuous. For all $q \in M$ and all $g \in G$ we have $\pi(g \cdot q) = g \pi(q)$, so by continuity of right translations we obtain $\pi(p \cdot q)= p \cdot \pi(q)$. 

Let $e \in M$ be an idempotent, and $p=\pi(e)$. For all $q \in M$ we have 
\[ \pi(q)= \pi(q \cdot e)= q \cdot \pi(e)= q \cdot p .\]
So $\pi$ is the right -translation by $p$.

Let $q$ be such that $q \cdot p = e$ (such a $q$ exists because $M \cdot p=M$), and set $\rho(x)=x \cdot q$. Then $\pi (\rho (x)) = \rho(x) \cdot p = x \cdot (q \cdot p)=x$. Hence $\rho$ is injective, and we already knew that it is surjective (by minimality of $M$, see the previous lemma). Hence $\rho$ is a bijection of $M$, and from $\pi \circ \rho= id$ we obtain that $\pi= \rho^{-1}$ is bijective.
\end{proof}

\begin{thm}
Let $G$ be a topological group. Up to isomorphism, there exists a unique universal minimal $G$-flow, which we denote by $M(G)$\index{$M(G)$}.
\end{thm}

\begin{proof}
We already established the existence of a universal minimal $G$-flow. So, let $M$ be a minimal subflow of $S(G)$, and $N$ another universal minimal $G$-flow. Applying the definition of universality, we obtain two $G$-equivariant continuous maps $\varphi \colon M \to N$ and $\psi \colon N \to M$.

Then $\psi \circ \varphi \colon M \to M$ is continuous and $G$-equivariant, hence bijective since $M$ is coalescent. Hence $\varphi$ is injective, so it is an isomorphism of $G$-flows.
\end{proof}

Note that it follows from our observations that any minimal subflow of $S(G)$ is isomorphic to $M(G)$.
We seize the opportunity to mention the following useful fact.

\begin{prop}\label{p:retraction}
Let $G$ be a topological group, and $M$ be a minimal subflow of $S(G)$. Then there exists a continuous $G$-equivariant retraction $r \colon S(G) \to M$.
\end{prop}

\begin{proof}
Let $\pi \colon S(G) \to M$ be $G$-equivariant. Then $\pi$ maps $M$ into itself, so by coalescence it induces an automorphism of $M$ which we denote by $\varphi$. Let $r= \varphi^{-1} \circ \pi$. 

Then $r$ is $G$-equivariant, continuous, maps $S(G)$ onto $M$, and for every $x \in M$ we have $r(x)= \varphi^{-1} (\pi(x))=x$. So $r$ is the desired retraction.
\end{proof}

Since we are going to consider universal minimal flows of subgroups of $\Sinf$ in the next chapter, we mention a result which implies one only needs to understand what happens for closed subgroups.

\begin{prop}
Let $G$ be a topological group, and assume that $H$ is a dense subgroup of $G$ acting continuously on a compact space $X$. 

Then the action of $H$ extends to a continuous action of $G$ on $X$.
\end{prop}

\begin{proof}
We saw in Proposition \ref{p:right_uniform_continuity_compact_flow} that each map $h \mapsto hx$ is right-uniformly continuous; by compactness of $X$ this map extends to the completion of $(H,\mcU_r)$, which coincides with the completion of $(G,\mcU_r)$ since $H$ is dense in $G$. In particular, for all $x$ we can continuously extend $h \mapsto hx$ to $G$. 

This gives us an action of $G$ on $X$; we want to prove that this action is by homeomorphisms, i.e., that each mapping $x \mapsto gx$ is continuous. Fix $g \in G$, $x \in X$ and $V$ a neighborhood of $gx$. By continuity of $H \actson X$ and compactness of $X$ , there is a neighborhood $U$ of $1$ in $H$ and an open $W$ containing $gx$ such that $\overline{UW} \subseteq V$.

Since $H$ is dense in $G$, $\overline{U}$ is a neighborhood of $1$ in $G$ and $k \mapsto kx$ is continuous there is $h \in H$ such that $g \in \overline{Uh}$ and $hx \in W$. Let $O \ni x$ be open and such that $hO \subseteq W$. Then for every $y \in O$ we have (using continuity of $k \mapsto ky$) that $g y \in \overline{Uh y} \subseteq \overline{UW} \subseteq V$.

Now that we know that our action is by homeomorphisms, fix $\varphi \in C(X)$ and $\varepsilon >0$. There is an open neighborhood $U$ of $1$ in $H$ such that $\|h \varphi - \varphi\| \le \varepsilon$ for each $h \in U$. Since each $g \mapsto gx$ is continuous, it follows that $\|g \varphi - \varphi\| \le \varepsilon$ for each $g \in \overline{U}$. This proves that $G \actson C(X)$ is continuous, and we are done.
\end{proof}

\begin{exo}\label{exo57}
Let $G$ be a topological group, and $H$ be a dense subgroup of $G$. Show that $M(G)=M(H)$, in the sense that both $G$ and $H$ act on the same space and the two 
$H$-actions coincide (hence the $G$-action on $M(G)=M(H)$ is the unique continuous extension of the $H$-action).
\end{exo}

\begin{defin}
A topological group is \emph{extremely amenable}\index{extremely amenable topological group} if every $G$-flow admits a fixed point.
\end{defin}

Note that this property is never satisfied by a nontrivial locally compact group: indeed, the existence of a Haar measure immediately implies that there exist nontrivial actions of $G$ on compact spaces (actually, irreducible unitary representations separate points, so there exist many such actions; for an even stronger fact, see Veech's theorem below).

\begin{exo}\label{exo57bis}
Let $G$ be a Polish group. Prove that $G$ is extremely amenable iff every continuous action of $G$ on a compact \emph{metrizable} space admits a fixed point.
\end{exo}

Since every $G$-flow contains a minimal subflow, and a minimal subflow with a fixed point is a singleton, $G$ is extremely amenable iff $M(G)$ is a singleton. 

In the next chapter, we discuss a combinatorial characterization of extremely amenable subgroups of $\Sinf$. 
The colorings that appear in the next chapter can be thought of as discretized versions of $1$-Lipschitz maps from $G$ to $[0,1]$; in hindsight, the following result explains why one may expect such colorings (hence, Ramsey theory) to be related to extreme amenability.

\begin{thm}\label{thm:finite_osc_stable}
Let $G$ be a Polish group and $d$ a compatible right-invariant metric on $G$. Then $G$ is extremely amenable if, and only if, the following property holds:

For every $1$-Lipschitz $\varphi \colon G \to [0,1]$, every $\varepsilon >0$ and every finite $A \subseteq G$ there exists $g \in G$ such that $|\varphi(ag)-\varphi(a'g)|< \varepsilon$ for all $a,a' \in A$.
\end{thm}

Let us provide some context for this condition: the space $L_G$ of all $1$-Lipschitz functions from $(G,d)$ to $[0,1]$, endowed with the pointwise convergence topology, is compact, and $G$ naturally acts on $L_G$ via $g \varphi(g')=\varphi(g'g)$. We prove below that the condition above amounts to the statement that for every $\varphi \in L_G$ there exists a $G$-fixed point in $\overline{G \varphi}$; thus this condition is clearly necessary for $G$ to be extremely amenable. When this condition is satisfied, one says that the action of $G$ on itself by left translation is \emph{finitely oscillation stable}; we refer to \cite{Pestov06} for more on this concept and its relevance when one is studying extreme amenability.

\begin{proof}
If $G$ is extremely amenable, then for every $\varphi$ there exists a fixed point in $\overline{G \varphi}$, which implies that the condition above is satisfied. Conversely, assume that that this condition holds and let $G \actson X$ be a $G$-flow. Pick $x_0 \in X$.

For every $\varphi \in C(X)$ the map $g \mapsto \varphi(gx_0)$ is right-uniformly continuous and bounded, hence is a uniform limit of a sequence $(\varphi_n)_{n< \omega}$ of bounded Lipschitz maps from $G$ to $\R$ (see for instance Theorem 6.8 of \cite{Heinonen2001}). Our assumption thus gives us that for every $\varepsilon >0$ and every finite $A \subseteq G$ there exists $g \in G$ such that $|\varphi(agx_0) - \varphi(a'gx_0)| \le \varepsilon$ for all $a,a' \in A$. Denote
\[F_{A,\varepsilon}= \lset x \in X : \forall a,a' \in A \ |\varphi(ax)-\varphi(a'x)| \le \varepsilon \rset . \]
These are closed subsets of $X$ with the finite intersection property, so by compactness there exists $x \in X$ such that $\varphi(gx)=\varphi(x)$ for every $g \in G$. The set of such $x$ is closed and $G$-invariant, whence by repeating this argument we obtain that for every finite subset $F$ of $C(X)$ there exists $x \in X$ such that $\varphi(gx)=\varphi(x)$ for all $g \in G$ and all $\varphi \in F$. Using compactness one last time, we conclude that there exists $x \in X$ such that $\varphi(gx)=\varphi(x)$ for all $g \in G$ and all $\varphi \in C(X)$, i.e., $x$ is a fixed point of the action.
\end{proof}

To finish this chapter, we establish a fundamental theorem due to Veech (we only prove it for Polish groups because it slightly simplifies the presentation, and all the important ideas of the proof are already present in that case). We first need two lemmas.

\begin{lem}
Let $\Gamma$ be a graph (i.e., an irreflexive symmetric binary relation) on a set $X$; assume that every element of $X$ has at most $N$ neighbors in $\Gamma$. Then one can color the vertices of $\Gamma$ using at most $N+1$ colors, that is, there exists a partition $X= \bigsqcup_{i=1}^{N+1} X_i$ such that no two elements of the same $X_i$ are adjacent in $\Gamma$.
\end{lem}

\begin{proof}
Fix an enumeration $(x_{\alpha})_{\alpha < \beta}$ of $X$. Define a map $c \colon X \to \{1,\ldots,N+1\}$ by letting 
\[c(x_\alpha) = \min \{ i \in \{1,\ldots,  N+1 \} : \textrm{ for each } \alpha' < \alpha \textrm{ such that } (x_\alpha,x_{\alpha}') \in \Gamma \textrm{ one has } c(x_{\alpha'}) \ne i \}. \]
Note that, since every element of $X$ has at most $N$ neighbors in $\Gamma$, this map is well defined; setting $X_i=c^{-1}(\{i\})$ satisfies the required condition.
\end{proof}

\begin{thm}[Veech]\index{Veech's theorem}
Let $G$ be a locally compact Polish group. Then $G$ acts freely on $S(G)$.
\end{thm}

\begin{proof}
Fix $g \ne 1_G$. We prove that there exist $\varphi_1, \ldots,\varphi_m \in \mathrm{RUC}_b(G)$ such that for each $h \in G$ there exists $i \in \{1,\ldots,m\}$ such that $|\varphi_i(h) - \varphi_i(gh)| \ge 1$. All those functions extend continuously to $S(G)$ and it follows that for each $x \in S(G)$ there exists $i \in \{1,\ldots,m\}$ such that $|\varphi_i(x) - \varphi_i(g \cdot x)| \ge 1$ (note that we use the fact that $m$ is finite here). It then follows immediately that $g$ does not fix any element of $S(G)$.

Applying Struble's theorem (see Exercise \ref{exo:Struble}), we pick a compatible, right-invariant metric $d$ on $G$ such that $9\le d(g,1_G)  \le 10$ and $K=\{h: d(h,1_G) \le 15 \}$ is compact. Replacing $d$ with $\min(d,20)$ if necessary, we may assume that $d$ is bounded.

Now let $A$ be a maximal subset of $G$ such that for any $a \ne a' \in A$ one has $B(a,1) \cap B(a',1)= \emptyset$. For every $g \in G$ there exists $a \in A$ such that $d(g,a) < 2$. We consider the graph $\Gamma$ on $A$ such that $(a,a') \in \Gamma$ iff $a \ne a'$ and $a' a^{-1} \in K^2$. 

Since $K^2$ is compact, it is covered by some finite number $N$ of open balls of radius $1$; the cardinality of any family of pairwise disjoint open balls of radius $1$ with centers belonging to $K^2$ must then be smaller than $N$. Note that if $a,a_1\ldots,a_{n}$ are pairwise distinct elements of $A$ and $a$ is adjacent to each $a_i$, then the balls $B(a_ia^{-1},1)$ are pairwise disjoint, so $n \le N$. Thus $\Gamma$ has degree $\le N$, hence it can be colored using at most $N+1$ colors.

Write $A = \bigsqcup_{i=1}^{m} A_i$ with $m \le N+1$ and each $A_i$ nonempty. Let $a \ne b \in A$ be such that $d(a,b) \le 15$. There exists $i \in \{1,\ldots,m\}$ such that $a \in A_i$; since $ba^{-1} \in K \subseteq K^2$ we have that $(a,b) \in \Gamma$ so $b \not \in A_i$. If some other $x \in A_i$ is such that $d(x,b) \le 15$ then we have both $ba^{-1} \in K$ and $xb^{-1} \in K$, which implies that $xa^{-1} \in K^2$, so that $x$ and $a$ both belong to $A_i$ and $(x,a) \in \Gamma$, a contradiction. 

Denote $\varphi_i(h)=d(h,A_i)$; it belongs to $\mathrm{RUC}_b(G)$ since it is $1$-Lipschitz from $(G,d)$ to $\R$. So far we have established that for every $a, b \in A$ such that $d(a,b) \le 15$ there exists $i \in \{ 1, \ldots, m \}$ such that $d(b,A_i)=|\varphi_i(a)-\varphi_i(b)| = d(a,b)$. 

Now pick any $h \in G$. There exist $a$, $b \in A$ such that $d(a,h)< 2$ and $d(b,gh)< 2$. Since $d(h,gh)=d(1_G,g)$ we have $d(a,b) \le 4+ d(g,1_G) \le 15$; also, $d(a,b)\ge d(1_G,g)-4 \ge 5$. Since for some $i$ we have $|\varphi_i(a)-\varphi_i(b)|= d(a,b)$ and $\varphi_i$ is $1$-Lipschitz, we obtain that $|\varphi_i(h)-\varphi_i(gh)| \ge d(a,b)-4  \ge 1$ and we are done.
\end{proof}

\begin{coro}
Let $G$ be Polish locally compact. Then $G$ admits a free continuous action on a compact metrizable space.
\end{coro}

\begin{proof}
For each $g \in G \setminus \{1\}$ there exists a finite set $\Phi \subseteq \mathrm{RUC}_b(G)$ such that $\sup_{\varphi \in \Phi}|\varphi(h) - \varphi(gh)| \ge 1 $ for each $h \in G$. Letting $A$ be the closed, unital subalgebra of $\mathrm{RUC}_b(G)$ generated by $\Phi$, the associated compactification $X_A$ is metrizable (since $A$ is separable) and the $G$-flow $G \actson X_A$ is such that $g$ fixes no element of $X_A$. By compactness of $X_A$ and continuity of the action, there exists a neighborhood $U$ of $g$ such that no element $h$ of $U$ fixes an element of $X_A$. 

So we may write $G \setminus\{1\} = \bigcup_{n < \omega} U_n$ in such a way that for each $n$ there exists a metrizable flow $G \actson X_n$ such that for every $g \in U_n$ and every $x \in X_n$ we have $g \cdot x \ne x$. Then $G \actson \prod_n X_n$ is a free, metrizable $G$-flow.
\end{proof}

\bigskip

\emph{Comments.} The book \cite{deVries} contains a wealth of related material. V. Pestov was the origin of much progress in our understanding of extreme amenability, and his survey \cite{Pestov99} was very influential. A lot of information and references  may be found in his book \cite{Pestov06}, which is compulsory reading for anyone interested in extreme amenability. The proof of Veech's theorem given above comes from \cite{Kechris2005} (actually the proof given there applies to any Hausdorff locally compact topological group. Using a variant of Struble's theorem with an adequate right-invariant pseudometric would make the proof presented here go through in that case). There is also in \cite{Kechris2005} a proof that the universal minimal flow of a non compact, locally compact group is never metrizable (a fact that we will later prove differently for Polish groups).

\chapter{A theorem of Kechris, Pestov and Todor\v{c}evic}
We now have all the tools to prove a theorem of Kechris, Pestov and Todor\v{c}evic, which is now widely known as the ``KPT correspondence''. This theorem reveals a connection between topological dynamics (extreme amenability of some automorphism groups) and combinatorics, specifically Ramsey theory.

\begin{defin}
A topological space $X$ is \emph{$0$-dimensional}\index{$0$-dimensional topological space} if every point of $X$ has a neighborhood basis consisting of clopen subsets.
\end{defin}

This is the case for instance for $\omega$, $2^\omega$, $\omega^\omega$ and $\beta \omega$.

Consider a subgroup $G \le \Sinf$. For every open subgroup $V$ of $G$, note that $\ell^\infty(V \rcoset G)$, the space of all bounded functions which are invariant when multiplied on the left by an element of $V$, is a subalgebra of $\mathrm{RUC}_b(G)$ (because $V$ is a neighborhood of $1$). And each of these functions takes only countably many values, hinting that $S(G)$ has many clopen sets (an intuition which will be confirmed shortly).

\begin{prop}\label{p:many_locally_constant}
Let $G$ be a subgroup of $\Sinf$, $\varepsilon >0$, and $\varphi \in \mathrm{RUC}_b(G)$. Then there exists an open subgroup $V$ of $G$ as well as $\psi \in \ell^\infty(V \rcoset G)$ which takes only finitely many values and is such that $\|\varphi - \psi\| \le \varepsilon$.
\end{prop}

\begin{proof}
Since $\varphi \in \mathrm{RUC}_b(G)$, there exists a neighborhood $V$ of $1_G$ such that $|\varphi(vg) - \varphi(g)| \le \varepsilon$ for every $g \in G$ and every $v \in V$, and we may as well assume that $V$ is an open subgroup of $G$ since $G$ is nonarchimedean.

Since $\varphi(G)$ is bounded in $\R$, there exist $z_1,\ldots,z_p \in \R$ such that $\varphi(G) \subseteq \bigcup_i B(z_i,\varepsilon)$.

Find a (at most) countable family $(g_i)_{i \in I}$ in $G$ such that $G= \bigsqcup_i V g_i$. For each $i \in I$ find some $k_i \in \lset 1,\ldots p\rset$ such that $|\varphi(g_i) - z_{k_i}| \le \varepsilon$, then set $\psi(vg_i)= z_{k_i}$ for every $v \in V$ and every $i \in I$. Then $\psi \in \ell^\infty(V \rcoset G)$. Given $g \in G$, write$g=v g_i$ for some $v \in V$ and $i \in I$ and observe that
\[|\varphi(g) - \psi(g)|=|\varphi(vg_i) - \psi(vg_i)| \le |\varphi(vg_i) - \varphi(g_i)| + |\varphi(g_i) - z_{k_i}|\le 2 \varepsilon .\]
Hence $\|\varphi - \psi\| \le 2 \varepsilon$ and we are done.
\end{proof}

One could give an alternative, and more conceptual, argument to prove the previous proposition, by applying the Stone--Weierstrass theorem to the Samuel compactification of $G$.

\begin{prop}[Pestov]
Let $G$ be a subgroup of $\Sinf$. Then $S(G)$ is $0$-dimensional (hence $M(G)$ is also $0$-dimensional).
\end{prop}

\begin{proof}
Let $p \ne q \in S(G)$. There exists $\varphi \in C(S(G))$ such that $\varphi(p) \ne \varphi(q)$.

Let $|\varphi(p) - \varphi(q)|=2 \varepsilon$. The restriction of $\varphi$ to $G$ belongs to $\mathrm{RUC}_b(G)$; by the previous lemma there exists $\psi \in \mathrm{RUC}_b(G)$ which takes only finitely many values on $G$ and such that $\|\varphi - \psi\|\le \varepsilon$.

Viewing $\psi$ as a continuous function on $S(G)$, it still takes only finitely many values, and $\psi(p) \ne \psi(q)$. Hence $\psi^{-1}(\lset \psi(p) \rset)$ and $\psi^{-1}(\lset \psi(q) \rset)$ are disjoint clopen sets containing $p$, $q$ respectively.

Now, let $U$ be an open set containing $p$. For every $q \in S(G) \setminus U$ we have some disjoint clopen $U_q \ni p$, $V_q \ni q$. By compactness of $S(G) \setminus U$ there exist $q_1,\ldots,q_n$ such that $\bigcup_i V_{q_i}$ contains $S(G) \setminus U$, so $\bigcap_i U_{q_i}$ is a clopen neighborhood of $p$ contained in $U$.
\end{proof}

Next, we identify a family of flows which suffice to understand whether $G$ is extremely amenable.

\begin{defin}
Let $G$ be a subgroup of $\Sinf$. For every open subgroup $V$ in $G$, we denote by $X_V$ the flow obtained by having $G$ act on $2^{V \rcoset G}$ by setting 
\[ g \cdot \varphi(Vh)= \varphi(Vhg). \]
\end{defin}

\begin{thm}[Kechris--Pestov--Todor\v{c}evic]\label{t:flows_giving_extreme_amenability}
Let $G$ be a subgroup of $\Sinf$. Then $G$ is extremely amenable iff any minimal subflow of any $X_V$ is trivial.
\end{thm}

\begin{proof}
One implication is immediate (if $G$ is extremely amenable any minimal flow is trivial).

Conversely, assume that $M(G)$ is nontrivial, and let $D$ be a nontrivial clopen subset of $M(G)$.

Then $V= \lset g \in G : g D = D\rset$ is a subgroup of $G$. For every $p \in D$ there exists a neighborhood $U_p$ of $1_G$ and an open $W_p \ni p$ such that $U_p W_p \subseteq D$. By compactness we have $D \subseteq \bigcup_{i=1}^n W_{p_i}$. Since $\bigcap_i U_{p_i}$ is a neighborhood of $1_G$ it contains an open subgroup $U$, and for every $u \in U$ we have $u D \subseteq D$, hence $uD=D$ for all $u \in U$. So $U \le V$, proving that $V$ is open.

Now we consider $\pi \colon M(G) \to 2^{V \rcoset G}$ defined by 
\[ \pi(p)(Vg) = 1 \Leftrightarrow g p \in D . \]
This map is well-defined since $VD=D$, and continuous since $D$ is clopen and $G \actson M(G)$ is continuous so each $g^{-1}D$ is clopen.

For every $g,h \in G$ and $p \in M(G)$ we have
\[
\pi(gp)(Vh)= 1 \Leftrightarrow h g p \in D \Leftrightarrow \pi(p)(Vhg) = 1 \Leftrightarrow (g \pi(p))(Vh)=1 .\]
Hence $\pi$ is $G$-equivariant, so $\pi(M(G))$ is a subflow of $X_V$. For $p \in D$ we have $\pi(p)(V)=1$ while $\pi(p)(V)=0$ for $p \not \in D$, hence $\pi(M(G))$ is nontrivial.
\end{proof}

It is probably worth pointing out that one can push this argument a little further to prove that, for a nonarchimedian group $G$, the Samuel compactification $S(G)$ is the inverse limit of the family of Stone-\v{C}ech compactifications $\beta(V \rcoset G)$, where $V$ belongs to the set of open subgroups of $G$, ordered by reverse inclusion. See \cite{Pestov98} and \cite{Zucker16}.

The only fixed points for $G \actson 2^{V \rcoset G}$ are the constant functions $0$ and $1$, so the above criterion amounts to saying that for every $V$ and every $c \in 2^{V \rcoset G}$ the closure of $Gc$ contains a constant function. Hence we obtain the following:

\begin{thm}[Kechris--Pestov--Todor\v{c}evic]\label{thm:KPT1}
Let $G$ be a subgroup of $\Sinf$. Then $G$ is extremely amenable iff for every open subgroup of $G$ and every $c \in 2^{V \rcoset G}$ the following holds:
\[\forall F \subseteq G \ \textrm{\emph{finite}} \ \exists g \in G \ \forall f_1, f_2 \in F \quad c(Vf_1g)= c(Vf_2g) .\]
\end{thm}

\begin{exo}\label{exo58}
Provide the details of the proof of the above statement.
\end{exo}

\begin{exo}\label{exo59}
Prove that in the theorem above one may equivalently consider $k^{V \rcoset G}$ for any $k \ge 2$, and that one also obtains an equivalent statement by letting $V$ run over some fixed family of open subgroups forming a basis of neighborhoods of $1_G$.
\end{exo}

Now we turn to a combinatorial interpretation of this property, when $G$ is the automorphism group of some relational ultrahomogeneous structure (we recall that every closed subgroup of $\Sinf$ is of this form). So, fix a relational Fraïssé class $\mcK$ with elements of arbitrarily high (finite) cardinality, denote by $\bF$ its limit and let $G=\Aut(\bF)$.

Given a finite $A \subseteq F$ we denote by $G_{A}$ the pointwise stabilizer of $A$. Using ultrahomogeneity we may identify $G/G_{A}$ with the set of all embeddings of $\bA$ into $\bF$.

\begin{defin}
Given $\bA$, $\bB \in \mcK$ we denote $\emb{A}{B}$\index{$\emb{A}{B}$} the (finite) set of embeddings of $\bA$ in $\bB$. For $\bA \in \mcK$ we similarly define $\emb{A}{F}$ as the (infinite) set of all embeddings from $\bA$ to $\bF$.
\end{defin}
Here we should warn about differing choices of notations in the literature: it is more common to denote by $\emb{A}{B}$ the set of all \emph{substructures} of $\bB$ isomorphic to $\bA$. As soon as some element of $\mcK$ has a nontrivial automorphism group, the two definitions differ. This distinction will be immaterial for us, see Proposition \ref{p:Ramsey_for_embeddings_implies_rigid} and Theorem \ref{thm:KPT}.

The group $G$ acts on $\emb{A}{F}$ by composition, via $g \cdot \alpha (a) = g(\alpha(a))$. The ultrahomogeneity of $\bF$ amounts to the statement that $\begin{cases} G \lcoset G_{A}  \to \emb{{A}}{F} & \\ g G_{A} \mapsto g_{|A} & \end{cases}$ is an isomorphism of $G$-spaces.

Note that the action $G \actson 2^{G \lcoset G_{A}}$ defined by $g \cdot x(h G_A)= x(g^{-1} h G_A)$ is isomorphic to $G \actson 2^{G_A \rcoset G}$ (by taking inverses), and we are going to work with the first action rather than the second.

\begin{defin}
Let $\bA \in \mcK$. A \emph{coloring}\index{coloring} of $\emb{A}{F}$ is a map $\gamma \colon \emb{A}{F} \to k$, where $k < \omega$. We similarly define colorings of $\emb{A}{B}$ for $\bA$, $\bB \in \mcK$.
\end{defin} 

So $2$-colorings of $\emb{A}{F}$ are simply elements of $2^{\emb{A}{F}}$, which we already met under the guise of elements of the flow $X_{G_A}$.

\begin{defin}
A class $\mcK$ of $\mcL$-structures has the \emph{Ramsey property (for embeddings)}\index{Ramsey property for embeddings} if for any $k$ and any $\bA$, $\bB \in \mcK$ there exists $\bC \in \mcK$ such that for any $k$-coloring $\gamma$ of $\emb{A}{C}$ there exists $\beta \in \emb{B}{C}$ such that $\gamma$ is constant on $\beta \circ \emb{A}{B}$.
\end{defin}

\begin{exo}\label{exo60}
Prove that the Ramsey property is equivalent to the property above where one only considers $2$-colorings.
\end{exo}

We could have similarly defined the Ramsey property for substructures, by asking that for any $\bA$, $\bB \in \mcK$ there exists $\bC \in \mcK$ such that any $k$-coloring of the copies of $\bA$ contained in $\bC$ is constant on some copy of $\bB$ (i.e., takes the same value on all copies of $\bA$ contained in some fixed copy of $\bB$).

\begin{prop}\label{p:Ramsey_for_embeddings_implies_rigid}
A class $\mcK$ has the Ramsey property for embeddings iff it has the Ramsey property for substructures and every $\bA \in \mcK$ has a trivial automorphism group (we say that elements of $\mcK$ are \emph{rigid}\index{rigid structure}).
\end{prop}

\begin{proof}
Note that if every element of $\mcK$ has a trivial automorphism group then coloring embeddings is the same thing as coloring substructures, since every embedding of some $\bA$ into $\bC$ is in that case uniquely characterized by its image.

So the above proposition really amounts to the statement that, if some $\bA \in \mcK$ has a nontrivial automorphism group, then $\mcK$ cannot have the Ramsey property for embeddings.

Assume that some $\bA \in \mcK$ has a nontrivial automorphism group then let $\bA = \bB$ in the definition of the Ramsey property (for embeddings). We want to prove that no $\bC$ in $\mcK$ can witness the Ramsey property for the pair
$(\bA,\bA)$, and a well-chosen set of colors which we introduce below. If $\bC$ does not contain a copy of $\bA$ this is immediate, so we may assume that $\bC$ contains $\bA$. 

For each $\bA' \le \bC$ isomorphic to $\bA$, we choose an isomorphism $g_{\bA'} \colon \bA' \to \bA$ and for $\alpha \in \emb{A}{C}$ we set $\gamma (\alpha)= g_{\alpha(\bA)} \circ \alpha$ (so our set of colors is $\Aut(\bA)$, which is finite and does not depend on $\bC$).

For every $\beta \in \emb{A}{C}$ we have $\beta \circ \emb{A}{A} = \beta \circ \Aut(\bA)$, and $\gamma$ takes all possible values on this set, so $\bC$ cannot witness the Ramsey property for embeddings, hence this property fails.
\end{proof}

The Ramsey property is related to the amalgamation property, as witnessed by the result of the next exercise.

\begin{exo}[Ne\v{s}et\v{r}il]\label{exo61}
Assume that $\mcK$ satisfies the joint embedding property as well as the Ramsey property for embeddings (so elements of $\mcK$ have trivial automorphism group).

Fix $\bA$, $\bB$, $\bC \in \mcK$ and embeddings $\alpha \colon \bA \to \bB$, $\beta \colon \bA \to \bC$. Find $\bE \in \mcK$ containing a copy of both $\bB$ and $\bC$ and view $\alpha$, $\beta$ as maps from $A$ to $E$.

Using an appropriate coloring with values in $2^{\{B,C\}}$, and some $\bD$ witnessing the Ramsey property for colorings of $\emb{A}{E}$ using $4$ colors, prove that $\mcK$ satisfies the amalgamation property.
\end{exo}

Next, we give an equivalent ``infinite'' formulation of the Ramsey property.

\begin{prop}
Let $\mcK$ be a Fraïssé class of rigid structures with limit $\bF$. Then $\mcK$ has the Ramsey property iff for any $\bA$, $\bB \in \mcK$, any integer $k$ and any $k$-coloring $\gamma$ of $\emb{A}{F}$ there exists $\beta \in \emb{B}{F}$ such that $\gamma$ is constant on $\beta \circ \emb{A}{B}$.
\end{prop}

\begin{proof}
Assume first that $\mcK$ has the Ramsey property. Fix $\bA$, $\bB \in \mcK$. If $\emb{A}{B}$ is empty we have nothing to prove, so we assume that $\bA \le \bB$. We find $\bC$ witnessing that the Ramsey property holds for $(\bA,\bB)$ and we realize $\bC$ as a substructure of $\bF$. 

Let $\gamma$ be a $k$-coloring of $\emb{A}{F}$. Since $C \subset F$, any element of $\emb{A}{C}$ can be seen as an element of $\emb{A}{F}$. So $\gamma$ induces a coloring of $\emb{A}{C}$, thus there exists $\beta \in \emb{B}{C}$ such that $\gamma$ is constant on $\beta \circ \emb{A}{B}$. Viewing $\beta$ as an element of $\emb{B}{F}$, we are done.

Conversely, assume that $\mcK$ does not have the Ramsey property and fix $\bA \le \bB$ witnessing that failure. We may assume that $\bA$, $\bB$ are substructures of $\bF$. 

For every substructure $\bC$ of $\bF$ in which $\bB$ embeds we fix a bad $2$-coloring $\gamma_{\bC}$ of $\emb{A}{C}$, i.e., a coloring which is not constant on $\beta \circ \emb{A}{B}$ for any $\beta \in \emb{B}{C}$. 
Next, fix an ultrafilter $\mcU$ on the set of finite subsets of $F$ which is such that for any finite $X$ the set $\lset Y : X \subseteq Y \rset$ belongs to $\mcU$.

We define a $2$-coloring $\gamma$ of $\emb{A}{F}$ by setting $\gamma(\alpha) = \lim_{\mcU} \gamma_{\bC}(\alpha)$.

Equivalently, $\gamma(\alpha)=\varepsilon \in \lset 0,1 \rset$ iff $\lset C : \gamma_{\bC}(\alpha)=\varepsilon \rset \in \mcU$; of course $\gamma_{\bC}(\alpha)$ is not defined for all $C$, but it is defined as soon as $B \cup \alpha(A) \subseteq C$, and $\lset C : B \cup \alpha(A)\subseteq C \rset \in \mcU$.

We now check that $\gamma$ is a bad coloring of $\emb{A}{F}$: let $\beta \in \emb{B}{F}$. Then $\lset C : \beta(B) \subseteq C \rset \in \mcU$, and for any $C$ in this set we have $\alpha_1$, $\alpha_2 \in \emb{A}{B}$ such that $\gamma_{\bC}(\beta \circ \alpha_1)=0$ and $\gamma_{\bC}(\beta \circ \alpha_2)=1$. 

Since $\emb{A}{B}$ is finite, there exist $\alpha_1$, $\alpha_2 \in \emb{A}{B}$ such that $\lset C : \gamma_{\bC}(\beta \circ \alpha_1)=0 \rset \in \mcU$ and $\lset C : \gamma_{\bC}(\beta \circ \alpha_2)=1 \rset \in \mcU$. Thus $\gamma(\beta \circ \alpha_1)=0 \ne \gamma(\beta \circ \alpha_2)$, and we are done.
\end{proof}

\begin{thm}[Kechris--Pestov--Todor\v{c}evic]\label{thm:KPT}
Let $\mcK$ be a Fraïssé class of relational structures with infinite Fraïssé limit $\bF$. Then $G=\Aut(\bF)$ is extremely amenable iff $\mcK$ has the Ramsey property (for embeddings).
\end{thm}

\begin{proof}
Translating Theorem \ref{thm:KPT1} in terms of $\emb{\bar a}{F}$ instead of $G_{\bar a} \rcoset G$, we obtain that $G$ is extremely amenable iff
\[\forall H \subseteq G \text{ finite } \forall A \subseteq F \text{ finite } \ \forall \gamma \colon \emb{A}{F} \to 2 \ \exists g \in G \text{ such that } \gamma \text{ is constant on } g \circ H_{|\bA} .\]

Assume $\mcK$ has the Ramsey property, then fix a finite $H \subseteq G$, a finite $A \subseteq F$ and $\gamma \colon \emb{A}{F} \to 2$. Let $B = \bigcup_{h \in H} h A$. Using the Ramsey property we obtain $\beta \in \emb{B}{F}$ such that $\gamma$ is constant on $\beta \circ \emb{A}{B}$. Extending $\beta$ to some $g \in \Aut(\bF)$ we are done, since $\emb{A}{B}$ contains $H_{|A}$.

Conversely, assume that $G$ satisfies the property above, and fix $\bA \le \bB \in \mcK$. We may assume that $A \subseteq B \subseteq F$. Then every element $\alpha$ of $\emb{A}{B}$ extends to some $g_{\alpha} \in G$. 

Let $H=\lset g_\alpha : \alpha \in \emb{A}{B} \rset$. Let $\gamma$ be a $2$-coloring of $\emb{A}{F}$; by assumption we obtain some $g \in G$ such that $\gamma$ is constant on $g \circ H_{|A}$. Letting $\beta= g_{|B}$ we obtain that $\gamma$ is constant on $\beta \circ \emb{A}{B}$, witnessing that $\mcK$ satisfies the Ramsey property.
\end{proof}

It is certainly time to discuss an example.

\begin{thm}[Ramsey] The class of all finite linear orders has the Ramsey property.
\end{thm}

\begin{proof}
We claim that it is enough to prove the following classical (infinite version of) Ramsey theorem: for any integer $n$ and any $2$-coloring $\gamma$ of the set $\omega^{[n]}$ of $n$-element subsets of $\omega$, there exists an infinite subset $I$ of $\omega$ such that $\gamma$ is constant on $I^{[n]}$.

Indeed, given a finite subset $A \subseteq \Q$ of cardinality $n$, a $2$-coloring of $\emb{A}{\Q}$ is the same thing as a coloring of $\Q^{[n]} \cong \omega^{[n]}$. Since an infinite $I$ obviously contains finite sets of any cardinality, the property above gives us an $I$ witnessing the Ramsey property. 

We now prove the classical Ramsey theorem by induction on $n$. 
For $n=1$ the statement above follows from the pigeonhole principle. Assume that we have proved the result up to   
some $n$; let $m=n+1$ and $\gamma$ be a $2$-coloring of $\omega^{[m]}$. 

For each $a < \omega$ we let $\gamma_a \colon (\omega \setminus \lset a \rset))^{[n]}  \to 2$ be defined by $\gamma_a(X)= \gamma(X \cup \lset a \rset)$. 

By the induction hypothesis, there exists an infinite $I_0 \subseteq \omega \setminus \lset 0\rset$ such that $\gamma_0$ is constant on $I_0^{[n]}$. Let $i_1= \min I_0$, and apply the same argument to find an infinite $I_1 \subseteq I_0 \setminus \lset i_1 \rset$ such that $\gamma_{i_1}$ is constant on $I_1^{[n]}$. Set $i_2= \min I_1$ and keep going.

For $i \in I$, denote by $\varepsilon_i$ the constant color taken by $\gamma_i$ on $I_{i}^{[n]}$. By the pigeonhole principle, there is $\varepsilon \in \{0,1\}$ and $J \subseteq I$ infinite such that $\varepsilon_j= \varepsilon$ for all $j \in J$.

Let $j_0 < \ldots < j_n$ enumerate a subset of $J$ of cardinality $m$. Then $\lset j_1,\ldots,j_n \rset$ is a $n$-element subset of $I_{j_0}$, whence $\gamma( \lset j_0,\ldots,j_m \rset) = \gamma_{j_0}(\lset j_1,\ldots,j_m \rset) = \varepsilon$. Thus $\gamma$ is monochromatic on $J^{[m]}$ and we are done.
\end{proof}

The next result (which was a precursor of the Kechris--Pestov--Todor\v{c}evic correspondence) is now an immediate consequence.

\begin{coro}[Pestov] $\Aut(\Q,<)$ is extremely amenable.
\end{coro}

\begin{exo}[Pestov]\label{exo62}
Let $H$ be the group of all homeomorphisms of $[0,1]$, endowed with the topology of uniform convergence.
\begin{enumerate}[(i)]
\item Let $H_+= \lset g \in H : g(0)=0 \rset$. Prove that $H_+$ is extremely amenable.

(Hint: map $\Aut(\Q,<)$ densely and continuously in $H_+$)
\item Use this to compute the universal minimal flow of $H$.
\end{enumerate}
\end{exo}

We conclude this chapter by explaining a general strategy to compute universal minimal flows of Polish groups, in the particular case when they happen to be metrizable (which never happens in the locally compact non compact case, see Exercise \ref{exo:non_metrizable_UMF_loc_cpact} and the comments at the end of this chapter).

\begin{defin}
Let $G$ be a topological group and $H$ be a subgroup of $G$. We endow $G \lcoset H$ with the uniformity $\mcU$ coming from the right uniformity, i.e., the uniformity generated by entourages of the form 
\[\lset (fH,ufH) : u \in U \rset\]
for $U$ an open neighborhood of $1_G$.

We say that $H$ is \emph{co-precompact}\index{co-precompact subgroup} in $G$ if the completion $\yhwidehat{(G \lcoset H, \mcU)}$ is compact.
\end{defin}

\begin{exo}\label{exo63}
Prove that $H$ is co-precompact in $G$ iff for every nonempty open $V$ there exists a finite $F$ such that $VFH=G$.
\end{exo}

\begin{exo}\label{exo:co-precompact=oligomorphic}
Let $H$ be a subgroup of $\Sinf$. Prove that $H$ is co-precompact in $\Sinf$ if, and only if, $H$ is oligomorphic.
\end{exo}

\begin{prop}
Let $G$ be a Polish group and $H$ a subgroup of $G$. The left translation action of $G$ on $G\lcoset H$ extends to a continuous action $G \actson \yhwidehat{G \lcoset H}$.
\end{prop}

We should probably recall here that the quotient uniformity on $G \lcoset H$ is metrizable (see Definition \ref{def:quotient_metric}) hence $\yhwidehat{G \lcoset H}$ is also metrizable. This enables us to use sequences in the argument below.

\begin{proof}
It follows from the proof of Theorem \ref{t:extending_group_operations} that if $(g_n)_{n < \omega}$ and $(k_n)_{n < \omega}$ are $\mcU_r$-Cauchy then $(g_n k_n)_{n < \omega}$ is $\mcU_r$-Cauchy. Assume $(g_n)_{n < \omega}$ is Cauchy in $(G,\mcU_r)$ and $(k_nH)_{n < \omega}$ is Cauchy in $(G \lcoset H, \mcU)$. 

Then there exists a subsequence $(k_{\varphi(n)})_{n< \omega}$ and $(h_n)_{n< \omega} \in H^\omega$ such that $(k_{\varphi(n)}h_n)_{n< \omega}$ is Cauchy in $(G,\mcU_r)$, whence $(g_{\varphi(n)} k_{\varphi(n)} h_n)_{n< \omega}$ is Cauchy in $(G,\mcU_r)$. This implies that $(g_{\varphi(n)} k_{\varphi(n)} H)_{n < \omega}$ is Cauchy in $(G\lcoset H,\mcU)$.

It follows that every subsequence of $(g_n k_n H)_{n < \omega}$ admits a Cauchy subsequence, which implies that $(g_n k_n H)_{n < \omega}$ is Cauchy. This proves that $(g,kH) \mapsto gkH$ extends continuously to $\yhwidehat{(G,\mcU_r)} \times \yhwidehat{(G \lcoset H, \mcU)}$, and that is more than we needed.
\end{proof}

Let us discuss an instructive example. Endow the space of linear orderings $\mathrm{LO}$ with the topology coming from viewing it as a subset of $2^{\omega \times \omega}$, and have $\Sinf$ act on $\mathrm{LO}$ via 
\[ i \, (\sigma \cdot \prec) \, j \Leftrightarrow \sigma^{-1}(i) \prec \sigma^{-1}(j) . \]
Note that $\mathrm{LO}$ is compact since it is closed in $2^{\omega \times \omega}$.

\begin{lem}
The flow $\Sinf \actson \mathrm{LO}$ is minimal.
\end{lem}

\begin{proof}
Let $U$ be a basic open subset of $\mathrm{LO}$, which we can assume to be made up of all orders $\prec$ such that $i_0 \prec i_ 1 \prec \ldots \prec i_n$ for some $(i_0,\ldots,i_n) \in \omega^{n+1}$. 

Let $I=\{i_0,\ldots,i_n\}$ and let $F$ be the finite set of all elements of $\Sinf$ whose support is contained in $I$.

For any ordering $\prec$ of $\omega$, there exists a bijection $\sigma$ of $I$ such that 
$i_0 (\sigma \cdot \prec) i_1 \ldots (\sigma \cdot \prec) i_n$ (map the smallest element of $(I,\prec)$ to $i_0$, the second smallest to $i_1$, and so on). We may extend $\sigma$ to an element of $F$, and we obtain that $\prec \in FU$. 

This proves that $\mathrm{LO}=FU$, so $\Sinf \actson \mathrm{LO}$ is minimal as promised.
\end{proof}

\begin{lem}\label{exo:univ_minimal_flow_Sinf}
The flow $\Sinf \actson \mathrm{LO}$ is isomorphic to $\Sinf \actson \yhwidehat{\Sinf \lcoset \Aut(\Q)}$.
\end{lem}

In particular, $\Aut(\Q)$ is co-precompact in $\Sinf$ (of course this can also be seen directly, for instance it is a consequence of the result of Exercise \ref{exo:co-precompact=oligomorphic}).

\begin{proof}
Fix an ordering $\prec$ of $\omega$ such that $(\omega,\prec)$ is dense and without endpoints. Denote by $H$ the stabilizer of $\prec$, which is conjugate to $\Aut(\Q)$ in $\Sinf$. 

Let $V$ be the pointwise stabilizer of a finite $F \subset \omega$. Let $<_1$, $<_2$ be two elements of $\Sinf \cdot \prec$. 

There exists $\sigma \in \Sinf$ and $\tau \in V$ such that $(<_1,<_2)=(\sigma \cdot \prec, \tau \sigma \cdot \prec)$ iff there exists $\tau \in V$ such that $<_1 = \tau \cdot <_2$, which implies that for any $i,j \in F$ we have $(i<_1 j) \Leftrightarrow (i <_2 j)$. 

Conversely, if ${<_1}_{|F}={<_2}_{|F}$ then, since $<_1$ and $<_2$ are dense and without endpoints there exists $\tau \in \Sinf$ which is the identity on $F$ and such that $<_1=\tau \cdot <_2$.

If we denote by $E_F$ the basic entourage $\lset (\sigma H, \tau \sigma H) : \tau \in V \rset $, and by $\Delta_F$ the neighborhood of the diagonal $\lset (<_1,<_2) : {<_1}_{|F}={<_2}_{|F} \rset$ we just proved that 
\[\lset (\sigma \cdot \prec, \tau \cdot \prec) : (\sigma H, \tau H) \in E_F \rset = \Delta_F .\]
Since the sets $E_F$ and $\Delta_F$ are fundamental systems for the two uniformities under consideration, this proves that $\sigma H \mapsto \sigma \cdot \prec$ is a uniform isomorphism from $\Sinf \lcoset H$ to $\Sinf \cdot \prec$. 

This uniform isomorphism extends to a uniform isomorphism from $\yhwidehat{\Sinf \lcoset H}$ to $\mathrm{LO}$, which is $\Sinf$-equivariant by continuity of the actions and $\Sinf$-equivariance on a dense subset.
\end{proof}

\begin{thm}
Assume that $G$ is a Polish group and $H$ is a closed extremely amenable subgroup. Then $G \actson \yhwidehat{G \lcoset H}$ factors onto every minimal $G$-flow.

In particular, if $H$ is co-precompact, extremely amenable and $G \actson \yhwidehat{G \lcoset H}$ is minimal, then it is the universal minimal flow $M(G)$.
\end{thm}

\begin{proof}
Let $G \actson X$ be a minimal $G$-flow. Then $H \actson X$ is an $H$-flow, hence has a fixed point $x_0$. Since $g \mapsto gx_0$ is right uniformly continuous, $gH \mapsto gx_0$ is $\mcU$-uniformly continuous, so it extends to a continuous map $\pi \colon \yhwidehat{G \lcoset H} \to X$.

For every $kH \in G \lcoset H$ and every $g \in G$ we have $\pi(gkH)= gk x_0 = g \pi(kH)$. By continuity we obtain $\pi(gy)= g \pi(y)$ for every $y \in \yhwidehat{G \lcoset H}$. Hence $\pi$ is $G$-equivariant, and it is surjective since $G \actson X$ is minimal.
\end{proof}

Note that it also follows from the previous result that if $H$ is co-precompact and extremely amenable then every minimal subflow of $\yhwidehat{G \lcoset H}$ is the universal minimal flow of $G$; and it turns out that if one can find such an $H$ then one can also find one such that $G \actson \yhwidehat{G \lcoset H}$ is minimal, enabling one to compute the universal minimal flow of $G$ (see theorem \ref{t:MNTBMT} below). 

The previous result, combined with Lemma \ref{exo:univ_minimal_flow_Sinf}, enables us to compute a new universal minimal flow.

\begin{thm}[Glasner--Weiss]
The universal minimal flow of $\Sinf$ is $\Sinf \actson \mathrm{LO}$.
\end{thm}

\begin{proof}
We know that $\Aut(\Q)$ is extremely amenable and co-precompact in $\Sinf$. Since $\Sinf \actson \yhwidehat{\Sinf \lcoset \Aut(\Q)}$ is isomorphic to $\Sinf \actson \mathrm{LO}$, which is minimal, we have identified the universal minimal flow of $\Sinf$. 
\end{proof}

In particular, the universal minimal flow of $\Sinf$ is metrizable. Note that it has a comeager orbit, consisting of all dense linear orders without endpoints. 
This result actually fits into a broader picture, as the following theorem (which we prove in the next chapter) shows.

\begin{thm}[Melleray--Nguyen van Thé--Tsankov; Ben Yaacov--Melleray--Tsankov]\label{t:MNTBMT}
Let $G$ be a Polish group. Then the following conditions are equivalent:
\begin{enumerate}[(i)]
\item\label{t:metrizable_UMF_1} The universal minimal flow $M(G)$ is metrizable.
\item\label{t:metrizable_UMF_2} There exists a co-precompact, extremely amenable subgroup $H \le G$ such that $M(G)=\yhwidehat{G \lcoset H}$.
\end{enumerate}
\end{thm}

We already discussed the implication \ref{t:metrizable_UMF_2} $\Rightarrow$ \ref{t:metrizable_UMF_1}. The new and interesting thing here is that the only way for a Polish group $G$ to have a metrizable universal minimal flow is for $G$ to contain a large extremely amenable subgroup. 

Replacing $H$ by its closure we may assume that $H$ above is closed; then $G\lcoset H$ is a dense $G_\delta$ orbit in $M(G)$ (since $G/H$ is a Polish subspace of $\yhwidehat{G/H}$). Thus a metrizable minimal flow of a Polish group always has a comeager orbit, and the stabilizer of a point in the comeager orbit must be extremely amenable.

\begin{exo}\label{exo:non_metrizable_UMF_loc_cpact}
Deduce from the previous theorem that the universal minimal flow of a non-compact locally compact Polish group is never metrizable.

(This result was known long before Theorem \ref{t:MNTBMT} was proved, and actually holds for every locally compact topological group, see \cite{Kechris2005})
\end{exo}

\bigskip

\emph{Comments.} The main theorem of this chapter comes from the paper \cite{Kechris2005} of Kechris--Pestov--Todor\v{c}evic. The last result is a combination of theorems obtained in two papers, \cite{MellerayNVTTsankov} and \cite{BenYaacovMellerayTsankov}. Another approach to this was proposed in \cite{Zucker17} by A. Zucker (who had obtained an earlier and different proof for subgroups of $\Sinf$, see \cite{Zucker16}). Recent work of Zucker and co-authors has brought significant progress, for instance the reader (and the author, if we are being honest) could do worse than spending some time with \cite{Zucker19}. 

We should probably note that Ramsey theory is not only concerned with finding monochromatic sets of copies of a given structure. For instance, given a Fraïssé class $\mcK$ with limit $\bF$, it may happen that for any $\bA \in \mcK$ there exists an integer $i_{\bA}$ such that for any coloring $c$ of $\emb{A}{F}$ there exists $\alpha \in \emb{F}{F}$ such that $c$ takes at most $i_{\bA}$ distinct values on $\alpha \circ \emb{A}{F}$. When that happens one says that $\mcK$ has \emph{finite big Ramsey degrees}; this is the focus of a lot of contemporary research, see \cite{Hubicka_Zucker} for a survey on this and related notions.

\chapter[Polish groups with a metrizable Universal minimal flow]{Polish groups with a metrizable universal minimal flow}

We turn to the proof of Theorem \ref{t:MNTBMT}. The proof splits in two parts: first, we prove that if $M(G)$ is metrizable then it has a comeager orbit; then we show that if $M(G)$ is metrizable and contains a comeager orbit there is a co-precompact, extremely amenable closed subgroup $H \le G$ such that $M(G) =\yhwidehat{G/H}$. 

\begin{defin}
  A \emph{compact topometric space}\index{topometric space} is a triple $(Z, \tau, \partial)$, where $Z$ is a set, $\tau$ is a compact Hausdorff topology on $Z$, and $\partial$ is a metric on $Z$ such that the following conditions are satisfied:
  \begin{itemize}
  \item the $\partial$-topology refines $\tau$;
  \item $\partial$ is $\tau$-lower semicontinuous, i.e., the set $\{(a,b) \in Z^2 : \partial(a,b) \le r \}$ is $\tau$-closed for every $r \geq 0$.
  \end{itemize}
\end{defin}

\begin{lem}
Let $(Z, \tau, \partial)$ be a compact topometric space. Then $\partial$ is complete.
\end{lem}
\begin{proof}
  Let $(z_n)$ be a Cauchy sequence; for all $n$, define $r_n= \sup \{\partial(z_n,z_m) : m \ge n\}$.
Then $r_n$ converges to $0$. Let $F_n$ denote the closed ball of radius $r_n$ centered at $z_n$. Each $F_n$ is $\tau$-closed, hence compact, and since $F_n$ contains $z_m$ for all $m \ge n$, this family has the finite intersection property. By compactness, $\bigcap_{n \in \bN} F_n$ is non-empty; it must be a singleton, which is the $\partial$-limit of the sequence $(z_n)$.
\end{proof}

We now define a topometric structure on the Samuel compactification $S(X)$ of a bounded metric space $(X,d)$ (i.e., the compactification associated to the algebra $\ucb(X)$): we endow $S(X)$ with its usual compact topology and define $\partial$ by:
\[
  \partial(a,b)= \sup \{|f(a)-f(b)| : f \in \cli(X) \},
\]
where $\cli(X)$ denotes the set of all bounded $1$-Lipschitz functions on $(X,d)$. To make sense of this, one must recall that functions in $\cli(X)$, being bounded and uniformly continuous, uniquely extend to continuous functions on $S(X)$.

Clearly, each $\{(a,b) : \partial(a,b) \le r\}$ is $\tau$-closed. To see why the $\partial$-topology refines $\tau$, recall that any function in $\ucb(X)$ is a uniform limit of Lipschitz maps, so $\tau$ has a basis of open sets of the form 
\[
  \{a \in S(X) : f_1(a) \in I_1, \ldots, f_n(a) \in I_n \},
\]
where each $f_j$ belongs to $\cli(X)$ and each $I_j$ is an open interval. Each set in this basis is $\partial$-open.

Note that $\partial$ and $d$ coincide on $X$, and that elements of $\cli(X)$ extend to maps on $S(X)$ which are both $\tau$-continuous and $\partial$-$1$-Lipschitz (the Lipschitz part is immediate from the definition of $\partial$). Also, if $d$ is the discrete $0$--$1$ metric on $X$, then $\partial$ is the discrete $0$--$1$ metric on $S(X)$.

The following result is the topometric analogue of the fact that for discrete $X$ the only convergent sequences in $\beta X$ are stationary.

\begin{thm}\label{t:caraccompact}
  Let $(X, d)$ be a bounded metric space. Then every $\tau$-convergent sequence in $S(X)$ is $\partial$-convergent.
\end{thm}

In order to prove this theorem, we first establish two lemmas. The first is the topometric analogue of the fact that in the discrete case, $S(X)$ is extremally disconnected (see proposition \ref{p:extremally_disconnected}). If $A$, $B$ are two subsets of a metric space $(Z, d)$, we denote
\begin{equation*}
  d(A, B) = \inf \{d(a, b) : a \in A, b \in B \}.
\end{equation*}
\begin{lem}\label{l:caraccompact1}
Let $U,V$ be nonempty $\tau$-open subsets of $S(X)$. Then we have
\[
  \partial(\cl[\tau]{U}, \cl[\tau]{V})= \partial(U, V)= d(U \cap X, V \cap X).
\]
\end{lem}
\begin{proof}
First note that, since $X$ is dense in $S(X)$, we have $\cl[\tau]{U}= \cl[\tau]{U \cap X}$ and $\cl[\tau]{V}= \cl[\tau]{V \cap X}$. Consider the function $f \in \cli(X)$ defined by $f(x) = d(x, U \cap X)$. Since it is $1$-Lipschitz, it extends to a $\tau$-continuous, $\partial$-$1$-Lipschitz map on $S(X)$, which we still denote by $f$. 

Since $f=0$ on $U \cap X$, we must also have $f=0$ on $\cl[\tau]{U}$ by continuity; similarly, $f \ge d(V \cap X, U \cap X)$ on $V \cap X$, so $f \ge d(V \cap X,U\cap X)$ on $\cl[\tau]{V}$.
Hence $f$ witnesses the fact that $\partial(\cl[\tau]{U}, \cl[\tau]{V}) \ge d(V \cap X, U \cap X)$; since $d(U \cap X,V \cap X)$ is equal to $\partial(U \cap X, V \cap X)$ by the definition of $\partial$, this inequality must in fact be an equality, and we are done.
\end{proof}

\begin{lem}\label{l:caraccompact2}
Assume that $(a_n)$ is a sequence in $S(X)$ and $\delta > 0$ is such that $\partial(a_n, a_m) > \delta$ for all $n \neq m$. Then, for every $\varepsilon < \delta/2$, there exist a subsequence $(b_n)$ of $(a_n)$ and $\tau$-open sets $(U_n)$ such that $b_n \in U_n$ and $\partial(U_n,U_m) \ge \varepsilon$ for all $n \ne m$.
\end{lem}
\begin{proof}
Let $f \in \cli(X)$ be such that $|f(a_0)-f(a_1)| > \delta$. The triangle inequality implies that, for all $n>1$, we have  $|f(a_0)-f(a_n)|> \frac{\delta}{2}$ or $|f(a_1)-f(a_n)| > \frac{\delta}{2}$. One of those cases happens infinitely many times. Thus we see that for any such sequence $(a_n)$, there exists $i_0 \in \{0,1\}$, an infinite subset $\{i_n\}_{n \ge 1} \subseteq \bN \setminus \{0,1\}$ and $f_0 \in \cli(X)$ such that $f_0(a_{i_0})=0$ and $f_0(a_{i_n}) > \frac{\delta}{2}$ for all $n \ge 1$. Repeating this infinitely many times, we build a subsequence $(b_n)$ of $(a_n)$ and a sequence of maps $f_n \in \cli(X)$ such that $f_n(b_n)=0$ for all $n$ and $f_n(b_m) > \frac{\delta}{2}$ for all $n < m$.

Set $U_n= \{a \in S(X) : f_n(a) < \frac{\delta}{2} - \varepsilon \text{ and } f_k(a) > \frac{\delta}{2} \text{ for all } k <n \}$.
We have $b_n \in U_n$ for all $n$, and the function $f_n$ witnesses the fact that $\partial(U_n,U_m) \ge \varepsilon$ for all $n<m$.
\end{proof}

\begin{proof}[Proof of Theorem \ref{t:caraccompact}]
  Let $(a_n)$ be a $\tau$-convergent sequence in $S(X)$ with limit $a$ and suppose that it does not admit a $\partial$-Cauchy subsequence. 
  
Passing to a subsequence if necessary, we may assume that there exists $\delta > 0$ such that $\partial(a_n, a_m) > \delta$ for all $n \neq m$ and we can apply Lemma~\ref{l:caraccompact2} to obtain a further subsequence $(b_n)$ of $(a_n)$ and $\tau$-open subsets $U_n$ of $S(X)$ such that $b_n \in U_n$ and $\partial(U_n,U_m) \ge \delta/2$ for all $n \ne m$. Let
\[
  U= \bigcup_{n} U_{2n} \quad \text{ and } \quad V = \bigcup_{n} U_{2n+1}.
\] 
Then we have both that $\partial(U,V) \ge \delta/2$ and $a \in \cl[\tau]{U} \cap \cl[\tau]{V}$, which contradicts Lemma~\ref{l:caraccompact1}.

Thus every $\tau$-convergent sequence admits a $\partial$-Cauchy subsequence, which combined with the facts that $\partial$ is complete and that the $\partial$-topology refines $\tau$ implies the statement of the theorem.
\end{proof}

\begin{coro}
  \label{c:metrizable-subsp}
  Let $K \subseteq S(X)$ be a subset such that $K$ equipped with the relative $\tau$-topology is a metrizable topological space. Then the $\partial$-topology and $\tau$ coincide on $K$ and in particular, if $K$ is $\tau$-closed, $(K, \partial)$ is compact.
\end{coro}
\begin{proof}
  We already know that the $\partial$-topology is finer than $\tau$. To see the converse, note that if $K$ is metrizable, its topology is determined by convergence of sequences and then apply Theorem~\ref{t:caraccompact}.
\end{proof}

Now we let $G$ be a Polish group. We fix a bounded, right-invariant metric $d$ on $G$. Then $(G, d)$ is a metric space and we construct $(S(G), \tau, \partial)$ as above.

\begin{proof}[Proof that if $M(G)$ is metrizable then it has a comeager orbit]$ $

We view $M(G)$ as a subflow of $S(G)$ and fix a $G$-equivariant retraction $r \colon S(G) \to M(G)$ (see Proposition \ref{p:retraction}). 

To show that there is a comeager orbit, we apply Rosendal's criterion (Lemma \ref{l:Rosendal}). Let $V \ni 1_G$ and $U \sub M(G)$ be given. We may assume that $V = \{g : d(g,1_G) < \varepsilon\}$ for some $\varepsilon > 0$. Since $M(G)$ is metrizable, we have that $\tau$ and the $\partial$-topology coincide on $M(G)$ by Corollary~\ref{c:metrizable-subsp} and we can find a non-empty $\tau$-open $U' \subseteq U$ of $\partial$-diameter $< \varepsilon$.

Let $W_1, W_2 \subseteq U'$ be non-empty, open. By choice of $U'$ we have $\partial(W_1,W_2) < \varepsilon$; since $W_1 \subseteq r^{-1}(W_1)$ and $W_2 \subseteq r^{-1}(W_2)$, we also have that $\partial(r^{-1}(W_1), r^{-1}(W_2))< \varepsilon$. Then Lemma~\ref{l:caraccompact1} tells us that $d(r^{-1}(W_1) \cap G, r^{-1}(W_2) \cap G) < \varepsilon$. So we can find $f_1 \in r^{-1}(W_1) \cap G$ and $f_2 \in r^{-1}(W_2) \cap G$ such that $d(f_1, f_2) < \varepsilon$, that is, $f_2 f_1^{-1} \in V$. Since $r(f_2) = f_2f_1^{-1}r(f_1) \in W_2 \cap f_2f_1^{-1}W_1$, Rosendal's criterion is verified and we are done.
\end{proof}

Now we know that $M(G)$, if metrizable, must have a comeager orbit. It remains to prove that the stabilizer $H$ of a point $x_0$ in this comeager orbit is co-precompact, extremely amenable and that $M(G)= \yhwidehat{G/H}$. 
We first establish that this stabilizer is co-precompact; the argument makes essential use of the coalescence of $M(G)$.

\begin{lem}
Let $G$ be a Polish group such that $M(G)$ is metrizable. Let $x_0$ be an element of the comeager orbit and $H$ be the stabilizer of $x_0$. Then $H$ is co-precompact in $G$, and $M(G)$ is isomorphic to $\yhwidehat{G/H}$.
\end{lem}

\begin{proof}
To avoid potentially confusing notation, denote by $y_0$ the point $H$ in $G/H$. Since $g \mapsto gx_0$ is right-uniformly continuous, the map $\psi \colon G/H \to M(G)$ defined by $\psi(gy_0)=gx_0$ is uniformly continuous hence, by compactness of $M(G)$, it extends to a continuous, $G$-equivariant $\psi \colon \yhwidehat{G/H} \to M(G)$.

By the universal property of $M(G)$, there exists a continuous, $G$-equivariant $\varphi \colon M(G) \to \yhwidehat{G/H}$. Then $\psi \circ \varphi= \tau$ is a $G$-equivariant map from $M(G)$ to itself; by coalescence, $\tau$ is a homeomorphism so, replacing $\varphi$ by $\varphi \circ \tau^{-1}$, we may as well assume that $\psi \circ \varphi= \mathrm{id}_{M(G)}$.

It remains to show that $\varphi$ is onto; for that, it is sufficient to prove that $y_0 \in \varphi(M(G))$ since $\overline{G y_0}= \yhwidehat{G/H}$. Denote $z_0 =\varphi(x_0)$; we now prove that $y_0=z_0$, which will conclude the proof.

Let $V$ be a neighborhood of $1_G$. Since $Gx_0$ is comeager, the Effros theorem tells us that there exists an open $U \ni 1_G$ such that for all $g \in G$ one has $g x_0 \in U x_0 \Rightarrow g \in VH$. Since $\psi(z_0)=x_0$, $\psi$ is continuous, and $U x_0$ is a neighborhood of $x_0$ in $Gx_0$, there exists an open $W \ni z_0$ such that for any $g$ one has $gy_0 \in W \Rightarrow \psi(gy_0)=gx_0 \in U	x_0$.

We have thus shown that, for any neighborhood $V$ of $1_G$, there exists an open $W \ni z_0$ such that for any $g$ one has $g y_0 \in W \Rightarrow g \in VH$. Pick a sequence $(g_n)_{n< \omega}$ of elements of $G$ such that $g_n y_0 \xrightarrow[n \to + \infty]{} z_0$. It follows from the previous observation that $g_n y_0 \xrightarrow[n \to +\infty]{} y_0$, whence $y_0=z_0$, as promised. Thus $\varphi$ is onto, and we are done.
\end{proof}

\begin{lem}\label{l:continuous_equivariant_nomalizer}
Let $G$ be a Polish group and $H$ be a co-precompact closed subgroup.

Assume that $\varphi \colon G/H \to G/H$ is continuous and $G$-equivariant. Then there exists $f$ in the normalizer of $H$ such that for all $g \in G$ one has $\varphi(gH)=gfH$.
\end{lem}

Before the proof, we note that the reason why this result is not immediate is that we do not know a priori that $\varphi$ is injective; proving that requires some work.

\begin{proof}
There is not much choice when it comes to choosing $f$: just pick some $f$ such that $\varphi(H)=fH$. We need to prove that $f$ belongs to the normalizer of $H$. Since $\varphi$ is $G$-equivariant, we have $hfH=fH$ for all $h \in H$, i.e., $f^{-1}Hf \subseteq H$.

Now, let $U$ be an open neighborhood of $1_G$ and note that for every $g_1,g_2 \in G$ we have 
\[g_1 \in Ug_2H \Rightarrow g_1f  \in Ug_2 Hf \subseteq Ug_2fH .\]

We conclude that $\varphi \colon G/H \to G/H$ is uniformly continuous (for the quotient uniformity on $G/H$) hence it extends to a continuous $\varphi \colon \yhwidehat{G/H} \to \yhwidehat{G/H}$. 

Denote for notational simplicity $Y=\yhwidehat{G/H}$. It follows from the computation above that there exists a fundamental system of entourages $O$ for the (unique compatible) uniformity of $Y$ such that $(\varphi \times \varphi)(O) \subseteq O$ for all $O$. We claim that this, associated to the compactness of $Y$, implies that $\varphi$ is injective. Granting this for the moment, consider $g$ such that $g \varphi(H)=\varphi(H)$. Then $gH=H$, whence $g \in H$. Since every element of $fHf^{-1}$ fixes $fH=\varphi(H)$ we conclude that $fHf^{-1} \subseteq H$ so that $fHf^{-1}=H$ and we are done.

We still have a claim to prove. Choose $y_1 \ne y_2 \in V$, let $V$ be a neighborhood of the diagonal of $Y$ such that $(y_1,y_2) \not \in V$, then choose a symmetric entourage $W$ such that $W^5 \subseteq V$ and $(\varphi \times \varphi)(W) \subseteq W$. Then consider a finite subset $A$ such that $W[A]=Y$ and such that, among all such subsets, $|\lset(a_1,a_2) : (a_1,a_2) \not \in W^3\rset|$ has the smallest cardinality (such a subset exists since $Y$ is compact hence totally bounded).

Now note that, by surjectivity of $\varphi$,  $W[\varphi(A)]=Y$ and for any $a_1,a_2 \in A$ we have $(a_1,a_2) \in W^3 \Rightarrow (\varphi(a_1),\varphi(a_2)) \in W^3$ (remember that $(\varphi \times \varphi)(W) \subseteq W$). The minimality condition on $A$ thus implies that if $(a_1,a_2) \not \in W^3$ then also $(\varphi(a_1),\varphi(a_2) \not \in W^3$.

There exists $a_1,a_2 \in A$ such that $(a_1,y_1) \in W$ and $(a_2,y_2) \in W$. If $(\varphi(y_1),\varphi(y_2)) \in W$ then $(\varphi(a_1),\varphi(a_2)) \in W^3$ whence $(a_1,a_2) \in  W^3$, which in turn implies $(y_1,y_2) \in W^5$. Since we know that $(y_1,y_2) \not \in W^5$ we thus have that $(\varphi(y_1),\varphi(y_2)) \not \in W$, in particular $\varphi(y_1) \ne \varphi(y_2)$, and $\varphi$ is injective as promised. 
\end{proof}

\begin{exo}\label{exo:surjective_Lipschitz-compact}
\begin{enumerate}
\item Prove that every surjective, $1$-Lipschitz map from a compact metric space to itself is an isometry (adapt the argument used in the proof of Lemma \ref{l:continuous_equivariant_nomalizer}).
\item Use this to reformulate the argument given in the proof above, by considering the quotient metric on $G/H$ associated to a right-invariant metric on $G$ and showing that $G$-equivariant maps are $1$-Lipschitz.
\end{enumerate}
\end{exo}

\begin{lem}
Let $G$ be a Polish group and $H$ a co-precompact, closed subgroup such that $G \actson \yhwidehat{G/H}$ is minimal. Let $z \in \yhwidehat{G/H}$ be such that $Hz=z$. Then there exists $f$ in the normalizer of $H$ such that $z=fH$.
\end{lem}

\begin{proof}
Again, to avoid confusion, denote $x_0= H \in \yhwidehat{G/H}$; for simplicity write $X=\yhwidehat{G/H}$.

Let $z$ be as above. Since $g \mapsto gz$ is right-uniformly continuous there exists a continuous, $G$-equivariant $\pi \colon X \to X$ such that $\pi(x_0)= z$. Let $U$ be open nonempty in $X$. By minimality and separability, there exists a countable sequence $(g_n)_{n< \omega}$ of elements of $G$ such that $\bigcup_n g_n U= X$, whence $\bigcup_{n< \omega} g_n \pi(U)=X$, so $\pi(U)$ is not meager. We conclude that $\pi$ is category-preserving. Since $Gx_0$ is comeager in $X$ we thus have that $\pi^{-1}(Gx_0)$ is comeager in $X$, whence $Gx_0 \cap \pi^{-1}(Gx_0) \ne \emptyset$.

We just proved that $\pi(Gx_0)=Gx_0$, so $\pi$ induces a $G$-equivariant map from $G/H$ to itself: applying the result of the previous lemma, we conclude that there exists $f$ in the normalizer of $H$ such that $\pi(gx_0)=gfx_0$ for all $g \in G$. In particular, $z=\pi(x_0)=fx_0$. 
\end{proof}

We are now ready to conclude the proof.

\begin{thm}
Let $G$ be a Polish group such that $M(G)$ is metrizable. Let $x_0$ be an element of the comeager orbit and $H$ be the stabilizer of $x_0$. Then $H$ is co-precompact, extremely amenable and $M(G) = \yhwidehat{G/H}$.
\end{thm}

\begin{proof}
The only remaining part is to prove that $H$ is extremely amenable. To that end, we use Theorem \ref{thm:finite_osc_stable}. Let $d$ be a right-invariant metric on $G$.

Denote by $L_H$ the space of all $1$-Lipschitz functions from $(H,d)$ to $[0,1]$ and endow it with the (compact) topology of pointwise convergence. Recall that $H$ acts on $L_H$ via $h \cdot \varphi(h_1)=\varphi(h_1h)$, and that we need to prove that for every $\varphi \in L$ there exists a fixed point in $\overline{H \varphi}$.

Similarly let $L_G$ denote the space of all $1$-Lipschitz functions from $(G,d)$ to $[0,1]$, with its compact topology.

Now, choose a $1$-Lipschitz $\varphi \colon (H,d) \to [0,1]$. We extend $\varphi$ to a $1$-Lipschitz $\varphi \colon (G,d) \to [0,1]$. 

Consider the diagonal action $G \actson L_G \times \yhwidehat{G/H}$; name $x_0$ the element $H$ of $G/H$.

Since $M(G)$ has a $H$-fixed point, every $G$-flow has a $H$-fixed point. In particular, there exists $(\psi,z) \in \overline{G (\varphi,x_0)}$ such that $H (\psi,z) = (\psi,z)$. We know that there exists $f$ in the normalizer of $H$ such that $z= fx_0$. Observe that $(f^{-1} \cdot \psi,x_0) \in \overline{G (\varphi,x_0)}$ is fixed by $H$. Denote $\psi'=f^{-1} \cdot \psi$.

There exists a sequence $(g_n)_{n< \omega}$ of elements of $G$ such that $g_n \cdot \varphi \xrightarrow[n \to + \infty]{} \psi'$ and $g_n x_0 \xrightarrow[n \to + \infty]{} x_0$. 

We may pick a sequence $(h_n)_{n< \omega}$ of elements of $H$ such that $h_n g_n^{-1} \xrightarrow[n \to + \infty]{} 1_G$. For every $n$ and every $f \in G$ we have 
\[|\varphi(fh_n)-\varphi(fg_n)| \le d(fh_n,fg_n)= d(f h_ng_n^{-1},f) \xrightarrow[n \to + \infty]{} 0 .\]
It follows that $h_n \cdot \varphi \xrightarrow[n \to + \infty]{} \psi'$, and we conclude that there exists an $H$-fixed point inside $\overline{H \varphi}$. Thus $H$ is extremely amenable and we are done.
\end{proof}

\medskip
\emph{Comments.} The proof that a metrizable UMF of a Polish group has a comeager orbit is lifted essentially verbatim from \cite{BenYaacovMellerayTsankov}; some of the ideas originate in work of Zucker \cite{Zucker16}. The second part of the argument is a slightly reformulated presentation of ideas from \cite{MellerayNVTTsankov}. 

Recent work of Basso and Zucker suggests that the existence of a comeager orbit in the universal minimal flow is a natural dividing line for general topological groups (they call such groups ``groups with tractable dynamics''; see \cite{Basso_Zucker}).

Topometric structures occur naturally in several contexts related to continuous logic (such as type spaces); on any Polish group $(G,\tau)$ it is interesting to consider the topometric structure $(G,\tau,\partial)$, where $\partial(g,h)= \sup \{d(gk,hk) : k \in G\}$ for some left-invariant metric $d$ on $G$. Another, related example is given by considering weak*-topologies on the dual of a normed vector space (or its unit ball), and the metric given by the operator norm. See for instance the papers \cite{BenYaacov2008} and \cite{Hanson2025}.

\epigraphhead{}
\clearpage
\epigraphhead{}
\pagestyle{plain}
\setlength{\epigraphwidth}{0.45\textwidth}

\part{Descriptive set theory}
\setcounter{chapter}{8}

\chapter{Polish spaces}

We recall that a subset $A$ of a topological space $X$ is $G_\delta$ if there exists a sequence of open sets $(O_n)_{n< \omega}$ such that $A= \bigcap_{n< \omega} O_n$; and $F_\sigma$ if $X\setminus A$ is $G_\delta$, that is, if there exists a sequence of closed sets $(F_n)_{n<\omega}$ such that $A=\bigcup_{n<\omega} F_n$.

Polish spaces have been defined in the first chapter; we cover some of their basic theory now.

\begin{exo}\label{exo:finite_union_g_delta}
Prove that a finite union of $G_\delta$ subsets of a topological space is still a $G_\delta$ subset.
\end{exo}

\begin{thm}[Alexandrov]\index{Alexandrov's theorem}\label{thm:Alexandrov}
Let $(X,d)$ be a complete metric space, and $Y \subseteq X$ be a $G_\delta$ subspace. Then there exists a compatible metric $\rho$ on $Y$ (for the topology induced by $d$) such that $(Y,\rho)$ is complete.
\end{thm}

\begin{proof}
Write $Y= \bigcap_{n< \omega} O_n$, where each $O_n$ is $G_\delta$. First consider for $n< \omega$ the metric $\rho_n$ on $Y$ defined by 
\[\rho_n(y,y')= d(y,y') + \left| \frac{1}{d(y,X \setminus O_n)} - \frac{1}{d(y,X \setminus O_n)} \right | .\]
This is a metric which induces the same topology as $d$ on $Y$ (the function $y \mapsto d(y,X \setminus O_n)$ is continuous and never $0$  on $Y$ since $Y \subseteq O_n$ and $O_n$ is open).

Let $(y_i)_{i < \omega}$ be a Cauchy sequence of elements of $(Y,\rho_n)$. Then $(y_i)_{i<\omega}$ is Cauchy in $(X,d)$, hence $y_i \xrightarrow[i \to + \infty]{} x \in X$. Furthermore, $\left( \frac{1}{d(y_i,X \setminus O_n)} \right)_{i< \omega}$ is also Cauchy in $\R$, hence bounded. So $x \in O_n$.

Now consider on $Y$ the metric 
\[\rho(y,y')  = \sum_{n=0}^{+ \infty} \frac{1}{2^n}\min(\rho_n(y,y'),1) .\] 
This is a compatible metric on $Y$, and a sequence $(y_i)_{i< \omega}$ is Cauchy for $\rho$ iff $(y_i)_{i< \omega}$ is Cauchy for each $\rho_n$, whence $(y_i)_{i < \omega}$ converges to some $x \in \bigcap_{n} O_n$. This proves that $(Y,\rho)$ is complete.
\end{proof}

\begin{defin}
Let $X,Y$ be two metric spaces, $A$ a nonempty subset of $X$ and $f \colon A \to Y$ a function. For every $x \in \overline{A}$ we define the \emph{oscillation}\index{oscillation of a function at a point} $w_f(x)$ of $f$ at $x$ by the formula 
\[\omega_f(x) = \inf_{\varepsilon >0} \sup \lset d(f(a),f(a')) : a,a' \in A \textrm{ and } a,a' \in B(x,\varepsilon) \rset .\]
\end{defin}

More or less by definition, one sees that $f$ is continuous at $x \in A$ iff $\omega_f(x)=0$; and, assuming that $Y$ is complete, that $f$ admits a limit at $x \in \overline{A} \setminus A$ iff $\omega_f(x)=0$.

\begin{exo}\label{exo:oscillation}
Let $X,Y$ be two metric spaces, $A$ a nonempty subset of $X$ and $f \colon A \to Y$ a function. Prove that for each $r >0$ the set $\lset x \in \overline{A} : \omega_f(x) < r \rset$ is open in $\overline{A}$.
\end{exo}

\begin{thm}[Kuratowski]
Let $X,Y$ be two metric spaces, $A$ a subset of $X$ and $f \colon A \to Y$ a continuous function. Assume that $Y$ is complete. Then there exists a $G_\delta$ subset $B$ such that $A \subseteq B \subseteq \overline{A}$ and $f$ extends to a continuous function from $B$ to $Y$.
\end{thm}

\begin{proof}
Granting the result of Exercise \ref{exo:oscillation}, the set $B=\lset x \in \overline{A} : \omega_f(x)=0 \rset$ 
contains $A$ and is $G_\delta$. Furthermore, $f$ is continuous on $A$ and admits a limit at every point of $B$ (because $Y$ is complete), so one can extend $f$ by continuity to $B$.
\end{proof}

\begin{thm}[Lavrentiev]\index{Lavrentiev's theorem}
Let $X$, $Y$ be two complete metric spaces. Assume that $A$, $B$ are subsets of $X$, $Y$ respectively and that $f \colon A \to B$ is a homeomorphism. Then there exist $G_\delta$ subsets $\widetilde A$, $\widetilde B$ of $X$, $Y$ respectively such that $f$ extends to a homeomorphism from $\widetilde A$ onto $\widetilde B$.
\end{thm}

\begin{proof}
Find a $G_\delta$ subset $A_1 $ containing $A$ and such that $f$ extends to a continuous function $f_1 \colon A_1 \to Y$; and a $G_\delta$ subset $B_1$ containing $B$ and such that $f^{-1}$ extends to a continuous function $g \colon B_1 \to X$. Then define 
\[\widetilde{A} = \lset x \in A_1 : f_1(x) \in B_1 \textrm{ and } g(f_1(x))=x \rset \, , \  \widetilde{B} = \lset y \in B_1 : g(y) \in A_1 \textrm{ and }  f_1(g(y))=y \rset .\]
Since $f_1$, $g$ are continuous, $\widetilde{A}$ is a closed subset of $f_1^{-1}(B_1) \cap A_1$, which is a $G_\delta$ subset of $X$. Thus $\widetilde{A}$ is $G_\delta$ in $X$. Similarly, $\widetilde{B}$ is $G_\delta$ in $Y$.

By definition, $A \subseteq \widetilde{A}$, $B \subseteq \widetilde{B}$. It is easy to check that $f_1(\widetilde{A})=\widetilde{B}$ and $g(\widetilde{B})=\widetilde{A}$, so we are done.
\end{proof}

\begin{coro}\label{coro:G_delta=Polish}
Let $X$ be a complete metric space, and $A$ a subset of $X$. There exists a compatible complete metric on $A$ iff $A$ is a $G_\delta$ subset of $X$.
\end{coro}

\begin{proof}
Simply consider $\mathrm{id} \colon A \to A$, where on the left $A$ is seen as a subset of $X$ and on the right it is a subset of itself (endowed with a compatible complete metric). By Lavrentiev's theorem there exists a $G_\delta$ subset $\widetilde{A}$ of $X$ such that $\mathrm{id}$ extends to a homeomorphism from $\widetilde{A}$ onto $A$. Of course one must have $\widetilde{A}=A$, so that $A$ is $G_\delta$ in $X$.
\end{proof}

\begin{exo}\label{exo:locally_compact_open_in_completion}
Let $X$ be a metrizable space, and $Y \subseteq X$ be a dense subset of $X$ which is locally compact for the induced topology. Prove that $Y$ is open in $X$.
\end{exo}

It follows from the result of the previous exercise that every locally compact metrizable space admits a compatible complete metric; in particular, every separable, metrizable locally compact space is Polish.

We now introduce two fundamental (and ubiquitous) examples.

\begin{defin}
Let $A$ be a discrete set. We denote $A^{< \omega}$\index{$A^{<  \omega}$} the set of finite sequences of elements of $A$; for $s=(s(0),\ldots,s(n))$ we let $|s|=n+1$ (and agree that the empty sequence is such that $|\emptyset|=0$). If $i=|t|=|s|+1$ and $t(n)=s(n)$ for all $n< |s|$ then we denote $t=s \smallfrown t(i-1)$.

We endow $A^{\omega}$ with the product topology, which has a basis of clopen sets of the form 
\[N_s = \lset \alpha \in A^\omega : \alpha_{|n}=s \rset \]
where $s \in A^{< \omega}$ and $|s|=n$.

The \emph{Cantor space}\index{Cantor space $\mcC$} is $\mcC=2^\omega$, and the \emph{Baire space}\index{Baire space $\mcN$} is $\mcN=\omega^\omega$.
\end{defin}

For future reference, we point out the following fact.

\begin{exo}\label{exo:copies_Baire_Cantor}
Show that each of $\mcC$ and $\mcN$ contains a homeomorphic copy of the other (to embed $\mcN$ in $\mcC$, consider the set of all binary sequences taking the value $1$ infinitely many times).
\end{exo}

\begin{lem}
Let $A$ be a discrete set and $X=A^\omega$ endowed with the product topology. Let $F \subseteq X$ be a nonempty closed subset. Then $F$ is a continuous retract of $X$, i.e., there exists a continuous surjection $r \colon X \to F$ such that $r \circ r= r$.
\end{lem}

\begin{proof}
First choose, for each $s$ such that $F \cap N_s \ne \emptyset$, some $\beta_s \in F \cap N_s$. Given $\alpha \in A^\omega$, using the fact that $F$ is closed we see that either $\alpha \in F$ or there exists some $n(\alpha) \in \omega$ such that $N_{\alpha_{|n(\alpha)}} \cap F \ne \emptyset$ and $N_{\alpha_{|n(\alpha)+1}} \cap F= \emptyset$. In the first case, set $r(\alpha)=\alpha$; in the second case let $s=\alpha_{|n(\alpha)}$ and set $r(\alpha)=\beta_s$.

Clearly, $r$ is a retraction of $X$ onto $F$. If $\alpha \not \in F$, then for every $\alpha'$ such that $\alpha'(i)=\alpha(i)$ for all $i \le n(\alpha)$ one has $r(\alpha)=r(\alpha')$, hence $r$ is constant on a neighborhood of $\alpha$. If $\alpha \in F$, then for any $n$ and any $\alpha'$ such that $\alpha'(i)=\alpha(i)$ for all $i \le n$ one has $r(\alpha')(i)=\alpha(i)$ for all $i \le n$. Hence in that case $r(\alpha')\xrightarrow[\alpha' \to \alpha]{} \alpha=r(\alpha)$, so $r$ is continuous at $\alpha$.
\end{proof}

\begin{prop}\label{prop:continuous_open_image_of_Baire}
Assume that $Y$ is a metrizable space, and that there exists a continuous, open surjection $\pi$ from $\omega^\omega$ onto $Y$. Then $Y$ is Polish.
\end{prop}

In the next chapter we will see that this theorem actually provides a characterization of Polish spaces. 

\begin{proof}
Fix a compatible metric $d$ on $Y$. We aim to prove that $Y$ is $G_\delta$ in its completion $(\yhwidehat{Y},d)$.

Denote $U_s= \pi(N_s)$. Since $\pi$ is open, there exists an open $V_s \subset \yhwidehat{Y}$ such that $U_s= V_s \cap Y$. Intersecting $V_s$ with $B(y,2\mathrm{diam}(U_s))$ for some $y \in U_s$, we ensure that $\mathrm{diam}(V_s) \le 4 \mathrm{diam}(U_s)$.

We would like to prove that $Y=\bigcap_{n \in \omega} \bigcup_{|s|=n} V_s$; this need not be true (inclusion from left to right is obvious, the converse may be false) so we need to work a little.

We can find a family of open sets $(W_s)_{s \in \omega^{< \omega}}$ such that :
\begin{itemize}
\item for all $s$, one has $W_s \subseteq V_s$.
\item For all $n$, $\bigcup_{|s|=n} W_s= \bigcup_{|s|=n} V_s$.
\item For all $n$ and all $z \in \yhwidehat{Y}$, the set $\lset s \in \omega^n : z \in W_s \rset$ is finite.
\end{itemize}
To see why this is indeed possible, see Exercise \ref{exo:loc_finite_subcover}.

Now we claim that $Y= \bigcap_{n \in \omega} \bigcup_{|s|=n} W_s$; granting this, we obtain as promised that $Y$ is $G_\delta$ in $\hat{Y}$, hence Polish.

We turn to the proof of the claim. Assume that $z \in \bigcap_{n \in \omega} \bigcup_{|s|=n} W_s$. Consider 
\[\mcT=\lset s \in \omega^{< \omega} : z \in W_s \rset .\]

This is an infinite subtree of $\omega^{< \omega}$ such that each node has finitely many successors (see Chapter \ref{chapter:universality} if necessary for more information concerning trees on $\omega$); hence by König's lemma there exists $\alpha \in \omega^\omega$ such that $z \in W_{\alpha_{|n}}$ for all $n< \omega$. But we also have $\pi(\alpha) \in V_{\alpha_{|n}}$ for all $n< \omega$, and the diameter of $V_{\alpha_{|n}}$ is vanishing since it is controlled by $4 \mathrm{diam}(U_{\alpha_{|n}})$. If follows that $z=\pi(\alpha) \in Y$, and we are done.
\end{proof}

In the proof above we used the following important fact about metric spaces (the demonstration of which is not immediate!).

\begin{exo}\label{exo:loc_finite_subcover}
Let $X$ be a metric space, and $(U_n)_{n < \omega}$ be a sequence of open subsets of $X$. Prove that there exists a sequence $(V_n)_{n< \omega}$ of open subsets of $X$ with the following properties:
\begin{enumerate}
\item For each $n$, $V_n \subseteq U_n$.
\item $\bigcup_{n< \omega} V_n= \bigcup_{n< \omega} U_n$.
\item For each $x \in X$, $\lset n : x \in V_n \rset$ is finite.
\end{enumerate}
\end{exo}

We turn to a structure theorem that, in practice, implies that many theorems about Polish spaces can be reduced to the case of \emph{perfect} Polish spaces (recall that a topological space is perfect \index{perfect topological space} if it has no isolated points). 

\begin{exo}\label{exo:F_sigma_not_G_delta}
Let $(X,d)$ be a perfect Polish space and $Q \subseteq X$ a countable dense subset. Prove that $Q$ is $F_\sigma$ but not $G_\delta$.
\end{exo}

\begin{thm}[Cantor--Bendixson]\index{Cantor--Bendixson theorem}
Let $X$ be a Polish space. Then one can write $X= P \sqcup D$, where $P$ is perfect and $D$ is open and (at most) countable. Moreover, such a decomposition is unique.
\end{thm}

\begin{proof}
Let $D= \lset x \in X : x \textrm{ has a countable neighborhood}\rset$. This set is open by definition (it can equivalently be described as the union of all countable open subsets of $X$), and countable because $X$ is second-countable hence has the Lindelöff property. 

Denote $P=X \setminus D$. Then every $x \in P$ is such that every neighborhood of $x$ has uncountable intersection with $P$, hence $P$ is perfect.

Note that if $U$ is a countable open subset of a Polish space then, since the Baire category theorem holds in $U$, $U$ must contain an isolated point. This implies that, in a perfect Polish space, every open subset is uncountable. Hence if one writes $X= P_1 \sqcup D_1$ with $P_1$ perfect and $D_1$ open countable, then $P_1 \cap D = \emptyset$; similarly $P \cap D_1 = \emptyset$, thus $P_1=P$ and this proves the uniqueness of the decomposition.
\end{proof}

The set $P$ above is called the \emph{perfect kernel}\index{perfect kernel} of $X$.

\begin{exo}\label{exo:derivative}
Let $X$ be a Polish space. For $A$ a subset of $X$ we define 
\[ D(A)= A \setminus \{a \in A : a \textrm{ is an isolated point of } A\}. \]
Then we define inductively subsets $X^{(\alpha)}$ for $\alpha$ an ordinal, as follows:
\begin{itemize}
\item $X^{(0)}=X$.
\item For every $\alpha$, $X^{(\alpha +1)}=D(X^{(\alpha)})$.
\item For every limit ordinal $\alpha$, $X^{(\alpha)}= \bigcap_{\beta < \alpha} X^{(\beta)}$.
\end{itemize}
Prove that there exists $\alpha < \omega_1$ such that $X^{(\alpha)}=X^{(\alpha+1)}=P$, where $P$ is the perfect kernel of $X$.
 
The smallest $\alpha$ such that  $X^{(\alpha)}=X^{(\alpha+1)}$ is called the \emph{Cantor--Bendixson rank}\index{Cantor--Bendixson rank} of $X$. 

The set-theoretically inclined reader may want to prove that, for every countable ordinal $\alpha$, the countable ordinal $\omega^\alpha +1$, endowed with the order topology, is Hausdorff, compact, countable (hence countable Polish) and has Cantor--Bendixson rank equal to $\alpha +1$.
\end{exo}

We conclude this chapter with another classical application of the Baire category theorem.

\begin{exo}\label{exo:Baire_class_1}
\begin{enumerate}
\item Let $X,Y$ be Polish spaces. Let $f \colon X \to Y$ be a function of \emph{Baire class $1$}, i.e, $f$ is a pointwise limit of a sequence of continuous functions from $X$ to $Y$.
\begin{enumerate}
\item Prove that for every closed subset $A \subseteq Y$ the set $f^{-1}(A)$ is $G_\delta$ in $X$.
\item Prove that the set of points of continuity of $f$ is dense $G_\delta$ in $X$.
\end{enumerate} 
\item Let $f \colon \R \to \R$ be a differentiable function. Prove that the set $\{x : f' \textrm{ is continuous at } x\}$ is dense $G_\delta$ in $\R$.
\end{enumerate}
\end{exo}

\bigskip
\emph{Comments.} We took for granted the Baire category theorem, which is covered in many introductory texts on the topology of metric spaces. The material presented here is very standard; many more details may be found in \cite{Kechris1995} or \cite{Srivastava}. The exposition in this chapter (and those that follow it) is closely influenced by \cite{Kechris1995}.

\chapter{Schemes and transfer theorems}

\begin{defin}
Let $(X,d)$ be a metric space. A \emph{Cantor scheme}\index{Cantor scheme} on $X$ is a family $(F_s)_{s \in 2^{< \omega}}$ of nonempty subsets of $X$ with the following properties:
\begin{enumerate}
\item For every $s \in 2^{< \omega}$ one has $F_{s \smallfrown 0} \cap F_{s \smallfrown 1} = \emptyset$.
\item For each $\varepsilon \in \{0,1\}$ and each $s \in 2^{< \omega}$ one has $F_{s \smallfrown \varepsilon} \subseteq F_s$.
\end{enumerate}
We say that the Cantor scheme is \emph{convergent}\index{convergent Cantor scheme} if it satisfies the following additional conditions:
\begin{enumerate}[resume]
\item For each $\varepsilon \in \{0,1\}$ and each $s$ one has $\overline{F_{s \smallfrown \varepsilon}} \subseteq F_s$.
\item For every $\alpha \in 2^\omega$, $\mathrm{diam}(F_{\alpha_{|n}}) \xrightarrow[n \to + \infty]{} 0 $.
\end{enumerate}
\end{defin}

If $(X,d)$ is complete, the four conditions above imply that for every $\alpha \in 2^\omega$ the intersection $\bigcap_n 
F_{\alpha_{|n}}$ is a singleton. 

\begin{exo}\label{exo:continuity_convergent_Cantor_scheme}
Prove that the function associated to a convergent Cantor scheme on a complete metric space is continuous and injective.
\end{exo}

\begin{thm}\label{t:embedding_a_Cantor}
Let $(X,d)$ be an uncountable Polish metric space. Then there exists a continuous injection from the Cantor space $\mcC$ into $X$.
\end{thm}

It follows that an uncountable Polish space must have cardinality greater than $2^{\aleph_0}$ (actually, equal to $2^{\aleph_0}$, see Exercise \ref{exo:CH_for_Polish} below).

\begin{proof}
Because of the Cantor--Bendixson theorem, it is enough to consider the case where $X$ is perfect; we assume that such is the case.

It is enough to build a convergent Cantor scheme on $X$. For this, it is sufficient (via an inductive process) to prove that for any $\varepsilon >0$ and every nonempty open $U \subseteq X$ there exist two nonempty disjoint open subsets $V,W$ of $U$ with $V \cap W = \emptyset$, $\overline{V}, \overline{W} \subseteq U$ and $\mathrm{diam}(V), \mathrm{diam}(W) \le \varepsilon$.

We may as well assume that $U= B(x,r)$ for some $r>0$. Since $X$ is perfect, we can find $y \ne x \in U$. Let $\delta >0$ be such that $2 \delta \le \varepsilon$ and $d(x,y)+ \delta < r$. Then $V=B(x,\delta)$ and $W=B(y,\delta)$ have all the desired properties.
\end{proof}

\begin{exo}\label{exo:CH_for_Polish}
Let $(X,d)$ be a separable metric space. Prove that the cardinality of $X$ is smaller than $2^{\aleph_0}$. 

It follows that for every uncountable Polish space $X$ there exists a bijection $f \colon X \to \R$: the continuum hypothesis is true for Polish spaces.
\end{exo}

We recall that a topological space is $0$-dimensional if every point admits a neighborhood basis consisting of clopen subsets. For compact Hausdorff spaces, this is equivalent to the requirement that every connected component is a singleton.

We may also use Cantor schemes to prove the following characterization of the Cantor space.

\begin{thm}[Brouwer]
Assume that $X$ is nonempty, compact, metrizable, $0$-dimensional and perfect. Then $X$ is homeomorphic to the Cantor space.
\end{thm}

\begin{proof}
It is enough to prove that there exists a convergent Cantor scheme $(F_s)_{s \in 2^{< \omega}}$ on $X$ with the additional properties that each $F_s$ is clopen, $F_\emptyset=X$ and for each $s$ one has $F_s= F_{s \smallfrown 0} \sqcup F_{s \smallfrown 1}$. Indeed the (necessarily injective and continuous) function associated to the scheme is then a bijection, hence a homeomorphism by compactness.

Fix a compatible metric $d$ on $X$. First, note that for every nonempty clopen subset $U$ of $X$ and every $\varepsilon >0$ there exists a partition $(U_i)_{i \le n}$ of $X$ with $n \ge 2$, each $U_i$ nonempty clopen and $\mathrm{diam}(U_i) \le \varepsilon$ for each $i \le n$. 

To see this, pick $x \ne y \in U$, $r>0$ such that $r < d(x,y)$ and cover $U$ by a finite family $(V_i)_{i \le m}$ of clopen subsets of diameter $\le \min(r,\varepsilon)$. Then replace each $V_i$ with $V_i \setminus \bigcup_{j<i} V_j$ to obtain the desired family (discarding the empty set if it occurs; there will be at least two nonempty subsets in the family since $x$, $y$ cannot belong to the same $V_i$).

Let us now produce the desired Cantor scheme. Let $F_{\emptyset}=X$. Applying the previous observation, we obtain $U_0,\ldots,U_n$ clopen nonempty and of diameter $\le 1$ such that $X= \bigsqcup_i U_i$. For $1 \le i \le n$, set $F_{0^i}= U_i \sqcup \ldots \sqcup U_n$, and for $0 \le i \le n-1$ set $F_{0^i \smallfrown 1}= U_{i}$. 

Then simply repeat the same construction within each $U_i$, this time with clopen subsets of diameter $\le \frac{1}{2}$, and keep going.
\end{proof}

A similar idea enables one to prove the following useful fact.
\begin{thm}\label{exo:perfect_set_property_analytic}
Let $X$, $Y$ be Polish spaces and $f \colon X \to Y$ a continuous map. Assume that $f(X)$ is uncountable. Then there exists a subset $K$ of $X$ homeomorphic to $\mcC$ and such that $f_{|K}$ is injective.
\end{thm}

\begin{proof}
We claim that one can build a convergent Cantor scheme $(F_s)_{s \in 2^{< \omega}}$ on $X$ with the property that $f(F_{s \smallfrown 0}) \cap f(F_{s \smallfrown 1}) = \emptyset$ for all $s \in 2^{< \omega}$. Denoting by $g \colon \mcC \to X$ the (injective) function associated to this Cantor scheme, we then obtain that $f \circ g \colon \mcC \to Y$ is injective, which establishes the desired result.

It is enough to show that if $F \subseteq X$ is closed and $f(F)$ is uncountable, then for every $r>0$ there exist two disjoint closed subsets $F_0$, $F_1$ of $F$ such that $f(F_0)$, $f(F_1)$ are uncountable, have empty intersection and $F_0$, $F_1$ have diameter $\le r$.  First, remove from $F$ the set of all $y \in F$ such that $f(B(y,\rho) \cap F)$ is countable for some $\rho >0$; we are removing an open set with countable image (because of the Lindelöff property), so we are left with a closed subset $F' \subseteq F$ which still has uncountable image. For any $x, y\in F'$ such that $f(x) \ne f(y)$, any sufficiently small closed balls around $x$, $y$ will satisfy the desired condition (by continuity of $f$).
\end{proof}

\begin{thm}[Mycielski]\label{thm:Mycielski}\index{Mycielski's theorem}
Let $X$ be a nonempty perfect Polish space and $R \subseteq X^2$ be comeager. Then there exists a continuous, injective map $f \colon \mcC \to X$ such that $(f(\alpha),f(\beta)) \in R$ for any $\alpha \ne \beta \in \mcC$.
\end{thm}

\begin{proof}
Let $d$ be a compatible metric on $X$. We may assume that $R= \bigcap_{k< \omega} O_k$, where each $O_k$ is dense open in $X^2$ and $O_{k+1} \subseteq O_k$ for all $k$.

We claim that one can define a convergent Cantor scheme $(U_s)_{s \in 2^{< \omega}}$ on $X$
with the following properties :
\begin{enumerate}
\item Each $U_s$ is open.
\item For each $k$ and each distinct $s,t \in 2^k$ one has $U_s \times U_t \subseteq O_k$.
\end{enumerate}
Assuming that this is feasible, we know that the continuous map $f \colon \mcC \to X$ induced by this Cantor scheme is injective, and the second condition above guarantees that  $(f(\alpha),f(\beta)) \in R$ for any $\alpha \ne \beta \in \mcC$.

Assume that this construction has been carried out up to some step $k$ (set $U_\emptyset=X$ to get started). Then notice that for any nonempty open subsets $V_0$, $V_1$ of $X$ and any $k$ there exist disjoint nonempty open subsets $W_0$, $W_1$ of $V_0$, $V_1$ such that $W_0 \times W_1 \subseteq O_{k+1}$ (because $O_{k+1}$ is dense open). Shrinking $W_0$, $W_1$ we may further assume that they have disjoint closures and have arbitrarily small diameter. Then inductively apply this fact (finitely many times, obtaining the desired open subsets as intersections of finitely many open subsets) to obtain the desired family $(U_t)_{t \in 2^{k+1}}$ from the family $(U_s)_{s \in 2^k}$. 
\end{proof}

Alternative arguments to prove theorems \ref{exo:perfect_set_property_analytic} and \ref{thm:Mycielski} are given in Chapter \ref{chapter:Effros}.

\begin{exo}\label{exo:meager_equivalence_relations}
Let $X$ be Polish and $R$ be a Baire measurable equivalence relation on $X$.
\begin{enumerate}
\item Prove that $R$ is meager iff each equivalence class of $R$ is meager.
\item Assume that $R$ is meager. Prove that there exists a perfect subset of $X$ intersecting each $R$-class in at most one point (in particular, there are continuum many distinct $R$-equivalence classes).
\end{enumerate}
\end{exo}

To solve the first question of this exercise it is convenient to use the Kuratowski--Ulam theorem, which was proved in the first chapter of these notes.

\begin{defin}
Let $(X,d)$ be a metric space. 
A \emph{Lusin scheme}\index{Lusin scheme} on $X$ is a family $(F_s)_{s \in \omega^{< \omega}}$ of sets such that 
\begin{enumerate}
\item For every $s \in \omega^{< \omega}$ and every $i \ne j \in \omega$ one has $F_{s \smallfrown i} \cap F_{s \smallfrown j} = \emptyset$.
\item For each $s \in \omega^{< \omega}$ and each $i< \omega$ one has $F_{s \smallfrown i} \subseteq F_s$.
\end{enumerate}
We say that a Lusin scheme is \emph{convergent} if it satisfies the following additional property.
\begin{enumerate}[resume]
\item For every $\alpha \in \omega^\omega$ one has $\mathrm{diam}(F_{\alpha_{|n}}) \xrightarrow[n \to + \infty]{} 0$ (by convention, $\mathrm{diam}(\emptyset)=0$).
\item For every $\alpha \in \omega^\omega$ one has $\bigcap_{k< \omega} F_{\alpha_{|k}} = \bigcap_{k< \omega} \overline{F_{\alpha_{|k}}}$.
\end{enumerate}
\end{defin}

Any injective continuous map $f \colon \mcN \to X$ induces a convergent Lusin scheme, by setting $F_s = \pi(N_s)$; conversely, assuming that $(X,d)$ is complete, every convergent Lusin scheme induces an injective continuous map $\pi \colon Y \to X$, where $Y=\lset \alpha \in \mcN : \forall n \ F_{\alpha_{|n}} \ne \emptyset \rset$ is a  closed subset of $\omega^\omega$ and $\{\pi(\alpha)\}= \bigcap_n F_{\alpha_{|n}}$ for every $\alpha \in Y$ (to see that $Y$ is closed, note that by definition it is an intersection of clopen sets). 

In practice, we are often going to ensure that the last condition above holds by asking that $\overline{F_{s\smallfrown i}} \subseteq F_s$ for every $s \in \omega^{< \omega}$ and every $i< \omega$.

Just as we used convergent Cantor schemes to build continuous, injective maps with domain $2^\omega$, one can use convergent Lusin schemes to build continuous, injective maps with domain a closed subset of $\omega^\omega$. Let us give an example.

\begin{thm}\label{t:Polish_continuous_image_of_Baire}
Let $X$ be a Polish space. Then there exists a closed subset $Y$ of $\omega^\omega$ and a continuous bijection $\pi \colon Y \to X$.
\end{thm}

\begin{proof}
We claim that one can build a convergent Lusin scheme with $F_\emptyset=X$ and such that for each $s$ one has $F_s= \bigcup_{n< \omega} F_{s \smallfrown n}$. Assume for now that this is possible, and (assuming that $F_s \ne \emptyset$) pick $x \in F_s$. Then there exists a unique $i$ such that $x \in F_{s \smallfrown i}$, then a unique $j$ such that $x \in F_{s \smallfrown i \smallfrown j}$, and so on. This produces $\alpha \in N_s$ such that $x= \pi(\alpha)$.

Denote $Y=\lset \alpha \colon \forall n \ F_{\alpha_{|n}} \ne \emptyset \rset$, which we know is a closed subset of $\omega^\omega$. We have shown that $\pi(N_s \cap Y)=F_s$ for each $s$, in particular $\pi$ is surjective. Since it is also continuous and injective by definition of a Lusin scheme, we are done.

It remains to justify that one can build the desired Lusin scheme; we can actually build one where each $F_s$ is an $F_\sigma$ subset of $X$. It is sufficient to prove that, given a $F_\sigma$ subset $F$ of $X$ and some $\varepsilon >0$, one can write $F=\bigsqcup_{n< \omega} F_n$ with $F_n$ an $F_{\sigma}$ subset of diameter $\le \varepsilon$ and $\overline{F_n} \subseteq F$ for all $n$.

Start from $F= \bigcup_n A_n$, where each $A_n$ is closed. Since each closed subset of a Polish space is a countable union of closed subsets with diameter $\le \varepsilon$, we may assume that $\mathrm{diam}(A_n) \le \varepsilon$ for all $n$.

For all $n$ let $B_n=A_n \setminus \bigcup_{i< n} A_i$. This is an $F_\sigma$ subset of $X$ since open sets are $F_\sigma$ and an intersection of two $F_\sigma$ subsets is still $F_\sigma$. Furthermore, we have by definition $\overline{B_n} \subseteq \overline{A_n} =A_n \subseteq F$ for all $n$. So we have as desired that each $B_n$ is $F_\sigma$, $\mathrm{diam}(B_n) \le \varepsilon$, $F=\bigsqcup_n B_n$ and that $\overline{B_n} \subseteq F$ for all $n$.
\end{proof}

\begin{exo}\label{exo:perfect_Polish_bijection_Baire}
Let $(X,d)$ be a Polish space. Prove that $X$ is perfect iff there exists a continuous bijection $\pi \colon \mcN \to X$.

(First, show that if $X$ is Polish perfect then one can write $X= \bigsqcup_{n< \omega} X_n$, where each $X_n$ is nonempty, perfect and $G_\delta$)
\end{exo}

It follows from the previous theorem that every Polish space is a continuous image of $\omega^\omega$ (recall that $\omega^\omega$ continuously retracts on each of its closed subsets). One can actually also ensure that $\pi$ is open.

\begin{thm}\label{t:Polish_open_image_of_Baire}
Let $(X,d)$ be a nonempty Polish space. Then there exists a continuous, open surjection $\pi \colon \omega^\omega \to X$.
\end{thm}

\begin{proof}
We use a small variation on the notion of Lusin scheme: we claim that there exists a family $(F_s)_{s \in \omega^{< \omega}}$ of nonempty open sets such that :

\begin{enumerate}
\item For every $s \in \omega^{< \omega}$ and every $i \in \omega$ one has $\overline{F_{s \smallfrown i}}
 \subseteq F_s$.
\item For every $s \in \omega^{< \omega}$ and every $i \in \omega$ one has $F_s = \bigcup_{i \in \omega} F_{s \smallfrown i}$.
\item For every $\alpha \in \omega^\omega$ one has $\mathrm{diam}(F_{\alpha_{|n}}) \xrightarrow[n \to + \infty]{} 0$.
\end{enumerate}

Assume for the moment that this is possible; then for each $\alpha$ one can set $\{\pi(\alpha)\}= \bigcap_n F_{\alpha_{|n}}$ and the same argument as in the proof of Theorem \ref{t:Polish_continuous_image_of_Baire} gives us that for all $s$ one has $F_s= \pi(N_s)$, which implies that $\pi$ is an open surjection.

To build a family of open sets satisfying the conditions from the list above, it is enough to observe that for every nonempty open $U \subseteq X$ and every $\varepsilon >0$ one can write $U= \bigcup_n U_n$, with $U_n$ nonempty open, $\overline{U_n} \subseteq U$ and $\mathrm{diam}(U_n) \le \varepsilon$. To achieve this, find for each $x \in U$ some $r_x >0$ such that $B(x,2r_x) \subseteq U$. Using the Lindelöff property, we find a sequence $(x_n)_{n < \omega}$ of elements of $U$ such that $U=\bigcup_n B(x_n,r_{x_n})$ (note that we allow repetitions in the sequence $(x_n)_{n< \omega}$). Then we obtain the desired family of open sets by setting $U_n=B(x_n,r_{x_n})$.
\end{proof}

\begin{exo}[Sierpinski]\label{exo:continuous_open_image_of_Polish}
Let $X$ be a Polish space, $Y$ a metrizable space and $f \colon X \to Y$ be a continuous, open surjection. Prove that $Y$ is Polish.
\end{exo}

\begin{exo}\label{exo:0_dim_Polish_in_Baire_and_Cantor}
Let $X$ be a $0$-dimensional Polish space. Prove that $X$ is homeomorphic to a closed subset of $\mcN$ and to a $G_\delta$ subset of $\mcC$.
\end{exo}

The following theorem is part of the reason why the Baire space $\mcN$ is ubiquitous in descriptive set theory.

\begin{thm}[Alexandrov--Urysohn]
The Baire space $\mcN$ is the unique (up to homeomorphism) nonempty $0$-dimensional Polish space in which every compact subset has empty interior.
\end{thm}

\begin{proof}
We know that $\mcN$ is $0$-dimensional and Polish. Denote by $\pi_n \colon \mcN \to \omega$ the projection $\alpha \mapsto \alpha(n)$. If $K \subseteq \omega$ is compact then $\pi_n(K)$ is bounded; by definition of the product topology, for any nonempty open subset $U\subseteq \mcN$ one must have $\pi_n(U)=\omega$ for all large enough $n$. It follows that any compact subset of $\mcN$ has empty interior.

Consider a Polish, perfect, $0$-dimensional space $X$ in which every compact subset has empty interior. Combining the ideas of the proofs of Theorems \ref{t:Polish_continuous_image_of_Baire} and \ref{t:Polish_open_image_of_Baire}, it is enough to show that for any nonempty clopen subset $U$ of $X$ and any $\varepsilon >0$ there exists a disjoint family $(U_n)_{n< \omega}$ of nonempty clopen subsets of $X$ such that $\mathrm{diam}(U_n) \le \varepsilon$ for all $n$ and $\bigsqcup_n U_n= U$. Indeed, this enables us to build a convergent Lusin scheme such that the associated map is defined on the whole of $\omega^\omega$ and is an open surjection.

So, fix a nonempty clopen $U \subseteq X$ and $\varepsilon >0$. Since $U$ is closed and not compact, it is not totally bounded so there exists $r \le \frac{\varepsilon}{2}$ such that $U$ cannot be covered by finitely many balls of radius $r$. Since for each ball $B(x,r)$ there exists a clopen $V$ such that $x \in V \subseteq B(x,r)$, we obtain (using the Lindelöff property) the existence of a sequence of clopen sets $U_n$ such that $U=\bigcup_n U_n$, $\mathrm{diam}(U_n) \le \varepsilon$ for all $n$ and for all $N$ one has $X \ne \bigcup_{n=0}^N U_n$. We may then set for all $n$ $V_n = U_n \setminus \bigcup_{i< n} U_i$. Each $V_n$ is clopen, $\bigsqcup_{n< \omega} V_n = X$ and $\{n : V_n \ne \emptyset\}$ is infinite.
\end{proof}

\begin{exo}\label{exo:irrationals_vs_Baire}

\begin{enumerate}
\item Assume that $X$ is a $G_\delta$ subset of $\R$ such that both $X$ and $\R \setminus X$ are dense in $\R$. Prove that $X$ is homeomorphic to $\mcN$.

Note that this proves in particular that the familiar space $\R \setminus \Q$ is homeomorphic to $\mcN$ (a fact that one could also prove using continued fraction expansions of irrational numbers).
\item Prove that the conclusion of the previous question still holds if one replaces $\R$ by a $0$-dimensional Polish space.
\end{enumerate}
\end{exo}

\begin{exo}\label{exo:dense_G_delta_copy_of_Baire_space}

\begin{enumerate}
\item Let $X$ be a Polish space. Prove that there exists a dense $G_\delta$ subset $Y$ of $X$ which is $0$-dimensional.
\item Let $X$ be an uncountable Polish space. Prove that there exists a $G_\delta$ subset $Y$ of $X$ which is homeomorphic to $\mcN$. Show that if $X$ is perfect then there exists such a subset $Y$ which is dense in $X$. 
\item Given a perfect Polish space $X$, denote $\mcA_{\mathrm{CAT}}(X)$ the Boolean algebra of all Baire measurable subsets of $X$, identified if their symmetric difference is meager. Prove that $\mcA_{\mathrm{CAT}}(X)$ and $\mcA_{\mathrm{CAT}}(\mcN)$ are isomorphic (this Boolean algebra is called the \emph{category algebra}\index{category algebra}). 
\end{enumerate}
\end{exo}

\begin{exo}[Sierpinski]\label{exo:topological_carac_rationals}
\begin{enumerate}
\item Prove that $\Q$ is, up to homeomorphism, the unique countable metrizable space which is perfect and nonempty.
\item Prove that $\Q^n$ is homeomorphic to $\Q$ for all nonzero $n< \omega$.
\end{enumerate}
\end{exo}

\emph{Comments.} Again this is very classical material, which is entirely covered in \cite{Kechris1995}.

\chapter{Borel and analytic subsets of Polish spaces}

We recall that a subset $A$ of a topological space $(X,\tau)$ is \emph{Borel}\index{Borel subset} if it belongs to the $\sigma$-algebra generated by the open subsets of $(X,\tau)$. We denote by $\mcB(X,\tau)$ (or simply $\mcB(X)$ or $\mcB(\tau)$ when there is no risk of confusion) the family of Borel subsets of $(X,\tau)$.

Borel subsets of a given topological space may be obtained ``from within'' by iterating the operations of countable union, countable intersection, and complementation. This leads to the following recursive definition.

\begin{defin}
Let $X$ be a Polish space. We define by transfinite induction on $1\le \xi< \omega_1$ two families of subsets $\mathbf{\Sigma}_\xi^0(X)$, $\mathbf{\Pi}_\xi^0(X)$ as follows:
\begin{enumerate}
\item $\mathbf{\Sigma}_1^0(X)$ is the family of open subsets of $X$. 
\item For each $\xi$, $\mathbf{\Pi}_\xi^0(X)$ is the family of subsets $A$ of $X$ such that $X \setminus A \in \mathbf{\Sigma}_{\xi}^0(X)$ (in particular, $\mathbf{\Pi}_1^0(X)$ is the family of closed subsets of $X$.)
\item For each $2 \le \xi< \omega_1$, $\mathbf{\Sigma}_{\xi}^0(X)$ is the family of all subsets $B$ of $X$ such that there exists a sequence $(A_n)_{n< \omega}$ of subsets of $X$ such that each $A_n$ belongs to $\mathbf{\Pi}_{\xi_n}^0(X)$ for some $\xi_n < \xi$ and $B= \bigcup A_n$.
\end{enumerate}
\end{defin}

One can also define the \emph{ambiguous pointclasses} $\mathbf{\Delta}_\xi^0(X)= \mathbf{\Sigma}_\xi^0(X) \cap \mathbf{\Pi}_\xi^0(X)$.

Here we should point out that the use of a bold font matters; in effective descriptive set theory one uses ``lightface'' notation such as $\Sigma^0_\xi(X)$ to denote other classes of subsets, see \cite{Moschovakis}.

\begin{exo}\label{exo:Borel_pointclasses}
Let $X$ be a Polish space and $\xi <\omega_1$. 
\begin{enumerate}
\item Prove that $\mathbf{\Sigma}_{\xi}^0(X)$ is stable under countable unions and finite intersections.
\item Prove that for all $\xi$ one has $\mathbf{\Sigma}_\xi^0(X) \subseteq \mathbf{\Sigma}_{\xi+1}^0(X)$ and $\mathbf{\Sigma}_\xi^0(X) \subseteq \mathbf{\Pi}_{\xi+1}^0(X)$.
\end{enumerate}
\end{exo}

\begin{exo}\label{exo:Borel_pointclasses_2}
Let $X$ be a Polish space. Prove that
\[ \mcB(X)= \bigcup_{\xi < \omega_1} \mathbf{\Sigma}_\xi^0(X) = \bigcup_{\xi < \omega_1} \mathbf{\Pi}_\xi^0(X) .\]
\end{exo}

Given a Borel subset $A$ of $X$, one can then think of the least $\xi$ such that $A \in \mathbf{\Sigma}_\xi^0(X)$ as a measure of the complexity of the definition of $A$; while this is a very fruitful idea we will not pursue this direction at all in these notes. Mostly, here we think of Borel subsets of a given Polish space $X$ as being ``nice'' or ``definable'', without a particular care for the precise complexity of the definition - with the notable exception of the $\mathbf{\Pi}_2^0$-level of this hierarchy, which corresponds to $G_\delta$ subsets and which we have already  abundantly discussed.

\begin{exo}[Lebesgue]\label{exo:Lebesgue}
Prove that the set of Borel functions from $\R$ to $\R$ is the smallest subset $A$ of $\R^\R$ with the following two properties:
\begin{itemize}
\item Every continuous function belongs to $A$.
\item If $(f_n)_{n < \omega}$ is a sequence of elements of $A$ which converges pointwise to $f$ then $f \in A$.
\end{itemize}
\end{exo}

We begin with a simple observation.
\begin{prop}
Let $X$, $Y$ be Polish spaces. Assume that $f \colon X \to Y$ is a Borel map. Then its graph $\Gamma_f$ is a Borel subset of $X \times Y$.
\end{prop}

Shortly we will prove that the converse of this result also holds, an important basic result of descriptive set theory.

\begin{proof}
Let $(V_n)_{n< \omega}$ be a countable basis for the topology of $Y$. Then for any $x \in X$, $y \in Y$ one has
$(f(x)=y) \Leftrightarrow (\forall n < \omega \ (f(x) \in V_n) \Leftrightarrow (y \in V_n))$.

This gives us the equality $\Gamma_f = \bigcap_{n< \omega} \left( f^{-1}(Y \setminus V_n) \times (Y \setminus V_n) \cup f^{-1}(V_n) \times V_n \right)$. 
\end{proof}

An important technique in descriptive set theory is to use refinements of topologies. Let us give an example, which will have many applications below.

\begin{thm}\label{thm:making_Borel_clopen}
Let $(X,\tau)$ be a Polish space. Then for each Borel subset $B$ of $X$ there exists a Polish topology $\tau_B$ on $X$ which refines $\tau$, has the same Borel subsets as $\tau$, and for which $B$ is a clopen subset.
\end{thm}

Here we should note that, soon, we will establish that if $\tau'$ is Polish and refines $\tau$ then necessarily $\mcB(X,\tau')=\mcB(X,\tau)$ so part of the conclusion above is actually redundant (and similarly for several other statements in this chapter).

\begin{proof}
Let us denote by $\mcA$ the family of subsets of $X$ which satisfy the conclusion of the theorem; clearly $\mcA$ contains $\emptyset$ and $X$ and is stable under complementation.

Let $F$ be a closed subset of $X$; we know that there exists a complete metric $d_1$ inducing the topology of $F$, and a complete metric $d_2$ inducing the topology of $O=X \setminus F$. We may also assume that $d_1$, $d_2$ are bounded by $1$, then define
\[d(x,x')= \begin{cases} d_1(x,x') & \textrm{ if } (x,x') \in F^2 \\
d_2(x,x') & \textrm{ if } (x,x') \in O^2 \\1 & \textrm{ otherwise} \end{cases}. \]
It is easy to check that $d$ generates a Polish topology $\tau_F$ on $X$ which refines $\tau$ and for which $F$ is clopen. 
Furthermore, each open ball for $d$ is clearly $\tau$-Borel, so every open subset of $(X,\tau_F)$ is $\tau$-Borel. Conversely, every $\tau$-open subset is $\tau_F$-open so $\mcB(X,\tau)=\mcB(X,\tau_F)$. This proves that $\mcA$ contains every $\tau$-closed subset as well as every $\tau$-open subset.

To see why $\mcA$ is closed under countable unions, it is enough to show that if $(\tau_n)_{n< \omega}$ is a sequence of Polish topologies refining $\tau$ and $\mcB(\tau_n)=\mcB(\tau)$ for all $n$, then the topology $\tau_\infty$ generated by $\bigcup_n \tau_n$ is a Polish topology such that $\mcB(\tau_\infty)= \mcB(\tau)$. Indeed, assume that this is proved and let $(A_n)_{n< \omega}$ be a sequence of elements of $\mcA$, as witnessed by a sequence of Polish topologies $(\tau_n)_{n< \omega}$. Then we obtain a Polish topology $\tau_\infty$ refining $\tau$, with $\mcB(\tau_\infty)=\mcB(\tau)$ and for which $A=\bigcup_n A_n$ is a union of clopen sets, hence an open set. Applying the same argument to $B=\bigcup_n (X \setminus A_n)$, we obtain a Polish topology $\tau_\infty'$ with the desired properties and for which $B$ is open. Then  the topology generated by $\tau_\infty \cup \tau_\infty'$ is a Polish topology refining $\tau$, with the same Borel subsets as $\tau$ and such that $A$ is clopen.

Now, consider a sequence $(\tau_n)_{n< \omega}$ of topologies on $X$ with the properties mentioned above. Each $\tau_\infty$-open subset is $\tau$-Borel, so that $\mcB(\tau_\infty) \subseteq \mcB(\tau)$, and the converse implication is immediate since each $\tau_n$ refines $\tau$.

To see that $\tau_\infty$ is a Polish topology, denote $Y= \prod_{n< \omega}(X,\tau_n)$, with the product topology. This is a Polish space, in which $X$ embeds via the map $\varphi \colon x \mapsto (x,x,\ldots,x,\ldots)$. This is a homeomorphism from $(X,\tau_\infty)$ onto $\varphi(X)$; we only need to prove that $\varphi(X)$ is $G_\delta$. To that end, pick $y=(x_i)_{i< \omega} \not \in  \varphi(X)$. Then there exists $n\ge 1$ such that $x_n \ne x_0$. Since $\tau$ is Polish, there exist two $\tau$-open subsets $U$, $V$ such that $x \in U$, $y \in V$ and $U \cap V= \emptyset$. Then $\lset(x'_n)_{n< \omega} : x_0' \in U \textrm{ and } x_n' \in V \rset$ is an open subset of $Y$ which contains $y$ and does not intersect $\varphi(X)$. Thus $\varphi(X)$ is actually closed, which is more than we need.
\end{proof}

\begin{exo}\label{exo:CH_for_Borel}
Let $(X,\tau)$ be a Polish space, and $B$ be an uncountable Borel subset of $X$. Prove that $|B|=2^{\aleph_0}$.
\end{exo}

\begin{exo}\label{exo:making_Borel_continuous}
Let $(X,\tau_X)$, $(Y,\tau_Y)$ be Polish spaces and $f \colon X \to Y$ a Borel map. Prove that there exists a Polish topology $\tau_f$ on $X$ which refines $\tau$, has the same Borel subsets as $\tau_X$, and such that $f \colon (X,\tau_f) \to (Y,\tau_Y)$ is continuous.
\end{exo}

\begin{exo}\label{exo:Borel_injective_image_of_closed_in_Baire}
Let $X$ be a Polish space. Prove that for every Borel subset $B \subseteq X$ there exists a closed $Y \subseteq \mcN$ and a continuous bijection $f \colon Y \to B$.
\end{exo}

A fundamental fact in descriptive set theory (indeed, in a way, the origin of the whole subject) is that, given a continuous map $f \colon X \to Y$ between Polish spaces $X$, $Y$, the image $f(X)$ need not be Borel. This leads us to the following definition.

\begin{defin}
Let $X$ be a Polish space, and $A$ be a subset of $X$. We say that $A$ is \emph{analytic}\index{analytic subset of a Polish space} if there exists a Polish space $Y$ and a continuous map $f \colon Y \to X$ such that $f(X)=A$.

We say that $A$ is \emph{coanalytic}\index{coanalytic subset of a Polish space} if $X \setminus A$ is analytic.
\end{defin}

One denotes by $\ana(X)$\index{$\ana$} the family of analytic subsets of $X$, and by $\coana(X)$\index{$\coana$} the family of coanalytic subsets. Note that every Borel subset of $X$ is both analytic and coanalytic, see Exercise \ref{exo:Borel_injective_image_of_closed_in_Baire} (so, often the notation $\mathbf{\Delta}_1^1(X)$ is used for Borel subsets of a Polish space $X$).

In the definition of an analytic set, one might equivalently require that $Y=\mcN$, since any Polish space is the image of $\mcN$ by a continuous map. It is also worth observing here that, by the result of Theorem \ref{exo:perfect_set_property_analytic}, every uncountable analytic set has cardinality $2^{\aleph_0}$. This fact need not be true for coanalytic subsets; see however Exercise \ref{exo:cardinalities_of_coanalytic sets} for a proof that an infinite coanalytic subset of a Polish space has three possible cardinalities: $\aleph_0$, $\aleph_1$, or $2^{\aleph_0}$.

\begin{exo}\label{exo:analytic_is_projection_of_closed}
Let $X$ be a Polish space and $A$ a subset of $X$. Prove that $A$ is analytic iff there exists a closed subset $F$ of $\mcN \times X$ such that $A=\pi(F)$ (where $\pi$ denotes the projection on the second coordinate).
\end{exo}

\begin{exo}\label{exo:unions_intersections_analytic}
Let $X$ be a Polish space, and $(A_n)_{n< \omega}$ be a sequence of analytic subsets of $X$. Prove that $\bigcup_{n< \omega} A_n$ and $\bigcap_{n< \omega} A_n$ are analytic subsets of $X$.
\end{exo}

\begin{exo}\label{exo:properties_pointclass_analytic}
Let $X$, $Y$ be Polish spaces. 
\begin{enumerate}
\item Assume that $A \subseteq X$ is analytic and $f \colon A \to Y$ is a continuous map. Prove that $f(A)$ is analytic.
\item Assume that $A \subseteq X$ is analytic and $f \colon X \to Y$ is Borel. Prove that $f(A)$ is analytic. 
\item Assume that $B \subseteq Y$ is analytic and $f \colon X \to Y$ is Borel.  Prove that $f^{-1}(B)$ is analytic.
\end{enumerate}
\end{exo}

Now we turn to the proof of the fact that for any uncountable Polish space $X$ there exists an analytic subset of $X$ which is not Borel. The proof goes though a coding argument. Note that it is enough to deal with the case $X=\mcN$. Indeed, assume that we know that $A \subseteq \mcN$ is analytic non Borel, and $X$ is Polish and uncountable. We proved above that $X$ contains a $G_\delta$ subset which is homeomorphic to $\mcN$, so we may as well consider that $A \subseteq \mcN \subseteq X$. Then $A$ is still an analytic subset of $X$, and if $A$ were Borel in $X$ then it would also be Borel in $\mcN$, which is not the case.

\begin{defin}
Let $\Gamma$ be a class of subsets of Polish spaces (e.g., open sets, closed sets, Borel sets, analytic sets...).

A subset $A \subseteq \mcN^2$ is \emph{universal for} $\Gamma(\mcN)$ if $A \in \Gamma(\mcN \times \mcN)$ and, for every $B \in \Gamma(\mcN)$, there exists $x \in \mcN$ such that $B=A_x=\{y \in \mcN : (x,y) \in A\}$.
\end{defin}

Observe that if $A \subseteq \mcN^2$ is universal for $\ana(\mcN)$ then $A$ is not a Borel subset of $\mcN^2$. Indeed, if $A$ were Borel then $B=\lset x \in \mcN : (x,x) \not \in A \rset$ would be a Borel, hence analytic subset of $\mcN$. Then there would exist $x_0 \in \mcN$ such that $B=A_{x_0}$, which would lead us to an instance of the classical barber paradox:
\[x_0 \in B = A_{x_0} \Leftrightarrow (x_0,x_0) \in A \Leftrightarrow x_0 \not \in B .\]

So, if we know that there exists a universal subset for $\ana(\mcN)$, then there exists an analytic subset of $\mcN \times \mcN$ which is not Borel. Since $\mcN^2$ is homeomorphic to $\mcN$, such a subset also exists in $\mcN$.

Establishing the next result becomes a pressing concern.

\begin{thm}
There exists a universal subset for $\ana(\mcN)$.
\end{thm}

\begin{proof}
Let us first show that there exists a universal subset for $\mathbf{\Sigma}_1^0(\mcN)$. To see this, fix a countable basis $(V_n)_{n< \omega}$ of the topology of $\mcN$, and consider 
\[A= \lset (\alpha, \beta) \in \mcN^2 : \beta \in \bigcup_{n< \omega} V_{\alpha(n)} \rset .\]
Observe that $A$ is open in $\mcN^2$ . To see this, pick $(\alpha,\beta) \in A$. There exists $n$ such that $\beta \in  V_{\alpha(n)}$. For any $\alpha'$ such that $\alpha(n)=\alpha'(n)$ and any $\beta' \in  V_{\alpha(n)}$ we also have $(\alpha',\beta') \in A$.

By definition, every open subset of $\mcN$ is of the form $\mcA_\alpha$ for some $\alpha \in \mcN$.

The complement $B$ of $A$ is then a $\mathbf{\Pi}_1^0$-universal subset of $\mcN$ : given a closed subset $F$ of $\mcN$ there exists $\alpha$ such that $\mcN \setminus F = A_\alpha$, whence $F= B_\alpha$.

It follows from all this that there exists a closed subset $P$ of $\mcN \times \mcN^2$ such that, for every closed subset $F$ of $\mcN^2$, there exists $\beta \in \mcN$ such that $F=P_\beta$. Then consider
\[C= \lset (\alpha,\beta) \in \mcN^2 : \exists \gamma \ (\alpha,\beta,\gamma) \in P \rset .\]
This set is analytic in $\mcN^2$, since it is the image of $P$ by the projection map $(\alpha,\beta,\gamma) \mapsto (\alpha,\beta)$.

Let $B \subseteq \mcN$ be analytic. Then $B$ is the continuous image of some function $f \colon \mcN \to \mcN$ hence (considering the graph of this function) the projection of some closed subset $F_B$ of $\mcN^2$. Since $P$ is universal, there exists $\alpha \in \mcN$ such that
\[ \forall (\beta,\gamma) \in \mcN^2 \ (\beta,\gamma) \in F_B \Leftrightarrow (\alpha,\beta,\gamma) \in P . \]
Fix such an $\alpha$. For each $\beta \in \mcN$ one has
\[
\beta \in B \Leftrightarrow \exists \gamma \ (\beta,\gamma) \in F_B 
        \Leftrightarrow  \exists \gamma \ (\alpha,\beta,\gamma) \in P 
        \Leftrightarrow  \beta \in C_\alpha . \qedhere
\]
\end{proof}

\begin{defin}
Let $X$ be a Polish space and $A,B$ be subsets of $X$. We say that $A,B$ are \emph{Borel-separable}\index{Borel-separable subsets} if there exists a Borel subset $C$ of $X$ such that $A \subseteq C$ and $B \cap C= \emptyset$.
\end{defin}

\begin{lem}
Let $X$ be a Polish space. Assume that $(A_n)_{n< \omega}$ and $(B_n)_{n< \omega}$ are sequences of subsets of $X$ such that $A_n$ and $B_m$ are Borel-separable for each $n,m$. Then $\bigcup_n A_n$ and $\bigcup_n B_n$ are Borel-separable.
\end{lem}

\begin{proof}
Choose for each $n,m$ some Borel subset $C_{n,m}$ of $X$ such that $A_n \subseteq C_{n,m}$ and $B_m \cap C_{n,m}=\emptyset$. Then $\bigcup_n A_n$ is contained in $C=\bigcup_n \bigcap_m C_{n,m}$, and $C \cap \bigcup_n B_n= \emptyset$. Clearly $C$ is Borel, so we are done.
\end{proof}

The following theorem is the key property of analytic subsets of Polish spaces (one says that the family of analytic sets has the \emph{separation property}\index{separation property of analytic sets}).
\begin{thm}[Lusin]\index{Lusin's theorem}
Let $X$ be a Polish space and $A$, $B$ be two disjoint analytic subsets of $X$. Then $A$ and $B$ are Borel-separable.
\end{thm}

\begin{proof}
Let $f \colon \mcN \to A$, $g \colon \mcN \to B$ be two continuous surjections. For each $n$ denote $A_s=f(N_s)$, $B_s= g(N_s)$. 

Assume that $A$ and $B$ are not Borel-separable. Since $A= \bigcup_{|s|=1} A_s$, $B= \bigcup_{|t|=1} B_t$, there exists some $s$, $t$ of length $1$ such that $A_s$ and $B_t$ are not Borel-separable. Since $A_s$ and $B_t$ are analytic, we similarly know that there exists $s'$ of length $2$ extending $s$, and $t'$ of length $2$ extending $t$, such that $A_{s'}$ and $B_{t'}$ are not Borel-separable.

Iterating this process, we obtain $\alpha \in \mcN$ and $\beta \in \mcN$ such that $A_{\alpha_{|n}}$ and $B_{\beta_{|n}}$ are not Borel-separable for any $n$. By continuity of $f$ and $g$, it follows that $f(\alpha)=g(\beta)$, hence $A \cap B \ne \emptyset$.
\end{proof}

\begin{coro}[Suslin]\index{Suslin's theorem}
Let $A$ be a subset of a Polish space $X$. Then $A$ is Borel if, and only if, both $A$ and $X \setminus A$ are analytic.
\end{coro}

\begin{proof}
If $A$ is Borel then both $A$ and $X \setminus A$ are Borel, hence analytic. Conversely, if $A$ and $X \setminus A$ are analytic then the separation property for analytic sets gives us some Borel subset $C$ such that $A \subseteq C$ and $C \cap (X \setminus A)= \emptyset$, which implies that $A=C$ is Borel.
\end{proof}

\begin{exo}\label{exo:separating_a_sequence_of_analytic_sets}
Let $X$ be a Polish space, and let $(A_n)_{n< \omega}$ be a sequence of analytic subsets of $X$ which are pairwise disjoint. Prove that there exists a sequence $(B_n)_{n < \omega}$ of Borel subsets of $X$ which are pairwise disjoint and such that $A_n \subseteq B_n$ for all $n$.
\end{exo}

\begin{thm}
Let $X$, $Y$ be two Polish spaces and $f \colon X \to Y$ a function. Then the following properties are equivalent:
\begin{enumerate}
\item $f$ is a Borel map.
\item The graph $\Gamma_f= \{(x,f(x)) : x \in X \}$ of $f$ is a Borel subset of $X \times Y$.
\item The graph of $f$ is an analytic subset of $X \times Y$.
\end{enumerate}
\end{thm}

\begin{proof}
We already saw the first implication, and the second is a triviality. So, assume that $\Gamma_f$ is analytic, and let $U$ be a Borel subset of $Y$.

For any $x \in X$ we have
\[ \left( x \in f^{-1}(U) \right) \Leftrightarrow \left( \exists y \ (x,y) \in \Gamma_f \textrm{ and } y \in U \right).\]
which proves that $f^{-1}(U)$ is analytic since it it the projection of an analytic subset of $X \times Y$.

We also have 
\[ \left( x \not \in f^{-1}(U) \right ) \Leftrightarrow \left( \exists y \ (x,y) \in \Gamma_f \textrm{ and } y \not \in U \right). \]
which shows that $X \setminus U$ is also an analytic subset of $X$. Lusin's theorem then ensures us that $f^{-1}(U)$ is Borel, and we are done.
\end{proof}

\begin{thm}[Lusin--Suslin]\label{thm:Lusin--Suslin}\index{Lusin--Suslin theorem}
Let $X,Y$ be Polish spaces, and $f \colon X \to Y$ an injective Borel function. Then $f(A)$ is Borel in $Y$ for any Borel subset $A$ of $X$. 
\end{thm}

Note that this immediately implies that $f^{-1} \colon f(A) \to X$ is a Borel map.

\begin{proof}
By the results above, the proof reduces to the case where $X=\mcN$, $A$ is a closed subset of $\mcN$ and $f$ is continuous. For $s \in \omega^{< \omega}$ denote $B_s= f(A \cap N_s)$.

Since $f$ is injective and continuous, $(B_s)_{s \in \omega^{< \omega}}$ is a convergent Lusin scheme. Furthermore, each $B_s$ is an analytic subset of $Y$. Since $B_s \cap B_t = \emptyset$ whenever $s \ne t$ and $|s|=|t|$, the separation theorem for analytic sets (and Exercise \ref{exo:separating_a_sequence_of_analytic_sets} ) allows us to upgrade this convergent Lusin scheme to another Lusin scheme $(B^*_s)_{s \in \omega^{< \omega}}$, where each $B^*_s$ is Borel and $B_s \subseteq B^*_s \subseteq \overline{B_s}$ for all $s$.

We now claim that $f(A)= \bigcap_{k < \omega} \bigcup_{s \in \omega^k} B^*_s$, which will yield the desired result. Since $B_s \subseteq B^*_s$ for all $s$, the inclusion from left to right is immediate. Conversely, choose 
$y\in \bigcap_{k < \omega} \bigcup_{s \in \omega^k} B^*_s $. Since for each $n$ and each $s \ne t$ with $|s|=|t|$ we have $B^*_s \cap B^*_t = \emptyset$, there exists for each $k$ a unique $s$ with $|s|=k$ and $y \in B^*_s$. Thus there exists $\alpha \in \mcN$ such that $y \in \bigcap_{k \in \omega} B^*_{\alpha_{|k}}$, and since $B^*_{\alpha_{|k}} \subseteq \overline{B_{\alpha_{|k}}}$ for all $k$ we obtain that $y=f(\alpha)$. Since $B_{\alpha_{|k}} \ne \emptyset$ for all $k$ we also have $A \cap N_{\alpha_{|k}} \ne \emptyset$ for all $k$, which means that $\alpha \in \overline{A}=A$, and we are done.
\end{proof}

\begin{exo}\label{exo:inclusion_topologies_equality_Borel}
Let $X$ be a set; assume that $\tau$, $\tau'$ are Polish topologies on $X$ such that $\tau \subseteq \tau'$. Prove that $\mcB(\tau)=\mcB(\tau')$.
\end{exo}

\begin{thm}\label{thm:Polish_same_cardinality_Borel_bijection}
Let $X$, $Y$ be two Polish spaces of the same cardinality. Then there exists a Borel bijection $f \colon X \to Y$.
\end{thm}

\begin{proof}
If $X$ and $Y$ are finite or countable then all of their subsets are Borel (because singletons are closed, hence Borel) so we have nothing to do. So we may as well assume that $X$ is uncountable and $Y=\mcC$. We know that there exists a continuous injection $f \colon \mcC \to X$. Furthermore, since there exists a continuous bijection from a closed subset of $\mcN$ onto $X$, there is a Borel injection of $X$ into $\mcN$, hence also into $\mcC$.

Thus it is enough to prove that if $X$, $Y$ are Polish and $f \colon X \to Y$ is a Borel injection, $g \colon Y \to X$ is a Borel injection then there exists a Borel bijection $h \colon X \to Y$. Since $g(Y)$ is a Borel subset of $X$ and $g \colon Y \to g(Y)$ is a Borel bijection, it is actually enough to deal with the case where $Y \subseteq X$.

So, let us consider a Borel injection $f \colon X \to X$ and let $Y$ be a Borel subset of $X$ such that $f(X) \subseteq Y$. For each $i< \omega$, define $X_i \subseteq X$ by setting $X_i = f^i(X \setminus Y)$. 

We now define $g \colon X \to Y$ by setting 
\[g(x)= \begin{cases} f(x) & \textrm{ if } x \in \bigcup_{i< \omega} X_i, \\ x & \textrm{ otherwise.} \end{cases} \]
Then $g$ is Borel (we are gluing together two Borel maps with disjoint domains). 

Since $g(\bigcup_{i < \omega} X_i) = \bigcup_{i< \omega} X_{i+1}$ and $f$ is injective, we see that $g$ is also injective. Since $g(X \setminus Y)= f(X) \setminus f(Y) \subseteq f(X)$ and $X_{i+1} \subseteq Y$ for every $i < \omega$, we have that $g(X) \subseteq Y$. Now, let $y \in Y$. If $y \not \in \bigcup X_i$ then $g(y)=y$ so $y \in g(X)$; else there is some $i \ge 1$ such that $y \in X_i$. Since $X_i=f(X_{i-1})$ we conclude that $y \in g(X)$. Finally we obtain as desired that $Y=g(X)$.
\end{proof}

\begin{defin}
A \emph{standard Borel space}\index{standard Borel space} is a measurable space $(X,\mcB)$ such that there exists a Polish space $Y$, and an isomorphism $f \colon (X,\mcB) \to (Y,\mcB(Y))$.
\end{defin}

We know that standard Borel spaces come in few flavors: if $X$ is at most countable, then $\mcB$ is the family of all subsets of $X$; and if $X$ is uncountable, then $(X,\mcB)$ is isomorphic to $\mcC$ endowed with its algebra of Borel subsets. This explains why standard Borel spaces are ubiquitous in analysis.

\begin{exo}\label{exo:bijection_Borel_subsets_Polish}
Let $X$, $Y$ be two Polish spaces and $A \subseteq X$, $B \subseteq Y$ two Borel sets. Prove that there exists a Borel bijection $f \colon X \to Y$ such that $f(A)=B$ iff $|A|=|B|$ and $|X \setminus A|=|X \setminus B|$.
\end{exo}

\begin{exo}\label{exo:extending_Kuratowski_theorem_to_Borel}
Let $(X,d)$ be a separable metric space. Prove that $(X,\mcB(X))$ is standard Borel iff $X$ is Borel in its completion iff $X$ is homeomorphic to a Borel subset of a Polish space.
\end{exo}

\emph{Comments.} The origins of the topic trace back to an error of Lebesgue, who mistakenly claimed that the projection on the first coordinate of a Borel subset of $\R^2$ is a Borel subset of $\R$; see \cite{Kanamori1995} for a brief overview of the beginning of the subject. Concerning this error, in a preface to a book of Lusin Lebesgue wrote the following famous words: `` La démonstration était simple, courte mais fausse. Fructueuse erreur ! Que je fus bien inspiré de la commettre!'' (which roughly translates to: ``The proof was simple, short, but false. A fruitful error! I was truly inspired to have made it!'').

Kanamori's notes also mention \emph{effective descriptive set theory}, which we do not discuss at all here; see \cite{Moschovakis}. Standard Borel spaces also feature prominently in ergodic theory, see e.g.~\cite{KerrLi}.

\chapter[Universal (co)analytic sets]{Universal (co)-analytic sets: trees and orderings on $\omega$}\label{chapter:universality}

We now discuss some important explicit examples of (co)analytic non Borel sets. 

For $s,t \in A^{< \omega}$, we write $s \sqsubseteq t$ when $|t| \ge |s|$ and for all $i< |s|$ one has $t(i)=s(i)$.

\begin{defin}
Let $A$ be a set. A \emph{tree}\index{tree} on $A$ is a subset $\mcT \subseteq A^{< \omega}$ such that for all $s,t$ such that $t \in \mcT$ and $s \sqsubseteq t$ one has $s \in \mcT$.

We say that $\mcT$ is \emph{pruned}\index{pruned tree} if every element of $\mcT$ has at least one immediate successor (in particular, $\mcT$ is infinite as soon as it is nonempty).

We say that $\alpha \in A^\omega$ is an \emph{infinite branch}\index{infinite branch in a tree} of $\mcT$ if $\alpha_{|n} \in \mcT$ for all $n< \omega$. We denote by $[\mcT]$ the set of infinite branches of $\mcT$; it is commonly called the \emph{body}\index{body of a tree} of $\mcT$.

We denote by $\mathrm{Tr}(A)$\index{$\mathrm{Tr}(A)$} the set of all trees on $A$.
\end{defin}

Note that any nonempty tree has a root, the empty sequence.
Thinking of $A$ as a discrete set, we may view $\mathrm{Tr}(A)$ as a subset of $2^{A^{< \omega}}$.

\begin{exo}\label{exo:trees_closed_subset}
Prove that $\mathrm{Tr}(A)$ is closed in $2^{A^{< \omega}}$.
\end{exo}

In particular, if $A$ is at most countable, then $\mathrm{Tr}(A)$ is a compact, metrizable space. 

The following simple observation explains why trees play an important role in classical descriptive set theory.

\begin{prop}\label{p:duality_trees_closed}
Let $F$ be a closed subset of the Baire space $\mcN$.

Then $\mcT_F=\{\alpha_{|n} : n < \omega \textrm{ and } \alpha \in F\}$ is a pruned tree on $\omega$ and $F=[\mcT_F]$. Conversely, for any $\mcT \in \mathrm{Tr}(\omega)$ the set $[\mcT]$ is closed in $\mcN$.
\end{prop}

\begin{proof}
It is clear that $\mcT$ defined above is a pruned tree, and that every $\alpha \in F$ is an infinite branch of $\mcT$. Assume that $\alpha \in [\mcT]$. Then for all $n < \omega$ there exists $\alpha_n \in F$ such that ${\alpha_n}_{|n}=\alpha_{|n}$; hence $\alpha_n \xrightarrow[n \to + \infty]{} \alpha$, so that $\alpha \in F$.

Conversely, if $(\alpha_n)_{n< \omega}$ is a sequence of elements of $[\mcT]$ which converges to some $\alpha$, then by definition for each $i < \omega$ there exists $n < \omega$ such that $\alpha_{|i}={\alpha_n}_{|i} \in \mcT$, whence $\alpha \in [\mcT]$.
\end{proof}

\begin{exo}\label{exo:finitely_branching_tree_compact_body}
Assume that $\mcT$ is a \emph{finitely branching} tree on $\omega$, i.e., every element of $\mcT$ has finitely many immediate successors. Prove that $[\mcT]$ is compact. Conversely, prove that every compact subset of $\omega^\omega$ is equal to $[\mcT]$ for some finitely branching, pruned tree $\mcT$ on $\omega$.
\end{exo}

\begin{exo}\label{exo:section_tree}
Let $\mcT$ be a pruned tree on $\omega \times \omega$. For $\alpha \in \mcN$ we define the \emph{section tree}\index{section tree} $\mcT_\alpha$ by setting 
\[s \in \mcT_\alpha \Leftrightarrow  |s|=n \textrm{ and} (\alpha_{|n},s) \in \mcT .\]
(in the definition above there is an implicit identification of $\omega^n \times \omega^n$ with $(\omega \times \omega)^n$)

Prove that $\alpha \mapsto \mcT_\alpha$ is a continuous map from $\mcN$ to $\mathrm{Tr}(\omega)$.
\end{exo}

\begin{defin}
We say that a tree $\mcT$ is \emph{well-founded}\index{well-founded tree} if $\mcT$ has no infinite branch, and \emph{ill-founded}\index{ill-founded tree} if $\mcT$ has an infinite branch. We denote by $\mathrm{IF}$ the set of all ill-founded trees on $\omega$.
\end{defin}

Note that $\mathrm{IF}$ is analytic, since $\{(\mcT,\alpha) \in \mathrm{Tr}(\omega) \times \omega^{ \omega} : \alpha_{|n} \in \mcT\}$ is clopen for all $n$, and 
\[ \mcT \in \mathrm{IF} \Leftrightarrow \exists \alpha \in \mcN \ \forall n \ \alpha_{|n} \in \mcT .\]

Much more is true: in some sense, $\mathrm{IF}$ is ``as complicated'' as an analytic set can be, as we explain now.

\begin{thm}
Let $X$ be a Polish space, and $A \subseteq X$ be an analytic set. Then there exists a Borel map $f \colon X \to \mathrm{Tr}(\omega)$ such that $A=f^{-1}(\mathrm{IF})$.
\end{thm}

One says that $\mathrm{IF}$ is a (Borel-)\emph{complete analytic} set\footnote{One could consider an \emph{ a priori} stronger notion, asking for a continuous $f$ whenever $X$ is $0$-dimensional. It turns out that both notions are equivalent; see 22B and 26B of \cite{Kechris1995}.}\index{complete analytic set}; the existence of analytic, non Borel subsets of $\mcN$ immediately implies that a complete analytic set cannot be Borel, since the inverse image of a Borel set by a Borel map is Borel.

\begin{proof}
We may assume that $X=\mcN$ since any Polish space is Borel isomorphic to a closed subset of $\mcN$. We also know that there exists a closed subset $F \subseteq \mcN \times \mcN$ such that for all $x \in \mcN$ one has 
\[x \in A \Leftrightarrow \exists \alpha \in \mcN \ (x,\alpha) \in F .\]
Choose some closed subset $F$ with this property, then fix a pruned tree $\mcT$ on $\omega \times \omega$ such that $F=[\mcT]$ and consider the continuous map $f \colon x \mapsto \mcT_x$ (we reuse the notation of exercise \ref{exo:section_tree}).

Then we have $x \in A \Leftrightarrow \exists \alpha \in \mcN \ (x,\alpha) \in F \Leftrightarrow \mcT_x(F) \in \mathrm{IF}$, and we are done.
\end{proof}

It follows that the set of well-founded trees on $\omega$ is a complete coanalytic set. When working with co-analytic sets, is is often fruitful to use another, closely related, concrete example of a complete coanalytic set.

\begin{defin}
Let $\mcT$ be a tree on $\omega$. We define an ordering $<_{\mcT}$ on $\mcT$ by setting $s<_{\mcT} t$ iff $s \sqsubset t$ or there exists some $i< |s|,|t|$ such that 
\[(\forall j < i \  s(j)=t(j)) \textrm{ and } s(i)< t(i) .\]
So, saying that $s$ is smaller than $t$ for this ordering means that either $t$ extends $s$, or that there is some $i$ such that $s(i) \ne t(i)$, and at the first such $i$ we have $s(i)< t(i)$.

This ordering is called the \emph{Kleene--Brouwer ordering}\index{Kleene--Brouwer ordering on a tree} on $\mcT$.
\end{defin}

\begin{exo}\label{exo:Kleene_Brouwer_ordering}
Let $\mcT$ be a tree on $\omega$. 
\begin{enumerate}
\item Prove that $<_{\mathrm{\mcT}}$ is a linear ordering on $\mcT$.
\item Prove that $\mcT$ is a well-founded tree iff $<_{\mathrm{\mcT}}$ is a well-ordering of $\mcT$.
\end{enumerate}
We recall that an order $<$ on a set $X$ is a well-order iff every nonempty subset of $X$ admits a minimal element for $<$ (equivalently, if there exists no infinite strictly decreasing sequence in $(X,<)$).
\end{exo}

Similarly to what we did for trees, we may view the set $\mathrm{LO}$ of all linear orderings on $\omega$ as a closed subset of $2^{\omega \times \omega}$, hence a compact metrizable set (we already considered this set in the first part of the book). Inside this set we single out the set $\mathrm{WO}$ of all well-orderings on $\omega$.

\begin{thm}
The set $\mathrm{WO}$ is complete coanalytic.
\end{thm}

\begin{proof}
First note that a linear ordering $\preceq$ on $\omega$ is not a well-ordering iff 
there exists some $\alpha \in \mcN$ such that for all $n < \omega$ one has $\alpha(n+1) \prec \alpha(n)$, an analytic condition since 
\[\{(\alpha, \preceq) : \forall n < \omega \ \alpha(n+1) \prec \alpha(n) \}\] is closed in $\mcN \times \mathrm{LO}$. This proves that $\mathrm{WO}$ is coanalytic.

We let $\mathrm{Tr}_\infty$ denote the subset of $\mathrm{Tr}(\omega)$ made up of all infinite trees; its complement is a countable subset of $\mathrm{Tr}(\omega)$.

Fix an enumeration $(s_n)_{n< \omega}$ of $\omega^{< \omega}$; then for each infinite $\mcT \subseteq \omega^{<\omega}$ we can set $n_0(\mcT)= \min(\{n : s_n \in \mcT\})$, $n_1(\mcT)= \min(\{n >n_0(\mcT): s_n \in \mcT\})$, and so on. For $t \in \mcT$ and $i < \omega$, we denote $i=n_\mcT(t)$ when $t=s_{n_i(\mcT)}$.

For each $\mcT \in \mathrm{Tr}_\infty$ define an ordering $\preceq_{\mcT}$ on $\omega$ by setting
\[\forall i,j < \omega \quad (i \preceq_{\mcT} j ) \Leftrightarrow (\exists t,t' \in \mcT \ i=n_\mcT(t) \textrm{ and } j=n_\mcT(t') \textrm{ and } t \le_\mcT t' ) .\]

Next, consider the map $f \colon \mcT \mapsto \preceq_{\mcT}$. 
It is straightforward to check that this map is continuous from $\mathrm{Tr}_\infty(\omega)$ to $\mathrm{LO}$. We conclude by observing that, for $\mcT \in \mathrm{Tr}_\infty$, $\mcT \not \in \mathrm{IF}$ iff 
$f(\mcT)  \in \mathrm{WO}$ and using the fact that $\mathrm{IF}$ is complete analytic (hence $\mathrm{Tr}_\infty \setminus \mathrm{IF}$ is complete coanalytic).
\end{proof}

This ability to encode coanalytic subsets of Polish spaces via maps into $\mathrm{WO}$ has many important consequences (it is in particular related to the notion of co-analytic rank, which we are not going to go into in these notes). The proof of the following result is an example.

\begin{thm}[The second separation theorem]
Let $X$ be a Polish space and $A,B$ be two analytic subsets of $X$. Then there exist two disjoint coanalytic subsets $C,D$ of $X$ such that $A \setminus B \subseteq C$ and $B \setminus A \subseteq D$.
\end{thm}

\begin{proof}
We may assume that $A$ and $B$ are both nonempty. We know that there exist Borel maps $f,g \colon X \to \mathrm{LO}$ such that $X \setminus A = f^{-1}(\mathrm{WO})$ and $X \setminus B=g^{-1}(\mathrm{WO})$. 

Given two orderings $\preceq$, $\preceq'$ on $\omega$, denote $(\preceq, \preceq') \in R$ if there exists a strictly increasing map $\varphi \colon (\omega,\preceq) \to (\omega,\preceq')$. By definition, $R$ is an analytic subset of $\mcN \times \mcN$. We also recall that if $\preceq$, $\preceq'$ both belong to $\mathrm{WO}$ then either $(\preceq, \preceq') \in R$ or $(\preceq', \preceq) \in R$.

We then define
\[C = X \setminus (B \cup \{x \in X : (f(x) ,g(x)) \in R \}) \textrm{ and } D = X \setminus (A \cup \{x \in X : (g(x),f(x)) \in R \}) . \]

These two sets are coanalytic since $A$, $B$ and $R$ are all analytic.

By definition, $C \cap D$ is contained in $(X \setminus B) \cap (X \setminus A)$; so for every $x \in C \cap D$ we have that both $f(x) \in \mathrm{WO}$ and $g(x) \in \mathrm{WO}$, whence $(f(x),g(x)) \in R$ or $(g(x) , f(x)) \in R$. Hence $C \cap D =\emptyset$.

Now, let $x \in A \setminus B$. Then $g(x)$ is a well-order on $\omega$ and $f(x)$ is not, hence we cannot have $(f(x), g(x)) \in R$. Thus $A \setminus B \subseteq C$, and similarly $B \setminus A \subseteq D$.
\end{proof}

\begin{exo}\label{exo:second_separation_sequence}
Let $(A_n)_{n< \omega}$ be a sequence of analytic subsets of a Polish space $X$. 

Prove that there exist pairwise disjoint coanalytic sets $(C_n)_{n< \omega}$ such that for each $n$ one has $A_n \setminus \bigcup_{m \ne n} A_m \subseteq C_n$.
\end{exo}

\begin{prop}
The set $\mathrm{UB}=\{\mcT \in \mathrm{Tr}(\omega): \mcT \textrm{ has a unique infinite branch} \}$ is coanalytic.
\end{prop}

\begin{proof}
For $s \in \omega^{< \omega}$ we let $A_s= \{\mcT \in \mathrm{Tr}(\omega): \exists \alpha \in [\mcT] \cap N_s \}$, $B_s = \bigcup_{|t|=|s|, t \ne s} A_s$. By definition, these are analytic subsets of $\mathrm{Tr}(\omega)$.

We have 
\[ \mathrm{UB}= \bigcup_{\alpha \in \mcN} \bigcap_{k< \omega} (A_{\alpha_{|k}} \setminus B_{\alpha_{|k}}) .\] 
By definition, for $s \ne t$ of the same length we have $(A_s \setminus B_s) \cap (A_t \setminus B_t)= \emptyset$. Furthermore, for each $s$ and $i$ we have $A_{s \smallfrown i} \subset A_s$ and $B_{s \smallfrown i} \supseteq B_s$. In particular, the sets  $(A_s \setminus B_s)_{s \in \omega^{< \omega}}$ form a Lusin scheme; but they are too complicated for our purposes since each $A_s \setminus B_s$ is the intersection of an analytic subset and a coanalytic subset.

Similarly to what we did to prove the Lusin--Suslin theorem, we would like to find a Lusin scheme made up of coanalytic sets $(C_s)_{s \in \omega^{< \omega}}$ such that for each $s \ne t$ of the same length $C_s \cap C_t= \emptyset$ and 
$\mathrm{UB}=\bigcup_{\alpha \in \mcN} \bigcap_{k< \omega} C_{\alpha_{|k}}$. Indeed, provided we can manage to do this, we  obtain that 
$ \mathrm{UB}= \bigcap_{k < \omega} \bigcup_{s \in \omega^k} C_s$
which implies that $\mathrm{UB}$ is coanalytic.

First, we apply the second separation theorem (or rather, its corollary established in Exercise \ref{exo:second_separation_sequence}), obtaining for each $k< \omega$ a countable family $(C_s)_{|s|=k}$ of coanalytic subsets of $\mathrm{Tr}(\omega)$ such that $A_s \setminus B_s \subseteq C_s$ for all $s$. Replacing $C_s$ by $\bigcap_{t \sqsubseteq s} C_t$, we may assume that $C_{s \smallfrown i} \subseteq C_s$ for each $s$ and $i$.

Next, for each $s$ we let $D_s= C_s \cap (\overline{A_s} \setminus B_s)$. Note that each $D_s$ is coanalytic; we have for each $s$ that 
\[A_s \setminus B_s \subseteq D_s \subseteq \overline{A_s} \setminus B_s .\] 
Since $C_s \cap C_t = \emptyset$ for $s \ne t$ with $|s|=|t|$, we also have $D_s \cap D_t= \emptyset$ for such $s$, $t$. 
It is also clear that $D_{s \smallfrown i} \subseteq D_s$ for each $i$ (recall that $A_{s \smallfrown i} \subset A_s$ and $B_{s \smallfrown i} \supset B_s$). So the family $(D_s)_{s \in \omega^{< \omega}}$ is a Lusin scheme.

We know that $\mathrm{UB} = \bigcup_{\alpha \in \mcN} \bigcap_{k< \omega} (A_{\alpha_{|k}} \setminus B_{\alpha_{|k}}) \subseteq \bigcup_{\alpha \in \mcN} \bigcap_{k < \omega} D_{\alpha_{|k}}$.

Conversely, fix $\alpha \in \mcN$ and let $\mcT$ belong to $\bigcap_{k < \omega} D_{\alpha_{|k}}$. Then, for each $k$, $\mcT$ belongs to $\overline{A_{\alpha_{|k}}} \setminus B_{\alpha_{|k}}$. This implies that $\alpha_{|k} \in \mcT$ for each $k$, so that $\alpha \in [\mcT]$ and $\mcT \in A_{\alpha_{|k}} \setminus B_{\alpha_{|k}}$ for all $k$. This finally proves that 
\[\mathrm{UB} = \bigcup_{\alpha \in \mcN} \bigcap_{k< \omega} D_{\alpha_{|k}} \]
and we are done.
\end{proof}

Our motivation for proving the previous theorem was to obtain the following fundamental result.

\begin{thm}\label{thm:sets_of_uniqueness}
Let $X$, $Y$ be two Polish spaces and $B \subseteq X \times Y$ be Borel. Then 
\[B_{!}=\{x : B_x \textrm{ is a singleton} \}\]
is coanalytic.
\end{thm}

\begin{proof}
As usual we may assume that $X=Y=\mcN$. Consider first the case where $B$ is closed. Then the map $\varphi \colon \alpha \mapsto \mcT_B(\alpha)$ (the section tree of $B$ at $\alpha$) is continuous and $\varphi^{-1}(\mathrm{UB})=B_{!}$, whence $B_{!}$ is coanalytic.

We now turn to the general case. Assume $B$ is Borel; find a closed $C \subseteq \mcN$ and a continuous injective map $f \colon C \to \mcN \times \mcN$ such that $f(C)=B$. Denote by $\pi \colon \mcN \times \mcN \to \mcN$ the projection on the first coordinate, and consider $D=\{(\pi(f(c)), c) : c \in C\}$.

Since $f$ is continuous, $D$ is closed; thus the result of the previous paragraph implies that $\{\alpha \in \mcN : \exists! c \in C \ \alpha=\pi(f(c))\}$ is coanalytic. This set is precisely $B_!$.
\end{proof}

\begin{exo}\label{exo:sections_fixed_cardinal_coanalytic}
Let $X$, $Y$ be Polish spaces and $B$ a Borel subset of $X \times Y$. Prove that for each $n< \omega$ the set $\{x \in X : |B_x|=n\}$ is coanalytic.
\end{exo}

\begin{exo}\label{exo:uniquely_branching_complete_coanalytic}
Prove that $\mathrm{UB}$ is complete coanalytic.
\end{exo}

\begin{exo}\label{exo:cardinalities_of_coanalytic sets}
\begin{enumerate}
\item Let $\alpha$ be a countable ordinal. Prove that the set of all orderings on $\omega$ which are isomorphic to an initial segment of $\alpha$ is a Borel subset of $\mathrm{LO}$.
\item Let $X$ be a Polish space and $A \subseteq X$ be a coanalytic set. Show that one can write $A= \bigcup_{\alpha < \omega_1} A_\alpha$, where each $A_\alpha$ is a Borel subset of $X$.
\end{enumerate}
It follows from this that infinite coanalytic subsets of Polish spaces have three possible cardinalities : $\aleph_0$, $\aleph_1$ and $2^{\aleph_0}$ (of course if one believes in the continuum hypothesis the last two possible cardinalities are equal).
\end{exo}

\begin{thm}\label{thm:Borel_projection_countable_fibers}
Let $X$, $Y$ be Polish spaces and $B \subseteq X \times Y$ be a Borel subset such that each $B_x$ is at most countable. Denote $\pi_X \colon X \times Y \to X$ the coordinate projection. Then $\pi_X(B)$ is Borel.
\end{thm}

\begin{proof}
We deal with the case where $B$ is closed in $X \times Y$ (as in the proof of Theorem \ref{thm:sets_of_uniqueness}, one can reduce to that case by using that each Borel subset of a Polish space is a continuous injective image of a closed subset of $\mcN$).

We know that for every $x \in \pi_X(B)$ the set $B_x$ is closed and (at most) countable, whence it has an isolated point. Letting $(U_n)_{n< \omega}$ be a countable basis of open sets for the topology of $Y$ and denoting $B^n= B \cap (X \times U_n)$, we obtain that $\pi_X(B) = \bigcup_{n< \omega}B^n_!$, whence $\pi_X(B)$ is coanalytic. Since it is also analytic, Suslin's separation theorem tells us that $\pi_X(B)$ is Borel.
\end{proof}

Later on we will establish a strengthening of Theorem \ref{thm:Borel_projection_countable_fibers}: the Lusin--Novikov theorem (Theorem \ref{thm:Lusin_Novikov}) which states that, under the assumptions above, one may write $B= \bigcup_{n< \omega} B_n$, where each $B_n$ is the graph of a partial Borel function.

\begin{exo}[Mackey] \label{exo:Mackey}$ $

\begin{enumerate}
\item Let $X$ be a standard Borel space and $(A_n)_{n< \omega}$ be a countable family of Borel subsets of $X$ such that for every $x,y \in X$ there exists $n$ such that $x \in A_n$ and $y \not \in A_n$. Prove that $(A_n)_{n< \omega}$ generates the Borel $\sigma$-algebra on $X$.
\item Assume that $(G,\tau)$ is a Polish group, and that there exists a countable family $(A_n)_{n < \omega}$ of subsets of $G$ satisfying the following two conditions:
\begin{enumerate}
\item For any $x \ne y \in G$ there exists $n$ such that $x \in A_n$ and $y \not \in A_n$.
\item For any Polish group topology $\tau'$, $A_n$ is Borel in $(G,\tau')$.
\end{enumerate} 
Show that $\tau$ is the unique Polish group topology on $G$.
\end{enumerate}
\end{exo}

\emph{Comments.} Above, we saw a classical method to prove that a given subset $X$ of a Polish space $P$ is complex (e.g., analytic non Borel): find a ``nice'' (for instance, Borel) map $f \colon P' \to P$, where $P'$ is Polish, and a complex subset $Y$ of $P'$ such that for all $z \in P'$ one has $z \in Y \Leftrightarrow f(z) \in X$. Intuitively, by applying $f$, one reduces the problem of determining whether a given $z$ belongs to $Y$ to the problem of determining whether $f(z)$ belongs to $X$, so the membership problem for $X$ is more complicated than the corresponding problem for $Y$. This idea of reduction appears throughout descriptive set theory, and has had a strong influence on the subject's development. This motivates the search for more sets sitting at a given complexity level (complete (co)analytic for instance); see \cite{Kechris1995} for many examples.

The last exercise of the chapter illustrates a branch of research that has been very active over the past few decades: given a Polish group $G$, how canonical is the topology on $G$? Is it the only Polish group topology, or even, the unique second-countable Hausdorff group topology? Going further, what conditions on a Polish group $G$ imply that every morphism from $G$ to a separable group is continuous? For a discussion of this and related topics, see Rosendal's survey on automatic continuity of group homomorphisms \cite{Rosendal09a}.

\chapter{Definability and Baire measurability}
We remind the reader that the notion of Baire measurability has been discussed in the first chapter; in particular, we refer to that chapter for the definition of Baire measurability (Definition \ref{def:Baire_measurable}) and the proof of the Kuratowski--Ulam theorem (see Theorem \ref{t:splitting_category}).

Baire measurability is in many ways parallel to the usual notion of measurability (for instance, the Kuratowski--Ulam theorem is an analogue of the Fubini theorem). We discuss here how it relates to some of the notions of ``definability'' that we encountered earlier. First, we note the following fact, which we have used without proof in the first part of the text.

\begin{thm}\label{thm:Baire_measurable_continuous_on_dense_Gdelta}
Let $X$, $Y$ be Polish spaces and $f \colon X \to Y$ be a Baire measurable map. Then there exists a comeager subset $\Omega$ of $X$ such that $f_{|\Omega}$ is continuous.
\end{thm}
This result applies in particular to any Borel map.

\begin{proof}
Let $(U_n)_{n< \omega}$ be a countable basis of open sets for the topology of $Y$. For each $n$ there exists an open subset $O_n$ and a meager subset $M_n$ of $X$ such that $f^{-1}(U_n) \mathrel{\Delta} O_n = M_n$. Denote $\Omega= X \setminus \bigcup_{n< \omega} M_n$; it is a comeager subset of $X$.

For every $n$, we have that $f^{-1}(U_n) \cap \Omega = O_n \cap \Omega$ is open in $\Omega$, hence $f_{|\Omega} \colon \Omega \to X$ is continuous.
\end{proof}

\begin{thm}\label{thm:separate_continuity}
Let $X,Y,Z$ be Polish spaces and suppose $f \colon X \times Y \to Z$ is separately continuous. Then there exists a comeager subset $\Omega$ of $X$ such that $f$ is continuous at every point of $\Omega$.
\end{thm}

\begin{proof}
Fix compatible metrics on $X$, $Y$, $Z$. Given $n,k < \omega$ define
\[F_{n,k} =\{(x,y) : \forall u,v \in B(y,2^{-k}) \ d(f(x,u),f(x,v)) \le 2^{-n} \}. \]
Since $y \mapsto f(x,y)$ is continuous for all $x$, we have $X \times Y = \bigcap_n \bigcup_k F_{n,k}$. 

Now we prove that each $F_{n,k}$ is closed: assume $((x_i,y_i))_{i< \omega}$ is a sequence of elements of $F_{n,k}$ which converges to $(x,y) \in X \times Y$. Then pick $u,v \in B(y,2^{-k})$. For $i$ large enough we have $u,v \in B(y_i,2^{-k})$, whence $d(f(x_i,u),f(x_i,v)) \le 2^{-n}$. We can now let $i$ go to $+ \infty$ and conclude that $d(f(x,u),f(x,v)) \le 2^{-n}$.

Now note that since for all $k$ we have $X \times Y = \bigcup_k F_{n,k}$, the Baire category theorem gives us that $\bigcup_k \Int{F_{n,k}}$ is dense (and of course open) so $\Omega= \bigcap_n \bigcup_k \Int{F_{n,k}}$ is comeager. We now prove that every $(x,y) \in \Omega$ is a point of continuity for $f$.

Fix $(x,y) \in \Omega$ and let $(x_i,y_i) \xrightarrow[i \to + \infty]{} (x,y)$. Fix $n< \omega$, then $k< \omega$ such that $(x,y) \in \Int{F_{n,k}}$; for $i$ large enough $(x_i,y)$ belongs to $F_{n,k}$ and $d(y_i,y) < 2^{-k}$.  In particular, for any $i$ large enough we have
\[d(f(x_i,y_i),f(x,y)) \le d(f(x_i,y_i),f(x_i,y)) + d(f(x_i,y),f(x,y)) \le 2^{-n} + d(f(x_i,y),f(x,y)). \]
Since $f(x_i,y) \xrightarrow[i \to + \infty]{} f(x,y)$ we are done.
\end{proof}

With a bit more care we could have proved that there is a subset $\Omega$ of $X \times Y$ such that each $\{x : (x,y) \in \Omega\}$ is comeager and $f$ is continuous at every point of $\Omega$; see Theorem 8.51 of \cite{Kechris1995} (the proof given above is directly taken from that argument).

\begin{exo}\label{exo:separate_continuity_group_operations}
Let $G$ be a group endowed with a Polish topology. Assume that $(g,h) \mapsto gh$ is separately continuous. Prove that $G$ is a Polish group.\end{exo}

\begin{exo}\label{exo:Baire_measurable_product}
Let $X$, $Y$ be two Polish spaces. Assume that $A \subseteq X$ and $B \subseteq Y$ are Baire measurable. Show that $A \times B$ is a Baire measurable subset of $X \times Y$. 
\end{exo}

\begin{exo}\label{exo:no_Baire_measurable_well_ordering}
Assume that $\preceq$ is a well-ordering on $\R$; denote $R= \{(x,y) : x \preceq y \}$. Assume for a contradiction that $R$ is Baire measurable.
\begin{enumerate}
\item Let $I \subseteq X$ be a non meager, Baire measurable initial segment. Prove that $R \cap I^2$ is not meager.
\item Let $I \subseteq X$ be a non meager, Baire measurable initial segment. Show that there exists $x \in I$ such that $\{y: y \preceq x\}$ is not meager and Baire measurable.
\item Derive a contradiction.
\end{enumerate}
\end{exo}

One could prove similarly that a non principal ultrafilter on $\omega$, viewed as a subset of $2^\omega$, can never be Baire measurable. We see that sets that are often considered to be ``non-constructible'' (like a well-ordering of the reals, or a nonprincipal ultrafilter on $\omega$) cannot be Baire measurable.

If we think that a meager set is ``negligible'', and a comeager set is ``full'', the definition of Baire measurability asserts that a non negligible, Baire measurable subset if actually full in some nonempty open subset; this is extremely useful in practice, hence it is important to enlarge our catalog of Baire measurable subsets. This is what we turn to now.

We recall that, for $X$ a Polish space and $A$ a subset of $X$, we denote
\[U(A) = \bigcup \lset O \text{ open in } X \colon O \setminus A \text{ is meager} \rset .\]
The set $U(A) \setminus A$ is always meager, and $A$ is Baire measurable iff $A \setminus U(A)$ is also meager.

\begin{prop}
Let $X$ be a Polish space and $A$ a subset of $X$. There exists a Borel (hence, Baire measurable) subset $B$ of $X$ such that $A \subseteq B$ and, for any Baire measurable $B'$ containing $A$, the set $B \setminus B'$ is meager. We call such a $B$ an \emph{envelope}\index{envelope} of $A$.
\end{prop}

Of course if $A$ is Baire measurable then one can simply take $B=U(A)$, since $U(A) \setminus A$ is meager. But the property above is true for \emph{any} $A$.

\begin{proof}
Let $(U_n)_{n< \omega}$ be a basis for the topology of $X$. Let $I = \{i < \omega : A \textrm{ is meager in } U_i \}$. For each $i \in I$ fix a sequence $U_{n,i}$ of dense open subsets of $U_i$ such that $A \cap \bigcap_{n < \omega} U_{n,i}= \emptyset$. Then consider
\[ B = \bigcap_{i \in I} \bigcup_{n < \omega} (X \setminus U_{n,i} ) .\]
This is a Borel subset of $X$ which contains $A$ by definition.

Assume that $C \subseteq B \setminus A$ is non meager and Baire measurable. Then there exists some $i < \omega$ such that $C$ is comeager in some nonempty $U_i$; since $A$ is disjoint from $C$ we have $i \in I$. But then, since $C \subseteq B$ we have that $C$ is disjoint from $\bigcap_{n} U_{n,i}$ hence $C$ is meager in $U_i$. This is a contradiction.
\end{proof}

\begin{thm}[Lusin--Sierpinski]\index{Lusin--Sierpinski theorem}
Let $X$ be a Polish space and $A \subseteq X$ be analytic. Then $A$ is Baire measurable.
\end{thm}

\begin{proof}
We choose a continuous map $f \colon \mcN \to X$ such that $A= f(\mcN)$; for $s \in \omega^{< \omega}$ we denote $A_s= f(N_s)$. 

For each $s \in \omega^{< \omega}$ we choose a Borel envelope $B_s \supseteq A_s$ and denote $C_s = B_s \cap \overline{A_s}$, a Borel (hence Baire measurable) subset of $X$. 

By definition (and continuity) of $f$ we have both
\[A = \bigcup_{\alpha \in \mcN} \bigcap_{k< \omega} A_{\alpha_{|k}} \quad \textrm{ and } \quad A = \bigcup_{\alpha \in \mcN} \bigcap_{k< \omega} \overline{A_{\alpha_{|k}}} \quad \textrm{ , hence} \quad A= \bigcup_{\alpha \in \mcN} \bigcap_{k< \omega} C_{\alpha_{|k}}. \]

Since $C_s$ contains $A_s$ for each $s$, we have
\[C_s \setminus \bigcup_{n< \omega} C_{s \smallfrown n} \subseteq C_s \setminus \bigcup_{n< \omega} A_{s \smallfrown n}  \subseteq B_s \setminus A_s .\]
In particular, $C_s \setminus \bigcup_{n< \omega} C_{s \smallfrown n}$ is meager for all $s$. Hence the set  
\[M = \bigcup_{s \in \omega^{< \omega}} (C_s \setminus \bigcup_{n< \omega} C_{s \smallfrown n}) \]
is also meager. Denote $\displaystyle C= \bigcap_{n< \omega} \bigcup_{s \in \omega^n} C_s$.
This is a Borel subset of $X$.

Clearly, $A$ is contained in $C$. Conversely, choose $x \in C$. If for each $s$ such that $x \in C_s$ we also have $x \in \bigcup_{n< \omega} C_{s \smallfrown n}$ then we can inductively build $\alpha \in \mcN$ such that $x \in C_{\alpha_{|k}}$ for all $k$, whence $x \in A$. This shows that $C \setminus A$ is contained in $M$, which is meager. Finally we conclude that $C \mathrel{\Delta} A$ is meager; since $C$ is Borel it follows that $A$ is Baire measurable.
\end{proof}

The definition of an envelope extends to arbitrary measurable spaces, as follows.

\begin{defin}
Let $(X,\mcS)$ be a measurable space and $A \subseteq X$. We say that $B \in \mcS$ is an \emph{$\mcS$-envelope}\index{$\mcS$-envelope} for $A$ if $A \subseteq B$ and, for every $B' \in \mcS$ such that $A \subseteq B'$, every subset of $B \setminus B'$ belongs to $\mcS$.
\end{defin}

To understand this definition: say that $A$ is $\mcS$-small if every subset of $A$ belongs to $S$. These sets are stable under countable unions. Saying that $B \in \mcS$ is an $\mcS$-envelope for $A$ amounts to saying that for every other $B' \in \mcS$ containing $A$ the set $B \setminus B'$ is $\mcS$-small. Obviously $\mcS$-envelopes are not unique in general, but the symmetric difference of two $\mcS$-envelopes is $\mcS$-small.

\begin{defin}
Let $X$ be a set. We call a family $(A_s)_{s \in \omega^{<\omega}}$ of subsets of $X$ a \emph{Suslin scheme}\index{Suslin scheme}. Given such a family, we can apply to it the \emph{Suslin operation}\index{Suslin operation} $\mcA$:
\[ \mcA((A_s)_{s \in \omega^{< \omega}}) = \bigcup_{\alpha \in \mcN} \bigcap_{k < \omega} A_{\alpha_{|k}}. \]
\end{defin}
We encountered particular cases of this operation when working with Lusin schemes and studying properties of analytic sets. 

\begin{exo}\label{exo:Suslin_schemes_measurability}
Let $X$ be a Polish space.
\begin{enumerate}
\item Prove that each analytic set is obtained by applying the Suslin operation to a family of closed subsets of $X$.
\item Let $(A_s)_{s \in \omega^{< \omega}}$ be a family of analytic sets. Prove that there exists an analytic subset $B \subseteq X \times \mcN$ such that for all $x \in X$ one has $(x \in \mcA((A_s)_{s \in \omega^{< \omega}})) \Leftrightarrow (\exists \alpha \in \mcN \ (x,\alpha) \in B)$.
\end{enumerate}
It follows from this that the Suslin operation, applied to a family of analytic sets, produces an analytic set. Actually a stronger fact holds: given a family $\mcS$ of subsets of a set $X$, denote by $\mcA(\mcS)$ the family of sets obtained by applying Suslin's operation to a Suslin scheme of elements of $\mcS$. Then $\mcA(\mcA(\mcS))=\mcA(\mcS)$. We will not prove that here, see for instance 25.6 of \cite{Kechris1995} for a proof.

The argument in the proof given above that analytic sets are always Baire measurable actually shows that applying the Suslin operation to a family of analytic sets yields a Baire measurable subset. We now turn to a generalization of this fact, due to Szpilrajn-Marczewski. Assume that $(X,\mcS)$ is a measurable space and that every $A \subseteq X$ admits an $\mcS$-envelope. Let $(A_s)_{s \in \omega^{< \omega}}$ be a Suslin scheme of elements of $\mcS$. We want to prove that applying the Suslin operation to this scheme results in an element of $\mcS$.
\begin{enumerate}[resume]
\item Consider $\widetilde{A}_s= \bigcap_{t \sqsubseteq s} A_t$. Prove that $\mcA((\widetilde{A}_s)_{s \in \omega^{< \omega}})=\mcA(({A}_s)_{s \in \omega^{< \omega}})$. So we now assume that $A_t \subseteq A_s$ whenever $s \sqsubseteq t$.
\item For $s \in \omega^{< \omega}$ denote $B_s=\mcA((A_t)_{s \sqsubseteq t})$. Show that there exists an $\mcS$-envelope for $B_s$ contained in $A_s$, and use this to further reduce to the case where each $A_s$ is an $\mcS$-envelope for $\mcA((A_t)_{s \sqsubseteq t})$.
\item Further reduce to the case where for each $s$ one has $A_s= \bigcup_{i<\omega} A_{s \smallfrown i}$ and conclude.
\end{enumerate}
In particular, it follows from this result that the family of Baire measurable subsets of a Polish space is stable under the Suslin operation (Nikodym).
\begin{enumerate}[resume]
\item Let $X$ be standard Borel and $\mu$ be a complete, $\sigma$-finite Borel measure on $X$. Use the Szpilrajn-Marczewski theorem to prove that every analytic subset of $X$ is $\mu$-measurable (one says that analytic sets are \emph{universally measurable}\index{universally measurable set}).
\end{enumerate}
\end{exo}

The next theorem is another illustration of the nice interplay between Baire measurability and definability.

\begin{thm}\label{thm:Montgomery--Novikov}
Let $X$, $Y$ be Polish spaces and $A \subseteq X \times Y$ be Borel. Then for any open set $U \subseteq Y$ the set $\{x \in X : A_x \textrm{ is meager in } U\}$ is Borel.
\end{thm}

\begin{proof}
Let $\mcA$ be the family of all Borel subsets of $X  \times Y$ for which the property in the statement of the theorem holds. Let us show that this family is stable under complementation; assume that $A \in \mcA$, denote $B= X \times Y \setminus A$ and let $U$ be nonempty open in $Y$. Then $B_x$ is not meager in $U$ iff $B_x$ is comeager in some nonempty open $V$ contained in $U$ iff $A_x$ is meager in some nonempty open $V$ contained in $U$. Thus $\{x: B_x \textrm{ is not meager in } U\}$ is Borel and we are done (we can restrict $V$ to run over some fixed countable basis of open sets for the topology of $Y$).

Clearly every open subset is contained in $\mcA$, and if $(A_n)_{n< \omega}$ is a sequence of elements of $\mcA$ and $U$ an open set then $(\bigcup_n A_n)_x$ is meager in $U$ iff each $(A_n)_x$ is meager in $U$, hence $\mcA$ is stable under countable unions.
\end{proof}

The proof above straightforwardly adapts to show that if $(X,\mcB)$ is a measurable space, $Y$ is Polish and $A \subseteq X \times Y$ is measurable for $\mcB \times \mcB(Y)$ then for any open $U$ the set $\{x \in X : A_x \textrm{ is meager in } U\}$ belongs to $\mcB$ (this is due to Montgomery and Novikov).

\begin{exo}\label{exo:Vaught}
Assume that $G$ is a Polish group, $X$ is a standard Borel space and $G$ acts on $X$ in such a way that $(g,x) \mapsto gx$ is Borel. Let $A$ be a subset of $X$.
We define the \emph{Vaught transforms} $A^*=\{x : \forall^* g \in G \ gx \in A\}$ and $A^\Delta= \{x : \exists^* g \in G \ gx \in A\}$.
\begin{enumerate}
\item Show that both $A^*$ and $A^\Delta$ are $G$-invariant, and that $A$ is invariant iff $A=A^*=A^\Delta$.
\item Assume that $A$ is Borel. Show that both $A^*$ and $A^\Delta$ are Borel.
\item Assume that $X$ is Polish, $G \actson X$ is continuous and $A$ is comeager in $X$. Show that $A^*$ is comeager in $X$.
\end{enumerate}
\end{exo}

\bigskip \emph{Comments.} For more about Baire measurability, and its interplay with Lebesgue measurability, I cannot recommend Oxtoby's book \cite{Oxtoby} highly enough. The wonderful book \cite{Godefroy2022} covers many of the applications of Baire's methods in functional analysis (and, as a bonus, the reader might have to learn French in order to read this book).
The very recent \cite{BinghamOstaszewski} provides interesting further material.
It is also worth noting that it is consistent with the axioms of ZF, along with dependent choice, that  every subset of a Polish space is Baire measurable, see Shelah's \cite{Shelah1984}, following earlier work of Solovay \cite{Solovay70}.

\chapter{The Effros Borel space}\label{chapter:Effros}
We now discuss an important class of examples of standard Borel spaces, which are frequently used in invariant descriptive set theory to encode classes of mathematical objects (see the exercise concerning separable Banach spaces at the end of this chapter).

\begin{defin}
Let $X$ be a topological space. We denote by $\mcK(X)$ the space of all compact subsets of $X$, and endow it with the \emph{Vietoris topology}\index{Vietoris topology}, which is the topology generated by sets of the form 
\[\lset K \in \mcK(X) : K \subseteq U \textrm{ and } \forall i < n \ K \cap U_i \ne \emptyset\rset \]
where $U$ and (if $n \ge 1$) $U_0,\ldots,U_{n-1}$ are open subsets of $X$.
\end{defin}

Note that, by definition, $\emptyset$ is an isolated point of $\mcK(X)$, so it is common to exclude it and only consider the subset $\mcK^*(X)$ of all \emph{nonempty} closed subsets of $X$.

\begin{exo} \label{exo: Vietoris_topology_Hausdorff_distance}
Assume that $X$ is a metrizable topological space, and that the topology of $X$ is induced by a metric $d$. Show that the Vietoris topology is induced by the \emph{Hausdorff metric}\index{Hausdorff metric}
\[d_H(K,L)= \inf \lset \varepsilon \ge 0 : K \subseteq L_\varepsilon \textrm{ and } L \subseteq K_\varepsilon \rset \]
where $A_\varepsilon= \lset x \in X : d(x,A) \le \varepsilon \rset$.

(to make sense of this formula when $K$ or $L$ is empty, assume that $d$ is bounded by $1$ and let $\inf(\emptyset)=1$).
\end{exo}

\begin{lem}
Let $X$ be a metrizable space. The Borel structure of $\mcK(X)$ is generated by sets of the form $\lset K :K \cap U \ne \emptyset \rset$, where $U$ runs over all open subsets of $X$.
\end{lem}

\begin{proof}
Fix an open subset $U$ of $X$, and denote $U_n =\lset x \in X : d(x,X \setminus U) < 2^{-n} \rset$. Then, for $K \in \mcK(X)$, we have 
\[K \subseteq U \Leftrightarrow \exists n < \omega \ K \cap U_n = \emptyset .\]
The implication from right to left above is immediate, and the converse implication comes from the fact that if $K$ is compact, $F$ is closed and $K \cap F= \emptyset$ then $d(K,F) >0$.
\end{proof}

\begin{exo} \label{exo:Vietoris_topology_complete_separable}
Assume that $X$ is a Polish space and that $d$ is a compatible metric on $X$. Prove that the following facts hold.
\begin{enumerate}
\item If $D$ is a dense subset of $X$, then the set of all finite subsets of $D$ is dense in $\mcK(X)$. In particular, $\mcK(X)$ is separable.
\item The Hausdorff metric $d_H$ associated to $d$ is complete (so that $\mcK(X)$, endowed with the Vietoris topology, is Polish).
\item If $X$ is compact then $\mcK(X)$ is compact.
\item The maps $(K,L) \mapsto K \cup L$ and $(K,L) \mapsto K \cap L$ are Borel.
\end{enumerate}
\end{exo}

In the next two exercises we recover the result of Theorem \ref{exo:perfect_set_property_analytic}  as well as a strengthening of Theorem \ref{thm:Mycielski}.

\begin{exo} \label{exo:alternative_cantor_scheme_Vietoris}
\begin{enumerate}
\item Let $X$ be a perfect Polish space. Prove that $\{K \in \mcK(X) : \mcK \textrm{ is perfect } \}$ is dense $G_\delta$ in $\mcK^*(X)$.
\item Let $X$, $Y$ be Polish spaces and $f \colon X \to Y$ a continuous map such that, for any nonempty open $U \subseteq X$, $f_{|U}$ is not constant. Assume also that $X$ is perfect. Prove that 
\[\{K \in \mcK^*(X) : f_{|K} \textrm{ is injective} \} \]
is dense $G_\delta$ in $\mcK^*(X)$.
\item Let $X$, $Y$ be Polish spaces and $f \colon X \to Y$ be a Borel map such that $f(X)$ is uncountable. Prove that there exists a homeomorphic copy $K \subseteq X$ of $\mcC$ such that $f_{|K}$ is injective\footnote{a fact that we already proved in Theorem \ref{exo:perfect_set_property_analytic} for continuous $f$.}. 
\end{enumerate}
\end{exo}

\begin{exo}\label{exo:Cantor_independent_over_Q} Let $X$ be a perfect Polish space.
\begin{enumerate}
\item Assume that for each $n < \omega$, $(R_n)_{n \ge 1}$ is a comeager subset of $X^n$. Prove that there exists a continuous map $f \colon \mcC \to X$ such that for any $n$ and any pairwise distinct $x_1,\ldots,x_n \in \mcC$ we have $(f(x_1),\ldots,f(x_n)) \in R_n$ (for $n=2$ we recover Mycielski's theorem).
\item Show that there is a Cantor subset $K$ of $\R$ which is linearly independent over $\Q$.
\end{enumerate}
\end{exo}

We turn to the main definition of this chapter.

\begin{defin}
Let $X$ be a topological space. We denote by $\mcF(X)$ the set of all closed subsets of $X$, and endow it with the \emph{Effros Borel structure}\index{Effros Borel structure}, which is the $\sigma$-algebra generated by sets of the form 
\[\lset F \in \mcF(X) : F \cap U \ne \emptyset \rset \]
where $U$ runs over all open subsets of $X$.
\end{defin}

\begin{exo}\label{exo:belonging_is_Borel_in_Effros}
Let $X$ be a Polish space. Prove that $\{(F,x): x \in F\}$ is a Borel subset of $\mcF(X) \times X$.
\end{exo}

Note that when $X$ is compact we have $\mcF(X)=\mcK(X)$ and if $X$ is moreover metrizable then the Effros Borel structure on $\mcK(X)$ is the Borel structure induced by the Vietoris topology on $\mcK(Y)$.

\begin{thm}
Let $X$ be a Polish space. Then $\mcF(X)$, endowed with the Effros Borel structure, is a standard Borel space.
\end{thm}

\begin{proof}
We begin by choosing some metrizable proper compactification $(Y,d)$ of $X$ and view $X$ as a subspace of $Y$. 

Then we claim that the map $\Phi \colon \mcF(X) \to \mcK(Y)$ defined by $\Phi(F)= \overline{F}$ is an isomorphism from $\mcF(X)$ endowed with the Effros Borel structure onto $\Phi(\mcF(X))$ endowed with the Borel structure induced by the Vietoris topology.

To see this, observe first that $\Phi$ is injective. Then note that if $U \subseteq Y$ is open then $V=U \cap X$ is open in $X$ and, for any $F \in \mcF(X)$, one has $\Phi(F) \cap U \ne \emptyset \Leftrightarrow F \cap V \ne \emptyset$. This shows that $\Phi$ is Borel. Similarly, given $V$ open in $X$ there exists $U$ open in $Y$ such that $V=U \cap X$ and then for any $F \in \mcF(X)$ we have $F \cap V \ne \emptyset \Leftrightarrow \Phi(F) \cap U \ne \emptyset$, whence $\Phi^{-1}$ is Borel.

To complete the proof, it remains to show that $\Phi(X) = \lset \overline{F} : F \in \mcF(X) \rset$ is Borel in $\mcK(Y)$. Note that, for $K \in \mcK(Y)$, $F=K \cap X$ is closed in $X$ and $K=\overline{F}$ iff $K \cap X$ is dense in $K$.  Since $X$ is Polish, it is a $G_\delta$ subset of $Y$, so there exists a sequence of open subsets $(U_n)_{n< \omega}$ of $Y$ such that $X= \bigcap_{n< \omega} U_n$. Applying the Baire category theorem, we see that, for any compact subset $K$ of $Y$, we have that $X \cap K$ is dense in $K$ iff $U_n \cap K$ is dense in $K$ for all $n< \omega$. Fix a basis $(V_n)_{n< \omega}$ of open sets for the topology of $Y$. It follows from our previous observation that for all $K \in \mcK(Y)$ one has
\[K \in \Phi(\mcF(X)) \Leftrightarrow \forall n, m < \omega \left((K \cap V_m) \ne \emptyset \Rightarrow (K \cap U_n \cap V_m  \ne \emptyset) \right). \]
For each fixed $n,m$ the subset $\lset K \in \mcK(X) : K \cap V_m = \emptyset \rset$ is closed in $\mcK(X)$, and the subset
$\lset K  \in \mcK(X) : K \cap U_n \cap V_m \ne \emptyset \rset$ is open; so their union is Borel, and we conclude that $\Phi(\mcF(X))$ is a Borel (actually, $G_\delta$) subset of $\mcK(X)$.
\end{proof}

It might be worth pointing out that it follows from the previous proof that the topology on $\mcF(X)$ induced by the topology of $\mcK(Y)$ is Polish.

\begin{prop}
Let $X$ be a Polish space. Then the map $(K,L) \mapsto K \cup L$ is a Borel map from $\mcF(X) \times \mcF(X)$ to itself.

The map $(K,L) \mapsto K \cap L$ is not Borel in general.
\end{prop}

\begin{proof}
Let $(Y,d)$ be a metrizable proper compactification of $X$. We know that $ \Phi \colon F \mapsto \overline{F}$ (where the closure is computed in $Y$) is a Borel isomorphism from $\mcF(X)$ onto its image. Let $U$ be an open subset of $X$, and $V$ be an open subset of $Y$ such that $U=V \cap Y$. Then for any two closed subsets $K$, $L$ of $X$ we have that 
\[ \left( (K \cup L) \cap U \ne \emptyset \right) \Leftrightarrow \left( \overline{K \cup L} \cap V \ne \emptyset \right) \Leftrightarrow \left( (\overline{K} \cup \overline{L}) \cap V \ne \emptyset \right)\]
and we conclude that $(K,L) \mapsto K \cup L$ is Borel by using that forming unions is a Borel map on $\mcK(Y) \times \mcK(Y)$.

To see that $(K,L) \mapsto K \cap L$ is not Borel in general, Let $B$ denote the closed subset of $\mcN$ defined by 
\[\alpha \in B \Leftrightarrow \forall k < \omega \  \alpha(2k)=0.\]
Then consider the map $\varphi \colon \omega^{< \omega} \to \omega^{< \omega}$ defined by $\varphi(s)=(0,s(0),\ldots,0,s(n-1))$ (where $n=|s|$); for a tree $\mcT$ on $\omega$ let $\pi(\mcT)= \{s \in \omega^{< \omega} : \exists t \in \mcT \ s \sqsubseteq\varphi( t) \}$. By definition $\pi(\mcT)$ is a tree, and it is straightforward to check that $\pi \colon \mathrm{Tr}(\omega) \to \mathrm{Tr}(\omega)$ is continuous.

Since $\mcT \mapsto [\mcT]$ is a Borel map from $\mathrm{Tr}(\omega)$ to $\mcF(\mcN)$, the map $\psi \colon \mcT \mapsto [\pi(\mcT)]$ is Borel. And by definition $\mcT \in \mathrm{IF}$ iff $\psi(\mcT) \cap B \ne \emptyset$. This proves that $\{F \in \mcF(\mcN) : F \cap B \ne \emptyset \}$ is not a Borel subset of $\mcF(\mcN)$, since otherwise $\mathrm{IF}$ would be Borel (this set is however clearly analytic, since $\{(F,x) : x \in F\}$ is a Borel subset of $\mcF(X) \times X$).
\end{proof}

\begin{exo}\label{exo:compact_Borel_in_Effros}
Let $X$ be a Polish space. Show that $\mcK(X)$ is a Borel subset of $\mcF(X)$.
\end{exo}

\begin{thm}[Kuratowski--Ryll-Nardzewski]\index{Kuratowski--Ryll-Nardzewski theorem}
Let $X$ be a Polish space. Then there exists a sequence of Borel functions $s_n \colon \mathcal{F}(X) \to X$ such that, for any nonempty $F \in \mcF(X)$, $\{s_n(F): n < \omega \}$ is a dense subset of $F$.
\end{thm}

\begin{proof}
As is often the case, the proof of this theorem morally reduces to the case where $X=\mcN$ via the use of an appropriate transfer theorem.
We assume that $X$ is nonempty and denote $\mcF^*(X)$ the set of all nonempty closed subsets of $X$. We begin by building a Borel map $s \colon \mcF^*(X) \to X$ such that $s(F) \in F$ for all $F \in \mcF^*(X)$.

We start by fixing a continuous, open surjection $\pi \colon \mcN \to X$ and denote $U_s = \pi(N_s)$ for $s \in \omega^{< \omega}$. 

For each $F \in \mcF^*(X)$ the set $\pi^{-1}(F)$ is closed nonempty in $\mcN$. We let $\alpha_F$ be the \emph{leftmost branch} of the corresponding pruned tree $\mcT_F$ on $\mcN$, i.e., the lexicographically smallest element of $[\mcT_F]$ (a precise definition is given below). We then want to set $s(F)=\pi(\alpha_F)$, and for that to suit our purposes we have to prove that this function $s$ is Borel. It is enough to prove that $F \mapsto \alpha_F$ is Borel, i.e., that its graph is a Borel subset of $\mcF^*(X) \times \mcN$. 

Given $s \ne t \in \omega^n$ say that $s <_{\textrm{lex}} t$ if $s(i)< t(i)$, where $i$ is smallest such that $s(i) \ne t(i)$. Then the leftmost branch $\alpha_F$ mentioned above is characterized as follows: for any $\alpha \in \mcN$ and any $F \in \mcF^*(X)$ we have
\[ (\alpha= \alpha_F) \Leftrightarrow \left( \forall n ((\forall s \in \omega^n \ s <_{\textrm{lex}} \alpha_{|n} \Rightarrow U_s \cap \pi^{-1}(F) = \emptyset) \textrm{ and } (U_{\alpha_{|n}} \cap \pi^{-1}(F) \ne \emptyset)) \right) .\]
It then follows from the definition of the Effros Borel structure on $X$ (and the fact that $\pi$ is open) that $\{(F,\alpha) : \alpha= \alpha_F\}$ is a Borel subset of $\mcF^*(X) \times X$. Thus $s$ is Borel; and for every $F \in \mcF^*(X)$, we have that $s(F)=\pi(\alpha_F) \in F$.

So far, we have described a way to choose, in a Borel way, an element in any nonempty closed subset of $X$. To obtain the general result, fix a countable basis $(V_n)_{n< \omega}$ for the topology of $X$. Then set, for any $n < \omega$,
\[s_n(F) = \begin{cases} s(F) & \textrm{ if } F \cap V_n = \emptyset , \\ s(\overline{F \cap V_n}) & \textrm{ if } F \cap V_n \ne \emptyset . \end{cases} \]
To see why this is Borel, note that by the map $\mcF(X) \to \mcF(V_n)$, $F \mapsto F \cap V_n$ is Borel, and (by the argument we used in the proof that $\mcF(X)$ is a standard Borel space when $X$ is Polish) $A \mapsto \overline{A}$ is a Borel map from $\mcF(V_n)$ to $\mcF(X)$.
\end{proof}

Let us give an example of application of the Kuratowski--Ryll-Nardzewski selection theorem.

\begin{thm}\label{thm:Burgess}
Let $X$ be a Polish space, and $E$ an equivalence relation on $X$ with closed equivalence classes and such that, for every open subset $U$ of $X$, the set  $U_E =\bigcup_{x \in U} [x]_E$ is Borel (where $[x]_E$ denotes the equivalence class of $x$). Then $E$ admits a \emph{Borel selector}\index{Borel selector}, i.e., there exists a Borel map $s \colon X \to X$ such that for $x, y \in X$ one has $s(x) \mathrel{E} x$ and $(x \mathrel{E} y ) \Rightarrow (s(x)=s(y))$.
\end{thm}

\begin{proof}
We first observe that the map $\pi \colon x \mapsto [x]_E$ is a Borel map from $X$ to $\mathcal F(X)$. Indeed, given a nonempty open subset $U$ of $x$, note that 
\[(U \cap [x]_E \ne \emptyset) \Leftrightarrow (x \in U_E),  \]
and our assumption is precisely that each $U_E$ is Borel. We then conclude by applying the Kuratowski--Ryll-Nardzewski selection theorem to obtain a Borel selector for $E$.
\end{proof}

\begin{exo}\label{exo:Burgess}
Let $G$ be a a Polish group, and $H$ be a closed subgroup of $G$. Prove that there exists a Borel subset $T$ of $G$ which intersects each left $H$-coset in exactly one point.
\end{exo}

Effros Borel spaces are commonly used to encode various classes of mathematical structures, for instance as a way to provide a setting to study isomorphism relations.

Let us simply give an example.

\begin{exo}\label{exo:standard_space_of_Banach_spaces}
\begin{enumerate}
\item Using the fact that the Cantor space $\mcC$ maps onto every nonempty compact metrizable space, show that every separable real Banach space embeds linearly isometrically into the space $E$ of all continuous functions from $\mcC$ to $\R$ (endowed as usual with the sup norm $\|\cdot\|_\infty$).
\item Show that $\mcB=\{F \in \mathcal F(E) : F \textrm{ is a vector subspace of } E\}$ is Borel in $\mcF(E)$.  We call $\mcB$ the Borel space of all separable Banach spaces.
\item Prove that 
\[\{F \in \mcB : F \textrm{ is isometric to a Hilbert space }\} \textrm{ and } \{F \in \mcB : F \textrm{ is finite dimensional}\}\]
are Borel in $\mcB$.
\end{enumerate}
Many interesting classes of separable Banach spaces turn out to be non Borel, this is for instance the case of Banach spaces with separable dual and reflexive separable Banach spaces (though both are coanalytic).
\end{exo}

\emph{Comments.} Given a complete metric $d$ on a Polish space $X$, Beer \Cite{Beer1991} has shown that the \emph{Wijsman topology}, i.e., the topology induced by the maps $F \mapsto d(x,F)$, is a Polish topology on 
the set $\mcF^*(X)$ of nonempty closed subsets of $X$, which induces the Effros Borel structure on $\mcF^*(X)$; see also \cite{Zsilinszky98}.

As we mentioned above, Effros Borel spaces are often used in invariant descriptive set theory, where one is concerned with structure theorems for various classes of equivalence relations on Polish spaces; see \cite{Gao2009a} and \cite{Hjorth2000}. In a different direction, the book \cite{Dodos2010} discusses some of the interactions between descriptive set theory and Banach space theory, in which Effros Borel spaces are used extensively (and we again mention \cite{Godefroy2022}, which discusses some related topics).

\chapter{$\mcG_0$-dichotomy and some applications}
We turn to a discussion of a fundamental dichotomy theorem proved by Kechris, Solecki and Todor\v{c}evic at the turn of the millennium. In this chapter we explain the dichotomy and discuss some of its corollaries; in the next chapter we present a relatively recent proof due to B. Miller.

\begin{defin}
Let $X$ be a set. A \emph{graph}\index{graph} on $X$ is a binary relation which is symmetric and irreflexive. If $X$ is a Polish space, we view $\mcG$ as a subset of $X \times X$ and say that it is a \emph{Borel graph} if $\mcG$ is a Borel subset of $X \times X$. We analogously define analytic graphs on Polish spaces.
\end{defin}

Given a graph $\mcG$ on a set $X$ and $A \subseteq X$, the associated \emph{induced subgraph}\index{induced subgraph} is $\mcG_{|A}= \mcG \cap (A \times A)$ (this corresponds to the notion of substructure we considered in Chapter \ref{chapter:Fraïssé}).

\begin{defin}
Let $X$ be a set, $\mcG$ a graph on $X$ and $A$ a subset of $X$. We say that $A$ is $\mcG$-\emph{independent}\index{$\mcG$-independent} if $\mcG_{|A}= \emptyset$ (that is, no two points of $A$ are adjacent in $\mcG$).

A \emph{coloring}\index{coloring of a graph} is a function $c \colon X \to Y$, where $Y$ is a set and $c^{-1}(\{y\})$ is $\mcG$-independent for each $y \in Y$ (that is, two elements of $\mcG$ which are connected by an edge must be assigned distinct colors).
\end{defin}

Of course we are concerned with definable versions of these combinatorial notions; we are going to focus exclusively on Borel $\omega$-colorings of graphs on Polish spaces, i.e., colorings $c \colon X \to \omega$ such that $X$ is Polish and $c$ is Borel.

Note that if $\mcG$ is a graph on a Polish space $X$, $\mcH$ is a graph on a Polish space $Y$, and $\varphi$ is a Borel homomorphism from $\mcG$ to $\mcH$ (i.e., $\varphi \colon X \to Y$ is Borel and maps edges in $\mcG$ to edges in $\mcH$) then any Borel $\omega$-coloring $c$ of $\mcH$ gives rise to a Borel $\omega$-coloring $c \circ \varphi$ of $\mcG$.

\begin{defin}\label{defin:G_s}
Let $S$ be a subset of  $2^{< \omega}$. We define the graph $\mcG_S \subset 2^\omega \times 2^\omega$ by saying that $(x,y) \in \mcG_S$ if, and only if, there exists some $s \in S$, $\varepsilon \in \{0,1\}$ and $z \in 2^\omega$ such that $x=s \smallfrown \varepsilon \smallfrown z$ and $y= s \smallfrown (1- \varepsilon) \smallfrown z$. 
\end{defin}

Note that $\mcG_S$ is by definition a Borel subset of $2^\omega \times 2^\omega$.

Given any graph $G$ on a set $X$, we can consider the transitive closure $E_G$ of $G$, which is an equivalence relation; two points $x,x'$  are in the same equivalence class for $E_G$ iff there is a path in $G$ linking $x$ to $x'$. The equivalence relation $E_0$ on $2^\omega$ is defined as follows:
\[ \forall x,y \in 2^\omega \quad ( (x,y) \in E_0) \Leftrightarrow (\exists n \ \forall m \ge n \ x(m)=y(m) ) .\]
For any $(x,y) \in \mcG_S$ there exists exactly one $n$ such that $x(n) \ne y(n)$ so $(x,y) \in E_0$.

\begin{exo}\label{exo:graph_generates_E_0}
Assume that $S \subseteq 2^{< \omega}$ is such that $|S \cap 2^n| \ge 1$ for all $n$. Show that $E_{\mcG_S}=E_0$.
\end{exo}

Intuitively, the relation $E_0$ is ``the simplest nontrivial Borel equivalence relation''. We chose not to go into the theory of Borel reducibility of definable equivalence relations in these notes and refer the reader to the references mentioned at the end of this chapter for more information on this topic, which is a cornerstone of modern descriptive set theory.

Note that each equivalence class of $E_0$ is dense.

\begin{exo}\label{exo:generic_ergodicity_E_0}
Let $X$ be a Polish space, and $f \colon 2^{\omega} \to X$ be a Baire measurable function such that for all $\alpha, \beta \in 2^\omega$ one has $(\alpha,\beta) \in E_0 \Rightarrow f(\alpha)=f(\beta)$. Prove that $f$ is constant on a comeager set.
\end{exo}

\begin{defin}
We say that a subset $S$ of $2^{< \omega}$ is \emph{dense} if for every $t \in 2^{< \omega}$ there exists $s \in S$ which extends $t$. We say that $S$ is \emph{sparse} if $|S \cap 2^n|=1$ for all $n$.
\end{defin}

\begin{lem}
Assume that $S$ is a dense subset of $2^{< \omega}$, and let $A$ be a Baire measurable, non meager subset of $2^\omega$. Then $A$ is not $\mcG_S$-independent.

In particular, if $S$ is dense then $\mcG_S$ does not admit a Borel $\omega$-coloring.
\end{lem}

\begin{proof}
Since $A$ is Baire measurable and non meager, there exists some $t \in 2^{< \omega}$ such that $A$ is comeager in $N_t$; and since $S$ is dense it follows that $A$ is comeager in $N_s$ for some $s \in S$. Denote $n=|s|$ and let $f \colon N_s \to N_s$ be the map from $N_s$ to itself defined by flipping $x(n)$, $i.e.$, $s\smallfrown \varepsilon \smallfrown z \mapsto s \smallfrown (1- \varepsilon) \smallfrown z$. This is a homeomorphism of $N_s$, so $A$, $f(A)$ are both comeager in $N_s$, hence $A \cap f(A) \ne \emptyset$. It follows that there exists $z \in 2^\omega$ such that $s \smallfrown 0 \smallfrown z$ and $s \smallfrown 1 \smallfrown z$ both belong to $A$, whence $A$ is not $\mcG_S$-independent.

If $c \colon 2^\omega \to \omega$ were a Borel (or even, Baire measurable) coloring, then for each $i \in \omega$ we would have that $c^{-1}(\{i\})$ is Baire measurable and $\mcG_S$-independent, hence meager. But then $2^\omega$ would be meager, a contradiction.
\end{proof}

\begin{exo}\label{exo:it_is_possible_to_be_dense_and_sparse}
Give a construction of a subset $S$ of $2^{< \omega}$ which is dense and sparse.
\end{exo}

\begin{thm}[Kechris--Solecki--Todor\v{c}evic \cite{KST}; the $\mcG_0$-dichotomy theorem]\index{$\mcG_0$-dichotomy theorem}
Let $S \subseteq 2^{< \omega}$ be dense and sparse, and $\mcG$ be an analytic graph on a Polish space $X$. Then exactly one of the following conditions holds:
\begin{enumerate}
\item The graph $\mcG$ admits a Borel $\omega$-coloring.
\item There exists a continuous homomorphism $\varphi \colon \mcG_S \to \mcG$.
\end{enumerate}
\end{thm}

Note that these conditions are mutually exclusive because $S$ is dense; actually, in the next chapter we are going to prove that, assuming only that $S$ is sparse, at least one of the two properties above holds whenever $G$ is an analytic graph on a Polish space $X$ (but they are not mutually exclusive in general). An interesting consequence of the statement above is that, if there exists a Borel homomorphism $\varphi \colon \mcG_S \to \mcG$ then there also exists a continuous homomorphism from $\mcG_S$ to $\mcG$.

In what follows, we fix a dense sparse set $\mcS$ and denote by $\mcG_0$ the associated graph $\mcG_S$.

Let us discuss some important theorems that can (as observed by B.Miller) be deduced from the $\mcG_0$-dichotomy theorem (see \Cite{Miller2012} for some related examples).

\begin{thm}[Lusin--Novikov]\index{Lusin--Novikov theorem}\label{thm:Lusin_Novikov}
Let $X, Y$ be two Polish spaces, and $A \subset X \times Y$ be Borel and such that for each $x$ the vertical section $A_x = \{y : (x,y) \in A\}$ is at most countable. Then $A=\bigcup_{n< \omega} A_n$, where each $A_n$ is the graph of a Borel function from a Borel subset of $X$ to $Y$.

It follows that the projection of $A$ on $X$ is Borel.
\end{thm}

We have already established that, under the above assumptions, the projection of $A$ on $X$ is Borel, see Theorem \ref{thm:Borel_projection_countable_fibers}. The stronger result we obtain here could also be deduced from \ref{thm:sets_of_uniqueness} with some additional work.

\begin{proof}
Consider the graph $\mcH$ on $X \times Y$ defined by setting 
\[((x,y) \mathrel{\mcH} (x',y')) \Leftrightarrow (x=x', y \ne y' \textrm{ and } y, y' \in A_x ). \]
This is a Borel graph on $X \times Y$, so the $\mcG_0$-dichotomy theorem applies. If there exists a Borel coloring $c \colon X \times Y \to \omega$ of $\mcH$ then for each $i < \omega$ the set $c^{-1}(\{i\})$ is $\mcH$-independent. In particular, the projection $\pi_X$ on the first coordinate is injective on $A_i= A \cap c^{-1}(\{i\})$. It follows that $X_i=\pi_X(A_i)$ is a Borel subset of $X$ and that $A_i$  is the graph of a Borel function with domain $X_i$. Since $A= \bigsqcup_{i< \omega} A_i$, we are done in that case; one can also observe that $\pi_X(A)= \bigcup_{i < \omega} X_i$ is Borel.

We now want to prove that there is no continuous homomorphism from $\mcG_0$ to $\mcH$. Assume for a contradiction that $\varphi$ is such a homomorphism. Denote $\varphi(\alpha)=(\varphi_X(\alpha),\varphi_Y(\alpha))$. By definition, 
\[\forall \alpha,\beta \in 2^\omega \ \alpha \mathrel{\mcG_0} \beta \Rightarrow \varphi_X(\alpha)=\varphi_X(\beta) .\]
Since $\varphi_X$ is continuous, and $E_0$ has dense  equivalence classes, it follows that $\varphi_X$ is constant on $2^\omega$; so there exists $x_0 \in X$ such that for all $(\alpha,\beta) \in 2^\omega \times 2^\omega$ one has $\varphi(\alpha,\beta)=(x_0,\varphi_Y(\beta))$.

But then $\varphi_Y$ is a continuous homomorphism from $\mcG_0$ to the graph $\mcH_{x_0}$ on $Y$ where
\[(y,y') \in \mcH_{x_0} \Leftrightarrow (y \ne y' \textrm{ and } y,y' \in A_{x_0}) \]
Since at most countably elements of $Y$ are adjacent in $\mcH_{x_0}$ to another element of $Y$ (only the elements of $A_{x_0}$ are connected by an edge to someone), it is immediate that $\mcH_{x_0}$ has a Borel $\omega$-coloring. This is a contradiction, which concludes the proof.
\end{proof}

\begin{exo}\label{exo:Lusin_Novikov_section_fixed_cardinality}
Let $X$, $Y$ be two Polish spaces, $A \subset X \times Y$ be Borel and such that for each $x$ the set $A_x = \{y : (x,y) \in A\}$ is at most countable. Use the Lusin--Novikov theorem to show that for each $n \le \omega$ the set $\{x : |A_x|=n\}$ is a Borel subset of $X$.
\end{exo}

\begin{coro}
Let $X$, $Y$ be two Polish spaces and $f \colon X \to Y$ a function such that each $f^{-1}(\{y\})$ is at most countable. Then $f(X)$ is a Borel subset of $Y$ and there exists a Borel map $s \colon s(X) \to X$ such that $f(s(y))=y$ for all $y \in s(X)$.
\end{coro}

\begin{proof}
Applying the Lusin--Novikov theorem to $A=\{(y,x) : y = f(x) \}$ we obtain that $f(X)$ is Borel. We can also write  $A= \bigcup_{i< \omega} A_i$, where each $A_i$ is the domain of a Borel function $s_i \colon Y_i \to X$, with each $X_i$ a Borel subset of $Y$.

We can now define $s$ by setting $s(y)=s_i(y)$ for every $y \in Y_i \setminus \bigcup_{j< i} Y_j$. We obtain as desired that $\mathrm{dom}(s)= \pi_Y(A)=f(X)$, and for each $y \in f(X)$ we have $f(s(y))=y$. 
\end{proof}

The argument used above shows that every Borel set $A \subseteq X \times Y$ with countable vertical sections admits a Borel \emph{uniformization}\index{uniformization}, i.e., there exists a Borel graph $A' \subseteq X \times Y$ contained in $A$ and such that $\pi_X(A')=\pi_X(A)$. Theorems justifying the existence of uniformizations are an important part of classical descriptive set theory which we will not tackle here. In general one cannot hope for the existence of a Borel uniformization; for instance there exists a closed subset $F \subseteq \mcN \times \mcN$ whose projection on the first coordinate is equal to $\mcN$, but which has no Borel uniformization (see e.g.~exercise 18.17 of \cite{Kechris1995}).

\begin{exo}\label{exo:Lusin_Novikov_injective}
Let $X$ be a Polish space, and $R \subseteq X \times X$ be a nonempty Borel subset of $X$ which is symmetric and such that $\{y : (x,y) \in R\}$ is countable for each $x \in X$. Prove that there exists a sequence of injective Borel functions $f_n \colon A_n \to B_n$, where $A_n$, $B_n$ are Borel subsets of $X$, such that $R= \bigsqcup_{n< \omega} \{(x,y) : y =f_n(x) \}$.
\end{exo}

\begin{thm}[Feldman--Moore]\index{Feldman--Moore theorem}
Let $X$ be a Polish space, and $E$ be a Borel equivalence relation on $X$ such that each equivalence class $[x]_E$ is at most countable. Then there exists a Borel action $G \actson X$ of a countable group $G$ such that $E$ is the equivalence relation induced by that action.
\end{thm}

Note that in the statement above we do not require the action to be free, which is in general impossible (even if all equivalence classes are infinite); so we might as well require that $G$ is the free group on countably many generators.

\begin{proof}
For each $x \in X$, the section $E_x$ is at most countable, so applying the Lusin--Novikov theorem (or rather, its corollary stated in Exercise \ref{exo:Lusin_Novikov_injective}) we find a sequence $(f_n)_{n< \omega}$ of injective partial Borel functions $f_n \colon A_n \to B_n$ whose graphs cover $E$. We may also assume that $f_0=\mathrm{id}$ and for each $n$ either $f_n= \mathrm{id}$ or $f_n(x) \ne x$ for all $x \in A_n$.

Let $A=\{n : f_n \ne \mathrm{id} \}$; fix a countable basis $(U_i)_{i< \omega}$ of the topology of $X$. For each $k \in A$ and each $i,j < \omega$ such that $U_i \cap U_j = \emptyset$ define $V^k_{i,j}=U_i \cap f_k^{-1}(U_j)$ and $W^k_{i,j} = f_k(V^k_{i,j})$. These two sets are Borel, disjoint (perhaps empty!) and we may then define a Borel involution $g^k_{i,j}$ by mapping each $x \in V^k_{i,j}$ to $f_k(x)$, and each $y \in W^k_{i,j}$ to $f_k^{-1}(y)$. 

Let $G$ denote the countable group of Borel automorphisms of $X$ generated by the involutions $g^k_{i,j}$. Then for each $x \in X$ and each $g \in G$ we have that $(x,g(x)) \in E$ since the graph of each $g^k_{i,j}$ is contained in $E$. Conversely, assume that $(x,y) \in E$ and $x \ne y$. Then there exist $k \in A$ such that $y= f_k(x)$, and $i,j < \omega$ such that $U_i \cap U_j = \emptyset$, $x \in U_i$ and $y \in U_j$. We thus have $y=f^k_{i,j}(x)$.
\end{proof}

One can also use the $\mcG_0$-dichotomy theorem to characterize cardinalities of quotient spaces associated to coanalytic equivalence relations on Polish spaces.

\begin{thm}[Silver]\index{Silver's theorem}
Let $X$ be a Polish space, and $E$ be a coanalytic equivalence relation on $X$. Then either $E$ has at most countably many classes or there exists a continuous map $f \colon 2^\omega \to X$ such that for any $\alpha \ne \beta$ one has $[f(\alpha)]_E \neq [f(\beta)]_E$ (in particular, the quotient space $X \lcoset E$ is either at most countable or of cardinality continuum).
\end{thm}

\begin{proof}
Let $\mcG= (X \times X) \setminus E$, which is an analytic graph on $X$. If $c \colon X \times X \to \omega$ is a Borel $\omega$-coloring of $\mcG$ then for each $x,y \in X$ such that $[x]_E \ne [y]_E$ we must have $c(x) \ne c(y)$ since $(x,y) \in \mcG$, whence $E$ has at most countably many classes.

So we may assume that there exists a continuous homomorphism $\varphi$ from $\mcG_0$ to $\mcG$. Consider 
$E'=(\varphi \times \varphi)^{-1}(E)$. This is a coanalytic (hence Baire measurable) equivalence relation on $2^\omega$.
By definition, for each $\alpha$ the set $E'_\alpha=\{\beta \in 2^\omega : (\alpha,\beta) \in E'\}$ is $\mcG_0$-independent, hence it must be meager. The Kuratowski--Ulam theorem then allows us to conclude that $E'$ is meager; and it follows, by applying Mycielski's theorem, that there exists a continuous map $\psi \colon 2^\omega \to 2^\omega$ such that for all $\alpha \ne \beta \in 2^\omega$ one has $[\psi(\alpha)]_{E'} \ne [\psi(\beta)]_{E'}$.

The map $f=\varphi \circ \psi$ is then a continuous map from $2^\omega$ to $X$ such that for each $\alpha \ne \beta \in 2^\omega$ one has $[f(\alpha)]_E \ne [f(\beta)]_E$, and we are done.
\end{proof}

Whether the conclusion of the above theorem holds true when $E$ is induced by a continuous action of a Polish group (in which case $E$ is an analytic equivalence relation with Borel classes) is a famous open problem known as the \emph{topological Vaught conjecture}\index{topological Vaught conjecture} (the original Vaught conjecture, also still open, is closely related to the case where the acting group is $\Sinf$).

\bigskip

\emph{Comments.}  Silver's original proof of his dichotomy theorem, as well as Kechris, Solecki and Todor\v{c}evic's proof of the $\mcG_0$-dichotomy theorem given in \cite{KST} used elaborate set-theoretic techniques. Groundbreaking work of B. Miller provided a new approach to these questions and helped usher in a plethora of new results in the area of Borel combinatorics, which has seen some dramatic progress in recent times. See \cite{Miller2012} and the survey \cite{KechrisMarks}. For more information on the Vaught conjecture and equivalence relations arising from Borel actions of Polish groups see \cite{Becker1996}.

\chapter{Miller's proof of the $\mcG_0$-dichotomy}

Les us now try to give some intuition for a proof of the $\mcG_0$-dichotomy due to B. Miller, and further formalized by A. Bernshteyn (our write-up here borrows heavily from Bernshteyn's).

We fix $S \subseteq 2^{< \omega}$ such that $|S \cap 2^n|=1$ for all $n< \omega$ (in particular, $S$ contains the empty sequence) and denote by $s_n$ the unique element of $S \cap 2^n$. 

Consider the graph $\mcG_S$ introduced in Definition \ref{defin:G_s}. It is fruitful to consider a sequence of finite approximations of $\mcG_S$. Namely, for each $n$, consider the graph $\mcG_s(n)$ on $2^{n}$ such that $u,v$ are adjacent in $\mcG_s(n)$ iff there exist some $s \in S$ and $\varepsilon \in \{0,1\}$ such that $u= s \smallfrown \varepsilon \smallfrown t$ and $v= s \frown (1-\varepsilon) \frown t$ for some $t$.

Forgetting about the labels for a moment, what graphs do we obtain? $\mcG_S(0)$ is simply a point; no matter what $S$ is, $\mcG_s(1)$ is a graph with two vertices linked by an edge, and $\mcG_S(2)$ is a path of length $4$. The situation becomes more complicated starting from $n=3$, where the particular choice of $S$ starts having an influence on the structure of $\mcG_S$ (for instance, $\mcG_S(3)$ may or may not have vertices of degree $\ge 3$). It becomes apparent that the graphs $\mcG_s(n)$ are obtained by iterating an operation on pointed graphs, starting from the trivial graph $\bullet$ on a singleton.

\begin{defin}
Let $H \subseteq X \times X$ be a graph and $u \in X$ be a vertex of $H$. We denote by $H \sg{u} H$ the graph defined as follows:
\begin{itemize}
\item Its set of vertices is $V(H \sg{u} H)= X_0 \sqcup X_1$, with $X_i= \{v \smallfrown i \colon v \in X \}$  for $i=0,1$.
\item $(v \smallfrown i ,w \smallfrown j) \in H \sg{u} H$ iff one of the following happens:
\begin{itemize}
\item $i=j$ and $(v,w) \in H$.
\item $i \ne j$ and $v=w=u$.
\end{itemize}
\end{itemize}
\end{defin}

Slightly less formally: form two vertex-disjoint copies of $H$, and add an edge between the two copies of $u$ to form $H \sg{u} H$. The graphs $\mcG_s(n)$ are obtained by iterating this operation, starting from $\mcG_S(0)=\bullet$ (formally, the domain of $\mcG_S(0)$ is $\{\emptyset\}$), then defining $\mcG_S(n+1)= \mcG_S(n) \sg{s_n} \mcG_S(n)$.

Given $n$, we can consider the graph $\widetilde{\mcG}_S(n)$ on $2^\omega$ whose edges are $(u \smallfrown x, v \smallfrown x)$ for $(u,v) \in \mcG_S(n)$ and $x \in 2^\omega$. Then $\mcG_S = \bigcup_n \widetilde{\mcG}_S(n)$, and each $\widetilde{\mcG_S}(n)$ has a natural epimorphism onto $\mcG_S(n)$ (mapping $u \smallfrown x$ to $u$) so we think of the finite graphs $\mcG_S(n)$ as providing a sequence of finite approximations of $\mcG_S$.

\begin{exo}\label{exo:G_s_is_connected_and_acyclic}
Prove that $\mcG_S(n)$ is a connected, acyclic graph for all $n$. Use this to show that $\mcG_S$ is acyclic.
\end{exo}

Now, assume that we want to build a homomorphism $\varphi \colon \mcG_S \to G$, where $G$ is some graph on a Polish space $X$. One natural strategy is to try to build $\varphi$ via a sequence of homomorphisms $\varphi_n \colon \mcG_S(n) \to G$ which converges to the desired homomorphism. Further, to build a homomorphism from $\mcG_S(n+1) = \mcG_S(n) \sg{s_n} \mcG_S(n)$ one can look for two homomorphisms $\varphi_0, \varphi_1 \colon \mcG_S(n) \to G$ such that $(\varphi_0(s_n),\varphi_1(s_n)) \in G$; indeed one can then set $\varphi(t \smallfrown i) = \varphi_i(t)$. This suggests an inductive construction; but it may happen that at some step we cannot find any suitable homomorphisms $\varphi_0$, $\varphi_1$ to move on to the next step. This is what motivates the introduction of what we call thwarting sets below.

After this somewhat handwavy discussion, we turn to the actual proof. We fix an analytic graph $G$ on a Polish space $X$, and a continuous function $\varphi \colon \mcN \to X \times X$ such that $\varphi(\mcN)=G$. We also choose a compatible complete metric $d$ on $X$ and a compatible complete metric $\rho$ on $\mcN$.

\begin{defin}
Let $H$ be a graph on a finite set $F$. A \emph{lifted homomorphism} from $H$ to $G$ is a pair of maps $\pi=(\pi_X,\pi_E)$ where $\pi_X \colon F \to X$, $\pi_E \colon H \to \mcN$ are such that for each $(h,h') \in H$ one has 
\[\varphi(\pi_E(h,h'))= (\pi_X(h),\pi_X(h')) .\]
We denote by $\mathrm{Hom}(H,G)$ the set of all lifted homomorphisms from $H$ to $G$, which we view as a closed subset of the Polish space $X^F \times \mcN^H$.
\end{defin}

By working with lifting of edges, which live in the complete metric space $(\mcN,\rho)$ rather than in the analytic space $G$, we can use completeness to build our desired homomorphism via a sequence of approximations.

\begin{defin}
Let $H$ be a graph on a finite set $F$, and $\mcH$ a subset of $\mathrm{Hom}(H,G)$. For $v \in F$ we denote by $\mcH \sg{v} \mcH$ the subset of $\mathrm{Hom}(H \sg{v}H,G)$ made up of all $(\tau_X,\tau_E)$ such that there exist $\pi^0, \pi^1 \in \mathrm{Hom}(H,G)$ with 
\[\forall v \in F \ \forall \varepsilon \in \{0,1\} \  \tau_X(v \smallfrown \varepsilon) = \pi^\varepsilon_X(v) ;\]
\[\forall v_1,v_2 \in F \ \forall \varepsilon \in \{0,1\} \ \tau_E(v_1 \smallfrown \varepsilon, v_2 \smallfrown \varepsilon) = \pi^\varepsilon_E(v_1,v_2) .\]
\end{defin}
These correspond to (lifted) homomorphisms from $H \sg{v} H$ to $G$ which are obtained by gluing together two (lifted) homomorphisms from $H$ to $G$. Note that we can glue together $\pi^0$, $\pi^1$ to form a (lifted) homomorphism from $H \sg{v} H$ to $G$ precisely if $(\pi^0_X(v),\pi^1_X(v)) \in G$.

\begin{defin}
Let $H$ be a graph on a finite set $F$ and $\mcH$ be a subset of $\mathrm{Hom}(H,G)$. For $v \in F$ we let $\mcH(v) = \{ \pi_X(v) : \pi \in \mcH \}$. We say that :
\begin{itemize}
\item $\mcH$ is \emph{thwarting} if there exists $v \in F$ such that $\mcH(v)$ is $G$-independent.
\item $\mcH$ is \emph{facilitating} if it cannot be covered by countably many thwarting Borel subsets.
\end{itemize}
\end{defin}

The terminology ``thwarting'' is meant to convey the idea that such a subset $\mcH$ prevents us from building an element of $\mathrm{Hom}(H\sg{v}H)$ by gluing two elements of $\mcH$ together, thus thwarting the inductive construction outlined above. Of course, a Borel facilitating subset is not thwarting; and the stronger condition we consider, namely not belonging to the $\sigma$-ideal generated by thwarting Borel subsets, allows us to prove the following lemma, which is the linchpin of the proof.

\begin{lem}
Let $H$ be a graph on a finite set $F$ and $v \in F$. Assume that a Borel subset $\mcH \subseteq \mathrm{Hom}(H,G)$ is facilitating. Then $\mcH \sg{v} \mcH$ is also facilitating.
\end{lem}

\begin{proof}
Fix a a countable sequence of thwarting Borel subsets $\mcA_n$ of $\mathrm{Hom}(H \sg{v} \mcH,G)$.

Since each $\mcA_n$ is thwarting, we can choose for each $n$ some $v_n \in F$ and $\varepsilon_n \in \{0,1\}$ such that $\mcA_n(v_n \smallfrown \varepsilon_n)$ is $G$-independent. Since $\mcA_n$ is Borel, $\mcA_n(v_n \smallfrown \varepsilon_n)$ is analytic, and we can find a Borel subset $B_n$ of $X$ which contains $\mcA_n(v_n \smallfrown \varepsilon_n)$ and is $G$-independent (see Lemma \ref{lem:Borel_G_independent}).

Denote $\mcH_n =\{\pi \in \mcH : \pi_X(v_n) \in B_n \}$. Since $B_n$ is $G$-independent, $\mcH_n$ is a thwarting Borel subset of $\mathrm{Hom}(H_n,G)$. Since $\mcH$ is facilitating, the Borel subset $\mcH' = \mcH \setminus \bigcup_{n< \omega} \mcH_n$ is not thwarting. In particular, $\mcH'(v)$ is not $G$-independent, so that $\mcH' \sg{v} \mcH'$ is nonempty. By construction, we have $(\mcH' \sg{v} \mcH') \cap \bigcup_{n< \omega} A_n = \emptyset$ : $\mcH \sg{v} \mcH$ is facilitating.
\end{proof}

Observe that if $H$ is a finite graph and $\mcH \subseteq \mathrm{Hom}(H,G)$ is Borel then $\mcH \sg{v} \mcH$ is also a Borel subset of $\mathrm{Hom}(H \sg{v} H,G)$, a fact that we will use without further reference below.

\begin{lem}\label{lem:linchpin_dichotomy}
Let $H$ be a graph on a finite set $F$. Assume that a Borel $\mcH \subseteq \mathrm{Hom}(H,G)$ is facilitating. Fix $r>0$. Then there exists a facilitating Borel subset $\mcH' \subseteq \mcH$ such that:
\begin{itemize}
\item For each $u \in F$ $\mathrm{diam}(\mcH'(u)) \le r$.
\item For each $(v,w) \in H$ $\mathrm{diam}(\{\pi_E(v,w) : \pi \in \mcH' \}) \le r$. 
\end{itemize}
\end{lem}

\begin{proof}
Write $X= \bigcup_{n< \omega} X_n$, $\mcN = \bigcup_{n< \omega} A_n$ where each $X_n$, $A_n$ is Borel and has diameter less than $r$. Fix $u \in F$ and $(v,w) \in H$. Denote 
\[\mcH_{n,m} = \{\pi \in \mcH : \pi_X(u) \in X_n \textrm{ and } \pi_E(v,w) \in A_m\} .\]
These are Borel subsets of $\mcH_{n,m}$ which cover $\mcH$. Since $\mcH$ is facilitating, some $\mcH_{n,m}$ is facilitating by Lemma \ref{lem:linchpin_dichotomy}. 

Thus, if $\mcH$ is facilitating then for any $u \in F$ and $(v,w) \in H$ there exists a Borel facilitating subset $\mcH_{u,v,w}$ of $\mcH$ such that one has $\mathrm{diam}(\mcH_{u,v,w}(u)) \le r$ and $\mathrm{diam}(\{\pi_E(v,w) : \pi \in \mcH_{u,v,w} \}) \le r$. 

To obtain the desired conclusion, enumerate the (finite) set $F \times H$ and repeat the above argument finitely many times to obtain the desired facilitating subset of $\mcH$.
\end{proof}

\begin{proof}[End of the proof of the $\mcG_0$-dichotomy theorem]

Recall that $\mcG_s(0)$ is the trivial graph on the singleton $\{\emptyset\}$; let $\mcH_0=\mathrm{Hom}(\mcG_S(0),G)$. Clearly $\mcH_0(\{\emptyset\})= V(G)$, so either $\mcH_0$ is facilitating or $X$ is covered by countably many $G$-discrete Borel subsets, i.e., $G$ admits a Borel $\omega$-coloring.

Assuming that $G$ does not admit a Borel $\omega$-coloring, we are thus in the situation where $\mcH_0$ is facilitating. Using our work above, we then obtain a sequence of facilitating Borel subsets $\mcH_n \subseteq \mathrm{Hom}(\mcG_S(n),G)$ such that :
\begin{itemize}
\item For each $n$, $\mcH_{n+1} \subseteq \mcH_n \sg{s_n} \mcH_n$.
\item For each $u \in \{0,1\}^n $ $\mathrm{diam}(\mcH_n(u)) \le 2^{-n}$.
\item For each $(u,v) \in \mcG_S(n)$ $\mathrm{diam}(\{ \pi_E(u,v) : \pi \in \mcH_n\}) \le 2^{-n}$.
\end{itemize}

Given $x \in 2^\omega$ we then let $f(x)$ denote the unique element of $\bigcap_{n< \omega} \overline{\mcH_n(x_{|n}})$; this makes sense because for each $n$, each $s \in 2^n$ and each $\varepsilon \in \{0,1\}$ we have $\mcH_{n+1}(s \smallfrown \varepsilon) \subseteq \mcH_n(s)$.

By construction, $f$ is continuous. Assume that $(x,y) \in \mcG_S$; then there exists $N < \omega$ and $s \in S$ of length $N$ such that $x=s \smallfrown \varepsilon \smallfrown z$, $y= s \smallfrown (1-\varepsilon) \smallfrown z$.

Let $\alpha \in \mcN$ be such that $\{\alpha\}= \bigcap_{n >N} \overline{\{\pi_E(x_{|n},y_{|n}) : \pi \in \mcH_n \}} $. Recall that for each $n$ and each $\pi \in \mcH_n$ we have $\varphi(\pi_E((x_{|n},y_{|n})))= (\pi_X(x_{|n}), \pi_X(y_{|n}))$. Hence, by definition of $f$ and by continuity of $\varphi$, we obtain that $\varphi(\alpha)=(f(x),f(y))$ so that $(f(x),f(y)) \in G$ as expected.
\end{proof}

We still have a bit of book-keeping to do for our work to be completely done, namely we need to explain why the following fact holds true.

\begin{lem}\label{lem:Borel_G_independent}
Let $G$ be a an analytic graph on a Polish space $X$. Assume that $A$ is a $G$-independent analytic subset. Then $A$ is contained in a $G$-independent Borel subset.
\end{lem}

\begin{proof}
Consider $B=\{x \in X : \exists a \in A \ (a,x) \in G\}$; this set is analytic and disjoint from $A$, so by the separation theorem we may find a Borel $C$ containing $A$ and such that $(a,c) \not \in G$ for every $a \in A$ and every $c \in C$.

Consider $D=\pi((C \times C) \cap G)$, where $\pi$ denotes the projection on the first coordinate. This is again an analytic set which is disjoint from $A$, so there exists a Borel subset $E$ containing $A$ and disjoint from $D$. Finally $C \cap E$ is Borel, $G$-independent, and contains $A$.
\end{proof}

We conclude these notes with an elegant observation about Borel chromatic numbers of locally finite Borel graphs.

\begin{exo}[Kechris--Solecki--Todor\v{c}evic]\label{exo:Bernshteyn}
Let $G$ be a Borel graph on a Polish space $X$ such that for each $x$ the set $\{y : (x,y) \in G \}$ is finite. Prove that $G$ admits a Borel $\omega$-coloring.

(Note: no need to use the $\mcG_0$-dichotomy here!)
\end{exo}

\emph{Comments.} Miller's proof of the $\mcG_0$-dichotomy is outlined in \cite{Miller2012} and fully developed in notes of A. Bernshteyn available online (see \cite{Bernshteyn}). Our write-up heavily borrows from Bershteyn's.

\epigraphhead{}
\clearpage
\epigraphhead{}
\pagestyle{plain}
\setlength{\epigraphwidth}{0.45\textwidth}
\part{Appendices}

\backmatter
\setcounter{chapter}{17}

\chapter{Basics of Set Theory}
We briefly discuss some basics of set theory (and make explicit the underlying assumptions used throughout these notes), discussing only what little is needed to understand these notes.

Ordinals have appeared here and there in the notes; we are going to recall their definition shortly. Before that, we need to discuss some classical lemmas on well-ordered sets.

\begin{defin}
Let $(X,\le)$ be a well-ordered set. We say that $A \subseteq X$ is an \emph{initial segment}\index{initial segment} if 
\[\forall a, x \in X \quad (a \in A \textrm{ and } x \le a) \Rightarrow x \in A\ .\]
Given $x \in X$ we denote by $X_x=\{x' \in X : x' < x\}$ the corresponding initial segment of $X$ and call it a \emph{strict} initial segment.
\end{defin}

Note that every initial segment $A$ which is not equal to $X$ is a strict initial segment (consider the smallest element of $X \setminus A$ and the associated initial segment).

\begin{lem}
Assume that $(X,\le)$ is a well-ordered set and that $f \colon X \to X$ is strictly increasing. Then $f(x) \ge x$ for all $x \in X$.
\end{lem}

\begin{proof}
Assume for a contradiction that $A=\{x \in X : f(x) < x \}$ is nonempty. Then $A$ has a smallest element $a$; by definition we have $f(a)< a$, whence $f(f(a)) < f(a)$, so $f(a) \in A$, which is impossible since $f(a)< a$ and $a$ is the smallest element of $A$.
\end{proof}

\begin{lem}
Let $X$ be a well-ordered set, $A,B$ two initial segments and $f \colon A\to B$ an isomorphism. Then $A=B$ and $f= \mathrm{id}$.  
\end{lem}

\begin{proof}
First observe that either $A \subseteq B$ or $B \subseteq A$. Indeed, if there exists $b \in B \setminus A$ then  $A$ cannot contain any element larger than $b$ since $A$ is an initial segment, so $A \subseteq \{x : x < b\} \subset B$.

Assume without loss of generality that $B \subseteq A$. Then $f \colon A \to B \subseteq A$ is strictly increasing, so by the previous lemma $a \le f(a)$ for all $a \in A$. Since $f(a) \in B$ and $B$ is an initial segment, we conclude that $a \in B$ for all $a \in A$, so that $A=B$. Now $f^{-1} \colon A \to A$ is also strictly increasing, so that $f^{-1}(a) \ge a$ for all $a \in A$, in other words we have both $a \ge f(a)$ and $f(a) \ge a$ for all $a \in A$: $f= \mathrm{id}$.
\end{proof}

\begin{thm}
Let $(X,\le_X)$ and $(Y,\le_Y)$ be two well-ordered sets. Then exactly one of the three following possibilities holds:
\begin{itemize}
\item $X$ is isomorphic to a strict initial segment of $Y$.
\item $Y$ is isomorphic to a strict initial segment of $X$.
\item $X$ and $Y$ are isomorphic.
\end{itemize}
\end{thm}

\begin{proof}
We already know that these three possibilities are mutually exclusive, since no well-ordered set is isomorphic to a strict initial segment of itself. Now consider the set 
\[ f= \{(x,y) \in X \times Y : X_x \textrm{ is isomorphic to } Y_y \} .\]
It follows from our previous lemmas that $f$ is (the graph of) a strictly increasing function from an initial segment $A_f$ of $X$ to an initial segment $B_f$ of $Y$. 

If either $A_f=X$ or $B_f=Y$ then we are done; assume for a contradiction that such is not the case. Then one can consider $x = \min(X \setminus A_f)$, $y= \min(Y \setminus B_f)$. But then $A_f= X_x$, $B_f= Y_y$ and $f \colon X_x \to Y_y$ is an isomorphism. Hence $(x,y) \in f$, a contradiction.
\end{proof}

Well-ordering a set $X$ gives a precise mathematical meaning to the intuitive notion of ``enumerating $X$''; if a set is infinite it can be enumerated in many different ways, but each of these ways has a uniquely defined order type. Ordinal numbers, which we define now, provide a complete set of representatives for the order types of well-orderings.

\begin{defin}
Let $\alpha$ be a set. We say that $\alpha$ is an \emph{ordinal}\index{ordinal} if $\alpha$ is strictly well-ordered by $\in$. In particular, every element of an element of $\alpha$ is an element of $\alpha$; and every nonempty subset of $\alpha$ admits a least element (for $\in$)
\end{defin}

\begin{thm}
Let $(X,\le)$ be a well-ordered set. Then there exists a unique ordinal $\alpha$ such that $(X,\le)$ is isomorphic to $(\alpha,\in)$.
\end{thm}

We do not give the proof of this fact here (we do know that such an $\alpha$, if it exists, is unique, but proving existence requires a little more work). 

\begin{defin}
We denote by $\omega$\index{$\omega$} the smallest infinite ordinal, and by $\omega_1$ the smallest uncountable ordinal (the union of the set of all countable ordinals).
\end{defin}
Throughout the text we have identified $\omega$ with the set of integers; in particular, $0=\emptyset$ (the smallest possible ordinal), $1= \{\emptyset\}$, $2=\{\emptyset,\{\emptyset\}\}$, and so on. 

For every two ordinal numbers $\alpha$, $\beta$, one is equal to an initial segment of the other; we then write $\alpha \le \beta$ or $\beta \le \alpha$. Note that then $\alpha \in \beta$ is equivalent to $\alpha < \beta$.
This explains why we often used the notation $n< \omega$ in these notes : finite ordinal numbers (i.e., natural integers) can equivalently be thought of as elements of $\omega$ or ordinals strictly smaller than $\omega$.

The particular coding we have used to define ordinals (via $\in$) does not really matter for us; but it is useful to have built this class of representatives of all possible order types of well orderings. 

Now we can do arithmetic with ordinals.

\begin{defin}
Let $(X,<_X)$ and $(Y,<_Y)$ be two well-ordered sets. Their sum $(X+Y,<)$ is defined as follows:
\begin{itemize}
\item The underlying set is the disjoint union $X \sqcup Y$.
\item The restriction of $<$ to $X$ coincides with $<_X$, the restriction of $<$ to $Y$ coincides with $<_Y$, and $x < y$ for every $x \in X$ and every $y \in Y$.
\end{itemize}
If $\alpha$, $\beta$ are two ordinals, then $\alpha+ \beta$ denotes the unique ordinal isomorphic to the sum of $(\alpha,\in)$ and $(\beta,\in)$.
\end{defin}

Intuitively, considering $\alpha + \beta$ coincides with adding a copy of $\beta$ after a copy of $\alpha$. Note that, for every $\alpha$, $\alpha +1$ is strictly larger than $\alpha$ since $\alpha$ is a strict initial segment of $\alpha +1$; but for instance $1+ \omega=\omega$ (in fact $1+\alpha=\alpha$ for every infinite ordinal $\alpha$). 

One could equivalently define this operation by \emph{transfinite induction}, as follows: fix an ordinal $\alpha$. Then agree that:
\begin{itemize}
\item If $\beta=\gamma \cup \{\gamma\}$ is a \emph{successor} ordinal (equivalently, $\beta$ admits a maximal element) then $\alpha + \beta= (\alpha+\gamma) \cup \{\alpha + \gamma\}$ (the smallest ordinal strictly larger than $\alpha + \beta$).
\item If $\beta= \bigcup_{\gamma < \beta} \gamma$ then $\alpha + \beta = \bigcup_{\gamma < \beta}  \alpha + \gamma$.
\end{itemize}
With these notations, the successor of an ordinal $\alpha$ is the ordinal $\alpha +1$.

One can similarly define ordinal multiplication, exponentiation, etc. Since we have not used these notions in the text we do not go further in this direction and conclude this short appendix by a discussion of cardinal numbers and the axiom of choice.

\begin{defin}
Let $X$ be a set. A \emph{choice function}\index{choice function} on $X$ is a function $\Phi \colon \mathcal{P}(X) \setminus\{\emptyset\} \to X$ such that for every nonempty $A \subseteq X$ one has $\Phi(A) \in A$.
\end{defin}

The \emph{axiom of choice}\index{axiom of choice} is the statement that every set admits a choice function; it has many other equivalent formulations, notably Zermelo's theorem\index{Zermelo's theorem} stating that every set can be well-ordered, or Zorn's lemma\index{Zorn's lemma} stating that every inductive, partially ordered set admits a maximal element (recall that a partially ordered set $(X,\le)$ is inductive if every totally ordered subset has an upper bound in $X$). 

In particular, Zorn's lemma is sufficient to prove that every filter on a set $X$ is contained in an ultrafilter, a fact that we have used to characterize compact Hausdorff spaces in terms of convergence of ultrafilters (though this ``ultrafilter axiom'' is weaker than the full axiom of choice). 

The axiom of choice is also equivalent to the statement that for any two sets $A$, $B$, either there is an injection from $A$ into $B$ or there is an injection from $B$ into $A$; one way to see that this property follows from the axiom of choice is to first notice that this is true for well-orderable sets, and then apply Zermelo's theorem (the converse requires a little more work). Then, given a set $X$, one can define its \emph{cardinal}: the least ordinal $\alpha$ such that $X$ injects into $\alpha$. The cardinality of the set of reals, often denoted $\mathfrak{c}$, is strictly bigger than the smallest infinite cardinal $\aleph_0$ which is simply the set $\omega$; one uses different notations for cardinals in order to avoid confusions (one should pay attention to the fact that cardinal arithmetic and ordinal arithmetic use the same symbols for different operations). The \emph{continuum hypothesis} is the statement that any cardinal $\kappa$ such that $\aleph_0 \le \kappa \le \mathfrak{c}$ is equal to either $\aleph_0$ or $\mathfrak{c}$. Using the $\aleph$ notation, this statement amounts to the equality $\aleph_1=2^{\aleph_0}$ (which may or may not be true depending on the axioms one is working with; see the bibliographical references on set theory mentioned at the end of this chapter). It is notable that analytic subsets of Polish spaces satisfy the conclusion of the continuum hypothesis (i.e., an uncountable analytic subset of a Polish space is of cardinality $\mathfrak{c}$), and possible cardinalities of coanalytic sets ($\aleph_0$, $\aleph_1$, $\mathfrak{c}$) are also fairly limited.

The axiom of choice is very useful when one wishes to establish the existence of universal objects; it has the defect of being highly nonconstructive, and it contradicts some other possible axioms one might want the universe of sets to satisfy. In these notes, it (or rather, the ultrafilter axiom, which is also highly nonconstructive) was integral in our proofs of structural theorems in topological dynamics; throughout the first part of the book we have freely used the axiom of choice whenever it was convenient. One could however recover, for instance, the Kechris--Pestov--Todor\v{c}evic correspondence (Theorem \ref{thm:KPT}) using only a weaker axiom, known as ``dependent choice''. Indeed, applying Proposition \ref{p:many_locally_constant} and an argument similar to the one we used to prove that finite oscillation stability implies extreme amenability, one can prove directly (i.e., without mentioning the universal minimal flow, whose existence is problematic without choice), using an argument similar to the one used to prove that finite oscillation stability implies extreme amenability, that a subgroup $G$ of $\Sinf$ is extremely amenable iff for every open subgroup $V$ of $G$ and every $\varphi \in \ell^\infty(V\rcoset G)$ taking only finitely many values there exists a fixed point in $\overline{G \varphi}$. This is the Ramsey property for embeddings in disguise.

Every theorem proved in the second part of these notes only requires dependent choice, which we state now. 

\begin{defin}
The \emph{axiom of dependent choice}\index{axiom of dependent choice} is the following statement: 

Assume that $X$ is a set, and that $R$ is a binary relation on $X$ with the property that for every $x \in X$ there exists $x' \in X$ such that $x \mathrel{R} x'$. Then there exists a sequence $(x_n)_{n < \omega}$ of elements of $X$ such that for all $n < \omega$ one has $x_n \mathrel{R} x_{n+1}$.
\end{defin}

This axiom plays a part in the proof of many theorems in analysis (for instance, the Baire category theorem) and infinite combinatorics, and is really fundamental for the topics discussed in these notes. It is however too weak to prove the existence of sets that are often considered ``pathological'', such as subsets of $\R$ which are not Lebesgue-measurable or not Baire measurable, or uncountable subsets of $\R$ which do not contain a perfect subset.

To see that the axiom of dependent choice is a consequence of the axiom of choice, fix a set $X$ and a binary relation $R$ on $X$ such that for each $x \in X$ there exists $y \in X$ satisfying $x \mathrel{R} y$ (in particular, $X$ is nonempty). Let $\Phi$ be a choice function on $X$; for $x \in X$ define $f(x)=\Phi(\{y : x \mathrel{R} y\})$ then pick $x_0 \in X$ and define $x_n = f^n(x_0)$ for all $n< \omega$. We then have, as expected, $x_n \mathrel{R} x_{n+1}$ for all $n< \omega$.

As a warning that the axiom of choice is sometimes used without it being obvious, we finish this appendix with the following easy, but (hopefully) thought-provoking, fact.

\begin{thm}
The axiom of choice is equivalent to the following statement: for any sets $X$, $Y$ and every surjective function $f \colon X \to Y$ there exists a function $g \colon Y \to X$ such that $f(g(y))=y$ for all $y \in Y$.
\end{thm}

\begin{proof}
Assume first that the axiom of choice holds, and fix a choice function $\Phi$ on $X$. Then $g \colon y \mapsto \Phi(f^{-1}(\{y\}))$ satisfies the desired condition.

Conversely, assume that the above statement holds and consider a set $X$, which we may assume to be nonempty. Denote by $\Omega$ the set of all nonempty subsets of $X$, and consider the disjoint union $Z=\bigcup_{A \in \Omega} A \times\{A\}$. Then for every $z \in Z$ there is a unique $A \in \Omega$ such that $z=(a,A)$ for some $a \in A$ and we may set $f(z)=A$. Then $f \colon Z \to \Omega$ is surjective, hence our assumption implies that there exists $g \colon \Omega \to Z$ such that $f(g(A))=A$ for all $A \in \Omega$. This $g$ (or rather, its first coordinate) induces a choice function on $X$.
\end{proof}

\bigskip
\emph{Comments.} As an introduction to set theory, Halmos' book \cite{Halmos} is quite concise and accessible; in a similar spirit Kaplansky's book \cite{Kaplansky} is a good introduction to set theory for analysts. One can also mention Ciesielski's book \cite{Ciesielski}, which provides a quick introduction to classical set-theoretic concepts before going into deeper topics such as forcing, while remaining accessible to non specialists. The reader who wants to go more deeply into set theory can start their journey by opening Jech's encyclopedic \cite{Jech}.

\chapter{Solutions of the exercises}

{\bf \ref{exo1}}. 1.
To prove that $(\sigma,\tau) \mapsto \sigma \tau$ is continuous, pick $(\sigma_0,\tau_0) \in \Sinf \times \Sinf$ and a neighborhood $U$ of $\sigma_0 \tau_0$. There exists a finite $A \subset \omega$ such that any element of $\Sinf$ which coincides with $\sigma_0 \tau_0$ on $A$ belongs to $U$. 

By definition,
\[V=\lset (\sigma,\tau) : \sigma_{|\tau_0(A)}= {\sigma_0}_{|\tau_0(A)} \textrm{ and } \tau_{|A}={\tau_0}_{|A} \rset \]
is a neighborhood of $(\sigma_0,\tau_0)$ and $\sigma \tau \in U$ whenever $(\sigma,\tau) \in V$.

To prove continuity of $\sigma \mapsto \sigma^{-1}$, fix $\sigma_0 \in \Sinf$ and a neighborhood $U$ of $\sigma_0^{-1}$. There exists a finite $A \subset \omega$ such that any $\tau \in \Sinf$ which coincides with $\sigma_0^{-1}$ on $A$ belongs to $U$. Note that $B=\sigma_0^{-1}(A)$ is finite; and any $\tau$ which coincides with $\sigma_0$ on $B$ is such that $\tau^{-1}$ coincides with $\sigma_0$ on $A$, hence such that $\tau^{-1}$ belongs to $U$.

2. By definition, for any $i$ we have $d(\sigma_i,\sigma_{i+1})=2^{-i-1}$, so the sequence $(\sigma_i)_{i< \omega}$ is Cauchy. In the product space $\omega^\omega$ this sequence converges to the shift map $n \mapsto n+1$, which does not belong to $\Sinf$. Hence $\Sinf$ is not closed in $(\omega^\omega,d)$ and is thus not complete.

3. Assume that $(\sigma_i)_{i < \omega}$ is Cauchy for $\rho$; note that $(\omega^\omega,d)$ is complete, so we obtain that $\sigma_i \xrightarrow[i \to + \infty]{} \sigma$ and $\sigma_{i}^{-1} \xrightarrow[i \to + \infty]{} \tau$, where $\sigma$, $\tau$ are two elements of $\omega^\omega$. Fix $n \in \omega$; there exists $I$ such that $\sigma_i(n)= \sigma(n)$ for all $i \ge I$, and $J$ such that $\sigma_i^{-1}(\sigma(n))= \tau(\sigma(n))$ for all $i \ge J$. Thus, for $i \ge \max(I,J)$ we have 
\[\tau(\sigma(n))= \sigma_{i}^{-1}(\sigma(n))= \sigma_{i}^{-1}(\sigma_i(n))=n .\]
One checks similarly that $\sigma \circ \tau= \mathrm{id}$. Consequently $\sigma \in \Sinf$, and $\sigma_i \xrightarrow[i \to + \infty]{} \sigma$; thus $(\Sinf,\rho)$ is complete.

4. Saying that $\sigma \colon \omega \to \omega$ is surjective is the statement 
\[\forall i < \omega \, \exists j < \omega \ \sigma(j)=i .\]
Thus the set of all surjective maps from $\omega$ to $\omega$ is equal to $\bigcap_{i< \omega} \bigcup_{j < \omega} \lset \sigma : \sigma(j)=i \rset$.
Each $\lset \sigma : \sigma(j)=i \rset$ is open, whence being surjective is a $G_\delta$ condition.

Similarly, the set of injective maps from $\omega$ to itself is equal to 
$ \bigcap_{i \ne j \in \omega} \lset \sigma : \sigma(i) \ne \sigma(j) \rset $
and so it is also $G_\delta$. 

Thus $\Sinf$ is the intersection of two $G_\delta$ subsets of $\omega^\omega$; hence $\Sinf$ is $G_\delta$ in $\omega^\omega$. So $\Sinf$ is Polish since it is a $G_\delta$ subset of a Polish space.

\bigskip

{\bf \ref{exo2}}.
For every $n < \omega$ the map $\pi_n \colon \sigma \mapsto \sigma(n)$ is a continuous map from $\Sinf$ to the discrete space $\omega$. Thus if $G$ is compact in $\Sinf$ then $\pi_n(G)$ is a compact, hence finite, subset of $\omega$. This shows that if $G$ is compact then each $G$-orbit is finite.

Conversely, assume that the $G$-orbit of each $n < \omega$ is finite. Let $(A_i)_{i < \omega}$ enumerate all the $G$-orbits for the action $G \actson \omega$; consider the natural map $\varphi \colon G \to \prod_{i< \omega} \mathfrak{S}(A_i)$. If we endow each $\mathfrak{S}(A_i)$ with the discrete topology, and $X=\prod_{i< \omega} \mathfrak{S}(A_i)$ with the product topology, then $\varphi$ is a homeomorphism onto its image. Furthermore, it is straightforward to check that $\varphi(G)$ is closed in $X$ since $G$ is closed in $\Sinf$. Since $X$ is compact, so is $\varphi(G)$, whence $G$ is compact.

\bigskip

{\bf \ref{exo3}}.
The fact that a product of groups (endowed with the obvious operations) is a group should not be a problem. Thus the only point is to check the topological conditions; in the case of a finite product this is straightforward.

In the case of a countable product of Polish groups $(G_i,d_i)$, where each $d_i$ is complete, begin by replacing $d_i$ with $\min(d_i,1)$, which is also complete and induces the same topology. Then consider the metric $d$ on $G$ defined by 
\[d((g(i))_{i< \omega},(h(i))_{i< \omega})= \sum_{i=0}^{+ \infty} d_i(g(i),h(i)) .\]
This metric has the following properties: a sequence $(g_n)_{n < \omega}$ of elements of $G$ converges iff each $(g_n(i))_{n< \omega}$ converges in $G_i$; and $(g_n)_{n < \omega}$ is Cauchy iff each $(g_n(i))_{n< \omega}$ is Cauchy in $(G_i,d_i)$. It is then clear that $(G,d)$ is complete.

Separability of $G$ is a standard fact about countable products of separable spaces (here one can simply use that sequences which eventually coincide with $(1_{G_i})_{i < \omega}$ are dense in $G$, and use separability of $G_i$).

The characterization of convergent sequences given above immediately implies that the group operations on $G$ are continuous as soon as they are continuous on each $G_i$, and we are done.

\bigskip
{\bf \ref{exo4}}.
Assume first that $G$ is a subgroup of $\Sinf$. Given $A$ a finite subset of $\omega$, let $G_A$ denote $\lset g \in G : \forall a \in A \ g(a)=a \rset$. By definition of the topology of $\Sinf$, the family of all $G_A$, as $A$ runs over finite subsets of $\omega$, is a basis of open neighborhoods of $1_G$, and clearly each $G_A$ is a subgroup.

The converse is more interesting. Let $(G_n)_{n< \omega}$ be a countable family of open subgroups which form a basis of neighborhoods of $1_G$ (we only consider the case where this family is infinite; otherwise $G$ is discrete and there is nothing to do). For each $n$ denote $X_n = G/G_n$; note that since $G_n$ is open in $G$ and $G$ is second-countable $G_n$ has at most countable index in $G$ (see exercise \ref{exo5} if necessary). Let $X = \bigsqcup_{n< \omega} X_n$; this is a countable set, and $G$ naturally embeds (as an abstract group) in the Polish group $H=\mathfrak{S}_X$ of all permutations of $X$ , via the map $\varphi$ defined by 
\[\varphi(g)(hG_n)= gh G_n . \]
Assume that $(g_i)_{i < \omega}$ converges to $g$ in $G$, and let $hG_n$ be an element of $X$. Since $G_n$ is an open subgroup, for all $n$ large enough we have $h^{-1}g_i ^{-1} g h \in G_n$ whence $g h G_n=g_i hG_n$. This proves that $(\varphi(g_i))_{i< \omega}$ converges to $\varphi(g)$; in other words, $\varphi$ is continuous.

Conversely, assume that $(\varphi(g_i))_{i < \omega}$ converges to some $\varphi(g)$ in $\mathfrak{S}_X$. Given $n< \omega$, we then have that for all $i$ large enough $g_i G_n = g G_n$, i.e., $g_i \in gG_n$. Since $(gG_n)_{n < \omega}$ is a basis of neighborhoods of $g$, we conclude that $(g_i)_{i < \omega}$ converges to $g$ in $G$.

Thus, $\varphi$ is a topological group isomorphism from $G$ to a subgroup of $\mathfrak{S}_X$; since $G$ is Polish, $\varphi(G)$ is a Polish subgroup of the Polish group $\mathfrak{S}_X$, hence a closed subgroup.

\bigskip
{\bf \ref{exo5}}.
Let $(V_n)_{n < \omega}$ be a countable basis for the topology of $X$, and consider $U= \bigcup_{i \in I} O_i$ an arbitrary union of open sets.

For each $i \in I$, one may choose some subset $A_i$ of $\omega$ such that $U_i= \bigcup_{a \in A_i} V_a$. Let $A= \bigcup_{i \in I} A_i$, which is a subset of $\omega$, hence (at most) countable. For each $a \in A$ choose some $j_a \in I$ such that $V_a \subseteq U_{j_a}$, and let $J = \lset j_a \colon a \in A\rset$; this is a countable subset of $I$. 

Pick $x \in \bigcup_{i \in I} U_i$, and let $i \in I$ be such that $x \in U_i$. Then there is some $a \in A_i$ such that $x \in V_a$; and $V_a \subseteq \bigcup_{j \in J} U_j$. Hence $x \in \bigcup_{j \in J} U_j$. This proves that $\bigcup_{i \in I} U_i \subseteq \bigcup_{j \in J} U_j$, and the converse inclusion is immediate.

\bigskip
{\bf \ref{exo6}}.
By Theorem \ref{thm:Banach}, we know that $\varphi$ is continuous. Now prove the following assertions:
\begin{itemize}
\item For each $x \in \R$ and each $n \in \Z$, $\varphi(nx)=n \varphi(x)$.
\item For each $x \in \R$ and each $q \in \Q$, $\varphi(qx)=q \varphi(x)$.
\item For each $x \in \R$ and each $t \in \R$, $\varphi(tx)=t \varphi(x)$.
\end{itemize}
(Only the last step requires continuity of $\varphi$). Denoting $\alpha=\varphi(1)$, we obtain that for all $t \in \R$ one has $\varphi(t)=t \alpha$.

\bigskip
{\bf \ref{exo7}}. 1. Under these assumptions on $H$, we know that $1$ belongs to the interior of 
$H H^{-1}=H$. Thus $H$ has nonempty interior hence (since $H$ is a subgroup) $H$ is open in $G$. But then $G \setminus H = \bigcup_{g \not \in H} gH$ is also open, thus $H$ is clopen.

2. Since $\varphi$ is Borel, it is Baire-measurable, hence actually continuous since it is a homomorphism of Polish groups. Similarly, $\varphi^{-1}$ is Borel (by the Lusin--Suslin theorem), hence continuous, and we are done.

3. Since $\tau_1 \subseteq \tau_2$, the identity map $\mathrm{id} \colon (G,\tau_2) \to (G,\tau_1)$ is a continuous isomorphism of Polish groups. Hence its inverse is also continuous: $\tau_2 \subseteq \tau_1$.

\bigskip
{\bf \ref{exo7bis}}.
1. Denote $V_A = \lset \sigma : \forall b \not \in A \ \sigma(b)=b \rset$. 

If $A$ has at least three elements, it is not hard to check that $\sigma \in U_A$ iff for all $\sigma' \in V_A$ one has $\sigma' \circ \sigma = \sigma \circ \sigma'$. In other words, $U_A = \bigcap_{\sigma' \in V_A} \lset \sigma : \sigma' \circ \sigma = \sigma \circ \sigma' \rset$ whence (by continuity of the group operations) $U_A$ is $\tau'$-closed. 

If $A$ has $1$ or $2$ elements then by the result of the previous paragraph $U_A$ contains, as a subgroup of countable index, a $\tau'$-closed subset, hence $U_A$ is $\tau'$-Borel (actually, $F_\sigma$).

2. It follows that each $\tau$-open subset is $\tau'$-Borel, so $\mathrm{id} \colon (\Sinf,\tau') \to (\Sinf,\tau)$ is Borel. Hence, by the result of Exercise \ref{exo7}, $\mathrm{id}$ is a topological group isomorphism, i.e., $\tau=\tau'$.

\bigskip
{\bf \ref{exo8}}.
Let $U$ be a nonempty open subset of $G$. Since $G$ is Polish, there is a sequence $(g_n)_{n< \omega}$ such that $G= \bigcup_{n< \omega} g_n U$. Hence $H= \bigcup_{n < \omega} \varphi(g_n) \varphi(U)$, so that $\varphi(U)$ is not meager in $H$. Furthermore, $\varphi(U)$ is analytic, hence Baire measurable .

Assume now that $U$ contains $1_G$, and let $V$ be open nonempty and such that $V V^{-1} \subseteq U$. By a corollary of Pettis' Lemma, since $\varphi(V)$ is not meager we have $1_H \in \mathrm{Int}(\varphi(V) \varphi(V)^{-1})=\mathrm{Int}(\varphi(V) \varphi(V^{-1})) \subseteq \mathrm{Int}(\varphi(U))$. Thus, as soon as $U$ is open and contains $1_G$, $1_H$ belongs to the interior of $\varphi(U)$. Since $\varphi$ is a homomorphism this is enough to prove that $\varphi$ is open.

\bigskip
{\bf \ref{exo9}}.
We clearly have $d(fH,gH)=d(gH,fH)$ and $d(fH,fH)=0$. Assume that $d(fH,gH)=0$; then there exists a sequence $(h_n)_{n < \omega}$ of elements of $H$ such that $\rho(fh_n,g)$ converges to $0$. Using right-invariance of $\rho$, this means that $f h_ng^{-1}$ converges to $1_G$, so (using continuity of left translation by $f^{-1}$ and right translation by $g$) $h_n \xrightarrow[n \to + \infty]{} f^{-1}g$. Since $H$ is closed, we then have that $f^{-1}g \in H$, in other words $fH=gH$.

The triangle inequality is straightforward to check: fix $\varepsilon >0$ and let $h_1,h_2 \in H$ be such that $\rho(f,gh_1) \le d(fH,gH)+ \varepsilon$ and $\rho(gh_1,kh_2) \le d(gH,kH)+ \varepsilon$. Then we have
\[d(fH,kH) \le \rho(f,kh_2) \le \rho(f,gh_1)+ \rho(gh_1,kh_2) \le d(fH,gH)+ d(gH,kH)+ 2 \varepsilon .\]

By definition, $\pi \colon g \mapsto gH$ is $1$-Lipschitz from $(G,\rho)$ to $(G/H,d)$, hence continuous; it remains to prove that $\pi$ is open. Pick $U$ open in $G$ and $g_0 \in U$; then find $\varepsilon >0$ such that $B(g_0,\varepsilon) \subseteq U$. Now note that for every $g \in G$ we have
\[d(g_0H,gH)< \varepsilon \Rightarrow \exists h \ \rho(g_0,gh) < \varepsilon \Rightarrow \exists h \ gh \in U \Rightarrow g \in UH \ .\]
Thus the ball of radius $\varepsilon$ centered at $g_0H$ is contained in $\pi(U)$. This shows that $\pi(U)$ is open, and we are done.

\bigskip
{\bf \ref{exo10}}.
Assume first that there exists $x \in X$ with a dense orbit, and let $U$, $V$ be nonempty open. By assumption there exists $f \in G$ such that $fx \in U$, and $h \in G$ such that $h x \in V$. Then $hf^{-1}U \cap V \ne \emptyset$.

Conversely, assume that for each nonempty open $U$, $V$ there exists $g \in G$ such that $gU \cap V \ne \emptyset$. Let $(V_n)_{n < \omega}$ be a basis of nonempty open sets for the topology of $X$; given $n< \omega$ denote 
\[U_n = \lset x : Gx \cap V_n \ne \emptyset \rset .\]
Since $V_n$ is open, each $U_n$ is also open; furthermore, our assumption amounts to the statement that any nonempty open subset of $X$ intersects every $U_n$, i.e., each $U_n$ is open dense.

Thus $A= \bigcap_{n< \omega} U_n$ is comeager in $X$; and $A$ precisely consists of the elements of $X$ which have a dense orbit.

\bigskip
{\bf \ref{exo11}}.
Assume that $A$ is Baire measurable; note that for any homeomorphism $g$ of $X$ we have $U(g A)=gU(A)$. Now, assume that $A$ is not meager; then $U(A)$ is nonempty, open, and $G$-invariant since $gA=A$ for all $g \in G$. By topological transitivity of $G \actson X$, this implies that $U(A)$ is dense, thus comeager. So $A$ is also comeager.

\bigskip
{\bf \ref{exo12}}. 1. Assume that $x$ belongs to $U_r$, and let $\varepsilon < r$ be such that $f(x) \in \mathrm{Int}\overline{f(B(x,\varepsilon))}$. Denote $V= \mathrm{Int}\overline{f(B(x,\varepsilon))}$; by continuity of $f$, $f^{-1}(V)$ is open and we know that it contains $x$. 

Let $\delta >0$ be such that $\varepsilon + \delta < r$; then, for any $y \in B(x,\delta) \cap f^{-1}(V)$ we have that $B(y,\varepsilon + \delta) \supseteq B(x,\varepsilon)$, whence
\[\mathrm{Int}\overline{f(B(y,\varepsilon + \delta))} \supseteq \mathrm{Int}\overline{f(B(x,\varepsilon))} = V \ni f(y) .\]
Thus $f^{-1}(V) \cap B(x,\delta)$ is an open neighborhood of $x$ which is contained in $U_r$.

2. This is obvious, since balls of vanishing radius form a neighborhood basis (and one only needs to consider rational $r$, so that the set of points of local density for $f$ is, as promised, a countable intersection of open sets).

\bigskip
{\bf \ref{exo13}}.
Note that $A$ is meager iff $Y \setminus A$ is comeager, and $f^{-1}(A)$ is meager iff $X\setminus f^{-1}(A)$ is comeager; so equivalence of the two first items is immediate.

If \ref{exo13_ii} holds and $A$ is dense open, then $f^{-1}(A)$ is open (because $f$ is continuous) and comeager, so $f^{-1}(A)$ is dense so that \ref{exo13_iii} also holds. The converse implication is also clear: if \ref{exo13_iii} holds and $A \subseteq Y$ is comeager, then $A$ contains a countable intersection $\bigcap_n O_n$ of dense open sets of $Y$, so $f^{-1}(A)$ contains $\bigcap_n f^{-1}(O_n)$, which is an intersection of dense open sets of $X$.

Assume that \ref{exo13_iv} is false and let $U$ be nonempty open such that $f(U)$ is meager. Then $f^{-1}(f(U))$ is not meager, since it contains $U$. So \ref{exo13_i} is false too.

Any set which is not meager is somewhere dense, so the implication \ref{exo13_iv} $\Rightarrow$ \ref{exo13_v} is immediate.

Assume that \ref{exo13_iii} does not hold; let $A$ be dense open such that $f^{-1}(A)$ is not dense. Choose $U$ open nonempty which does not intersect $f^{-1}(A)$. Then $f(U)$ is nowhere dense. In other words, \ref{exo13_v} implies \ref{exo13_iii}, and we are done.

\bigskip
{\bf \ref{exo14}}.
Since $f$ is continuous, we can find a dense $G_\delta$ subset $A$ of $Y$ such that $f \colon f^{-1}(A) \to A$ is open; and since $f$ is category-preserving, $f^{-1}(A)$ is comeager in $X$. Let $X'=f^{-1}(A)$ (a Polish space since it is $G_\delta$ in $X$).

For any $\Omega \subseteq X$, $\Omega$ is comeager in $X$ iff $\Omega \cap X'$ is comeager in $X'$ (dense open subsets of $X'$ are precisely intersections with $X'$ of dense open subsets of $X$); similarly a set is comeager in $Y$ iff its intersection with $A$ is comeager in $A$. Finally, $f \colon X' \to A$ is open, hence category-preserving.

Assume that $\Omega$ is comeager in $X$. Then $\Omega \cap X'$ is comeager in $X'$, whence (using the fact that the result of Theorem \ref{t:splitting_category} holds for continuous open maps) we obtain that
\[\lset a \in A : \Omega \cap f^{-1}(\{a\}) \textrm{ is comeager in } f^{-1}(\{a\}) \rset\] 
is comeager in $A$, hence in $Y$.

Conversely, if $\lset y \in Y : \Omega \cap f^{-1}(\{y\}) \textrm{ is comeager in } f^{-1}(\{y\}) \rset$
is comeager in $Y$, then 
\[\lset a \in A : \Omega \cap f^{-1}(\{a\}) \textrm{ is comeager in } f^{-1}(\{a\}) \rset\]
is comeager in $A$ so (using the fact that $\Omega \cap X'$ is Baire measurable in $X'$) we obtain that $\Omega \cap X'$ is comeager in $X'$, equivalently $\Omega$ is comeager in $X$. 

\bigskip

{\bf \ref{exo15}}.
The language for pure sets consists of only the equality symbol $=$. A finitely generated substructure of a set $X$ is simply a finite set; thus an isomorphism between finite structures $N_1$, $N_2$ is simply a bijection $f \colon N_1 \to N_2$. One can extend $f$ to a bijection of $N_1 \cup N_2$ onto itself, then to a bijection of $X$ by mapping each $x \not \in N_1 \cup N_2$ to itself.

Ultrahomogeneity of $(\Q,<)$ is classically proved using a back-and-forth argument (this goes back to Cantor). Start from an order preserving bijection $f \colon A \to B$ between two finite subsets $A,B$ of $\Q$. Enumerate $\Q=\{q_n : n \in \omega\}$.
Using the fact that the order of $\Q$ is dense and without endpoints, we may build increasing families of finite sets $A_i$, $B_i$ and order preserving bijections $f_i \colon A_i \to B_i$ with the following properties:
\begin{itemize}
\item $A_0=A$, $B_0=B$, $f_0=f$.
\item For each $n$ $f_{n+1}$ extends $f_n$.
\item For every $k$, $q_k \in A_{2k+2}$.
\item For every $k$, $q_k \in B_{2k+1}$.
\end{itemize}
Then the sequence $(f_n)_{n< \omega}$ induces an order-preserving bijection $f_\infty \colon \bigcup_n A_n \to \bigcup_n B_n$, which extends $f$. Since we ensured that $\bigcup_n A_n=\bigcup_n B_n=\Q$, we are done.

For the case of vector fields over a countable field $\bK$ (the assumption of countability is just there to ensure that we work with countable structures), form a language $\mcL$ which contains a constant $0$ (to be interpreted as the origin of the vector field), a binary function $+$ (interpreted as addition of vectors) as well as a function symbol $\cdot_\lambda$ for each $\lambda \in \bK$ (interpreted as multiplication by $\lambda$). Then finitely generated substructures are precisely finite-dimensional vector subspaces, and a linear isomorphism between subspaces of a vector space extends to an isomorphism of the whole space (modulo some use of the axiom of choice, which we allow in these notes).

Finally, the language of Boolean algebras contains (besides $=$) two constant symbols $0$, $1$ as well as binary functions $\vee$ (union) and $\wedge$ (intersection) and a unary function ${}^c$ (complementation). Any finitely generated boolean algebra is finite, so a partial isomorphism between finitely generated substructures corresponds to an isomorphism between finite subalgebras. Let $\bB$ be a countable, infinite, atomless Boolean algebra and let $\varphi \colon \bM \to \bN$ be an isomorphism between finite subalgebras of $\bB$. Choose $a \in \bB$ and let $\bM'$ be the boolean subalgebra generated by $\bM$ and $a$; using a back-and-forth argument, it is enough to show that $\varphi$ extends to a partial isomorphism with domain $\bM'$. Since $\bB$ has no atoms, there exists $c \in \bB$ such that for all $m \in \bM$ one has $c \wedge \varphi(m) \ne 0 \Leftrightarrow a \wedge m \ne 0$. Setting $\varphi(a)=c$ provides the desired extension (we leave the verification to the reader).

\bigskip
{\bf \ref{exo16}}.
Fix a closed subgroup $G$ of $\Sinf$, and consider the structure $\bM$ built in \ref{p:auts_countable_structure}. We know that $G=\Aut(\bM)$ and need to show that $\bM$ is ultrahomogeneous.
Let $A$, $B$ be finite subsets of $M$ and $f \colon \bA \to \bB$ an isomorphism. In particular, writing $A=\{a_1,\ldots,a_n \}$ and $B=\{f(a_1),\ldots,f(a_n)\}$ we know that $(a_1,\ldots,a_n)$ and $(f(a_1),\ldots,f(a_n))$ belong to the same $G$-orbit (for the diagonal action $G \actson \omega^n$). This amounts to saying that there exists $g \in G=\Aut(\bM)$ such that $g(a_i)=f(a_i)$ for all $i$, so that $f$ indeed extends to an automorphism of $\bM$.

\bigskip
{\bf \ref{ex:left_completion}}.
Assume that the universe of $\bM$ is equal to $\omega$, and view $\Aut(\bM)$ as a subset of $\omega^\omega$. 
For $f,g \in \omega^\omega$ we let $d(f,g)= \inf \lset 2^{-n} \colon \forall i < n \ f(i)=g(i) \rset$; then $d$ induces the product topology on $\omega^\omega$, $(\omega^\omega,d)$ is complete and the restriction of $d$ to $\Aut(\bM)$ is a compatible left-invariant metric.

Denote by $E$ the set of all embeddings from $\bM$ to itself; it is straightforward to check that $E$ is a closed subset of $\omega^\omega$. Hence $(E,d)$ is complete and it remains to prove that $\Aut(\bM)$ is dense in $E$ when $\bM$ is ultrahomogeneous. Choose $f \in E$, $A$ a finite subset of $\omega$ and a basic neighborhood $U=\lset g \in E \colon \forall a \in A \ g(a)=f(a) \rset$. Then $U$ consists precisely of elements of $E$ which coincide with $f$ on the structure $\langle A \rangle$ generated by $A$. Since $\langle A \rangle$ is a finitely generated substructure and $\bM$ is ultrahomogeneous, there exists $g \in \Aut(\bM)$ which coincides with $f$ on $\langle A \rangle$. Hence $\Aut(\bM) \cap U \ne \emptyset$, proving that $\Aut(\bM)$ is dense in $E$.

\bigskip
{\bf \ref{exo18}}.
Assume that ${\mathbf \Gamma}$, ${\mathbf \Gamma_1}$, ${\mathbf \Gamma_2}$ are three finite graphs and that $\alpha \colon {\mathbf \Gamma } \to {\mathbf \Gamma_1}$, $\beta \colon {\mathbf \Gamma} \to {\mathbf \Gamma_2}$ are embeddings. Without loss of generality we may assume that $\Gamma_1 \cap \Gamma_2=\Gamma$ and $\alpha$, $\beta$ are induced by the inclusions of $\Gamma$ into $\Gamma_1$, $\Gamma_2$ respectively. Then form a new graph ${\mathbf \Lambda}$, whose set of vertices is equal to $\Gamma_1 \cup \Gamma_2$ and such that there is an edge connecting $\lambda$ and $\lambda'$ iff they both belong to $\Gamma_1$ and there is an edge connecting them in ${\mathbf \Gamma_1}$, or they both belong to $\Gamma_2$ and there is an edge connecting them in ${\mathbf \Gamma_2}$. This provides the desired amalgam of ${\mathbf \Gamma_1}$, ${\mathbf \Gamma_2}$ over ${\mathbf \Gamma}$.

For linear orders, we can similarly reduce to the case where $A$ a subset of both $B$ and $C$, $B \cap C = A$, and the order on $A$ coincides with the order induced by $B$ and with the order induced by $C$. Then we define a linear order on $D=B \cup C$ by declaring that all elements of $B \setminus A$ are strictly smaller than all elements of $C \setminus A$. This witnesses that the class of finite linear orders satisfies the amalgamation property.

Finally, consider the class of all finite groups; let $G$, $H$, $K$ be finite groups and $\alpha \colon G \to H$, $\beta \colon G \to K$ be two group embeddings. Let $L$ denote the corresponding amalgamated free product; there exists embeddings $i \colon H \to L$ and $j \colon K \to L$ such that $i \circ \alpha= j \circ \beta$. The group $L$ is not finite, but it is residually finite, so there exists a finite group $\tilde L$ and a group homomorphism $\varphi \colon L \to \tilde L$ which is injective on $i(H) \cup j(K)$. Then $\varphi \circ i$ is a group embedding from $H$ into $\tilde L$, $\varphi \circ j$ is a group embedding from $K$ into $\tilde L$, and $(\varphi \circ i) \circ \alpha= (\varphi \circ j) \circ \beta$. 

\bigskip
{\bf \ref{exo19}}.
Let $\bM$ be an ultrahomogeneous structure; fix a finitely generated substructure $\bA$ of $\bM$, $\bB \in \age(\bM)$ and an embedding $\alpha \colon \bA \to \bB$. Since $\bB \in \age(\bM)$ there exists an embedding $\beta \colon \bB \to \bM$. Then $\beta \circ \alpha \colon \bA \to \bM$ is an embedding, so it extends by ultrahomogeneity to $g \in \Aut(\bM)$. It remains to notice that $g^{-1} \circ \beta$ is an embedding of $\bB$ into $\bM$ and
\[\forall a \in A \quad (g^{-1} \circ \beta) (\alpha(a))= g^{-1}(g(a)) =a .\]

\bigskip
{\bf \ref{exo20}}.
The hereditary property is clear (a subspace of a metric space, with the induced metric, is itself a metric space). 

To prove that the joint embedding property holds, take two finite rational metric spaces $X$, $Y$, let $L$ be a rational number larger that $\mathrm{diam}(X)$ and $\mathrm{diam}(Y)$. Consider $Z= X \sqcup Y$; define a distance function on $Z$ by setting $d(x,x')=d_X(x,x')$ for any $x,x' \in X$; $d(y,y')=d_Y(y,y')$ for any $y,y' \in Y$; and $d(x,y)=L$ for any $x \in X$ and $y \in Y$.

It remains to prove amalgamation over a nonempty subset; this reduces to the case where $X$ is a nonempty rational metric space which is a common subspace of two finite rational metric spaces $Y$, $Z$ such that $Y \cap Z=X$. We want to define a metric on $Y \cup Z$ which takes rational values and extends the metrics on $Y$, $Z$. One may for instance use the shortest-path metric: for any $y \in Y \setminus X$, $z \in Z \setminus X$, set
\[d(y,z) = \min \lset d(y,x)+d(x,z) : x \in X \rset.\] 
(We leave the verification to the reader).

\bigskip

{\bf \ref{exo21}}.
1. Let $\bA, \bB$ be finite subgraphs of $\bR$. Then $\bA \cup \bB$ is a finite graph, and one can obtain a finite graph $\bC$ with vertices $A \cup B \cup \{*\}$ by declaring that $*$ is connected to all elements of $A$ and to no element of $B$.

Since $\bR$ has the extension property, there is an embedding $g$ of $\bC$ into $\bR$ which fixes $A \cup B$ pointwise; $g(*)$ is then an element of $R$ which is connected to each element of $A$ and to no element of $B$.

Conversely, assume that $G$ satisfies such a condition. It is then not hard to prove by induction on $n$ that any graph on $n$ vertices embeds in $G$, so $G$ is universal for finite graphs; similarly, the fact that $G$ has the extension property is straightforwardly prove by induction on $|B| - |A|$ (the condition of the exercise amounts to the extension property for one-point extensions).

2. Letting $E=\{(i,j) \in \omega \times \omega : i< j\}$, we identify the set $\mcG$ of graphs with vertex set $\omega$ with $2^E$. In this exercise, $\mcG$ is endowed with the Bernoulli measure $\mu$ such that 
\[\mu(\lset x \colon x(e_1)=\varepsilon_1,\ldots,x(e_k)=\varepsilon_k\rset)= \frac{1}{2^k} \] for all $k$ and all pairwise distinct $e_1,\ldots,e_k$. We want to prove that the set of all elements of $\mcG$ which are isomorphic to $\bR$ has measure $1$ for $\mu$.

To that end, fix two disjoint finite subsets $A$, $B$ of $\omega$. Let $N=|A \sqcup B|$. For each $i > \max(A,B)$, the set 
\[\mcG_i= \lset G \in \mcG : \textrm{ there is an edge from } i \textrm{ to every } a \in A \textrm{ and to no } b \in B \rset\]
is such that $\mu(\mcG_i)= \frac{1}{2^N}$. Furthermore, the $\mcG_i$ form an independent family; thus $\mu (\bigcup_i \mcG_i)=1$. Since there are only countably many pairs of disjoint finite subsets of $\omega$, we conclude that, almost surely, an element $G$ of $\mcG$ is such that for any two disjoint finite $A,B \subset \omega$ there is an element of $G$ which is adjacent to every element of $A$ and to no element of $B$. It follows from the result of the first question that $G$ is almost surely isomorphic to $\bR$.

\bigskip
{\bf \ref{exo22}}.
This is very similar to the case of graphs. Let ${\mathbf G}$ be a finite subtournament of two finite tournaments ${\mathbf G_1}$, ${\mathbf G_2}$ such that $G=G_1 \cap G_2$. Then one can simply add an oriented edge from every element of $G_1$ to every element of $G_2$ to obtain a finite tournament which witnesses the amalgamation property.

\bigskip 
{\bf \ref{exo:Hall}}.
1. By construction, $\bH$ is the increasing union of a family of finite groups, so it is locally finite. Let $h_n=|G_n|$. We have $h_0=3$ and $h_{n+1}=h_n !$ for all $n$, so $h_n \xrightarrow[n \to + \infty]{} + \infty$. 

For any finite group $G$, and any $n \ge |G|$, $G$ embeds into the symmetric group on $n$ elements, thus $G$ embeds in $G_i$ for every big enough $i$; in particular $\bH$ is universal for the class of finite groups.

To prove ultrahomogeneity, fix some $i$ and let $G$, $G'$ be two subgroups of $G_i$ and $f \colon G \to G'$ be an isomorphism.  Take coset representatives $g_1,\ldots,g_n$ in $G_i \lcoset G$ and $g_1',\ldots,g_n'$ in $G_i  \lcoset G'$ (note that since $|G|= |G'|$ we have $|G_i \lcoset G|= |G_i \lcoset G'|$). Then define a bijection $\sigma$ of $G_i$ by setting $\sigma(g g_i) = f(g) g_i'$ for each $g \in G$ and $i \in \{1,\ldots,n\}$.

We view $\sigma$ as an element of $G_{i+1}= \mathfrak{S}(G_i)$. Let $\lambda \colon G_i \to G_{i+1}$ be the embedding given by the left translation action of $G_i$ on itself. Note that for every $g \in G$ we have
\[\lambda(g)(g_i)=gg_i = \sigma^{-1}(f(g)g_i')= \sigma^{-1}\lambda(f(g))\sigma(g_i).\]
From this it is not hard to check that $\lambda(g)= \sigma^{-1}\lambda(f(g)) \sigma$, equivalently $\lambda(f(g))= \sigma \lambda(g) \sigma^{-1}$. We conclude that conjugation by $\sigma$ is an inner automorphism of $G_{i+1}$ which extends $f$. This is more than what we need to conclude that $\bH$ is ultrahomogeneous (since any inner automorphism of $G_{i+1}$ extends to an inner automorphism of $\bH$).

2. We proved the existence of a countable, ultrahomogeneous group whose age is the class of all finite groups. Thus this class satisfies the amalgamation property (we had seen a proof of that fact earlier using properties of amalgamated free products).

3. Let $g,h \in \bH$ have the same order. Then there exists an isomorphism $\varphi \colon \langle g \rangle \to \langle h \rangle$ mapping $g$ to $h$, and when dealing with the first question we have proved that $\varphi$ extends to an inner automorphism of $\bH$. Thus $g$ and $h$ are conjugate.

4. Let $N$ be a nontrivial normal subgroup. If $N$ contains an element of some given order $p$, then (by the result of the previous question) $N$ contains every element of order $p$. Transpositions generate every $G_i$ (and are conjugate in $G_i$) so if $N$ contains an element or order $2$ then $N=\bH$. If $N$ contains an element $g$ of even order, then by Cauchy's theorem applied to $\langle g \rangle$, $N$ contains an element of order $2$ and we are done by the previous argument. 

To conclude, note that for any odd integer $p \ge 3$ the product of the $p$-cycles $(1,\ldots,p)$ and $(2, \ldots,p+1)$ in $\mathfrak{S}_{p+1}$ has even order. Thus if $N$ contains an element of nontrivial odd order it must also contain an element of even order, and we are done.

\bigskip
{\bf \ref{exo24}}. 
We prove by induction that these two conditions imply that an ordering which is dense and without endpoints satisfies the extension property. Let us only deal with the first step. Assume that $A$ is a finite subset of $M$, $|B|=|A|+1$, $B$ is a finite linearly ordered set and $\alpha \colon A \to B$ is an embedding. Let $b$ be such that $\{b\}= B \setminus \alpha(A)$. 

If $b$ is smaller than all elements of $\alpha(A)$, then (using that $\bM$ is without endpoints) pick some $m \in M$ such that $m < a$ for all $a \in A$, and define $\varphi \colon B \to M$ by setting $\varphi(b)=m$ and $\varphi(\alpha(a))=a$ for all $a \in A$. The case where $b$ is larger than all elements of $\alpha(A)$ is dealt with similarly.

The remaining case is when both $A_-=\{a \in A : \alpha(a) < b \}$ and $A_+=\{a \in A : \alpha(A) >b \}$ are nonempty; note that $a < a'$ for any $a \in A_-$ and $a' \in A_+$. Since $A_-$ is a finite ordered set it has a maximum $x$, and similarly $A_+$ has a minimal element $y$. Since $x<y$ and the ordering of $\bM$ is dense we can pick $m \in M$ such that $x < m < y$ and then set $\varphi(\alpha(a))=a$ for all $a \in A$ and $\varphi(b)=m$ to obtain the desired embedding of $B$ in $(M,<)$.

This proves the desired result. Clearly $(\Q,<)$ is dense and without endpoints, so it is the Fraïssé limit of finite linear orders; and since the extension property characterizes ultrahomogeneous structures we obtain that any countable linear order which is dense and without endpoints is isomorphic to $\Q$.

\bigskip

{\bf \ref{exo25}}.
The fact that $d(x,x)=0$ for all $x$ yields the first property; symmetry of $d$ gives the second. Stability under finite intersection comes from the fact that $\min(\varepsilon_1,\varepsilon_2)>0$ whenever $\varepsilon_1$, $\varepsilon_2 >0$. The fourth condition is satisfied by definition; and the last one follows from the triangle inequality.

Now, assume that $d_1$, $d_2$ are uniformly equivalent; let $U \in \mcU_{d_1}$. Then there exists $\varepsilon >0$ such that $\lset (x,y) \in X^2 : d_1(x,y) < \varepsilon \rset \subseteq U$. Find $\delta >0$ such that for all $x,y \in X$ one has $d_2(x,y) < \delta \Rightarrow d_1(x,y)< \varepsilon$. Then $\lset (x,y) \in X^2 : d_2(x,y) < \delta\rset \subseteq U$, proving that $U \in \mcU_{d_2}$. The other inclusion is similar.

Conversely, assume that $\mcU_{d_1}=\mcU_{d_2}$. Fix $\varepsilon >0$. Then $\lset (x,y) \in X^2 : d_1(x,y) < \varepsilon \rset \in \mcU_{d_1}=\mcU_{d_2}$ hence there exists $\delta >0$ such that $\lset (x,y) \in X^2 : d_2(x,y) < \delta \rset \subseteq \lset (x,y) \in X^2 : d_1(x,y) < \varepsilon \rset $. In other words, $d_2(x,y)< \delta \Rightarrow d_1(x,y) < \varepsilon$.
 
\bigskip
{\bf \ref{exo26}}.
Again this is mostly immediate. Assume that $O$ is open for the topology induced by $d$; then for any $x \in O$ there exists $r >0$ such that for all $y \in X$ one has $d(x,y)< r \Rightarrow y \in O$. Letting $U=\lset (a,b) \colon d(x,y) < r \rset$, we have that $U \in \mcU_d$ and $U[x] \subseteq O$. This proves that $O$ is open for the topology induced by $\mcU_d$ and the converse implication is similar. 

\bigskip
{\bf \ref{exo27}}.
Fix $U \in \mcU$ and $x \in X$. Consider $O= \lset y : \exists V \in \mcU\ V[y] \subseteq U[x] \rset$. We clearly have that $x \in O \subseteq U[x]$; thus to prove that $U[x]$ is a neighborhood of $x$ it is enough to prove that $O$ is open (actually $O$ is the interior of $U[x]$). Fix $y \in O$ and find $V \in \mcU$ such that $V[y] \subseteq U[x]$; then pick a symmetric $W  \in \mcU$ such that $W \circ W \subseteq V$. Then for every $z \in W[y]$ we have that $W[z] \subseteq V[y] \subseteq U[x]$, proving that $W[y]$ is contained in $O$. Hence $O$ is open: $U[x]$ is a neighborhood of $x$.

Conversely, let $A$ be a neighborhood of $x$; then there exists $U \in \mcU$ such that $U[x] \subseteq A$. Let $V= U \cup \{(x,y) \colon y \in A \}$. Then $V \in \mcU$ since it contains $U$; and $V[x]= U[x] \cup A=A$.

\bigskip
{\bf \ref{exo28}}.
Fix $U \in \mcU$, then choose some symmetric $V \in \mcU$ such that $V \circ V \circ V \subseteq U$. We claim that $V \subseteq \Int{U}$, which (once proved) implies that $\Int{U} \in \mcU$.

Take $(x,y) \in V$; then $V[x] \times V[y]$ is a neighborhood of $(x,y)$ and for any $(x',y') \in V[x] \times V[y]$ we have $(x',x) \in V$ and $(x,y) \in V$ and $(y,y') \in V$, so that $(x',y') \in V \circ V \circ V \subseteq U$.

\bigskip
{\bf \ref{exo29}}.
This follows immediately from Theorem \ref{t:uniform_continuity_compact}; indeed, if $\mcU$ and $\mcV$ are two uniformities on $X$ which induce the same compact topology, then $\mathrm{id}$ is uniformly continuous from $(X,\mcU)$ to $(X,\mcV)$ and from $(X,\mcV)$ to $(X,\mcU)$, in other words $\mcU=\mcV$. We know that neighborhoods of the diagonal form a compatible uniform structure (recall that all our spaces are Hausdorff), so we are done.

\bigskip
{\bf \ref{exo30}}.
Let $\mcF$ be a filter on $X$; assume first that for some $A \subseteq X$ we have neither $A \in \mcF$ nor $X \setminus A \in \mcF$. Then for any $F_1,\ldots,F_n \in \mcF$ we must have $A \cap F_1 \ldots \cap F_n \ne \emptyset$, otherwise we would have $\bigcap_{i=1}^n F_i \subseteq X \setminus A$, which would imply $X \setminus A \in \mcF$. Then $\lset B \subseteq X : \exists F \in \mcF \ B \supseteq F \cap A \rset$ is a filter which strictly contains $\mcF$, whence $\mcF$ is not an ultrafilter.

Next, assume that $\mcF$ is a filter which is not an ultrafilter; then there exists a filter $\mcF'$ which contains $\mcF$ strictly, and one can choose some $A \in \mcF' \setminus \mcF$. If $X \setminus A \in \mcF$ then it also belongs to $\mcF'$, which is not possible since $A \cap (X \setminus A)= \emptyset$. Hence neither $A$ nor $X \setminus A$ belong to $\mcF$.

\bigskip
{\bf \ref{exo31}}.
If $x \ne y$ then there exists a neighborhood $U$ of $x$ and a neighborhood $V$ of $y$ such that $U \cap V = \emptyset$. If $\mcF$ converges to $x$, then $U \in \mcF$, whence $V \not \in \mcF$ since $U \cap V= \emptyset$. Hence $\mcF$ cannot converge to $y$, so that any filter converges to at most one point.

\bigskip
{\bf \ref{exo32}}.
Saying that $f$ is continuous at $x$ is equivalent to the statement that for any neighborhood $V$ of $f(x)$ there is a neighborhood $U$ of $x$ such that $f(U) \subseteq V$. In particular, the neighborhood filter of $f(x)$ is contained in the image of the neighborhood filter of $x$, which proves that the image of any filter converging to $x$ is a filter converging to $f(x)$.

Conversely, if the image of the neighborhood filter at $x$ converges to $f(x)$, then for any neighborhood $V$ of $x$ we have that $f^{-1}(V)$ is a neighborhood of $x$, i.e., $f$ is continuous at $x$.

\bigskip
{\bf \ref{exo33}}.
Again this is a direct application of the definitions: we have that
\begin{align*}
(x_n)_{n < \omega} \textrm{ converges to } x & \Leftrightarrow  \forall V \textrm{ neighborhood of } x  \ \exists N < \omega \ \forall n \ge N \ x_n \in V \\
&\Leftrightarrow  \forall V \textrm{ neighborhood of } x, \   \lset n : x_n \not \in V \rset \textrm{ is finite} \\
&\Leftrightarrow  \forall V \textrm{ neighborhood of } x, \ \varphi^{-1}(V) \textrm{ is cofinite} \\
&\Leftrightarrow  \forall V \textrm{ neighborhood of } x, \ V \in \mcF \\
&\Leftrightarrow \mcF \textrm{ converges to } x.
\end{align*}

\bigskip
{\bf \ref{exo34}}.
We already know one implication. To see the converse, let $X$ be a space such that any ultrafilter on $X$ is convergent; then consider a family $\mcF$ of closed sets with the property that the intersection of any finite subfamily of $\mcF$ is nonempty. Next (using the axiom of choice, or more precisely the ultrafilter axiom) find an ultrafilter $\mcU$ containing $\mcF$.

By assumption on $X$, $\mcU$ converges to some $x$. Let $V$ be a neighborhood of $x$ and $F \in \mcF$. Since $V \cap F \in \mcU$ we must have $V \cap F \ne \emptyset$; this proves that $x \in \overline{F}=F$. Hence $x \in F$ for all $F \in \mcF$. This proves that $\cap \mcF \ne \emptyset$, so $X$ is compact.

To prove that a product of compact spaces is compact, first show that a filter $\mcF$ on a product $X=\prod_{i \in I} X_i$ of nonempty topological spaces is convergent iff each image $\pi(\mcF)$ is convergent (where $\pi_i \colon X \to X_i$ is the coordinate projection); this follows straightforwardly from the definition of the product topology (incidentally, the full axiom of choice would be needed here if we allowed non-Hausdorff spaces). Then consider an ultrafilter $\mcU$ on a product $X=\prod_{i \in I} X_i$ of compact spaces: since each $\pi_i(\mcU)$ is an ultrafilter on $X_i$ and $X_i$ is compact we obtain that $\pi_i(\mcU)$ is convergent, thus $\mcU$ is convergent, thus $X$ is compact.

\bigskip
{\bf \ref{exo35}}.
Again one simply needs to unravel the definitions. Denote the Fréchet filter by $\mcG$ and the metric uniformity by $\mcU_d$. Then we have, for any sequence $(x_n)_{n< \omega}$:
\begin{align*}
(x_n)_{n < \omega} \textrm{ is a Cauchy sequence} & \Leftrightarrow \forall \varepsilon >0 \ \exists N < \omega \ \forall n,m \ge N \ d(x_n,x_m) < \varepsilon \\
&\Leftrightarrow  \forall \varepsilon >0 \ \exists G \in \mcG \ \forall (n,m) \in G \times G \ d(\varphi(n),\varphi(m)) < \varepsilon \\
& \Leftrightarrow \forall U \in \mcU_d \ \exists G \in \mcG \ \varphi(G) \times \varphi(G) \subseteq U \\
& \Leftrightarrow \forall U \in \mcU_d \ \exists F \in \mcF \ F \times F \subseteq U .
\end{align*}
The last equivalence holds because $\mcF$ precisely consists of all $F$ which contain some $\varphi(G)$ for some $G \in \mcG$.

\bigskip
{\bf \ref{exo36}}.
1. Pick $x \in X$ and let $\mcV_x$ be the neighborhood filter of $x$. Fix $U \in \mcU$ then find a symmetric $V \in \mcU$ such that $V \circ V \subseteq U$. We know that $V[x] \in \mcV_x$; furthermore for any $(y,z) \in V[x] \times V[x]$ we have $(y,x) \in V$ and $(x,z) \in V$, so that $(y,z) \in V \circ V \subseteq U$. We proved that $V[x] \times V[x] \subseteq U$, which is enough to show that $\mcV_x$ is a Cauchy filter.

2. Assume first that $(X,\mcU)$ is complete and let $(x_n)_{n < \omega}$ be a Cauchy sequence. Let $\mcF$ denote the image of the Fréchet filter on $\omega$ under $\varphi$. Then we know that $(x_n)_{n< \omega}$ is Cauchy iff $\mcF$ is a Cauchy filter, and $(x_n)_{n< \omega}$ is convergent iff $\mcF$ is a convergent filter. Since $(X,\mcU)$ is complete we conclude that $(x_n)_{n< \omega}$ is convergent: $(X,d)$ is complete.

Conversely, assume that $(X,d)$ is complete and let $\mcF$ be a Cauchy filter. For all $n< \omega$ choose $F_n \in \mcF$ such that $F_n \times F_n \subseteq \lset (x,y) : d(x,y) < 2^{-n}\rset$ (so that $\mathrm{diam}(F_n) \le 2^{-n}$); we may also ensure that $F_{n+1} \subseteq F_n$ for all $n$ since $\mcF$ is stable under finite intersections. Then $(F_n)_{n < \omega}$ is a decreasing sequence of nonempty sets with vanishing diameter, so that $\bigcap_{n< \omega} \overline{F_n} =\{x\}$ for some $x \in X$ since $(X,d)$ is complete. Choose some $\varepsilon >0$, then find $n$ such that $2^{-n} < \varepsilon$. Since $x \in \overline{F_n}$ we have $d(x,y) \le 2^{-n}$ for all $y \in F_n$, thus $F_n \subseteq \lset y \colon d(x,y) < \varepsilon \rset$. We conclude that $\mcF$ contains the neighborhood filter of $x$, in other words $\mcF$ converges to $x$ and $(X,\mcU)$ is complete.

\bigskip
{\bf \ref{exo37}}.
Assume that $\mcF$ is a Cauchy filter with an adherent point $x$. Let $O$ be a neighborhood of $x$, which we may assume to be of the form $U[x]$ for some $U \in \mcU$. Find a symmetric $V \in \mcU$ such that $V \circ V \subseteq U$, then $F \in \mcF$ such that $F \times F \subseteq V$.

Choose $y \in F$. By assumption, $F \cap V[x] \ne \emptyset$ so we can find $z \in F$ such that $(x,z) \in V$; we also have $(z,y) \in V$ hence $(x,y) \in U$. Thus $F \subseteq U[x]$, proving that $\mcF$ contains the neighborhood filter of $x$: $\mcF$ converges to $x$.

\bigskip

{\bf \ref{exo38}}.
We already know that any compact uniform space is complete. Assume that $(X,\mcU)$ is a compact uniform space, and choose $U \in \mcU$. Then $X=\bigcup_{x \in X} U[x]$ and $x$ belongs to the interior of $U[x]$ for all $x$ so, by compactness, there exists a finite $A \subseteq X$ such that $X= \bigcup_{a \in A} U[a]=U[A]$. Hence $(X,\mcU)$ is totally bounded.

Assume now that $(X,\mcU)$ is totally bounded. Let us prove that every ultrafilter on a totally bounded space is Cauchy; using the characterization of compactness via convergence of ultrafilters will then give the desired result. Fix an ultrafilter $\mcF$ on $X$ and consider $U \in \mcU$ and a symmetric $V \in \mcU$ such that $V \circ V \subseteq U$. By total boundedness, there exists a finite $A \subseteq X$ such that $X=V[A]$ whence, since $\mcF$ is an ultrafilter, there is some $a$ such that $V[a] \in \mcF$. Since $V[a] \times V[a] \subseteq U$, we obtain as desired that $\mcF$ is a Cauchy filter.

\bigskip
{\bf \ref{exo39}}.
Let $(Y,\mcV)$ be a complete Hausdorff uniform space, and let $f \colon (X,\mcU) \to (Y,\mcV)$ be uniformly continuous.

First, we prove that the image by $f$ of a Cauchy filter $\mcF$ on $X$ is a Cauchy filter on $Y$. To that end, choose $V \in \mcV$; since $f$ is uniformly continuous there exists $U \in \mcU$ such that $(f(x),f(y)) \in V$ for any $(x,y) \in U$. Since $\mcF$ is a Cauchy filter, there exists $F \in \mcF$ such that $F \times F \subseteq U$, whence $(f(x),f(y)) \in V$ for any $(x,y) \in F \times F$. Equivalently, $f(F) \times f(F) \subseteq V$; since $f(F) \in f(\mcF)$ this proves that $f(\mcF)$ is a Cauchy filter.

Since we assumed that $(Y,\mcV)$ is complete, we can then map any minimal Cauchy filter $\mcF$ on $X$ to the limit of $f(\mcF)$; denote this map $\hat{f}$. By definition we have $\hat{f} \circ i=f$ and it remains to prove that $\hat{f}$ is uniformly continuous.

Pick $V \in \mcV$, and choose a symmetric $W \in \mcV$ such that $W \circ W \circ W \subseteq V$. We know that $U=f^{-1}(W) \times f^{-1}(W) \in \mcU$. Consider now the set 
\[C(U)=\lset (\mcF,\mcG) \in \yhwidehat{X} \times \yhwidehat{X} : \exists A \in \mcF \cap \mcG \ A \times A \subseteq U \rset \] 
and choose $(\mcF,\mcG) \in C(U)$. We have by definition $(\hat{f}(\mcF),\hat{f}(\mcG))= (\lim f(\mcF),\lim f(\mcG))$. Denote $y=\lim f(\mcF)$, $y'= \lim f(\mcG)$. To establish uniform continuity of $\hat{f}$, it is enough to prove that $(y,y') \in V$.

By definition of $C(U)$, we can fix some $A \in \mcF \cap \mcG$ such that $A \times A$ is contained in $U$, so that $f(A) \times f(A) \subseteq W$. Since $f(A) \in f(\mcF)$, which converges to $y$, we know that $f(A) \cap W[y] \ne \emptyset$; similarly, $f(A) \cap W[y'] \ne \emptyset$. Choose $z \in f(A) \cap W[y]$ and $z' \in f(A) \cap W[y']$. We then have
\[ (z,z') \in W \textrm{ and } (z,y) \in W \textrm{ and } (z',y') \in W \]
so that $(y,y') \in W \circ W \circ W \subseteq V$. This proves that $\hat{f}$ is a uniformly continuous extension of $f$ to $\yhwidehat{X}$ (and such an extension is obviously unique since $X$ is dense in $\yhwidehat{X}$, so there is at most one continuous extension to $\yhwidehat{X}$ of any continuous function on $X$).

This universal property characterizes $(\yhwidehat{X},\yhwidehat{\mcU})$ in the following sense: assume that $(X',\mcU')$ is a complete uniform space with a uniformly continuous embedding $j \colon X \to X'$ satisfying the same uniform property. Then $j^{-1} \circ i$ extends to a uniform isomorphism from $(X',\mcU')$ to $(\yhwidehat{X},\yhwidehat{\mcU})$. 

\bigskip
{\bf \ref{exo40}}.
Denote by $\mcU_Y$ the uniformity on $Y$ induced by $\mcU$ and let $i_Y \colon Y \to X$ and $i \colon X \to \yhwidehat{X}$ be the embedding maps.

Let $f \colon (Y,\mcU_Y) \to (Z,\mcV)$ be a uniformly continuous map from $(Y,\mcU_Y)$ to some complete uniform space $(Z,\mcV)$. Then, by density of $Y$ in $X$ and an argument similar to the one given in the previous exercise, there exists a unique uniformly continuous map $f_X \colon (X,\mcU) \to (Z,\mcV)$ such that $f_X \circ i_Y=f$. By the universal property of $(\yhwidehat{X},\yhwidehat{\mcU})$, there exists a unique uniformly continuous map $\hat{f} \colon (\yhwidehat{X},\yhwidehat{\mcU}) \to (Z,\nu)$ such that $\hat{f} \circ i= f_X$, hence such that $\hat{f} \circ i \circ i_Y=f$. This proves that the uniform embedding $i \circ i_Y \colon Y \to \yhwidehat{X}$ satisfies the universal property characterizing the completion of $Y$, whence $(\yhwidehat{X},\yhwidehat{\mcU})$ is (uniformly isomorphic to) the completion of $(Y,\mcU_Y)$.

\bigskip
{\bf \ref{exo41}}.
Let us prove that $(X,\mcU)$ is totally bounded iff $(\yhwidehat{X},\yhwidehat{\mcU})$ is totally bounded. We identify $X$ with $i(X) \subseteq \yhwidehat{X}$.

Assume first that $(X,\mcU)$ is totally bounded, and fix $U \in \yhwidehat{\mcU}$. Choose a symmetric $V \in \yhwidehat{\mcU}$ such that $V \circ V \subseteq U$; note that $V \cap (X \times X) \in \mcU$, so by total boundedness there exists a finite $A \subseteq X$ such that $X \subseteq \bigcup_{a \in A} V[a]$. By density of $X$, it follows that $\yhwidehat{X} =  \bigcup_{a \in A} \overline{V[a]}$. Since by definition of the topology induced by a uniform structure we have $\overline{V[a]} \subseteq (V \circ V)[a] \subseteq U[a]$, we are done.

Conversely, assume that $(\yhwidehat{X},\yhwidehat{\mcU})$ is totally bounded and pick $U \in \mcU$. Choose some $\yhwidehat{U}$ such that $\yhwidehat{U} \cap (X \times X)=U$. Then find a 
symmetric $\yhwidehat{V} \in \yhwidehat{\mcU}$ such that $\yhwidehat{V} \circ \yhwidehat{V} \subseteq \yhwidehat{U}$; by total boundedness of $\yhwidehat{X}$ there exists a finite $A$ subset $\yhwidehat{X}$ such that $X \subseteq \bigcup_{a \in A} \yhwidehat{V}[a]$. Each $\yhwidehat{V}[a]$ is open nonempty in $\yhwidehat{X}$ and $X$ is dense, so we can pick $x_a \in X \cap \yhwidehat{V}[a]$. Now pick $x \in X$; there exists $a \in A$ such that $(a,x) \in \yhwidehat{V}$, and by choice of $x_a$ we also have $(a,x_a) \in \yhwidehat{V}$. Finally, $(x_a,x) \in (\yhwidehat{V} \circ \yhwidehat{V}) \cap (X \times X) \subseteq U$ and we conclude that $X = \bigcup_{a \in A} U[x_a]$.

\bigskip

{\bf \ref{exo:right_translations_left_UC}}.
Fix $g \in G$, and let $U$ be an entourage of $(G,\mcU_l)$; we may assume that there exists an open neighborhood $O$ of $1_G$ such that $U=\lset(f,h) : f^{-1}h \in O\rset$. By continuity of the group operations on $G$, there exists an open neighborhood $O'$ of $1_G$ such that $g^{-1}O'g \subseteq O$.

Now consider $V= \lset (f,h): f^{-1}h \in O' \rset$, which is a basic entourage for $\mcU_l$. For any $(f,h) \in V$ we have $(fg)^{-1}(hg)= g^{-1}(f^{-1}h)g \in g^{-1}O'g \subseteq O$. In other words, $(f,h) \in V \Rightarrow (fg,hg) \in U$, proving that $h \mapsto hg$ is a uniform isomorphism of $(G,\mcU_l)$.

The corresponding result for left translations follows from the fact that $h \mapsto h^{-1}$ is a uniform isomorphism from $(G,\mcU_r)$ to $(G,\mcU_l)$. The fact that left and right translations are uniform isomorphisms for $\mcU_+$ is immediately visible from the fundamental system given in the definition.

For $\mcU_{Roelcke}$, one may simply note that $h \mapsto hg$ is a uniform isomorphism of both $(G,\mcU_l)$ and $(G,\mcU_r)$, hence it is a uniform isomorphism of their meet $\mcU_{Roelcke}$ (or directly adapt the proof given above).

\bigskip
{\bf \ref{exo43}}.
Denote by $\mcU_d$ the uniformity induced by $d$, and let $U=\{(g,h) :d(g,h) < r \}$ be a basic entourage. By left-invariance we have $U=\{(g,h) : d(g^{-1}h,1_G)<r \} = \{(g,h) : g^{-1}h \in B(1_G,r)\}$ so $U \in \mcU_l$.

Conversely, pick a basic entourage $U$ for $\mcU_l$, of the form $\{(g,h) : g^{-1}h \in O\}$ for some open neighborhood $O$ of $1_G$. Since $d$ induces the topology of $G$ there is $r>0$ such that $B(1_G,r) \subseteq O$, and we conclude by the same computation as above that $\{(g,h) : d(g,h) < r \} \subseteq U$. 

This proves that $\mcU_d=\mcU_l$. It follows that any two left-invariant compatible metrics on $G$ induce the same uniformity, so they are uniformly equivalent and in particular have the same Cauchy sequences.

\bigskip
{\bf \ref{exo:Struble}}.
1. Fix a left-invariant metric $d$ on $G$. Since $G$ is locally compact there exists a compact neighborhood $K$ of $1_G$, and we can pick $r$ such that the closed ball of center $1_G$ and radius $r$ is contained in $K$, hence is compact. For any $g \in G$ we have $\overline{B}(g,r)= g \overline{B}(1_G,r)$, so $\overline{B}(g,r)$ is compact.

Consider a Cauchy sequence $(g_n)_{n< \omega}$ in $G$. For some $M < \omega$ we have $d(g_n,g_M) \le r$ for all $n \ge M$, so that $(g_n)_{n \ge M}$ is a Cauchy sequence of elements of a compact metric space, hence a convergent sequence. We have proved that $(G,d)$ is complete.

3. First observe that $G$ is $\sigma$-compact since it is second-countable and locally compact. Fix a countable sequence $(K_n)_{n< \omega}$ of compact sets which covers $G$. 
Let $U$ be a relatively compact open neighborhood of $1_G$. Then define $V_0 = U \cup K_0$, $V_{n+1} = V_n^3 \cup K_n$ for all $n \in \omega$. By continuity of the group operations each $V_n$ is relatively compact, and $\bigcup_{n < \omega} V_n = G$.

Since $G$ is first countable we may, using continuity of group operations, build a sequence $(V_n)_{n \in \Z, n <0}$ of identity neighborhoods such that $V_{n}^3 \subseteq V_{n+1}$ for all $n$ and $\bigcap_{n \in \Z} V_n= \{1_G\}$. As in the proof of Theorem \ref{t:metrizable_uniformity} we can then set 
\[\rho(g,h) = \inf \{2^n : g^{-1}h \in V_n \} \]
and 
\[d(g,h) = \inf \left\{\sum_{i=0}^n \rho(g_i,g_{i+1}) : g_0 =g \ , g_{n+1}=h \right \}.\] 
The same computation as in the proof of Theorem \ref{t:metrizable_uniformity} yields $\frac{1}{2} \rho \le d \le \rho$, and it follows that $d$ is a compatible metric on $G$. Left-invariance is built into the definition, and compactness of closed balls comes from the fact that for each $n$ the closed ball $\{g: d(g,1_G) \le 2^n \}$ is contained in $V_{n+1}$.

\bigskip
{\bf \ref{exo45}}.
Assume that $G$ admits a compatible bi-invariant metric $d$. Then for each $g,h \in G$ and each $r >0$ we have 
\[d(g,1_G)< r \Leftrightarrow d(hg,h)< r \Leftrightarrow d(hgh^{-1},1_G)<r \]
whence $1_G$ admits a basis of neighborhoods which are conjugacy invariant.

Assume that the converse holds; pick a compatible left -invariant metric $d$ on $X$ then define 
\[\forall g,h \in G \quad \rho(g,h) = \sup_{k \in G} d(gk,hk). \]
This is a bi-invariant metric on $G$ and $d(g,h) \le \rho(g,h)$ for all $g,h$. Fix $r>0$; by assumption on $G$ we may find $\varepsilon >0$ such that for any $g\in G$ one has 
\[ d(g,1_G)< \varepsilon \Rightarrow \forall k \in G \ d(k^{-1}gk,1_G) < r .\]
(pick a conjugacy-invariant neighborhood of $1_G$ inside $B(1_G,r)$, then $\varepsilon >0$ so that $B(1_G,\varepsilon)$ is contained in that neighborhood) 

For any $g,h,k \in G$ we then have
$d(g,h)< \varepsilon \Rightarrow d(gk,hk)= d(k^{-1}gk,k^{-1}hk) < r$; in particular, $d(g,h) < \varepsilon \Rightarrow \rho(g,h) \le r$. This proves that $d$ and $\rho$ are uniformly equivalent, so that $G$ admits a compatible bi-invariant metric.

\bigskip
{\bf \ref{exo46}}.
Assume that $G$ admits a compatible bi-invariant metric $d$; then for any $g_1,g_2,h_1,h_2$ one has $ d(g_1h_1,g_2h_2) \le d(h_1,h_2) +d(g_1,g_2)$ (by bi-invariance) so, since $d$ induces $\mcU_l$, we conclude that $(g,h)\mapsto gh$ is left-uniformly continuous. The same argument proves that $(g,h) \mapsto gh$ is right-uniformly continuous in that case.

Conversely, assume that $(g,h) \mapsto gh$ is left-uniformly continuous. Let $U$ be an open neighborhood of $1_G$; there exists an open neighborhood $V$ of $1_G$ such that for all $g_1,g_2,h_1,h_2 \in G$ one has
\[(g_1^{-1} g_2 \in V \textrm{ and  } h_1^{-1} h_2 \in V) \Rightarrow (g_1h_1)^{-1}(g_2h_2) \in U .\]
Apply this with $g_1=1_G$, $g_2=g$, $h_1=h_2=h$ to obtain that for any $g \in V$ and any $h \in G$ one has $h^{-1}g h \in U$. It follows that $\bigcup_{h \in G} h V h^{-1} \subseteq U$, so that $1_G$ admits a basis of neighborhoods which are conjugacy invariant and it follows from the result of the previous exercise that $G$ admits a compatible bi-invariant metric. 

We already knew that for any topological group each right translation is a uniform isomorphism of $(G,\mcU_l)$, and of course the same is true for each left translation. We see now that it is however rarely the case that $(g,h) \mapsto gh$ is left-uniformly continuous.

\bigskip
{\bf \ref{exo47}}.
Fix a compatible left-invariant metric $d$ on $G$. Then $\mcU_+$ is induced by the metric $\rho$ defined by $\rho(g,h)= d(g,h)+ d(g^{-1},h^{-1})$, which is complete by a corollary of Theorem \ref{t:extending_group_operations}. Hence $(G,\mcU_+)$ is complete.

If $\mcU_l$ and $\mcU_r$ coincide then they both coincide with $\mcU_+$, which is complete.

\bigskip
{\bf \ref{exo48bis}}.
Let $G$ be the inverse limit of a surjective inverse system $(G_i)_{i \in I}$ and denote by $\pi_i \colon G \to G_i$ the associated projection. Let $U$ be an open neighborhood of $1_G$ in $G$; we may assume that there exists $i \in I$ such that $U=\pi_i^{-1}(U_i)$ with $U_i$ an open neighborhood of $1$ in $G_i$. Choose $V \ni 1_{G_i}$ open such that $V\cdot V \subseteq U_i$. Since $G_i$ is Roelcke precompact, there exists a finite $F_i \subseteq G_i$ such that $VF_iV= G_i$; we can pick a finite $F \subseteq G$ such that $\pi_i(F)=F_i$.

Let $g \in G$. Then there exists $f \in F$ and $g_1, g_2 \in \pi_i^{-1}(V)$ such that $\pi_i(g) = \pi_i(g_1 f g_2)$. In particular, there exists $h \in \ker(\pi_i) \subseteq \pi_i^{-1}(V)$ such that $g = g_1 f g_2 h$, so that $g \in \pi_i^{-1}(V) F \pi_i^{-1}(V \cdot V) \subseteq U F U$. We conclude that $UFU=G$.

\bigskip
{\bf \ref{exo48}}.
To prove the first implication, proceed by induction on $n$, using that (for any group action) if $G \actson X$ and $G \actson Y$ are transitive, $G_x$ is the stabilizer of some $x \in X$ and $G_y$ is the stabilizer of some $y \in Y$ then (by the same argument as in the proof of Theorem \ref{thm:oligomorphic_is_Roelcke}) the number of orbits for the diagonal action $G \actson X \times Y$ is equal to the cardinality of the double coset space $G_x \rcoset G \lcoset G_y$. Assuming that $G \le \Sinf $ is Roelcke precompact, this cardinality is finite since $G_x$ and $G_y$ are both open subgroups of $G$ so $G_x \cap G_y$ is also an open subgroup hence there exists a finite set $F$ such that $G= (G_x \cap G_y) F (G_x \cap G_y)$.

To prove the second implication, assume that $G \le \Sinf$, and let $(x_i)_{i \in I}$ be representatives of the orbits of $G \actson \omega$. For each finite $F \subseteq I$, let $X_F = \bigcup_{i \in F} G \cdot x_i$, $\pi_F \colon G \to \mathfrak{S}(X_F)$ the restriction map and $G_F= \pi_F(G)$. The assumption on $G$ implies that each $G_F$ is oligomorphic, and $G$ is the inverse limit of the family $(G_F)_{F \textrm{ finite } \subset I}$ (where $F_1 \le F_2$ iff $F_1 \subseteq F_2$; this is a surjective inverse system). 

The last implication follows immediately from the result of the previous exercise, along with the fact that an oligomorphic subgroup of $\Sinf$ is Roelcke precompact.

\bigskip
{\bf \ref{exo49}}.
Let $G$ be a Roelcke precompact nonarchimedean Polish group and $V$ be an open subgroup of $G$. Then there exists a finite subset $F$ such that $VFV=G$, so that there are finitely many subgroups $H$ of $G$ such that $V \subseteq H \subseteq G$, since each such subgroup must be a disjoint union of double cosets $VhV$ and there are finitely many double cosets of $V$ in $G$. Since $G$ admits a countable basis of neighborhoods $(V_n)_{n< \omega}$ consisting of open subgroups, each open subgroup of $G$ must contain some $V_n$, so that $G$ has only at most countably many open subgroups.

\bigskip

{\bf \ref{exo50}}.
Denote by $i \colon \omega \to \beta \omega$ the map which assigns to each $n \in \omega$ the corresponding principal ultrafilter.
Let $j \colon \omega \to L$ be an injection with dense image in a compact space $L$, and assume that for every compact $X$ and every $f \colon \omega \to X$ there exists a continuous $\hat{f} \colon L \to X$ such that $\hat{f} \circ j= f$. We claim that there exists a homeomorphism $g \colon \beta \omega \to L$ such that $g \circ i=j$.

Indeed, applying the universal properties of $\beta \omega$ and $L$, we obtain the existence of continuous maps $g \colon \beta \omega \to L$ and $h \colon L \to \beta \omega$ such that $g \circ i = j$ and $h \circ j=i$.  Given $p=i(n) \in i(\omega)$ we then have $g \circ h (j(n))= g(i(n))=j(n)$; by density we conclude that $g \circ h = \mathrm{id}_L$. The same argument yields $h \circ g = \mathrm{id}_{\beta \omega}$.

\bigskip
{\bf \ref{exo:product_unif_continuous}}.
Let $M>0$ be such that $\|f\|_\infty \le M$ and $\|g\|_\infty \le M$. Note that, for any $x,x' \in X$ one has 
\begin{align*}
|f(x)g(x)- f(x')g(x')| &\le |f(x)(g(x)-g(x'))|+ |(f(x)-f(x'))g(x')|) \\
&\le M (|f(x)-f(x')|+ |g(x)-g(x')|) .
\end{align*}

Fix $\varepsilon >0$. Since $f$ and $g$ are uniformly continuous, there exists $U_f \in \mcU$ such that $|f(x)-f(x')| \le \frac{\varepsilon}{2M}$ whenever $(x,x') \in U_f$ and $U_g \in \mcU$ such that $|g(x)-g(x')| \le \frac{\varepsilon}{2M}$ whenever $(x,x') \in U_g$. Then $U= U_f \cap U_g \in \mcU$ and $|g(x)f(x)-g(x')f(x')| \le \varepsilon$ whenever $(x,x') \in U$.

\bigskip
{\bf \ref{exo52}}.
We need to prove that $\mathrm{RUC}_b(G)$ separates points and closed sets; pick $x \in G$ and a closed $F \subset G$ such that $x \not \in F$. Fix a bounded right-invariant metric $d$, and let $\varphi(g)=d(g,F)$. We have $\varphi(x) >0$ and $\varphi$ is constant equal to $0$ on $F$; furthermore, $\varphi$ is $1$-Lipschitz with respect to $d$, hence right-uniformly continuous.

\bigskip
{\bf \ref{exo53}}.
Let $g,h \in G$ and $y \in S(G)$. For every $\varphi \in \mathrm{RUC}_b(G)$ we have
\[ ((gh) \cdot y) (\varphi)= y((gh)^{-1}\varphi)= y(h^{-1}g^{-1} \varphi) \textrm{ and }\]
\[(g \cdot (h \cdot y))(\varphi)= (h \cdot y)( g^{-1}\varphi) = y (h^{-1} g^{-1} \varphi) .\]
So the formula indeed defines an action. Given $g,h \in G$ and $\varphi \in \mathrm{RUC}_b(G)$ we have
\[i(gh)(\varphi)= \varphi(gh) \textrm{ and } (g \cdot i(h)) (\varphi)= i(h)(g^{-1}\varphi)= (g^{-1} \varphi)(h)=\varphi(gh) \ .\]

\bigskip
{\bf \ref{exo54}}.
Use the same definition as in Exercise \ref{exo53} to define an action on $G$ on $X_A$ (the only difference is that in the product $\varphi$ belongs to $A$ instead of $\in \mathrm{RUC}_b(G)$), and the same verification to check that the left translation of $G$ on itself extends to an action of $G$ on $X_A$. Since $A \subseteq \mathrm{RUC}_b(G)$ the action of $G$ on $A$ is continuous, whence $G \actson X_A$ is continuous.

\bigskip
{\bf \ref{exo55}}.
In the case where $G$ is discrete (and countable since that is the context in which we discussed Stone-\v{C}ech compactifications; but that plays no role here), $\mathrm{RUC}_b(G) = \ell^\infty(G)$, the algebra of all bounded functions on $G$. Denote as usual by $i$ the embedding from $G$ to $S(G)$.  We know that for any compact space $X$ and any $f \colon G \to X$ there is a continuous $\hat{f}$ from $S(G)$ to $X$ such that $\hat{f} \circ i=f$. That is the universal property which characterizes $\beta G$, so that $\beta G= S(G)$.

The action defined by $A \in g \cdot p \Leftrightarrow g^{-1}A \in p$ extends the left translation of $G$ on itself, so by density of $G$ and continuity of $G \actson S(G)$ we simply have to check that this action is continuous. Since $G$ is discrete, it is enough to show that for each $g \in G$ the map $p \mapsto g \cdot p$ is continuous; this is clear since it maps each basic clopen set $[A]$ to the basic clopen set $[gA]$.

\bigskip

{\bf \ref{exo56}}.
(i) Assume that $G \actson X$ is minimal and $F \subseteq X$ is closed, $G$-invariant and nonempty. Then for any $x \in F$ we have $X= \overline{G \cdot x} \subseteq F$, so $F=X$. Conversely, if $G \actson X$ is not minimal, then one can find $x \in X$ such that $F=\overline{G \cdot x} \ne X$, and $F$ is a nontrivial $G$-invariant closed subset. 

The corresponding result for open subsets follows by considering complements. 

To prove the last equivalence, note that if for any nonempty open $U$ there exists a finite $F$ such that $F \cdot U =X$ then the only $G$-invariant open set is $X$, whence $G \actson X$ is minimal. Conversely, assume that $G \actson X$ is minimal and let $U$ be a nonempty open set. Then $G \cdot U$ is open and $G$-invariant, so by minimality $G \cdot U=X$. By compactness, since each $g \cdot U$ is open and these sets cover $X$ there exists a finite $F \subseteq G$ such that $F \cdot U=X$.

(ii) The intersection of an arbitrary decreasing family of nonempty $G$-invariant closed subsets is closed, nonempty (by compactness) and $G$-invariant. So $\supseteq$ defines an inductive ordering of the set of nonempty $G$-invariant closed subsets, and Zorn's lemma then gives us that there exists a $G$-invariant closed subset $Y$ which is maximal for $\supseteq$, i.e., minimal for inclusion among all nonempty $G$-invariant closed subsets. This implies that $Y$ cannot have a proper, $G$-invariant closed subset, so that $G \actson Y$ is minimal.

\bigskip

{\bf \ref{exo57}}. 
First, note that the action $H \actson M(H)$ extends to a continuous action $G \actson M(H)$, which is minimal since $H \actson M(H)$ is minimal.

Next, consider a minimal $G$-flow $G \actson X$; by density of $H$ in $G$ the $H$-flow $H \actson X$ is also minimal (for each $x$ and each $g \in G$, $g \cdot x \in \overline{H \cdot x}$ so  $\overline{G \cdot x}= \overline{H \cdot x}$), so by the universal property of $M(H)$ there exists a continuous, $H$-equivariant $\pi \colon M(H) \to X$. By density of $H$ in $G$ and continuity of $\pi$ we obtain that $\pi$ is also $G$-equivariant, and we are done.

\bigskip
{\bf \ref{exo57bis}}.
Assume that every action of $G$ on a compact metrizable space has a fixed point. Consider the action $G \actson S(G)$. Given a separable, closed, $G$-invariant subalgebra $A$ of $\mathrm{RUC}_b(G)$, denote by $X_A$ the corresponding (metrizable) compactification; by assumption there exists $x_A \in X_A$ such that $g \cdot x_A=x_A$ for all $g \in G$. This implies that, for every finite subset $F$ of $\mathrm{RUC}_b(G)$, the set 
\[\Sigma_F = \{p \in S(G): \forall f \in F \, \forall g \in G \ f(g \cdot p)=f(p)\} \]
is nonempty. Since each of these sets is closed, the compactness of $S(G)$ allows us to conclude that the intersection of all these sets is nonempty, i.e., there exists $p \in S(G)$ such that $f(g \cdot p)=f(p)$ for all $g \in G$. Since any minimal subflow of $S(G)$ is isomorphic to $M(G)$, we conclude that $M(G)$ is a singleton: $G$ is extremely amenable.

\bigskip

{\bf \ref{exo58}}.
Fix $c \in 2^{V \rcoset G}$ and assume that there exists a fixed point $x$ in $\overline{Gc}$. Given a finite $F \subseteq G$ there exists $g \in G$ such that $gc_{|Vf}=x_{|Vf}$, in other words $f \mapsto c(Vfg)$ is constant on $F$.

Conversely, assume that the displayed condition in the statement of Theorem \ref{thm:KPT1} holds. Then for $\varepsilon=0$ or $\varepsilon=1$ it must happen that for every $F \subseteq G$ finite there exists $g \in G$ such that $f \mapsto c(Vfg)$ is constant equal to $\varepsilon$ on $F$ (if this condition is not satisfied for $\varepsilon=0$, then there is a finite set $F$ witnessing this failure, and then the condition must be satisfied for $\varepsilon=1$ for any finite set containing $F$). Then the constant function equal to $\varepsilon$ belongs to $\overline{Gc}$.

\bigskip

{\bf \ref{exo59}}.
If $G$ is extremely amenable then, just as in the case $k=2$, one obtains that for every $c \in k^{V \rcoset G}$ there exists a fixed point in $\overline{G.c}$; the only fixed points are constant functions and we obtain the same conclusion as in the case $k=2$.

Conversely, if the condition holds for $k \ge 2$ it clearly holds for $k=2$ since $2^{V \rcoset G}$ may be seen as a subflow of $k^{V \rcoset G}$; and it follows from Theorem \ref{thm:KPT1} that $G$ is extremely amenable. 

To see that it is enough to let $V$ run over a fixed family of open subgroups forming a neighborhood basis at $1_G$, simply use the fact that if $V \le W$ are open subgroups then $k^{V \rcoset G}$ equivariantly maps onto $k^{W \rcoset G}$. So if there exists a fixed point in the closure of the orbit of every $c \in k^{V \rcoset G}$
then there exists a fixed point in the closure of the orbit of every $c \in k^{W \rcoset G}$.

\bigskip

{\bf \ref{exo60}}.
One implication is obvious. Assume that the Ramsey property for embeddings holds for $2$-colorings, and let us prove by induction that it holds for $k$-colorings for any $k \ge 2$. Assume that we have established the desired result up to some $k \ge 2$.

Pick $\bA, \bB \in \mcK$; we may assume $\bA$ embeds in $\bB$ since otherwise there is nothing to do. Find $\bC$ witnessing that the Ramsey property holds for $k$-colorings of $\emb{A}{B}$. Fix $\beta \in \emb{B}{C}$.
Then find $\bD$ witnessing that the Ramsey property holds for $2$-colorings of $\emb{A}{C}$. We claim that $\bD$ witnesses that the Ramsey property holds for $(k+1)$-colorings of $\emb{A}{B}$.

To see this, let $c \colon \emb{A}{D} \to k+1$ be a $(k+1)$-coloring. For $\alpha \in \emb{A}{D}$, define 
\[c'(\alpha)=\begin{cases} 1 & \textrm{ if } c(\alpha)=k \\ 0 & \textrm{otherwise} \end{cases} .\]

This is a $2$-coloring of $\emb{A}{D}$, so we may find $\gamma \in \emb{C}{D}$ such that $c'$ is constant on $\gamma \circ \emb{A}{C}$. If this constant is equal to $1$ then $c$ is constant on $\gamma \circ \beta \circ \emb{A}{B}$ and we are done. Else, $c$ induces a $k$-coloring $c''$ of $\emb{A}{C}$ defined by $c''(\alpha)= c(\gamma \circ \alpha)$. By assumption on $\bC$, there exists $\delta \in \emb{B}{C}$ such that $c''$ is constant on $\delta \circ \emb{A}{B}$, whence $c$ is constant on $(\gamma \circ \delta) \circ \emb{A}{B}$ and we are done.

\bigskip
{\bf \ref{exo61}}.
Assume that $\mcK$ satisfies the joint embedding property as well as the Ramsey property for embeddings. Then fix $\bA, \bB, \bC \in \mcK$ and embeddings $\alpha \colon \bA \to \bB$, $\beta \colon \bA \to \bC$. We also fix $\bE \in \mcK$ containing both $\bB$ and $\bC$ (thanks to the joint embedding property) and view $\alpha$, $\beta$ as elements of $\emb{A}{E}$.

Now pick $\bD \in \mcK$ witnessing that the Ramsey property holds for colorings of $\emb{A}{E}$ with values in $2^{\{B,C\}}$ (the subsets of $\{B,C\}$, which form a set with $4$ elements) and define a coloring of $\emb{A}{D}$ by setting: 
\[B \in c(f) \Leftrightarrow \text{ There exists an embedding } i \colon \bB \to \bD \text{ with } i \circ \alpha(A)= f(A) \]
and defining similarly when $C \in c(f)$.

We then find $\delta \in \emb{E}{D}$ such that $c$ is constant on $\delta \circ \emb{A}{E}$. Note that by definition $B \in c(\delta \circ \alpha)$, and similarly $C \in c(\delta \circ \beta)$. So $c$ is constant equal to $\{B,C\}$ on $\delta \circ \emb{A}{E}$. 

Since we know that $B \in c(\delta \circ \beta)$, it follows that there exists $i \in \emb{B}{D}$ such that $i \circ \alpha(A)= \delta \circ \beta(A)$. Since $\mcK$ has the Ramsey property for embeddings, it follows that $i \circ \alpha(a)= \delta \circ \beta(a)$ for all $a \in A$, so the embeddings $i \colon \bB \to \bD$ and $\delta_{|C} \colon \bC \to \bD$ show that $\mcK$ has the amalgamation property.

\bigskip
{\bf \ref{exo62}}.
(i) Let $A= \Q \cap \left]0,1 \right[$, with its usual ordering. It is order isomorphic to $(\Q,<)$ so we identify $G=\Aut(A,<)$ with $\Aut(\Q,<)$.

Given $g \in \Aut(A,<)$ we may extend $g$ to $\pi(g) \colon [0,1] \to [0,1]$ by setting $\pi(g)(0)=0$, $\pi(g)(1)=1$ and, for any $t \in \left]0,1 \right[$, $\pi(g)(t)= \sup\{g(q) \colon q \in A \textrm{ and } q \le t \}$. 
It is straightforward to check that $\pi(g)$ is a homeomorphism of $[0,1]$ and that $\pi \colon G \to \mathrm{Homeo}([0,1])$ is a homomorphism. 

To see that $\pi$ is continuous at $1_G$, fix $\varepsilon >0$ and find a positive integer $n$ such that $\frac{1}{n} \le \varepsilon$. Then for any $g\in G$ such that $g\left(\frac{k}{n}\right)= \frac{k}{n}$ for all $k \in \{0,\ldots,n\}$ we have that $\|\pi(g)- \mathrm{id}\|_{\infty} \le \frac{1}{n} \le \varepsilon$. Since those $g$ form an open subgroup, we conclude that $\pi$ is continuous at $1_G$, hence everywhere since it is a homomorphism.

To check that $\pi$ has a dense image, simply observe that any increasing homeomorphism may be approximated by a sequence of piecewise affine homeomorphisms with rational slopes and rational junction points, and those maps are in the image of $\pi$.

Since $G$ is extremely amenable, so is $\pi(G)$. Once all this book-keeping is complete, we conclude that $M(H_+)=M(\pi(G))$ is a singleton : $H_+$ is extremely amenable.

(ii) Let $H \actson X$ be a minimal flow. We know that $H_+ \actson X$ has a fixed point; since $H_+$ has index $2$ in $H$ ($H_+$ is the stabilizer of $0$, and $H \cdot 0= \{0,1\}$) it follows that $X$ has cardinality $1$ or $2$. Since $H \lcoset H_+$ has cardinality $2$, we conclude that $M(H)$ is equal to $H \lcoset H_+$ (with the left translation action of $H$ on $H \lcoset H_+$).

\bigskip
{\bf \ref{exo63}}.
We know that $H$ is co-precompact in $G$ iff $G/H$ is totally bounded. Basic entourages for the uniform structure on $G/H$ are of the form $\{(fH,ufH) : f \in G,  u \in U\}$ for some neighborhood $U$ of $1_G$. 

It follows from this that $H$ is co-precompact in $G$ iff for every open $U \ni 1_G$ there exists a finite $F \subseteq G$ such that $G/H= \bigcup_{f \in F} \{gH : \exists u \in U \ gH= ufH \}$, equivalently iff for every open $U \ni 1_G$ there exists a finite $F \subseteq G$ such that $G= UFH$.

\bigskip
{\bf \ref{exo:co-precompact=oligomorphic}}.It follows from the result of the previous exercise that $H$ is co-precompact in $\Sinf$ iff for every nonempty open $U$ in $\Sinf$ there exists a finite $F$ such that $\Sinf=HFU$. Furthermore, by definition of the topology of $\Sinf$ it is enough to consider the case where $U=\{\sigma : \forall i< k \ \sigma(i) =i \}$, for $k$ an integer.

Having said this, assume that $H$ is oligomorphic, and let $U=\{\sigma : \forall i < k \ \sigma(i)=i\}$ for some $k< \omega$, which we may assume to be $>0$. Denote by $\Omega_k$ the set of elements of $\omega^k$ with pairwise distinct coordinates; by assumption, there exist $x_1,\ldots,x_p \in \Omega_k$ such that for any $y \in \Omega_k$ one has $y \in H x_j$ for some $j$. Denote $x=(0,1,\ldots,k-1)$ and choose $\sigma_1,\ldots,\sigma_p \in \Sinf$ such that $\sigma_i(x)=x_i$ (where $\sigma_i(x)=(\sigma_i(0),\ldots,\sigma_i(k-1))$). We now have that for any $\sigma \in \Sinf$ there exists some $j$ such that $\sigma(x) \in H \sigma_j(x)$, equivalently, $\sigma \in H \sigma_j U$. This proves that $H$ is co-precompact.

\medskip Conversely, assume that $H$ is co-precompact, choose some nonzero integer $k$ and let $U=\{\sigma : \forall i < k \ \sigma(i)=i\}$. By assumption there exists a finite $F \subset \Sinf$ such that $\Sinf = HFU$, whence for every $\sigma \in \Sinf$ there exists $\tau \in F$ and some $h \in H$ such that $(\sigma(0),\ldots,\sigma(k-1))=(h \tau(0),\ldots,h \tau(k-1))$. This proves that the diagonal action of $H$ on $\Omega_k$ has finitely many orbits, whence $H$ is oligomorphic.

\bigskip
{\bf \ref{exo:non_metrizable_UMF_loc_cpact}}.
Let $G$ be a locally compact Polish group with metrizable universal minimal flow $M(G)$. Consider the co-precompact, extremely amenable subgroup $H$ provided by the theorem. Then $H$ is locally compact and extremely amenable, so $H$ is a singleton (one can apply Veech's theorem here, or note that the existence of a Haar measure on $H$ yields nontrivial actions on compact spaces). Hence $G$ is precompact, and we know that since $G$ is Polish this implies that $G$ is compact (see Proposition \ref{p:Solecki}). 

\bigskip

{\bf \ref{exo:surjective_Lipschitz-compact}}.
1. Let $(X,d)$ be compact metric and $\varphi \colon X \to X$ be $1$-Lipschitz and surjective.

We adapt the argument used in the proof of Lemma \ref{l:continuous_equivariant_nomalizer}. Choose $x \ne x' \in X$, let $d=d(x,y)$. Given $0< r < d(x,y)$ let $A$ be a finite $\frac{r}{4}$-covering of $(X,d)$ such that $|\{(a_1,a_2) : d(a_1,a_2) >d(x,y)- \frac{r}{2}\}|$ is minimal among all finite $\frac{r}{4}$-coverings.

Then $\varphi(A)$ is also a $\frac{r}{4}$-covering and the cardinality considered above can only decrease when passing from $A$ to $\phi(A)$. Hence for every $(a_1,a_2) \in A$ such that $d(a_1,a_2) > d(x,y)- \frac{r}{2}$ we must also have $d(\varphi(a_1),\varphi(a_2)) > d(x,y)- \frac{r}{2}$.

Pick $a_1$, $a_2 \in A$ such that $d(a_1,x) < \frac{r}{4}$ and $d(a_2,y)< \frac{r}{4}$. Then $d(a_1,a_2) > d(x,y)- \frac{r}{2}$ hence $d(\varphi(a_1),\varphi(a_2)) > d(x,y)- \frac{r}{2}$. Since $d(\varphi(x),\varphi(a_1))< \frac{r}{4}$ and $d(\varphi(y),\varphi(a_2))< \frac{r}{4}$ we conclude that $d(\varphi(x),\varphi(y)) > d(x,y) - r$. Letting $r$ go to $0$ we conclude that $d(\varphi(x),\varphi(y)) \ge d(x,y)$, hence $d(\varphi(x),\varphi(y))= d(x,y)$.

2. Let $d$ be a right-invariant compatible metric on $G$. One obtains a compatible metric on $G/H$ by setting 
\[d(g_1H,g_2H)= \inf \{d(g_1h_1,g_2h_2): h_1,h_2 \in H\} .\]
Reusing the notations of the proof of Lemma \ref{l:continuous_equivariant_nomalizer}, the argument used in the beginning of that proof implies that $\varphi \colon G \lcoset H \to G \lcoset H$ is $1$-Lipschitz; hence it extends to a $1$-Lipschitz map from $\yhwidehat{G \lcoset H}$ to $\yhwidehat{G \lcoset H}$. Since the image of $\varphi$ contains $G/H$ it is dense, hence its extension is surjective. By the result of the first question $\varphi$ is an isometry, so in particular it is injective.

\bigskip

{\bf \ref{exo:finite_union_g_delta}}.
Let $A$, $B$ be $G_\delta$ subsets; write $A= \bigcap_{i< \omega} A_i$ and $B=\bigcap_{j< \omega} B_j$.

 We then have
$A \cup B = \bigcap_{(i,j) \in \omega^2} (A_i \cup B_j)$; $\omega^2$ is countable 
and each $A_i \cup B_j$ is open as soon as $A_i$ and $B_j$ are open.

\bigskip

{\bf \ref{exo:oscillation}}.
We may assume $\overline{A}=X$. Let $x$ be such that $\omega_f(x)<r$ for some $r>0$. Then for $\varepsilon >0$ small enough we have $\sup\{d(f(a),f(a')): a,a' \in B(x,\varepsilon) \cap A \} < r$. Let $y$ belong to $B(x,\varepsilon)$ and $\delta$ be such that $d(x,y)+ \delta \le \varepsilon$. Then $B(y,\delta) \subseteq B(x,\varepsilon)$ so $\omega_f(y) < r$ and we are done.

\bigskip

{\bf \ref{exo:locally_compact_open_in_completion}}.
Fix $y \in Y$. Let $V \ni  y$ be open in $Y$ and such that the closure of $V$ in $Y$ is a compact subset (say, $K$) of $Y$. Then choose an open subset $U$ of $X$ such that $V= U \cap Y$. Since $Y$ is dense in $X$ and $U$ is open, we have $\overline{U}=\overline{U \cap Y}= \overline{V}=K$ (here the closure is taken in $X$, and the last equality holds by compactness of $K$). Hence $U$ is actually contained in $Y$, and we have shown that $Y$ is open in $X$.

\bigskip

{\bf \ref{exo:copies_Baire_Cantor}}.
It is clear that there is a copy of the Cantor space $\mcC= \{0,1\}^\omega$ contained in the Baire space $\mcN=\omega^\omega$. To see the converse, consider the map $f \colon \mcN \to \mcC$ defined by $f(\alpha)= 0^{\alpha(0)}10^{\alpha(1)}1 \ldots$.

If $\alpha_{|n}=\beta_{|n}$ then also $f(\alpha)_{|n}=f(\beta)_{|n}$ so clearly $f$ is continuous.
The image of $f$ is the subset $Y$ of $\mcC$ consisting of all $\beta$ such that $\{n : \beta(n)=1\}$ is infinite. 

Assume that $f(\alpha_n) \xrightarrow[n \to + \infty]{} f(\alpha)$ and fix $i< \omega$. Let $j$ be such that $f(\alpha)_{|j}$ takes at least $i$ times the value $1$. Then for $n$ large enough we have $f(\alpha_n)_{|j}=f(\alpha)_{|j}$, and this implies that ${\alpha_n}_{|i}=\alpha_{|i}$: $\alpha_n \xrightarrow[n \to + \infty]{} \alpha$, in other words $f$ is a homeomorphism onto its image.

\bigskip
{\bf \ref{exo:loc_finite_subcover}}.
Since $(X,d)$ is metric, each $U_n$ is a $F_\sigma$ subset of $X$. Namely, $U_n= \bigcup_{q \in \Q^*_+} F_q$, where $F_q = \lset x : d(x,X \setminus U_n) \ge q \rset$). So we write $U_n= \bigcup_{i< \omega} F_{n,i}$, where each $F_{n,i}$ is closed and $F_{n,i} \subseteq F_{n,i+1}$ for all $n$ and $i$.

Then set $V_n = U_n \setminus \bigcup_{m,i< n} F_{m,i}$. Clearly this is an open subset of $U_n$.

If $x \in \bigcup_n U_n$, then let $n= \min \lset k : x \in U_k \rset$. By definition we have $x \in V_n$, since $F_{m,i} \subseteq U_m$ for each $i$ and $m$. So $\bigcup_{n< \omega} U_n= \bigcup_{n < \omega} V_n$.

To see why the third point holds, pick $x \in \bigcup_n U_n$. Then for some $m$ and $i$ we have $x \in F_{m,i}$, so that for each $n> \max(m,i)$ we have $x \not \in V_m$.

\bigskip
{\bf \ref{exo:F_sigma_not_G_delta}}.
Write $Q= \{q_n : n < \omega \}$. Then each $\{q_n\}$ is closed so $Q= \bigcup_{n< \omega} \{q_n\}$ is $F_\sigma$. Since $X$ is perfect, each $\{q_n\}$ has empty interior, equivalently, $X \setminus \{q_n\}$ is a dense open subset of $X$. It follows that $X \setminus Q$ is dense $G_\delta$. If $Q$ were $G_\delta$ then the two dense $G_\delta$ subsets $Q$ and $X \setminus Q$ would have empty intersection, which would contradict the Baire category theorem.

\bigskip
{\bf \ref{exo:derivative}}.
Note that for any $A$ $D(A)$ is closed in $A$; since an intersection of closed sets is closed this implies that $X^{(\alpha)}$ is closed for every ordinal $\alpha$. By the Lindelöff property, there exists an ordinal $\alpha < \omega_1$ such that $\bigcap_{\beta \le \alpha} X^{(\beta)}= \bigcap_{\beta < \omega_1} X^{(\beta)}$. In particular one must have $X^{(\alpha+1)}=X^{(\alpha)}$ (so $X^{(\beta)}= X^{(\alpha)}$ for all $\beta \ge \alpha$).

Let us prove by induction that $P \subseteq X^{(\alpha)}$ for all $\alpha < \omega_1$. This is obvious if $P=\emptyset$ so we may assume that $P \ne \emptyset$. We clearly have $P \subseteq X^{(0)}$ and, since no point of $P$ is isolated in $P$, the implication $P \subseteq X^{(\alpha)} \Rightarrow P \subseteq X^{(\alpha +1)}$ is also clear. The case of limit ordinals is immediate.

Let $\alpha$ be such that $X^{(\alpha)}=X^{(\alpha +1)}$; we know that $P \subseteq X^{(\alpha)}$. If the inclusion is strict, then some element of $X^{(\alpha)}$ has a countable neighborhood in $X$, hence also in $X^{(\alpha)}$. But then, since $X^{(\alpha)}$ is Polish, it has an isolated point, so $X^{(\alpha+1)} \ne X^{(\alpha)}$, a contradiction. Hence $X^{(\alpha)}=P$.

\bigskip
{\bf \ref{exo:Baire_class_1}}. 
1.a) Write $A= \bigcap_{k< \omega} O_k$, with $O_k=\{x : d(x,A) < 2^{-k}\}$.
 
Then each $f_n^{-1}(O_k)$ is open since each $f_n$ is continuous, and we claim that 
\[ f^{-1}(A)= \bigcap_{k< \omega} \bigcup_{n\ge k} f_n^{-1}(O_k). \] 
One inclusion is immediate : if $x \in f^{-1}(A)$ then for all $k$ we have that $f(x) \in O_k$, whence $f_n(x) \in O_k$ for large enough $n$ since $O_k$ is open and $f_n(x) \xrightarrow[n \to + \infty]{} f(x)$ . 

Conversely, if $x \in \bigcap_{k< \omega} \bigcup_{n\ge k} f_n^{-1}(O_k)$ then for all $k$ we have that $f_{n_k}(x) \in O_k$ for some $n_k \ge k$, i.e., $d(f_{n_k}(x),A) < 2^{-k}$. Letting $k$ go to $+ \infty$ we obtain that $n_k \xrightarrow[k \to + \infty]{} + \infty$ so $d(f(x),A)=0$, i.e., $f(x) \in A$ since $A$ is closed.

b) Let $(V_n)_{n< \omega}$ be a basis of open sets for the topology of $Y$. Note that for every $x$, $f$ is continuous at $x$ iff for every neighborhood $V$ of $f(x)$ the set $f^{-1}(V)$ is a neighborhood of $x$. This means that the set of points at which $f$ is not continuous is equal to
\[\bigcup_{n < \omega} f^{-1}(V_n) \setminus \LongInt{f^{-1}(V_n)} . \]
In the previous question we showed that $f^{-1}(A)$ is $G_\delta$ for every closed $A$; taking complements, it follows that $f^{-1}(U)$ is $F_\sigma$ for every open $U$.

In particular, each $f^{-1}(V_n)$ is $F_\sigma$, whence $ f^{-1}(V_n) \setminus \LongInt{f^{-1}(V_n)}$ is $F_\sigma$ (it is the intersection of two $F_\sigma$ subsets) and has empty interior, hence is meager. Hence the set of points at which $f$ is not continuous is meager since it is a countable union of meager subsets of $X$.

2. For all $n$ define $f_n(x)= 2^n (f(x+2^{-n})-f(x))$. Each $f_n$ is continuous and $f_n(x) \xrightarrow[n\to + \infty]{} f'(x)$ for all $x \in \R$, whence $f'$ is a function of Baire class $1$ and we may apply the result of the previous question.

\bigskip

{\bf \ref{exo:continuity_convergent_Cantor_scheme}}.
Assume that $\alpha \ne \beta \in 2^\omega$ and let $n$ be the smallest integer such that $\alpha(n) \ne \beta(n)$. Let $s= \alpha_{|n}=\beta_{|n}$. Then we have that $f(\alpha) \in F_{s \smallfrown \alpha(n)}$ and $f(\beta) \in F_{s \smallfrown \beta(n)}$. Since $\alpha(n) \ne \beta(n)$ it follows from the definition of a Cantor scheme that $f(\alpha) \ne f(\beta)$. This shows that $f$ is injective.

To see that $f$ is continuous, fix $\alpha \in \mcC$ and $\varepsilon >0$. For $n$ large enough we have $\mathrm{diam}(F_{\alpha_{|n}}) \le \varepsilon$, so for every $\beta \in \mcC$ we have $\beta_{|n}=\alpha_{|n} \Rightarrow d(f(\alpha),f(\beta)) \le \varepsilon$.

\bigskip
{\bf \ref{exo:CH_for_Polish}}.
Let $(U_n)_{n< \omega}$ be a countable basis of open sets for the topology of $X$. Consider the map $\Phi \colon X \to \mathcal{P}(\omega)$ defined by 
\[\Phi(x)  = \{n < \omega : x \in U_n \}.\]
Since $X$ is Hausdorff we see that if $\Phi(x)=\Phi(x')$ then $x=x'$, i.e., $\Phi$ is injective.

\bigskip
{\bf \ref{exo:meager_equivalence_relations}}.
1. If every equivalence class of $R$ is meager then the Kuratowski--Ulam theorem implies that $R$ is meager in $X \times X$. Conversely, assume that $R$ is meager. Then, again by Kuratowski--Ulam, the set $\Omega=\{x : [x]_R \textrm{ is meager} \}$ is comeager in $X$. If some $[x]_R$ were not meager  then $[x]_R \cap \Omega \ne \emptyset$, which is impossible since $[x]_R$ would be both meager and  not meager in $X$ and $X$ is Polish. So every equivalence class is meager.

2. This follows directly from Mycielski's theorem applied to the complement of $R$.

\bigskip
{\bf \ref{exo:perfect_Polish_bijection_Baire}}. 
If there exists a continuous bijection from $\mcN$ to $X$ then clearly $X$ is perfect since $\mcN$ is perfect.

Conversely, assume that $X$ is perfect Polish and nonempty; fix a compatible metric $d$ and $r>0$. As a first step, we want to write $X$ as the disjoint union of a sequence of nonempty perfect $G_\delta$ subsets with small diameter.

Assume first that $X$ is not compact; then one can find a sequence of balls $(B_i)_{i< \omega}$ of diameter $\le r$  such that $X= \bigcup_{i< \omega} \overline{B_i}$ and for all $I$ one has $\bigcup_{i=0}^I \overline{B_i} \ne X$. 
Now for each $i$ let $C_i= \overline{B_i} \setminus \overline{\bigcup_{j<i} B_j}$; each $C_i$ is a perfect $G_\delta$ subset, $X= \bigsqcup_{i< \omega} C_i$ and $C_i$ must be nonempty for infinitely many $i$. So we are done in that case.

If $X$ is compact, first build an open subset $U$ of diameter $\le r$ such that $\overline{U}$ is not open. Then $\overline U$ is perfect $G_\delta$ and has small diameter, and $X \setminus \overline{U}$ is perfect and not compact, so we can apply the reasoning of the previous paragraph to write $X \setminus \overline{U} = \bigsqcup_{i< \omega} U_i$ with each $U_i$ perfect $G_\delta$.

We conclude that one may write $X= \bigsqcup_{n< \omega} U_n$, where each $U_n$ is $G_\delta$, nonempty, perfect, and has diameter less than $r$. Note that by Alexandrov's theorem each $U_n$ admits a compatible complete metric $d_n$, and we may additionally require that $d_n(u,u') \ge d(u,u')$ for every $u,u' \in U_n$.

We may use this to inductively build a Lusin scheme $(U_s)_{s \in \omega^{< \omega}}$ and compatible complete metrics $d_s$ on each $U_s$ such that:
each $U_s$ is a $G_\delta$ subset of $X$; $U_\emptyset=X$; for each $s$, $\bigsqcup_{i} U_{s \smallfrown i}=U_s$; $d_s \ge {d_t}_{|U_s}$ whenever $t \sqsubseteq s$; each $U_{s \smallfrown i}$ has diameter less than $2^{-|s|}$ for $d_s$.

Fix $\alpha \in \mcN$ and some integer $m$. Any sequence $(x_i)_{i \ge m}$ such that $x_i \in U_{\alpha_{|i}}$ for all $i$ is Cauchy in $(U_{\alpha_{|m}},d_{\alpha_{|m}})$, so it converges to some $x$ in $U_{\alpha_{|m}}$. This shows that $\bigcap_{m}U_{\alpha_{|m}}$ is nonempty; because of the conditions on diameters this intersection must then be a singleton. 

It follows that our Lusin scheme induces a continuous map $\pi \colon \mcN \to X$, and from the properties of the scheme we deduce that $\pi$ is bijective.

\bigskip
{\bf \ref{exo:continuous_open_image_of_Polish}}.
We know that there exists a continuous, open surjection of $\omega^\omega$ onto $X$, so we may as well assume that $X=\omega^\omega$. Then we can simply apply Proposition \ref{prop:continuous_open_image_of_Baire}.

\bigskip
{\bf \ref{exo:0_dim_Polish_in_Baire_and_Cantor}}.
Since $X$ is $0$-dimensional, one can build a convergent Lusin scheme $(F_s)_{s \in \omega^{< \omega}}$ with $F_\emptyset=X$, $F_s = \bigsqcup_{n} F_{s \smallfrown n}$ and each $F_s$ is clopen (possibly empty). Indeed, it suffices to notice that for every clopen $U$ and every $r<0$ one can write $U= \bigsqcup_n A_n$ where each $A_n$ is clopen (possibly empty) and of diameter less than $r$. 

Then we know that $Y=\{\alpha : \forall n \ F_{\alpha_{|n}} \ne \emptyset\}$ is closed and that the map $\pi \colon Y \to X$ associated to our convergent Lusin scheme is continuous. Furthermore, since each $F_s$ is open we obtain that $\pi$ is an open map, so $\pi \colon Y \to X$ is a homeomorphism. 

So far, we have shown that $X$ is homeomorphic to a closed subset of $\mcN$; since $\mcN$ is homeomorphic to a $G_\delta$ subset of $\mcC$, it follows that $X$ is also homeomorphic to a $G_\delta$ subset of $\mcC$.

\bigskip

{\bf \ref{exo:irrationals_vs_Baire}}.
1. We know that $X$ (with the induced topology) is Polish. Since $X$ is dense in $\R$ it is perfect. Given $x \in X$ and $\varepsilon >0$, there exists $a,b \not \in X$ such that $x- \varepsilon< a< x < b < x + \varepsilon$. Then $\left]a,b \right[ \cap X=[a,b] \cap X$ is a clopen neighborhood of $X$ contained in $B(x,\varepsilon)$; thus every element of $X$ admits a basis of clopen neighborhoods, so that $X$ is $0$-dimensional.

Let $K$ be a compact subset of $X$. Assume that the interior of $K$ is nonempty; then there exists $x \in X$ and $r>0$ such that $]x-r,x+r[ \cap X \subseteq K$. But then, by density of $X$ in $\R$, it follows that $\left]x-r,x+r \right[ \subseteq \overline{K}=K$. This contradicts the density of $\R \setminus X$. We conclude that every compact subset of $X$ has empty interior. 

All the conditions of the Alexandrov--Urysohn theorem are thus satisfied by $X$: $X$ is homeomorphic to $\mcN$.

2. Assume now that $Y$ is nonempty, $0$-dimensional and that $X$ is a $G_\delta$ subset such that both $X$ and $Y \setminus X$ are dense in $Y$. Observe that $Y$ must then be perfect, so  that $X$ is perfect too. Clearly $X$ is $0$-dimensional since it is a subset of a $0$-dimensional space. The proof that compact subsets of $X$ have empty interior is the same as above (replacing $]x-r,x+r[$ by the open ball $B(x,r)$) and we conclude similarly.

\bigskip

{\bf \ref{exo:dense_G_delta_copy_of_Baire_space}}.
1. Let $(U_n)_{n< \omega}$ be a basis of open sets for the topology of $X$. For every $n< \omega$ let $F_n = \overline{U_n} \setminus U_n$, which is closed and has empty interior. Hence $Y= X \setminus \bigcup_{n} F_n$ is dense $G_\delta$ in $X$.

For every $n$, we have that $U_n \cap Y = \overline{U_n} \cap Y$; so $U_n \cap Y$ is clopen in $Y$. Thus the topology of $Y$ admits a basis consisting of clopen sets, which implies that $Y$ is $0$-dimensional.

2. First, writing $X=P \sqcup D$ where $P$ is closed perfect and $D$ is open and (at most) countable we reduce to the case where $X$ is perfect. Then the previous result gives us a dense $G_\delta$ subset $Y$ of $X$ which is $0$-dimensional, and perfect since $X$ is perfect and $Y$ is dense. Given the result of Exercise \ref{exo:irrationals_vs_Baire}, it is thus enough to show that there exists a dense $G_\delta$ subset $Z$ of $Y$ such that $Y \setminus Z$ is also dense in $Y$ (this will also prove that every perfect Polish space contains a dense $G_\delta$ subset homeomorphic to $\mcN$).

Let $(U_n)_{n< \omega}$ be a basis of nonempty open sets for the topology of $X$. For each $n$, choose $y_n \in Y \cap U_n$ then set $Z= \bigcap_{n< \omega} Y \setminus \{y_n\}$. Then $Z$ is dense $G_\delta$ and its complement is dense in $Y$, and we are done.

3. Let $Y$ be a dense $G_\delta$ subset of $X$ which is homeomorphic to $\mcN$; abusing notation we assume $Y=\mcN$. Then consider the map $\Phi \colon \mcP(X) \to \mcP(\mcN)$ defined by $\Phi(A)=A \cap \mcN$. Note that $A$ is meager in $X$ iff $\Phi(A)$ is meager in $\mcN$; $A$ is Baire measurable in $X$ iff $\Phi(A)$ is Baire measurable in $\mcN$ (open subsets of $\mcN$ differ from open subsets of $X$ by a meager set); and $U \mathrel{\Delta} V$ is meager in $X$ iff $\Phi(U) \mathrel{\Delta} \Phi(V)$ is meager in $\mcN$. Thus $\Phi$ induces an isomorphism of the category algebras.

\bigskip

{\bf \ref{exo:topological_carac_rationals}}.
1. Let $X =\{x_n : n < \omega\}$ be countable, metrizable, perfect, and nonempty (we assume that $n \mapsto x_n$ is injective). Let $d$ be a compatible metric on $X$. Since $d$ takes only countably many values, for every $x$ the set of all $r>0$ such that the closed and open balls of radius $r$ centered at $x$ coincide is dense in $]0,+ \infty[$; this shows that $X$ is $0$-dimensional.

Let $A=\{\alpha \in \omega^\omega : \exists N \ \forall n \ge N \  \alpha(n)=0 \}$. This set is countable, perfect, $0$-dimensional and nonempty; if we prove that $X$ is homeomorphic to $A$ then we are done.

Let $U$ be a nonempty clopen subset of $X$, $x \in X$ and $r>0$. Since $U$ is perfect, we can write $U= \bigsqcup_{n< \omega} U_n$ where each $U_n$ is nonempty clopen, $x \in U_0$ and $\mathrm{diam}(U_n) \le r$ for all $n$.

This allows us to build a convergent Lusin scheme $(U_s)_{s \in \omega^{< \omega}}$ with $U_\emptyset=X$ and each $U_s$ nonempty clopen in $X$. Along with that scheme we let $y_s= x_{n(s)}$ where $n(s)$ is the smallest $n$ such that $x_n \in U_s$, and can arrange that $y_s \in U_{s \smallfrown 0}$ for each $s$ (in particular, $y_0=x_0$).

Let $Y= \{\alpha : \bigcap_{n} U_{\alpha_{|n}} \ne \emptyset\}$ and define $\pi \colon Y \to X$ by setting $\{\pi(\alpha)\}= \bigcap_{n} U_{\alpha_{|n}}$. By construction, $Y$ contains $A$, and if $\alpha= s \smallfrown 0^\infty$ then $\pi(\alpha)=y_s$. By definition of a Lusin scheme, $\pi$ is continuous and injective.

We prove by induction that for each $n$ there exists $s$ such that $n=\min \{i : x_n \in U_s\}$. This property is true for $n=0$. Assume that it holds for $0,\ldots,n$. Then for each $i \le n$ there exists $s_i \in \omega^{< \omega}$ such that $\pi(s_i \smallfrown 0^\infty)=x_i$. Let $\alpha_i= s_i \smallfrown 0^\infty$. For $N$ large enough, for each $i \le n$ the diameter of $U_{{\alpha_i}_{|n}}$ is strictly smaller than $d(x_i,x_{n+1})$. Let $s$ be such that $|s|=N$ and $x_{n+1} \in U_s$; then none of $x_0,\ldots,x_n$ belongs to $U_s$ so that $n+1=\min \{i : x_n \in U_s\}$.

It follows that for every $x \in X$ there exists $\alpha \in A$ such that $\pi(\alpha)=x$. By injectivity of our Lusin scheme, we thus have $Y=A$. 

So far, we have shown that $\pi \colon A \to X$ is a continuous bijection. By construction, we have for all $s \in \omega^{< \omega}$ that $\pi(A \cap N_s)=U_s$; since each $U_s$ is clopen it follows that $\pi$ is an open map, and we are done.

2. Clearly $\Q^n$ is countable, metrizable, perfect, and nonempty.

\bigskip

{\bf \ref{exo:Borel_pointclasses}}.
1. The stability of ${\mathbf \Sigma}_1^0(X)$ under countable unions is immediate since a union of open sets is open. For $\xi \ge 2$ the stability under countable unions is built into the definition. For finite intersections, the case $n=1$ is also immediate; and the case of limit ordinals follows immediately from the successor case.

So, assume that $\xi=\lambda+1$ and let $A= \bigcup_{n< \omega}A_n$, $B=\bigcup_{n< \omega} B_n$ where each $A_n$, $B_n$ belongs to $\mathbf{\Pi}_{\lambda}^0(X)$. Then we have
\[A \cap B  = \bigcup_{n < \omega} \bigcup_{k< \omega} (A_n \cap B_k) \]
and each $A_n \cap B_k$ belongs to ${\mathbf \Pi}_{\lambda}^0(X)$ since $\mathbf{\Pi}_\lambda^0(X)$ is stable under finite intersections. 

We conclude by transfinite induction.

2. Let us again prove by transfinite induction that these two inclusions hold true. Note that every open subset of a metrizable space if $F_\sigma$, which shows that ${\mathbf \Sigma}_{1}^0(X) \subseteq {\mathbf \Sigma}_{2}^0(X)$. For $\xi \ge 2$ the inclusion ${\mathbf \Sigma}_{\xi}^0(X) \subseteq {\mathbf \Sigma}_{\xi+1}^0(X)$ holds by definition. 

The other inclusion is even simpler : for any $\xi$ and any $A \in {\mathbf \Sigma}_{\xi}^0(X)$, taking $A_n=A$ for all $n$ shows that $A=\bigcap_n A_n$ belongs to ${\mathbf \Pi}_{\xi+1}^0(X)$.

\bigskip
{\bf \ref{exo:Borel_pointclasses_2}}.
Let us prove by transfinite induction that for all $\xi$ every element of $\mathbf{\Sigma}_\xi^0(X)$ is Borel:
\begin{itemize}
\item This is true for $\xi=1$ by definition.
\item If this property holds for $\xi$ it also holds for $\xi +1$ since the complement of a Borel subset as well as a countable union of Borel subsets are Borel.
\item If $\xi$ is limit and the property holds for every $\lambda < \xi$ then every element of ${\mathbf  \Sigma}^0_\xi(X)$ is a countable union of Borel subsets of $X$, hence ${\mathbf \Sigma}^0_\xi(X) \subset \mathcal{B}(X)$.
\end{itemize}
This shows that $\bigcup_{\xi< \omega_1} \mathbf{\Sigma}_\xi^0(X) \subseteq \mcB(X)$. 

Conversely, each $\mathbf{\Sigma}_\xi^0(X)$ is stable under unions, and complements of elements of $\mathbf{\Sigma}_\xi^0(X)$ belong to $\mathbf{\Sigma}_{\xi +1}^0(X)$ (these two facts follow from the result of the previous exercise). Hence $\bigcup_{\xi< \omega_1} \mathbf{\Sigma}_\xi^0(X)$ is a $\sigma$-algebra which contains all open sets, so $\bigcup_{\xi< \omega_1} \mathbf{\Sigma}_\lambda^0(X) \supseteq  \mcB(X)$.

The analogous result for $\bigcup_{\xi< \omega_1} \mathbf{\Pi}_\xi^0(X)$ follows immediately (complementation maps $\mcB(X)$ onto $\mcB(X)$) 

\bigskip

{\bf \ref{exo:Lebesgue}}.
It is well-known (and not hard to prove) that a pointwise limit of Borel functions is Borel so the set of Borel functions has both properties under consideration. Next, let $A \subseteq \R^\R$ satisfy both conditions. Define $A_0=C(X,\R)$. Then, by transfinite induction, let $A_{\alpha +1}$ be the set of functions which are a pointwise limit of elements of $A_\alpha$, and for limit $\alpha$ set $A_\alpha = \bigcup_{\beta < \alpha} A_\beta$. We have $A=\bigcup_{\alpha < \omega_1} A_\alpha$ (this set is contained in $A$, contains all the continuous functions and is stable under taking pointwise limits). By transfinite induction it is then straightforward to check that for any $f,g \in A$ one has both that $fg \in A$ and that $f-g \in A$.

We want to prove that every Borel function from $\R$ to $\R$ belongs to $A$. It is enough to prove this for characteristic functions of Borel subsets (since $A$ is stable under finite sums and pointwise limits). 

We consider the family $\Sigma$ of all subsets of $\R$ the characteristic function of which is a pointwise limit of a sequence of elements of $A$. 

Assume first that $I$ is an open interval; then it is not hard to write the characteristic function of $I$ as a pointwise limit of continuous functions. Since every open subset is a countable, disjoint union of open intervals and $A$ is stable under finite sums and pointwise limits, we see that open sets belong to $\Sigma$.

Observe that if $F$, $G$ belong to $\Sigma$ then so does $F \setminus G$ : if $(f_n)_n$, $(g_n)_n$ converge pointwise to $\chi_F$, $\chi_G$ respectively then $f_n(f_n-g_n)$ converges pointwise to $\chi_{F \setminus G}$. Now if $(F_n)_{n < \omega}$ is any sequence of elements of $\Sigma$ we observe that the characteristic function $\varphi_n$ of $F_n \setminus \bigcup_{m< n} F_m$ belongs to $A$. Hence the characteristic function of $\bigcup_n F_n$, which is equal to $\sum_n \varphi_n$, belongs to $A$.

This proves that $\Sigma$ is a $\sigma$-algebra which contains open sets: it contains all Borel subsets of $X$, and we are done.

\bigskip
{\bf \ref{exo:CH_for_Borel}}.
We know that there exists a Polish topology on $X$ for which $B$ is a clopen subset; in particular, there exists a topology on $B$ which turns $B$ into a Polish space. Combining the result of Exercise \ref{exo:CH_for_Polish} and Theorem \ref{t:embedding_a_Cantor}, we know that an uncountable Polish space has cardinality $2^{\aleph_0}$ and we are done.

\bigskip
{\bf \ref{exo:making_Borel_continuous}}.
Let $(U_n)_{n< \omega}$ be a basis of open sets for $(Y,\tau_Y)$. Since each $f^{-1}(U_n)$ is Borel there exists a Polish topology $\tau_n$ on $X$ refining $\tau$ such that $f^{-1}(U_n)$ is clopen in $X$. Let $\tau_\infty$ be the topology generated by $\bigcup_n \tau_n$; by the argument used in the proof of Theorem \ref{thm:making_Borel_clopen} we know that $\tau_\infty$ is Polish. By construction $f^{-1}(U)$ is open in $(X,\tau_\infty)$ for every open subset $U$ of $Y$ hence $f \colon (X,\tau_\infty) \to (Y,\tau_Y)$ is continuous.

\bigskip
{\bf \ref{exo:Borel_injective_image_of_closed_in_Baire}}.
Let $B$ be a Borel subset of $X$. There exists a Polish topology $\tau_B$ on $X$ refining the topology of $X$ such that $B$ is clopen in $(X,\tau_B)$. Hence $(B,\tau_B)$ is Polish and $\mathrm{id} \colon (B,\tau_B) \to (X,\tau)$ is continuous. Since $(B,\tau_B)$ is Polish there exists a continuous bijection $f$ from a closed subset $Y$ of $\mcN$ to $(B,\tau_B)$; then $f$ is also continuous from $Y$ to $(X,\tau)$ and we are done.

\bigskip
{\bf \ref{exo:analytic_is_projection_of_closed}}.
Since $\mcN \times X$ is Polish, any closed $F \subseteq \mcN \times X$ is also Polish; then it follows from the definition of an analytic set (and continuity of $\pi$) that $\pi(F)$ is analytic.

Conversely, assume that $A \in \ana(X)$. Then there exists a Polish space $Y$ and a continuous map $f \colon Y \to X$ such that $f(Y)=A$. Since there exists a continuous surjection $g \colon \mcN \to Y$, we may as well assume that $Y=\mcN$. Then consider $F=\{(\alpha,f(\alpha)): \alpha \in \mcN\}$. This is a closed subset of $\mcN \times X$ since $f$ is continuous, and $\pi(F)= \{f(\alpha): \alpha \in \mcN\}=f(\mcN)=A$.

\bigskip
{\bf \ref{exo:unions_intersections_analytic}}.
For each $n$ let $(Y_n,\tau_n)$ be Polish and $f_n \colon Y_n \to X$ be continuous such that $f(Y_n)=A_n$. 

Denote $Y=\prod_n Y_n$, which is a Polish space. Then $F=\{y: \forall n,m < \omega \ f_n(y(n))=f_m(y(m))\}$ is closed in $Y$. The map $f \colon y \mapsto f_0(y(0))$ is continuous on $Y$, and we have $f(F)=\bigcap_n A_n$. This proves that $\bigcap_n A_n$ is analytic.

Denote $Z= \bigsqcup_n Y_n$, and let $\tau$ be the topology on $Z$ in which each $Y_n$ is clopen and the restriction of $\tau$ to $Y_n$ coincides with $\tau_n$. Then $Z$ is Polish; and the map $f \colon Z \to X$ which is equal to $f_n$ on $Y_n$ is continuous. Clearly $f(Z)=\bigcup_{n} A_n$, hence $\bigcup_n A_n$ is analytic.

\bigskip
{\bf \ref{exo:properties_pointclass_analytic}}.
1. This is immediate since a composition of continuous maps is continuous.

2. There exists a continuous map $g \colon \mcN \to X$ such that $g(\mcN)=A$. Since $f \circ g \colon \mcN \to Y$ is Borel, there exists a finer Polish topology $\tau'$ on $\mcN$ such that $f \circ g \colon (\mcN,\tau') \to Y$ is continuous. Thus $f(A)$ is the image of a Polish space by a continuous map, i.e., $f(A)$ is analytic.

3. Denote $\pi \colon X \times Y \to X$ the projection on the first coordinate, and $\Gamma_f$ the graph of $f$, which is a Borel (hence analytic) subset of $X \times Y$. Observe that
\[ f^{-1}(B)= \pi\left( \Gamma_f \cap (X \times B) \right) .\]
Clearly a product of two analytic sets is analytic, hence $X \times B$ is analytic; $\Gamma_f$ is Borel hence analytic, and the intersection of two analytic sets is analytic by the result of the previous exercise. Hence $f^{-1}(B)$ is the image of an analytic set by a continuous map: it is analytic.

\bigskip
{\bf \ref{exo:separating_a_sequence_of_analytic_sets}}.
For each $n \ne m$ there exists a Borel subset $C_{n,m}$ such that $A_n \subseteq C_{n,m}$ and $A_m \cap C_{n,m}=\emptyset$.  For each $n$ define $B_n = \bigcap_{m \ne n} (C_{n,m} \cap (X \setminus C_{m,n}))$.

By definition, each $B_n$ is Borel and $A_n \subseteq B_n$ for all $n$. If $n \ne m$ then we have $B_n \subseteq X \setminus C_{m,n}$ and $B_m \subseteq C_{m,n}$ so that $B_n \cap B_m=\emptyset$.

\bigskip
{\bf \ref{exo:inclusion_topologies_equality_Borel}}.
Obviously, since $\tau \subseteq \tau'$, we have $\mcB(\tau) \subseteq \mcB(\tau')$. Note that $\mathrm{id} \colon (X,\tau') \to (X,\tau)$ is continuous, hence Borel, so its inverse $\mathrm{id} \colon (X,\tau) \to (X,\tau')$ is Borel by the Lusin--Suslin theorem. This amounts to the statement that $\mcB(\tau') \subseteq \mcB(\tau)$ and we are done.

\bigskip

{\bf \ref{exo:bijection_Borel_subsets_Polish}}.
If there exists a bijection $f \colon X \to Y$ such that $f(A)=B$ then we also have $f(X \setminus A)=Y \setminus B$ so both $|A|=|B|$ and $|X \setminus A|=|Y \setminus B|$.

Conversely, assume that $A$, $B$ are Borel, $|A|=|B|$ and $|X \setminus A|=|Y \setminus B|$. Since one can refine the topologies of $X$, $Y$ to finer Polish topologies so that $A$ and $B$ become clopen, and the Borel $\sigma$-algebras are unaffected by this refinement, Theorem \ref{thm:Polish_same_cardinality_Borel_bijection} implies that there exists a Borel bijection $f \colon A \to B$ (where of course $A$, $B$ are endowed with the induced topologies from $X$, $Y$). Similarly there exists a Borel bijection $g \colon X \setminus A \to Y \setminus B$. One obtains the desired Borel bijection by mapping each $x \in A$ to $f(x)$ and each $x \in X \setminus A$ to $g(x)$.

\bigskip

{\bf \ref{exo:extending_Kuratowski_theorem_to_Borel}}.
Certainly if $X$ is Borel in its completion then $X$ is homeomorphic to a Borel subset of a Polish space.

Next we prove that Borel subsets of Polish spaces (with the induced Borel $\sigma$-algebra) are standard Borel: let $X$ be a Borel subset of a Polish space $Y$; refining the topology of $Y$ without changing the Borel subsets we may assume that $X$ is clopen, hence Polish. Thus $(X,\mcB(X))$ is standard Borel.

Finally, assume that $(X,\mcB(X))$ is standard Borel. Let $Z$ be the completion of $(X,d)$. By assumption, there exists a Polish space $Y$ and a Borel isomorphism $f \colon (X,\mcB(X)) \to (Y,\mcB(Y))$. Let $(U_n)_{n< \omega}$ be a basis for the topology of $Y$;  for each $n$ there exists a Borel set $B_n$ of $Z$ such that $f^{-1}(U_n)=B_n \cap X$. Pick a Polish topology $\tau'$ on $Z$, refining the original topology, and such that each $B_n$ is open. If we endow $X$ with the topology induced by $\tau'$, $f \colon X \to Y$ is continuous; of course $f^{-1} \colon (Y,\tau_Y) \to (Z,\tau')$ may not be continuous, but it is still Borel so we can refine the topology of $Y$ so that it becomes continuous, then refine the topology on $Z$ again, and so on (all these refinements do not change the Borel structure on either $Z$ or $Y$). After $\omega$ steps we obtain Polish topologies on $Z$, $Y$ respectively which refine the original Polish topologies (hence have the same Borel subsets) and are such that $f \colon X \to Y$ is a homeomorphism. Since $Y$ is Polish this means that $X$ is $G_\delta$ in $Z$ for this finer topology, hence it is a Borel subset of $Z$. 

\bigskip

{\bf \ref{exo:trees_closed_subset}}. 
By definition, we have that $\mcT \in 2^{A^{< \omega}}$ is a tree on $A$ iff for every $s,t \in A^{< \omega}$ such that $s \sqsubseteq t$ one has $\mcT(t)=1 \Rightarrow \mcT(s)=1$. For fixed $s,t$ the subset $\{\mcT : \mcT(t)=1 \Rightarrow \mcT(s)=1\}$ is clopen in $2^{A^{< \omega}}$ so $\mathrm{Tr}(A)$ is closed in $2^{A^{< \omega}}$ since it is an intersection of clopen subsets.

\bigskip

{\bf  \ref{exo:finitely_branching_tree_compact_body}}.
Assume that $\mcT$ is finitely branching. Then for each $n$ there are finitely many elements of $\mcT$ of length $n$. This means that $[\mcT]$ is precompact in $\mcN$ (with the metric given by $d(\alpha,\beta)= \inf\{2^{-n}: \alpha_{|n}=\beta_{|n}\}$. Since $[\mcT]$ is closed in $\mcN$ and $(\mcN,d)$ is complete we conclude that $[\mcT]$ is compact.

Conversely, assume that $[\mcT]$ is compact. For every $n$, the map $[\mcT] \to \omega^n$ defined by $\alpha \mapsto \alpha_{|n}$ is continuous, hence has finite range. So there are finitely many elements in $[\mcT] \cap \omega^n$ for all $n$, which amounts to saying that $\mcT$ is finitely branching.

\bigskip
{\bf \ref{exo:section_tree}}.
This amounts to proving that for any $s \in \omega^{<\omega}$ the set $\{\alpha : s \in \mcT_\alpha\}$ is clopen in $\mcN$. Let $n=|s|$; then for any $\beta \in \mcN$ such that $\beta_{|n}=\alpha_{|n}$ we have $s \in \mcT_{\alpha} \Leftrightarrow s \in \mcT_\beta$, which gives the desired result. 

\bigskip
{\bf \ref{exo:Kleene_Brouwer_ordering}}.
1. By definition $<_\mcT$ is irreflexive, and any two distinct elements are comparable - if neither $s$ extends $t$ nor $t$ extends $s$ then for some $i < |s|, |t|$ we have $s(i) \ne t(i)$ whence $s <_{\mcT} t$ or $t <_{\mcT} s$.

Now assume that $s <_{\mcT} t <_{\mcT} u$. If $u$ extends $t$ then  clearly $s <_{\mcT} u$. Else let $i$ be smallest such that $t(i) \ne u(i)$. If $i<|s|$ then if $s(i) \ne t(i)$ either there is some $j < i$ such that $s(j) < t(j)=u(j)$, and $s <_{\mcT} u$; or $s$ and $t$ coincide up to $i$ hence $s(i) < t(i)< u(i)$ and again $s <_{\mcT} u$. If $i \ge |s|$ then either both $t$ and $u$ extend $s$, or there exists $j< |s|$ such that $s(j) < t(j)=u(j)$. In both cases we obtain that $s <_{\mcT} u$.

2. Assume that $(s_n)_{n< \omega}$ is a strictly decreasing sequence of elements of $\mcT$. Then no $s_n$ is empty, and for every $n$ we have $s_{n+1}(0) \le s_{n}(0)$. Thus there is some integer $\alpha(0)$ and an integer $N_0$ such that for every $n \ge N_0$ we have $s_n(0)=\alpha(0)$. Then for every $n \ge N_0+1$ we have $|s_n|>1$ and $s_{n+1}(1) \le s_n(1)$. So there exists $N_1 >N_0$ and some integer $\alpha(1)$ such that for all $n \ge N_1$ we have $s_n(0)=\alpha(0)$ and $s_n(1)= \alpha(1)$. Repeating this process, we inductively build $\alpha \in \mcN$ such that for all $i$ and all $n$ larger than some $N_i$ we have ${s_n}_{|i}= \alpha_{|i}$. In particular $\alpha \in [\mcT]$, whence $\mcT$ is ill-founded.

Conversely, assume that $\mcT$ is ill-founded and let $\alpha$ be an infinite branch. Then $(\alpha_{|n})_{n < \omega}$ is a strictly decreasing sequence of elements of $(\mcT,<_{\mcT})$.

\bigskip
{\bf \ref{exo:second_separation_sequence}}.
Recall that a countable union of analytic sets is analytic. Thus, applying the second separation theorem, we obtain that for each $n$ there exist disjoint coanalytic sets $C_n$, $D_n$ such that $A_n \setminus \bigcup_{m \ne n} A_m \subseteq C_n$ and $\bigcup_{m \ne n} A_m \setminus A_n \subseteq D_n$.

Now define $\widetilde{C}_n = C_n \cap \bigcap_{m \ne n} D_m$. These sets are coanalytic since countable intersections of coanalytic sets are coanalytic. By definition, $A_n \setminus A_m \subseteq D_m$ for all $m \ne n$, so $A_n \setminus \bigcup_{m \ne n} A_m$ is contained in $\bigcap_{m \ne n} D_m$. This shows that $A_n \setminus \bigcup_{m \ne n} A_m \subseteq \widetilde{C}_n$ for all $n$. Furthermore, if $n \ne m$ then $\widetilde{C}_n \subseteq C_n$ and $\widetilde{C}_m \subseteq D_n$ whence $\widetilde{C}_n \cap \widetilde{C}_m = \emptyset$.

\bigskip 
{\bf \ref{exo:sections_fixed_cardinal_coanalytic}}.
Denote $A_n= \{x \in X : |B_x|=n\}$. By definition $A_0$ is the complement of the projection of $B$ on $X$, which is analytic since $B$ is Borel. The fact that $A_1$ is coanalytic is the content of Theorem \ref{thm:sets_of_uniqueness}.

Now let $n$ be any integer $\ge 2$, fix a Borel linear ordering $\preceq$ on $Y$ (using the fact that $Y$ is Borel isomorphic to a Borel subset of $\R$) and let 
\[B_n =\{(x,y_1,\ldots,y_n) \in X \times Y^n : \forall i \ (x,y_i) \in B \textrm{ and } y_i \prec y_j \textrm{ for all } i < j \}. \]
Observe that $B_n$ is Borel, and we have that $x \in A_n$ iff there exists a unique $(y_1,\ldots,y_n) \in Y^n$ such that $(x,y_1,\ldots,y_n) \in B_n$. Hence Theorem \ref{thm:sets_of_uniqueness} again yields that $B_n$ is coanalytic.

\bigskip

{\bf \ref{exo:uniquely_branching_complete_coanalytic}}.
We already know that $\mathrm{UB}$ is coanalytic. Denote $\omega^*= \omega \setminus\{0\}$. We also know that the set of well-founded trees on $\omega^*$ is a complete coanalytic subset of $\mathrm{Tr}(\omega^*)$. Now, for each $\mcT \in \mathrm{Tr}(\omega^*)$ we define a tree $f(\mcT)$ on $\omega$ by setting that 
\[s \in f(\mcT) \Leftrightarrow (\forall i < |s| \ s(i)= 0 ) \textrm{ or } s \in \mcT .\]
Then $f$ is a continuous mapping and for each $\mcT \in \mathrm{Tr}(\omega^*)$ we have that $\mcT$ is well-founded iff $f(\mcT)$ has a unique infinite branch (namely, $0^\infty$). Hence $\mathrm{UB}$ is complete coanalytic.

\bigskip

{\bf \ref{exo:cardinalities_of_coanalytic sets}}.
1. First, fix some $\alpha \in \mathrm{LO}$ and consider 
\[\Sigma_\alpha = \lset (\beta,\varphi) : \beta \in \mathrm{LO} \textrm{ and } \varphi \textrm{ is an isomorphism from } (\omega,\alpha) \textrm{ to } (\omega,\beta) \rset .\]
This set is $G_\delta$ in $\mathrm{LO} \times \mcN$; the set of all elements of $\mathrm{LO}$ which are isomorphic to $\alpha$ is equal to $\pi(\Sigma)$, where $\pi$ is the projection on the first coordinate, hence it is analytic.

Now assume that $\alpha$ is a countable ordinal. Then $\beta \in \Sigma_\alpha \Leftrightarrow \exists! \varphi \in \mcN \ (\beta,\varphi) \in \Sigma_\alpha$, whence $\Sigma_\alpha$ is coanalytic. Thus $\Sigma_\alpha$ is both analytic and coanalytic: it is Borel. The set under consideration is the union of $\Sigma_\beta$ for $\beta< \alpha$, hence it is Borel as a countable union of Borel subsets.

Note that one could also prove that $\Sigma_\alpha$ is Borel for every $\alpha \in \mathrm{LO}$, as follows: have $\Sinf$ act continuously on $\mathrm{LO}$ (with the usual action), note that $\Sigma_\alpha$ is the orbit of $\alpha$ under this action and use the fact that orbits of continuous actions of Polish groups on Polish spaces are always Borel.

2. We know that there is a Borel map $f \colon X \to \mathrm{LO}$ such that for all $x$ one has $x \in A \Leftrightarrow f(x) \in \mathrm{WO}$. We also know that $\mathrm{WO}= \bigcup_{\alpha < \omega_1} B_\alpha$, where $B_\alpha$ is the set of all orderings isomorphic to some initial segment of $\alpha$; each $B_\alpha$ is Borel by the result of the first question. Hence $A= \bigcup_{\alpha <\omega_1} f^{-1}(B_\alpha)$ and each $f^{-1}(B_\alpha)$ is Borel. 

\bigskip
{\bf \ref{exo:Mackey}}.
1. Consider the map $\Phi \colon X \to 2^\omega$ defined by $\Phi(x)(n)=1 \Leftrightarrow x \in A_n$. This map is Borel since each $A_n$ is Borel, and injective by assumption on the sequence $(A_n)_{n< \omega}$. Thus by the Lusin--Suslin theorem $\Phi(X)$ is Borel in $2^\omega$ and $\Phi \colon X \to \Phi(X)$ is a Borel isomorphism. Since every Borel subset of $2^\omega$ belongs to the $\sigma$-algebra generated by $\{\alpha : \alpha(n)=1\}$ (where $n$ varies over $\omega$) it follows that every Borel subset of $X$ belongs to the $\sigma$-algebra generated by $(A_n)_{n< \omega}$.

2. Let $\tau'$ be another Polish group topology on $G$. Our two assumptions, along with the result of the first question, imply that every Borel subset in $(G,\tau')$ belongs to the $\sigma$-algebra generated by $(A_n)_{n< \omega}$, so every $\tau'$-Borel subset is $\tau$-Borel. Hence $\tau$ and $\tau'$ have the same Borel sets, so $\tau=\tau'$ since a Borel isomorphism between two Polish groups is a homeomorphism.

\bigskip

{\bf \ref{exo:separate_continuity_group_operations}}.
Denote $\varphi(g,h)= gh$. By Theorem \ref{thm:separate_continuity}, there exists a dense $G_\delta$ subset $\Omega$ of $G \times G$ such that $\varphi$ is continuous at every point of $\Omega$. By Kuratowski--Ulam, the set 
\[ \Sigma= \lset h : \forall^* g \in G \ \varphi \textrm{ is continuous at } (g,h) \rset\]
is comeager in $G$. If $h \in \Sigma$ and $h' \in G$ then since $\varphi(g,h'h)=\varphi(gh',h)$ and $g \mapsto gh'$ is a homeomorphism of $G$ we obtain that also $h'h \in \Sigma$. Thus $\Sigma$ is nonempty and invariant under left translation: $\Sigma=G$. 

Now fix $(g,h) \in G$ and let $(g_n,h_n)$ converge to $(g,h)$. Let $f$ be such that $\varphi$ is continuous at $(f,h)$. Then observe that $g_n h_n = gf^{-1} (fg^{-1}g_n h_n)$; by continuity of left translation by $fg^{-1}$ we have that $fg^{-1}g_n \xrightarrow[n \to + \infty]{} f$ so by assumption on $f$ we obtain $fg^{-1}g_n h_n \xrightarrow[n \to + \infty]{} fh$. Using again continuity of left translations (this time by $gf^{-1}$) we conclude that $g_n h_n$ converges to $gh$, proving that $(g,h) \mapsto gh$ is continuous everywhere.

It remains to show that $\psi : g \mapsto g^{-1}$ is continuous. By continuity of the product, the graph of $\psi$ is closed in $G \times G$, hence $\psi$ is Borel. This implies that there exists a dense $G_\delta$ subset $\Omega'$ such that the restriction of $\psi$ to $\Omega'$ is continuous. Now assume that $(g_n)_{n< \omega}$ converges to $g$. For each $n$ the set $\{k : kg_n \in \Omega'\}$ is comeager, similarly $\{k : k g \in \Omega'\}$ is comeager. Hence we may pick $k$ such that $kg_n \in \Omega'$ for all $n$ and $kg \in \Omega'$. It follows that $(kg_n)^{-1} \xrightarrow[ n\to + \infty]{} (kg)^{-1}$, i.e., $g_n^{-1}k^{-1} \xrightarrow[n \to + \infty]{} gk^{-1}$. Using continuity of right translation by $k$ we conclude that $g_n^{-1} \xrightarrow[n \to + \infty]{} g^{-1}$.

\bigskip
{\bf \ref{exo:Baire_measurable_product}}. 
There exists an open set $U$ and a meager set $M$ such that $A \mathrel{\Delta} U= M$, as well as an open set $V$ and a meager set $N$ such that $B \mathrel{\Delta} V=N$.

Then $(A \times B) \mathrel{\Delta} (U \times V) \subseteq (M \times Y) \cup (X \times N)$; and both $M \times Y$ and $X \times N$ are meager whereas $U \times V$ is open. Hence $A \times B$ is Baire measurable.

\bigskip
{\bf \ref{exo:no_Baire_measurable_well_ordering}}.
1. Certainly $I^2$ is non meager (apply Kuratowski--Ulam). Note that $I^2= \{(x,y) \in I^2: x \prec y \} \sqcup \{(x,y) \in I^2: y \prec x \} \sqcup \{(x,y) \in I^2: x=y \}$. 

The diagonal is meager since $\R$ has no isolated points. Furthermore, $\{(x,y) \in I^2: x \prec y \}$ is meager iff $\{(x,y) \in I^2: y \prec x \}$ is meager (one is obtained from the other by flipping coordinates); so either they are both meager and $I$ is meager (a contradiction) or $\{(x,y) \in I^2: x \prec y \} \subset R \cap I^2$ is not meager.

2. By Kuratoswki--Ulam there exists a non meager set $\Sigma$ such that for every $x \in \Sigma$ the set $\{y : (y,x) \in R \cap I^2\}$ is both Baire measurable and not meager. In particular $\Sigma$ is contained in $I$ and nonempty, so there exists $x \in I$ such that the associated initial segment is not meager.

3. Note that $\R$ itself is a Baire measurable, non meager initial segment; thus the set of all $x$ such that $\{y : y \preceq x\}$ is Baire measurable and not meager is nonempty.
Since $\preceq$ is a well-ordering, we can consider the smallest such $x$. Note that $I_x= \{y: y \prec x\}$ is then a non meager, Baire measurable initial segment; using the result of the previous question again, there exists $x' \in I_x$ such that $\{y : y \preceq x'\}$ is Baire measurable and non meager. This contradicts the definition of $x$.

\bigskip

{\bf \ref{exo:Suslin_schemes_measurability}}.
1. Let $A$ be an analytic subset of $X$. There exists a continuous $f \colon \mcN \to X$ such that $A=f(\mcN)$. For $s \in \omega^\omega$ let $A_s=\overline{f(N_s)}$; then each $A_s$ is closed and for each $\alpha \in \mcN$ we have $\{f(\alpha)\}= \bigcap_{k< \omega} A_{\alpha_{|k}}$. It follows that $A= \mcA((A_s)_{s \in \omega^{< \omega}})$.

2. Define a subset $B$ of $X \times \mcN$ by stipulating that 
\[(x,\alpha) \in B \Leftrightarrow \forall n < \omega \ \exists s \in \omega^n \ \alpha_{|n}=s \textrm{ and } x \in A_s .\]
It follows immediately from the definition of a Suslin scheme that $\mcA((A_s)_{s \in \omega^{< \omega}})$ is the projection of $B$. Since each $A_s$ is analytic, and countable intersections and unions of analytic sets are analytic, $B$ is analytic as well.

3. For each $\alpha \in \mcN$ we clearly have $\bigcap_{k< \omega} A_{\alpha_{|k}} = \bigcap_{k< \omega} \widetilde{A}_{\alpha_{|k}}$ so there is essentially nothing to prove.

4. We know that $B_s \subseteq A_s$, and that $A_s \in \mcS$. Let $C$ be an $\mcS$-envelope for $B_s$, then set $D=C \cap A_s$; this is an element of $\mcS$. Assume that $D' \in \mcS$ is such that $B_s \subseteq D'$. Then $D \setminus D' \subseteq C \setminus D'$; since $C$ is an $\mcS$-envelope for $B_s$ and $D' \in \mcS$ this implies that $D \setminus D' \in \mcS$. So $D$ is an $\mcS$-envelope for $B_s$ contained in $A_s$. Now, if we inductively replace each $A_s$ by an $\mcS$-envelope for $\mcA((A_t)_{s \sqsubseteq t})$ and apply the Suslin operation to this new family, the result is a set that differs from $\mcA((A_s)_{s \in \omega^{< \omega}})$ by a countable union of $\mcS$-small sets, hence by a $\mcS$-small set. Hence we may assume that each $A_s$ is an $\mcS$-envelope for $\mcA((A_t)_{s \sqsubseteq t})$.

5. Denote again $B_s=\mcA((A_t)_{s \sqsubseteq t})$. Then $B_{s}= \bigcup_{i< \omega} B_{s \smallfrown i}$, so $\bigcup_{i <\omega} A_{s \smallfrown i}$ is an $\mcS$-envelope for $B_s$. This implies that the symmetric difference of $\bigcup_{i< \omega} A_{s \smallfrown i}$ and $A_s$ is $\mcS$-small. Removing this countable family of $\mcS$-small sets, we are left with a family $(C_s)_{s < \omega^\omega}$ such that for each $s$ $C_s= \bigcup_{i< \omega} C_{s \smallfrown i}$, and the symmetric difference between $\mcA((C_s)_{s \in \omega^{< \omega}})$ and $\mcA((A_s)_{s \in \omega^{< \omega}})$ is $\mcS$-small.

By definition of a Suslin scheme, we have $\mcA((C_s)_{s \in \omega^{< \omega}}) = C_\emptyset$: one inclusion is always true, and for the other since $C_s = \bigcup_{i} C_{s \smallfrown i}$ we may build, for each $x \in C_s$, some $\alpha \in \mcN$ such that $x \in C_{\alpha_{|k}}$ for all $k < \omega$. So $\mcA((A_s)_{s \in \omega^{< \omega}})$ differs from an element of $\mcS$ by an $\mcS$-small set, and we are done.

6. In order to apply the result we just obtained, we need to show that every subset of $X$ admits an $\mcS$-envelope, where $\mcS$ is the $\sigma$-algebra of all $\mu$-measurable subsets of $X$. Let us first deal with the case where $\mu(X)$ is finite.

For any $A \subseteq X$, we may then consider $\mu^*(A)=\inf \{\mu(B): B \in \mcS \textrm{ and } A \subseteq B \} \in [0,1].$ Fix a family $(B_n)_{n< \omega}$ of subsets of $X$ containing $A$ and whose measures converge to $\mu^*(A)$, then let $B= \bigcap_n B_n$. This set is $\mu$-measurable and contains $A$. If $B'$ is $\mu$-measurable and contains $A$ then $\mu(B_n \setminus B') \le \mu(B_n \setminus A) \xrightarrow[n \to + \infty]{} 0$. So $\mu(B \setminus B')= \mu (\bigcap_n (B_n \setminus B'))=0$ (recall that we are assuming that $\mu(X)< + \infty$). Since $\mu$ is assumed to be complete we conclude that $B \setminus B' \in \mcS$ and we have shown that $A$ admits an $\mcS$-envelope.

Now assume that $\mu$ is $\sigma$-finite and let $A_n$ be a sequence of subsets such that $\mu(A_n)< + \infty$ for all $n$ and $\bigsqcup_{n< \omega} A_n=X$. For each $n$, define a measure $\tilde \mu_n$ on $\mcS$ by setting that $\tilde \mu_n(B)= 2^{-(n+1)} \mu(B \cap A_n)$ then set $\mu= \sum_{n=0}^{+ \infty} \tilde \mu_n$. This is a probability measure on $\mcS$, and for any $B \in \mcS$ we have $\mu(B)=0 \Leftrightarrow \tilde \mu(B)=0$. We conclude by applying our reasoning above to $\tilde \mu$ instead of $\mu$.

\bigskip
{\bf \ref{exo:Vaught}}. 1. Assume that $x \in A^*$ and let $h \in G$. We know that $\{g : gx \in A\}$ is comeager in $G$, hence $\{g: gh x \in A\}$ is also comeager in $G$ since $g \mapsto gh$ is a homeomorphism of $G$. So $hx \in A^*$ and we are done. The case of $A^\Delta$ is similar, replacing ``comeager'' by ``non meager'' in the previous argument.

If $A=A^*$ or $A=A^\Delta$ then $A$ is thus $G$-invariant. Conversely, if $A$ is $G$-invariant then for any $x \in A$ and any $g \in G$ we have $gx \in A$, so that $A \subseteq A^* \subseteq A^\Delta$. And if $x \in A^\Delta$ then (by the Baire category theorem applied in $G$) there exists $g \in G$ such that $gx \in A$, so $x \in A$ since $A$ is $G$-invariant.

2. If $A$ is Borel then $B=\{(x,g): gx  \in A\}$ is Borel. And $x \in A^*$ iff $B_x$ is comeager in $G$; this is a Borel condition by Theorem \ref{thm:Montgomery--Novikov}. Similarly $x \in A^\Delta$ iff there exists a nonempty open $U \subseteq G$ such that $B_x$ is comeager in $U$, again a Borel condition.

3. Since each $x \mapsto gx$ is a homeomorphism, we have that $\forall g \in G \ \forall^*x \in X \ gx \in A$. Applying the Kuratowski--Ulam theorem we obtain $\forall^* x \in X \ \forall^* g \in G \ gx \in A$, equivalently, $A^*$ is comeager in $X$.

\bigskip

{\bf \ref{exo: Vietoris_topology_Hausdorff_distance}}.
In both topologies the empty set is an isolated point so we only work with nonempty compact subsets below.

Let $\Omega = \{K: K \subseteq U \textrm{ and } \forall i < n \  K \cap U_i \ne \emptyset \}$ be a basic open subset for the Vietoris topology (with $U$, $U_0,\ldots,U_{n-1}$ all open and nonempty). Pick $K \in \Omega$. In particular $d(K,X \setminus U)=r >0$ (if $U=X$ set for instance $r=1$) since $K$ is compact, $U$ is open and $K$ is contained in $U$. Pick $x_i \in K \cap U_i$ and let $r_i$ be such that $B(x_i,r_i) \subseteq U_i$; then define $\rho = \min(\frac{r}{2},r_0,\ldots,r_{n-1}) >0$. For every $L$ such that $d(K,L) < \rho$ we have that $L \subseteq K_{\frac{r}{2}} \subseteq U$. Furthermore, every $x_i$ must belong to $L_\rho$ so for every $i$ there exists $y_i \in L$ such that $d(x_i,y_i) < r_i$, showing that $L \cap U_i \ne \emptyset$. We conclude that $\Omega$ is open for the topology induced by the Hausdorff metric.

Conversely, pick $r>0$ and $K$ a nonempty compact subset of $X$. Since $K$ is compact it is totally bounded, so there exist $x_0,\ldots,x_{n-1}$ such that $K \subseteq \bigcup_{i=0}^{n-1} B(x_i,r)$. Let $U=\{x : d(x,K) < r\}$ and $U_i= B(x_i,r)$ for all $i<n$. Then if $L \subseteq U$ we have $L \subseteq K_r$; and if $L \cap U_i \ne \emptyset$ for all $i< n$ then $K \subseteq L_r$. So the open ball centered at $K$ of radius $r$ is open for the Vietoris topology.

\bigskip

{\bf \ref{exo:Vietoris_topology_complete_separable}}.
1. Let $\Omega = \{K: K \subseteq U \textrm{ and } \forall i < n \ K \cap U_i \ne \emptyset \}$ be a basic open subset for the Vietoris topology (with $U$, $U_0,\ldots,U_{n-1}$ all open and nonempty). Since $D$ is dense there exists some $x_i \in U \cap U_i$ for all $i < n$. Then $\{x_0,\ldots,x_{n_1}\}$ is a finite subset of $D$ and belongs to $\Omega$.

2. Assume that $(K_n)_{n< \omega}$ is a Cauchy sequence in $(\mcK(X),d_H)$. We may assume that each $K_n$ is nonempty and (passing to a subsequence, and using that a Cauchy sequence with a convergent subsequence is convergent) that $d_H(K_i,K_{i+1}) < 2^{-i}$ for all $i$.

Let $K$ denote the set of all $x$ for which there exists a sequence $(x_n)_n$ which converges to $x$ such that for all $n$ there exists $m \ge n$ with $x_n \in K_m$. Assume that $(k_i)_{i < \omega}$ is a sequence of elements of $K$ which converges to some $x$. For each $i$ one can pick $n(i) \ge i$ and $x_i \in K_{n(i)}$ such that $d(x_i,k_i) < 2^{-i}$. This implies that $x \in K$: $K$ is closed.

Fix $r>0$ then find $N$ such that for all $n \ge N$ $d_H(K_n,K_N) < r$. Pick a finite subset $F$ of $K_N$ such that $\bigcup_{x \in F} B(x,r)$ covers $K_N$. Then $\bigcup_{x \in F} B(x,3r)$ covers $K$. Thus $K$ is totally bounded, so we obtain that $K$ is compact.

Fix again $r>0$ and find $N$ such that $d_H(K_n,K_m) \le r$ for all $n,m \ge N$. Given $x \in K$ there exists a sequence $(x_i)_{i< \omega}$ and $m_i \ge N$ such that $x_i \in K_{m_i}$ for all $i$ and $x_i \xrightarrow[i \to + \infty]{} x$. It follows that $d(x,K_n) \le r$ for $n \ge N$. Thus for any $n$ large enough $K \subseteq (K_n)_{r}$.

To obtain the converse inclusion, fix $N'$ such that $\sum_{i \ge N'} 2^{-i} \le r$ and pick $y \in K_n$. One may then find $(x_i)_{i\ge n}$ with $x_i \in K_i$, $x_n=y$, and $d(x_i,x_{i+1}) < 2^{-i}$ for all $i$. In particular the sequence $(x_i)_{i \ge n}$ is Cauchy in $(X,d)$, hence it converges to some $x \in K$. Furthermore we have $d(x,y) \le r$ so $K_n \subseteq K_r$ for every $n \ge N'$ and we are done.

3. Assume that $X$ is compact, and let $r>0$. Pick a finite set $F$ such that $X \subseteq F_{r}$. Given $K \in \mcK(X)$, let $F'=\{x \in F : B(x,r) \cap K \ne \emptyset\}$. By definition $K \subseteq F'_r$ and $F' \subseteq K_r$. This proves that $(\mcK(X),d_H)$ is totally bounded ($F$ has only finitely many subsets, and any element of $K$ is at distance $\le r$ of one of those subsets). We already know that $(\mcK(X),d_H)$ is complete so $\mcK(X)$ is compact.

4. Let $U$ be open. Then for any $K,L \in \mcK(X)$ we have that $(K \cup L) \cap U \ne \emptyset \Leftrightarrow (K \cap U \ne \emptyset) \textrm{ or } (L \cap U \ne \emptyset)$; so these sets form an open subset in $\mcK(X)^2$. Similarly $(K \cup L) \subseteq U$ iff both $K \subseteq U$ and $L \subseteq U$, so we see that $(K,L) \mapsto K \cup L$ is actually continuous. 

Let us deal now with the intersection map. The set of all $(K,L)$ such that $K \cap L \ne \emptyset$ is easily seen to be closed. Since every open subset of $X$ is an increasing, countable union of closed subsets, this implies that $\{(K,L) : K \cap L \cap U \ne \emptyset\}$ is $F_\sigma$ for every open subset $U$. And the set $\{(K,L) : K \cap L \subseteq U\}$ is open. This shows that $(K,L) \mapsto K \cap L$ is Borel (but not continuous in general, even if $X$ is assumed to be compact).

\bigskip
{\bf \ref{exo:alternative_cantor_scheme_Vietoris}}.
1. Let $(U_n)_{n< \omega}$ be a countable basis for the topology of $X$. Note that $K$ is perfect iff for every $n$, either $K \cap U_n = \emptyset$ (a closed condition) or there exist $i,j$ such that $U_i$, $U_j$ are disjoint, contained in $U_n$ and both $K \cap U_i \ne \emptyset$ and $K \cap U_j \ne \emptyset$ (a union of open conditions, hence an open condition). So being perfect is a $G_\delta$ condition.

To prove that perfect subsets are dense in $\mcK^*(X)$ one could show that each of the $G_\delta$ subsets considered above is dense and apply the Baire category theorem. Alternatively, pick a nonempty $K \in \mcK(X)$, $r>0$ and let $L= \overline{\bigcup_{x \in K} B(x,r)}$. Then $d_H(K,L) \le r$ and $L$ is perfect.

2. Let $D$ be a nonempty finite subset of $X$ and $r>0$. If $f_{|D}$ is not injective, by using the fact that $f$ is not constant on any $B(d,r)$ we can produce $D'$ finite such that $d_H(D,D') \le r$ and $f_{|D'}$ is injective. So $\{K \in \mcK^*(X) : f_{|K} \textrm{ is injective} \}$ is dense in $\mcK^*(X)$ (recall that finite subsets are dense in $\mcK^*(X)$), and it remains to show that it is $G_\delta$.

Let $(U_n)_{n<\omega}$ be a countable basis for the topology of $X$. It is not hard to check that $f_{|K}$ is injective iff 
\[\forall n,m \ (K \cap \overline{U_n} \ne \emptyset \textrm{ and } K \cap \overline{U_m} \ne \emptyset \textrm{ and } K \cap \overline{U_n} \cap \overline{U_m}= \emptyset) \Rightarrow (f(K) \cap f(\overline{U_n}) \cap f(\overline{U_m}) = \emptyset ). \]

Note that $\{K : K \cap F \ne \emptyset \}$ is closed for any closed $F \subseteq X$. Since $K \mapsto f(K)$ is continuous (because $f$ is uniformly continuous), it follows that $\{K : f_{|K} \textrm{ is injective} \}$ is $G_\delta$.

3. First, apply a transfer theorem to reduce to the case where $X$ is $0$-dimensional and $f$ is continuous. Then remove from $X$ the union of all open subsets $U$ such that $f(U)$ is countable to obtain a closed subset $X_1$; by the Lindelöff property, $f(X_1) \setminus f(X)$ is countable so for every nonempty open subset $U$ of $X_1$ we have that $f(U)$ is uncountable. Applying the Cantor--Bendixson theorem to $X_1$, we may as well assume that $X_1$ is perfect. Then we know that a generic element of $\mcK^*(X_1)$ is perfect and $0$-dimensional: it is homeomorphic to $\mcC$. And we also know that for a generic element $K$ of $\mcK^*(X_1)$ the restriction $f_{|K}$ is injective. We conclude by applying the Baire category theorem in $\mcK^*(X_1)$.

\bigskip

{\bf \ref{exo:Cantor_independent_over_Q}}.
1. We consider the space $E$ of all continuous functions from $\mcC$ to $X$, endowed with the metric $d_\infty(f,g)= \sup\{d(f(x),g(x)): x \in X\}$, which turns $E$ into a Polish space.

As a warm-up, let us show that injective functions form a dense $G_\delta$ subset of $E$. Indeed, fix $n$ and consider $\Omega_n=\{f: \exists (x,y) \ d(x,y) \ge 2^{-n} \textrm{ and } f(x)=f(y) \}$. By compactness of $\mcC$ $\Omega_n$ is closed, hence its complement $\Sigma_n$ is open. To show that $\Sigma_n$ is dense, pick a continuous function $f$ and some $\varepsilon >0$; since $\mcC$ is compact and $0$-dimensional, we may write $\mcC = \bigsqcup_{i=0}^N A_i$, where each $A_i$ is clopen, has diameter $< 2^{-n}$ and the oscillation of $f$ on each $A_i$ is smaller than $\varepsilon$. Using the fact that $X$ is perfect, we may then define a function $g$ which is constant on each $A_i$ and some $a_i \in A_i$ such that $d(g(a_i),f(a_i)) \le \varepsilon$ for all $i$ and $g(a_i) \ne g(a_j)$ for each $i \ne j$. In particular, $g \in \Sigma_n$ and $d_\infty(f,g) \le 2 \varepsilon$. Thus $\Sigma_n$ is dense open, so $\Sigma=\bigcap_n \Sigma_n = \{f \in E: f \textrm{ is injective }\}$ is dense $G_\delta$.

Now fix $n \ge 1$, fix $U$ dense $G_\delta$ in $X^n$ and consider 
\[\Delta_U = \{f : (\forall i \ne j \ x_i \ne x_j) \Rightarrow (f(x_1),\ldots,f(x_n)) \in U \}. \]

Assume first that $U$ is dense open in $X^n$.  An argument similar to the one above shows that $\Delta_U$ is dense $G_\delta$ in $E$ : consider functions such that 
\[ \forall i \ne j \ d(x_i,x_j) \ge 2^{-p} \Rightarrow (f(x_1),\ldots,f(x_n)) \in U  \]
for some fixed $p$, show that these functions form a dense open subset of $E$ using a variation on the argument above, and take the intersection. If $U$ is dense $G_\delta$, then $U= \bigcap_{i< \omega} U_i$ with $U_i$ dense open and $\Delta_U = \bigcap_{i< \omega} \Delta_{U_i}$ so $\Delta_U$ is still dense $G_\delta$ in $E$.

We conclude, by applying the Baire category theorem in $E$, that the set of all $f$ satisfying the desired conclusion is dense $G_\delta$ in $E$ hence, in particular, nonempty.

2. Let $X=\R$ and $R_n =\{(x_1,\ldots,x_n) \in \R^n : x_1,\ldots,x_n \textrm{ are linearly independent over } \Q\}$. The complement of $R_n$ is a countable union of sets of the form $\{(x_1,\ldots,x_n) : \sum_{i=1}^n \lambda_i x_i =0\}$ where $(\lambda_1,\ldots,\lambda_n)$ is a nonzero element of $\Q^n$, so the complement of $R_n$ is meager, i.e., $R_n$ is comeager. The result of the previous question then gives us the desired Cantor subset of $\R$.

\bigskip

{\bf \ref{exo:belonging_is_Borel_in_Effros}}.
Let $(U_n)_{n< \omega}$ be a countable basis for the topology of $X$. 
For a closed $F$ and $x \in X$ we have that 
\[ (x \in F ) \Leftrightarrow (\forall n  \ (x \in U_n \Rightarrow F \cap U_n \ne \emptyset)) .\]
For each fixed $n$, the set $\{(F,x): x \not \in U_n \textrm{ or } F \cap U_n \ne \emptyset \}$ is a union of two Borel subsets of $\mcF(X) \times X$, hence is Borel.

\bigskip
{\bf \ref{exo:compact_Borel_in_Effros}}.
Pick a metrizable proper compactification $(Y,d)$ of $X$. Consider the map $\Phi \colon K \mapsto K$, from $\mcK(X)$ to $\mcK(Y)$. If $U$ is open in $Y$ then $V  =U\cap X$ is open in $X$ and a compact subset $K$ of $X$ is contained in $U$ iff it is contained in $V$; similarly $K \cap U \ne \emptyset$ iff $K \cap V \ne \emptyset$. This shows that $\Phi$ is continuous; one shows similarly that $\Phi$ is a homeomorphism onto its image.  

Since the Borel structure on $\mcF(X)$ comes from the map $F \mapsto \overline{F}$ from $\mcF(X)$ to $\mcK(Y)$, this immediately implies that $\mcK(X)$ is Borel in $\mcF(X)$.

\bigskip
{\bf \ref{exo:Burgess}}.
Let us consider the equivalence relation $E$ on $G$ such that $(g,g') \in E$ iff $g' \in gH$; $E$ has closed equivalence classes and for every open $U$ we have (using the notations of Theorem \ref{thm:Burgess}) $U_E=UH$, which is open. Hence there exists a Borel map $s \colon G \to G$ such  that for all $g$ one has $s(g) \in gH$ and $gH=g'H \Rightarrow s(g)=s(g')$. This map induces an injective, Borel map $s \colon G \lcoset H \to G$. Then $T=s(G\lcoset H)$ is Borel since $s$ is Borel and injective, and $T$ intersects each left $H$-coset in exactly one point.

\bigskip
{\bf \ref{exo:standard_space_of_Banach_spaces}}.
1. Let $X$ be a separable Banach space. Let $B$ be the closed unit ball of $X^*$ endowed with the weak* topology (i.e., the topology generated by the maps $\delta_x \colon \varphi \mapsto \varphi(x)$ for $x \in X$). Then $B$ is compact by the Banach--Alaoglu theorem, and metrizable because $X$ is separable. Hence there exists a surjective, continuous map $\pi \colon \mcC \to B$.

For every $x \in X$ the map $\delta_x$ is continuous, $x \mapsto \delta_x$ is a linear map from $X$ to the set of all continuous functions on $B$, and $\|\delta_x\|_\infty  = \sup \{|\varphi(x)| : \varphi \in B \}= \|x\|$ by Hahn--Banach. Hence $X$ embeds linearly and isometrically into the space of continuous functions from $B$ to $\R$.

Similarly, the map $f \mapsto \pi \circ f$ is a linear isometry from the space of continuous functions on $B$ into the space of continuous functions from $\mcC$ to $\R$ endowed with the sup norm. We obtain the desired result by composing these two linear isometries.

2. Let $(s_n)_{n< \omega}$ be a sequence of Borel selectors from $\mcF(E)$ to $E$ provided by the Kuratowski--Ryll-Nardzewski theorem. By continuity of addition, we see that $F \in \mcF(E)$ is stable under addition iff 
\[\forall n,m < \omega \quad  s_n(F)+ s_m(F) \in F\]
and this is a Borel condition. Similarly for multiplication by a scalar, and we know that $\{F: 0 \in F\}$ is also Borel. So $\mcB$ is Borel.

3. A Banach space is isometric to a Hilbert space iff it satisfies the parallelogram identity, which by continuity is equivalent to 
\[\forall n,m < \omega \quad \|s_n(F)+s_m(F)\|^2 + \|s_n(F)-s_m(F)\|^2= 2 \|s_n(F)\|^2 + 2 \|s_m(F)\|^2 .\]
This is a Borel condition.

To deal with finite dimensional Banach spaces, recall that $F$ is finite dimensional iff the unit sphere of $F$ is precompact. Thus $F$ is finite dimensional iff for every $\varepsilon \in \Q^+$ there exist $n_1,\ldots,n_k \in \omega$ such that 
\[\forall n < \omega \quad (s_n(F) \ne 0) \Rightarrow \left(\exists i \in \{1,\ldots,k\} \ \left \|\frac{s_n(F)}{\|s_n(F)\|} - \frac{s_{n_i}(F)}{\|s_{n_i}(F)\|} \right \| < \varepsilon  \right).\]

\bigskip

{\bf \ref{exo:graph_generates_E_0}}.
Let us prove by induction on $n$ that if $x(i)=y(i)$ for all $i \ge n$ then $(x,y)$ belongs to $E_{\mcG_S}$. This is obvious if $n=0$ since then $x=y$. Assume that the property holds for $n$ and consider $x,y$ such that $x(i)=y(i)$ for all $i \ge n+1$. We may assume that $x(n) \ne y(n)$ otherwise there is nothing to do. So there exists $u,v \in 2^n$, $z \in 2^\omega$ and $\varepsilon \in \{0,1\}$ such that $x= u \smallfrown \varepsilon \smallfrown z$ and $y= v \smallfrown (1- \varepsilon) \smallfrown z$. Choose $s \in S \cap 2^n$. By definition of $G_S$, $(s \smallfrown \varepsilon \smallfrown z, s\smallfrown (1- \varepsilon) \smallfrown z) \in \mcG_S$. By the inductive hypothesis, there is a path in $\mcG_S$ from $u \smallfrown \varepsilon \smallfrown z$ to $s \smallfrown \varepsilon \smallfrown z$, as well as a path from $s\smallfrown (1- \varepsilon) \smallfrown z$ to $v\smallfrown (1- \varepsilon) \smallfrown z$. Hence there is a path in $\mcG_S$ from $x$ to $y$.

\bigskip

{\bf \ref{exo:generic_ergodicity_E_0}}.
Note that every subset of $2^\omega$ which is Baire measurable and $E_0$-invariant is either meager or comeager (see Exercise \ref{exo11}); indeed $E_0$ is induced by an action of a group acting by homeomorphisms on $2^\omega$, one can for instance take $G$ to be the countable group of all maps of the form $s \smallfrown \varepsilon \smallfrown z \mapsto   s \smallfrown (1- \varepsilon) \smallfrown z$ (There is actually a continuous $\Z$-action on $2^\omega$ whose induced equivalence relation is exactly $E_0$, but that is quite a bit harder to see). 

In particular, for every Borel subset $U$ of $X$ the subset $f^{-1}(U)$ is either meager or comeager. Fix a compatible complete metric $d$ on $X$. Covering $U$ by countably many balls of radius $ \le 1$, we see that there exists an open ball $B_0$ of radius less than $1$ such that $f^{-1}(B_0)$ is comeager. Repeating this within $B_0$ produces some ball $B_1$ of radius less than $\frac{1}{2}$ whose closure is contained in $B_0$ and such that $f^{-1}(B_1)$ is comeager. Inductively applying this process yields a sequence of open balls of vanishing diameter such that for all $n$ $\overline{B_{n+1}} \subseteq B_n$ and $f^{-1}(B_n)$ is comeager for all $n$. Then $\bigcap_n B_n = \{x\}$ for some $x \in X$, and $f^{-1}(\{x\})= \bigcap_{n< \omega} f^{-1}(B_n)$ is comeager.

\bigskip
{\bf \ref{exo:it_is_possible_to_be_dense_and_sparse}}.
Enumerate $2^{<\omega} = \{t_n : n <\omega\}$ with $|t_n| \le n$ for all $n$ (for instance $t_0=\emptyset$, $t_1=0$, $t_2=1$, $t_3=00$, $t_4=01$...). Then for each $n$ pick any $s_n$ of length $n$ which extends $t_n$ and let $S=\{s_n : n < \omega\}$.

\bigskip
{\bf \ref{exo:Lusin_Novikov_section_fixed_cardinality}}.
For $n=0$ we have that $\{x : |A_x|=0\}$ is the complement of the projection of $A$ on $X$, hence it is Borel by the Lusin--Novikov theorem.

Now write $A= \bigcup_n A_n$, where $A_n$ is the graph of a Borel function $s_n$ from some Borel $X_n \subseteq X$ to $Y$; let $Z=\bigcup X_n$ be the projection of $A$ on $X$. We may as well assume (by taking $s_n$ to be an appropriate constant on $Z \setminus X_n$) that each $s_n$ is defined on $Z$.

Fix $1 \le n < \omega$. Then for any $x \in X$ we have that $|A_x|=n$ iff there exists $i_1,\ldots,i_n$ such that 
\[\forall k \ne l \in \{1,\ldots,n\} \ s_{i_k}(x) \ne s_{i_l}(x) \textrm{ and } \forall j \ \exists k \in \{1,\ldots,n\} \ s_j(x)=s_{i_k}(x) .\]
This shows that for each $n< \omega$ the set $Z_n= \{x : |A_x|=n\}$ is Borel, and $Z \setminus \bigcup_{n< \omega} Z_n$ is precisely $\{x: |A_x|= \omega\}$, which is thus also Borel.

\bigskip
{\bf \ref{exo:Lusin_Novikov_injective}}.
Applying the Lusin--Novikov theorem, we obtain a sequence of partial Borel functions $g_n$ such that $R$ is the union of the graphs $G_n=\{(x,g_n(x)): x \in \mathrm{dom}(g_n)\}$. Replacing $G_n$ by $G_n \setminus \bigcup_{i<n} G_i$, we may assume that $R= \bigsqcup_{n< \omega} G_n$ (of course it may happen that $G_n$ is empty for some $n$).

There is no reason why $g_n$ should be injective; but since $R$ is symmetric we also have $R= \bigsqcup_{n< \omega} G_n^{-1}= \bigsqcup \{(g_n(x),x): x \in \mathrm{dom}(g_n)\}$. Now note that for every $n,m$ $\Gamma_{n,m}=G_n \cap G_m^{-1}$ is the graph of a partial Borel function $f_{n,m}$, which is injective by definition. It only remains to observe that 
\[A = \bigsqcup_{n < \omega} G_n = \bigsqcup_{n,m < \omega} \Gamma_{n,m} .\]

\bigskip

{\bf \ref{exo:G_s_is_connected_and_acyclic}}.
To prove that $\mcG_S(n)$ is connected and acyclic for all $n$, it is enough to prove that if $H$ is connected and acyclic, and $u$ is a vertex of $H$, then $H \sg{u} H$ is still connected and acyclic. Denoting by $n$ the number of vertices of $H$, we know that $H$ has $n-1$ edges, and $H \sg{u} H$ has $2n$ vertices and $2(n-1)+1=2n-1$ edges. So it is enough to prove that $H \sg{u} H$ is connected. Let $v$, $w$ be two vertices in $H \sg{u} H$. Either they belong to the same copy of $H$, in which case we can simply use a path in $H$ to connect them; or we have $v=v'\smallfrown \varepsilon$, $w=w' \smallfrown (1-\varepsilon)$ for some $\varepsilon \in \{0,1\}$. Then there is a path in $H \sg{u} H$ from $v$ to $u\smallfrown \varepsilon$ (because $H$ is connected), a path in $H \sg{u} H$ from $w$ to $u \smallfrown (1-\varepsilon)$ (again because $H$ is connected) and an edge in $H \sg{u} H$ between $u\smallfrown \varepsilon$ and $u \smallfrown (1-\varepsilon)$. This shows that $H \sg{u} H$ is connected.

In particular, the graphs $\mcG_S(n)$ are acyclic and connected for all $n$. If there were a nontrivial cycle in $\mcG_S$, it would be visible in some $\widetilde{\mcG}_S(n)$, and this would produce a nontrivial cycle in $\mcG_S(n)$, which is impossible. 

\bigskip
{\bf \ref{exo:Bernshteyn}}.
Let $(U_n)_{n< \omega}$ be a countable basis for the topology of $X$. Since every $x$ has only finitely many neighbors in $G$, there exists $i$ such that $x \in U_i$ and every $y \ne x$ such that $(x,y) \in G$ is outside of $U_i$. So we may define 
\[n(x)= \min \{i : x \in U_i \textrm{ and every neighbor of } x \textrm{ in } G \textrm{ is outside of } U_i\}. \]
Obviously if $(x,y) \in G$ then $n(x) \ne n(y)$, i.e., $n$ is an $\omega$-coloring. 

We still have to show that $n$ is Borel. Given the definition of $n$, it is enough to prove that for all open $U$ the set 
 $\{y : \textrm{there exists } z \in U \textrm{ such that } (y,z) \in G \}$
is Borel. This follows from the Lusin--Novikov theorem, since this set is the projection on $X$ of the Borel set $(X \times U) \cap G$, which is Borel and has finite vertical sections.

\pagestyle{empty}
\let\oldaddcontentsline\addcontentsline
\renewcommand{\addcontentsline}[3]{}
\let\addcontentsline\oldaddcontentsline
\printbibliography
\printindex
\end{document}